\documentclass[a4paper,12pt,bibliography=totocnumbered,numbers=noenddot]{scrreprt}

\usepackage{amsmath}
\usepackage{amssymb}
\usepackage{amsthm}
\usepackage[ansinew]{inputenc} 
\usepackage{enumerate}
\usepackage[ngerman]{babel}
\usepackage{bm} 
\usepackage{setspace} 
\usepackage[arrow, matrix, curve]{xy}
\usepackage{bibgerm}
\usepackage{makeidx} 
\usepackage[bookmarks,raiselinks,pageanchor,hyperindex,colorlinks,citecolor=black,linkcolor=black,urlcolor=black,filecolor=black,menucolor=black]{hyperref}
\usepackage{glossaries}

\allowdisplaybreaks 

\makeindex 

\makeglossaries

\glossarystyle{super}

\newtheoremstyle{dotless}{}{}{\itshape}{}{\bfseries}{}{ }{}
\theoremstyle{dotless} 

\newtheorem{lem}{Lemma}[chapter]
\newtheorem{defi}{Definition}[chapter] 
\newtheorem{bem}{Bemerkung}[chapter] 
\newtheorem{satz}{Satz}[chapter]
 
\newtheorem{fol}{Folgerung}[chapter]

\newtheorem{haupt}{Hauptsatz}[chapter]
\newtheorem*{satzohnenum}{Satz}
\newtheorem*{hauptohnenum}{Hauptsatz}

\newenvironment{bew}{\begin{proof}[Beweis]}{\end{proof}}

\usepackage{fancyhdr}
\begin{document}  

\renewcommand{\glossaryname}{Notation} 

\begin{center}
\vspace*{1,5cm}
\begin{Large}
\doublespacing
\begin{scshape}
Isometrische Einbettung von geschlossenen Riemannschen Mannigfaltigkeiten
\end{scshape}
\singlespacing
\end{Large}
\vspace{\fill}
An der Fakultät für Mathematik\\
der Otto-von-Guericke-Universität Magdeburg\\
zur Erlangung des akademischen Grades \\

Diplom-Mathematiker
\\
angefertigte\\
\vspace{\fill}
\begin{Large}

\textbf{Diplomarbeit}  
\\
\end{Large}
\vspace{\fill}

vorgelegt von\\
\textsc{Norman Zergänge}\\
geboren am 30.04.1988 in Staßfurt,\\
Studiengang Mathematik.\\ [2ex]
19. November 2013\\
\vspace{\fill}
Betreut am Institut für Analysis und Numerik von\\
\begin{scshape}
Prof. Dr. habil. Miles Simon
\end{scshape}
\end{center}

\thispagestyle{empty}

\onehalfspacing

\tableofcontents
\thispagestyle{empty}

\newglossaryentry{int(S)}
{name=\ensuremath{int(S)},
description={Inneres von S},
sort=Inneres von S}

\newglossaryentry{overline{S}}
{name=\ensuremath{\overline{S}},
description={Abschluss von S},
sort=Abschluss}

\newglossaryentry{supp(f)}
{name=\ensuremath{supp(f)},
description={Tr\"ager von f},
sort=Tr\"ager von f}

\newglossaryentry{(X,mathcal{T})}
{name=\ensuremath{(X,\mathcal{T})},
description={topologischer Raum},
sort=topologischer Raum}

\newglossaryentry{randS}
{name=\ensuremath{\partial S},
description={Rand von S},
sort=Rand von S}

\newglossaryentry{Metrik}
{name=\ensuremath{d},
description={Metrik},
sort=Metrik}

\newglossaryentry{metrischer Raum}
{name=\ensuremath{(M,d)},
description={metrischer Raum},
sort=metrischer Raum}

\newglossaryentry{offener Ball}
{name=\ensuremath{B_R(x)},
description={offener Ball um $x$ mit dem Radius $R$},
sort=offener Ball um x mit dem Radius R}

\newglossaryentry{euk}
{name=\ensuremath{|x|_{\mathbb{R}^n}},
description={euklidische Norm},
sort=euklidische Norm}

\newglossaryentry{dim(M)}
{name=\ensuremath{dim(M)},
description={Dimension einer topologischen Mannigfaltigkeit $M$},
sort=Dimension einer topologischen Mannigfaltigkeit M}

\newglossaryentry{M}
{name=\ensuremath{M},
description={glatte Mannigfaltigkeit},
sort=glatte Mannigfaltigkeit}

\newglossaryentry{Atlas}
{name=\ensuremath{\mathcal{A}},
description={Atlas auf $M$},
sort=Atlas auf M}

\newglossaryentry{prod_{i=1}^k {M_i}}
{name=\ensuremath{\prod_{i=1}^k {M_i}},
description={kartesisches Produkt der Mengen $M_1,...,M_k$},
sort=kartesische Produkt}

\newglossaryentry{id_S}
{name=\ensuremath{id_S},
description={Identit\"at auf S},
sort=Identit\"at auf S}

\newglossaryentry{C^{infty}(M,N)}
{name=\ensuremath{C^{\infty}(M,N)},
description={Menge aller glatten Abbildungen zwischen $M$ und $N$},
sort=Menge aller glatten Abbildungen zwischen M und N}

\newglossaryentry{C^{infty}(M)}
{name=\ensuremath{C^{\infty}(M)},
description={Raum aller glatten reellen Abbildungen auf $M$},
sort=Menge aller glatten Abbildungen zwischen M und R}

\newglossaryentry{T_p M}
{name=\ensuremath{T_p M},
description={Tangentialraum von $M$ im am Punkt $p$},
sort=Tangentialraum von M im am Punkt p}

\newglossaryentry{^{varphi}_{psi}F}
{name=\ensuremath{^{\varphi}_{\psi}F},
description={Koordinatendarstellung von F},
sort=Koordinatendarstellung von F}

\newglossaryentry{^{varphi}F}
{name=\ensuremath{^{\varphi}F},
description={Koordinatendarstellung von F},
sort=Koordinatendarstellung von F}

\newglossaryentry{phi i Koord}
{name=\ensuremath{\left. ^{\varphi}\frac{\partial}{\partial x^i}\right|_{p}},
description={i-ter Koordinatenvektor},
sort=i-ter Koordinatenvektor}

\newglossaryentry{i Koord}
{name=\ensuremath{\left. \frac{\partial}{\partial x^i}\right|_{p}},
description={i-ter Koordinatenvektor},
sort=i-ter Koordinatenvektor}

\newglossaryentry{Dualraum}
{name=\ensuremath{T_p^{\ast} M},
description={Dualraum zu $T_p M$},
sort=Dualraum zu $T_p M$}

\newglossaryentry{i Dual}
{name=\ensuremath{\left. ^{\varphi}dx^i\right|_{p}},
description={i-ter dualer Basisvektor},
sort=i-ter dualer Basisvektor}

\newglossaryentry{Vektorbundel}
{name=\ensuremath{(E,M,\pi)},
description={Vektorb\"undel},
sort=Vektorb\"undel}

\newglossaryentry{Raum kov 2}
{name=\ensuremath{T^2(T_p M)},
description={Vektorraum aller kovarianten 2- Tensoren auf $T_p M$},
sort=Vektorraum aller kovarianten 2- Tensoren auf T_p M}

\newglossaryentry{disVerein}
{name=\ensuremath{\coprod_{\alpha\in A}{X_{\alpha}}},
description={disjunkte Vereinigung},
sort=disjunkte Vereinigung}

\newglossaryentry{BunKov}
{name=\ensuremath{T^2(M)},
description={B\"undel von kovarianten 2-Tensoren},
sort=B\"undel von kovarianten 2-Tensoren}

\newglossaryentry{lok kov 2}
{name=\ensuremath{^{\varphi}T_{ij}},
description={lokale Darstellung eines kovarianten 2-Tensors T},
sort=lokale Darstellung eines kovarianten 2-Tensors T}

\newglossaryentry{mathcal{T}2M}
{name=\ensuremath{\mathcal{T}^2(M)},
description={Vektorraum aller glatten kovarianten 2-Tensorfelder in $T^2(M)$},
sort=Vektorraum aller glatten kovarianten 2-Tensorfelder in $T^2(M)$}

\newglossaryentry{g}
{name=\ensuremath{g},
description={Riemannsche Metrik},
sort=Riemannsche Metrik}

\newglossaryentry{(M,g)}
{name=\ensuremath{(M,g)},
description={Riemannsche Mannigfaltigkeit},
sort=Riemannsche Mannigfaltigkeit}

\newglossaryentry{Riem Norm}
{name=\ensuremath{|\cdot|_{g}},
description={Norm auf $T_p M$},
sort=Norm auf $T_p M$}

\newglossaryentry{Pullback}
{name=\ensuremath{F^{\ast}},
description={Pullback},
sort=Pullback}

\newglossaryentry{Sigma^2(T_p M)}
{name=\ensuremath{\Sigma^2(T_p M)},
description={Vektorraum aller symmetrischen kovarianten 2-Tensoren auf $T_p M$},
sort=Vektorraum aller symmetrischen kovarianten 2-Tensoren T_p M}

\newglossaryentry{TenProd}
{name=\ensuremath{\omega\otimes\eta},
description={Tensorprodukt},
sort=Tensorprodukt}

\newglossaryentry{symmTenProd}
{name=\ensuremath{\omega\eta},
description={symmetrisiertes Tensorprodukt},
sort=symmetrisiertes Tensorprodukt}

\newglossaryentry{g^{can}}
{name=\ensuremath{g^{can}},
description={Standard-Metrik auf $\mathbb{R}^n$},
sort=Standard-Metrik auf \mathbb{R}^n}

\newglossaryentry{mathbb{B}}
{name=\ensuremath{\mathbb{B}},
description={Einheitsball in $\mathbb{R}^n$},
sort=Einheitsball in \mathbb{R}^n}

\newglossaryentry{Raum2Karte}
{name=\ensuremath{^\varphi D^2_p(F)},
description={Raum aller Ableitungen bis zur zweiten Ordnung bzgl. Karte},
sort=Raum aller Ableitungen bis zur zweiten Ordnung bzgl. Karte}

\newglossaryentry{Raum2}
{name=\ensuremath{D^2_p(F)},
description={Raum aller Ableitungen bis zur zweiten Ordnung},
sort=Raum aller Ableitungen bis zur zweiten Ordnung}

\newglossaryentry{Sq}
{name=\ensuremath{\mathbb{S}^{q-1}},
description={Einheitssph\"are in $\mathbb{R}^q$},
sort=Einheitssph\"are R^q}

\newglossaryentry{TM}
{name=\ensuremath{TM},
description={Tangentialb\"undel von $M$},
sort=Tangentialb\"undel von M}

\newglossaryentry{UM}
{name=\ensuremath{UM},
description={Einheitstangentialb\"undel von $M$},
sort=Einheitstangentialb\"undel von M}

\newglossaryentry{ccone(X)}
{name=\ensuremath{ccone(X)},
description={konvex-konische H\"ulle von $X$},
sort=konvex-konische H\"ulle von $X$}

\newglossaryentry{conv(X)}
{name=\ensuremath{conv(X)},
description={konvexe H\"ulle von $X$},
sort=konvexe H\"ulle von $X$}

\newglossaryentry{Matrixnorm}
{name=\ensuremath{\left\Vert\cdot\right\Vert_{\mathbb{R}^{q\times q}}},
description={Matrixnorm, die von $|\cdot|_{\mathbb{R}^q}$ induziert wird},
sort=Matrixnorm}

\newglossaryentry{C^{m}(Omega,Rq)}
{name=\ensuremath{C^{m}(\Omega,\mathbb{R}^q)},
description={Raum aller $m$-fach stetig differenzierbaren vektorwertigen Funktionen auf einer offenen Menge $\Omega$},
sort=Raum aller $m$-fach stetig differenzierbaren Funktionen auf einer offenen Menge}

\newglossaryentry{C^{m}_0(Omega,Y)}
{name=\ensuremath{C^{m}_0(\Omega,Y)},
description={Raum aller Funktionen aus $C^{m}(\Omega,\mathbb{R}^q)$ deren Tr\"ager kompakt in $\Omega$ liegt},
sort=Raum aller Funktionen aus Cm}

\newglossaryentry{C^{m}(overline{Omega},Rq)}
{name=\ensuremath{C^{m}(\overline{\Omega},\mathbb{R}^q)},
description={Raum aller Funktionen aus $m$-fach stetig differenzierbaren Funktionen auf einer offenen, beschr\"ankten Menge $\Omega$, deren Ableitungen sich stetig bis zum Rand fortsetzen lassen},
sort=Raum aller Funktionen aus $m$-fach stetig differenzierbaren Funktionen auf einer offenen beschr\"ankten Menge beschrankten deren Ableitungen sich stetig bis zum Rand fortsetzen lassen}

\newglossaryentry{NormaufCm}
{name=\ensuremath{\left\Vert\cdot\right\Vert_{C^{m}(\overline{\Omega},\mathbb{R}^q)}},
description={Norm auf $C^{m}(\overline{\Omega},\mathbb{R}^q)$},
sort=Norm auf Cm}

\newglossaryentry{Holderkonstante}
{name=\ensuremath{[u]_{\alpha,S}},
description={H\"olderkonstante},
sort=H\"olderkonstante}

\newglossaryentry{C0alpha}
{name=\ensuremath{C^{0,\alpha}(\overline{\mathbb{B}})},
description={Raum aller $\alpha$-h\"olderstetigen reellen Funktionen auf $\mathbb{B}$},
sort=Raum aller alpha-h\"olderstetigen Funktionen auf B}

\newglossaryentry{NormaufC0alpha}
{name=\ensuremath{\left\Vert\cdot\right\Vert_{C^{0,\alpha}(\overline{\mathbb{B}})}},
description={Norm auf $C^{0,\alpha}(\overline{\mathbb{B}})$},
sort=Norm auf C0alpha}

\newglossaryentry{CmalphaB}
{name=\ensuremath{C^{m,\alpha}(\overline{\mathbb{B}})},
description={Raum aller $m$-fach stetig differenzierbaren reellen Funktionen auf $\mathbb{B}$ deren Ableitungen $\alpha$-h\"olderstetig sind},
sort=Raum aller m-fach stetig differenzierbaren reellenFunktionen deren Ableitungen alpha-h\"olderstetig sind}

\newglossaryentry{Norm1}
{name=\ensuremath{\left|u \right|_{C^{m,\alpha}(\overline{\mathbb{B}})}},
description={Norm auf $C^{m,\alpha}(\overline{\mathbb{B}})$},
sort=Norm auf C_0malphaB}

\newglossaryentry{Norm2}
{name=\ensuremath{\left\Vert u \right\Vert_{C^{m,\alpha}(\overline{\mathbb{B}})}},
description={Norm auf $C^{m,\alpha}(\overline{\mathbb{B}})$},
sort=Norm auf C_0malphaB}

\newglossaryentry{CmalphaBRq}
{name=\ensuremath{C^{m,\alpha}(\overline{\mathbb{B}},\mathbb{R}^q)},
description={Raum aller $m$-fach stetig differenzierbaren Funktionen auf $\mathbb{B}$ deren Ableitungen $\alpha$-h\"olderstetig sind},
sort=Raum aller m-fach stetig differenzierbaren Funktionen deren Ableitungen alpha-h\"olderstetig sind}

\newglossaryentry{NormaufCmalphaBRq}
{name=\ensuremath{\left|u \right|_{C^{m,\alpha}(\overline{\mathbb{B}},\mathbb{R}^q)}},
description={Norm auf $C^{m,\alpha}(\overline{\mathbb{B}},\mathbb{R}^q)$},
sort=Norm auf CmalphaBRq}

\newglossaryentry{SobolevRaum}
{name=\ensuremath{W^{1,p}_0(\mathbb{B})},
description={Sobolevraum},
sort=Sobolevraum}

\newglossaryentry{LRaum}
{name=\ensuremath{L^p},
description={Lebesgue-Raum},
sort=Lebesgue-Raum}

\newglossaryentry{mathbb{S}}
{name=\ensuremath{\mathbb{S}},
description={$[-\pi,\pi]\subseteq\mathbb{R}$},
sort=S}

\newglossaryentry{I}
{name=\ensuremath{I},
description={$[0,1]\subseteq\mathbb{R}$},
sort=I}

\newglossaryentry{Lambda}
{name=\ensuremath{\Lambda(\epsilon^k)},
description={$\epsilon^k\, h$ mit $h\in C^{\infty}(I\times\mathbb{S}\times\overline{\mathbb{B}})$ und $supp(h)\subseteq I\times\mathbb{S}\times\mathbb{B}$},
sort=Lambda}

\newglossaryentry{durchmesser}
{name=\ensuremath{diam(S)},
description={Durchmesser von $S$},
sort=Durchmesser}

\newglossaryentry{RaumGlatt}
{name=\ensuremath{C^{\infty}(\Omega,\mathbb{R}^q)},
description={Raum aller glatten Abbildungen auf einer offenen Menge $\Omega$},
sort=Raum aller glatten Abbildungen auf einer offenen Menge \Omega}

\newglossaryentry{RaumGlattFort}
{name=\ensuremath{C^{\infty}(\overline{\Omega},\mathbb{R}^q)},
description={Raum aller glatten Abbildungen auf einer offenen Menge $\Omega$, deren Ableitungen stetig bis zum Rand fortsetzbar sind},
sort=Raum aller glatten Abbildungen auf einer offenen Menge \Omega deren Ableitungen stetig bis zum Rand fortsetzbar sind}

\newglossaryentry{sobo}
{name=\ensuremath{W^{m,p}(\Omega)},
description={Sobolev-Raum},
sort=Sobolev-Raum}

\newglossaryentry{sobonull}
{name=\ensuremath{W^{m,p}_0(\Omega)},
description={Sobolev-Raum mit Nullrandbedingung},
sort=Sobolev-Raum mit Nullrandbedingung}

\newglossaryentry{sobonorm}
{name=\ensuremath{\left\Vert\cdot \right\Vert_{W^{m,p}(\Omega)}},
description={Norm auf $W^{m,p}(\Omega)$},
sort=Norm auf Wmp}

\newglossaryentry{differ}
{name=\ensuremath{\partial_{t,h}u},
description={Differenzenquotient},
sort=Differenzenquotient}

\chapter{Einleitung}
\thispagestyle{fancy}

\section{Darstellung des Problems}

In der vorliegenden Diplomarbeit wird gezeigt, dass f\"ur eine gegebene $n$-dimensionale geschlossene Riemannsche Mannigfaltigkeit $(M,g)$ eine Raumdimension $q(n)\in\mathbb{N}$ existiert, so dass sich $(M,g)$ isometrisch in den Raum $\mathbb{R}^{q(n)}$ einbetten l\"asst. Das bedeutet, es existiert eine glatte Einbettung $F\in C^{\infty}(M,\mathbb{R}^q)$, so dass $F^{\ast}(g^{can})=g$ gilt. Diese Fragestellung wurde erstmalig von J. Nash in der Arbeit \cite{nash1956imbedding} beantwortet. J. Nash hat gezeigt, dass f\"ur jede kompakte Mannigfaltigkeit, mit einer gegebenen $C^k$-Metrik f\"ur $k\geq 3$ eine isometrische Einbettung in den Raum $\mathbb{R}^{q(n)}$ f\"ur $q(n)=\frac{n}{2}(3n+11)$ existiert, wobei f\"ur die Regularit\"at der Einbettung mindestens $C^k$ gew\"ahrleistet wird. Der von J. Nash gef\"uhrte Beweis gilt als technisch sehr anspruchsvoll. Sp\"ater zeigten M. L. Gromov und V.A. Rokhlin in der Arbeit \cite{gromov1970embeddings}, dass sich die Wahl von $q(n)$ auf $\frac{1}{2}(n^2+9n+10)$ reduzieren l\"asst, sofern die Metrik glatt ist, wobei sich die Regularit\"at der Metrik auf die Einbettung \"ubertr\"agt. In den sp\"aten 80er-Jahren des 20. Jahrhunderts fand Matthias G\"unther einen neuen Weg, den Einbettungssatz zu beweisen, hierf\"ur sei auf die Arbeiten \cite{gunther1989perturbation}, \cite{gunther1989einbettungssatz} und \cite{gunther1991isometric} verwiesen, welche die Grundlage der vorliegenden Diplomarbeit bilden. Es wird gezeigt, dass sich jede $n$-dimensionale geschlossene Riemannsche Mannigfaltigkeit $(M,g)$ isometrisch in den Raum $\mathbb{R}^{q(n)}$ f\"ur $q(n):=\max\left\{\frac{n}{2}(n+5),\frac{n}{2}(n+3)+5 \right\}$ einbetten l\"asst. Hierbei wird sowohl von der Mannigfaltigkeit $M$, als auch von der Metrik $g$ Glattheit vorausgesetzt, welche sich auf die Einbettung \"ubertragen wird.\\
Der Vollst\"andigkeit halber sei erw\"ahnt, dass die zugrunde liegende Fragestellung, im lokalen Sinne, bereits im sp\"ateren 19. Jahrhundert von L. Schlaefi in \cite{schlaefli1871nota}, und in den 20er-Jahren des 20. Jahrhunderts von M. Janet \cite{janet1926} und \'E. Cartan \cite{cartan1927} diskutiert worden sind. F\"ur weitere historische Hintergrundinformationen sei auf \cite{hong2006isometric} und \cite{andrews2002notes} verwiesen. \\
Im Anschluss an die Darstellung der Beweisf\"uhrung von Matthias G\"unther wird eine M\"oglichkeit beschrieben, ein glatte Familie von Riemannschen
Mannigfaltigkeiten, wie zum Beispiel den in der Arbeit \cite{hamilton1982three} beschriebenen Ricci-Fluss, f\"ur eine kurze Zeit isometrisch in den Vektorraum $\mathbb{R}^{q(n)}$ einzubetten. 

\section{Beschreibung der Vorgehensweise}

In der gesamten Arbeit wird ausschlie\ss{}lich eine \textbf{geschlossene Riemannsche Mannigfaltigkeit}\index{Mannigfaltigkeit>Riemannsche $\sim$>geschlossene $\sim$} $(M,g)$ betrachtet. Das ist eine kompakte glatte Mannigfaltigkeit $M$ ohne Rand, gem\"a\ss{} Definition \ref{defi:B6}, im Unterschied zu \cite[Chapter 1, Manifolds with Boundary]{lee2003introduction}, zusammen mit einer Riemannschen Metrik $g\in\mathcal{T}^2(M)$. Aus dem folgenden Satz wird sich die Existenz einer isometrischen Einbettung ergeben, der hierbei verwendete Begriff der freien Einbettung wird in \autoref{sec:Ab2.1} exakt definiert:

\begin{haupt}\label{haupt1}
Gegeben sei eine $n$-dimensionale geschlossene Riemannsche Mannigfaltigkeit $(M,g)$, und eine freie Einbettung $F_0\in C^{\infty}(M,\mathbb{R}^q)$, mit $q\geq \frac{n}{2}(n+3)+5$, so dass $g-F_0^{\ast}(g^{can})\in\mathcal{T}^2(M)$ eine Riemannsche Metrik ist. Ferner sei ein $\delta\in\mathbb{R}_{>0}$ gegeben, dann existiert eine freie Einbettung $F\in C^{\infty}(M,\mathbb{R}^q)$, so dass $F^{\ast}(g^{can})=g$ und:
\begin{equation*}
\max_{x\in M}{|F(x)-F_0(x)|_{\mathbb{R}^q}}\leq \delta
\end{equation*}

gilt.

\end{haupt}

F\"ur sp\"atere Zwecke sei erw\"ahnt, dass mit \eqref{eq:B6} die Gleichung $F^{\ast}(g^{can})=g$ f\"ur eine Karte $(U,\varphi)\in\mathcal{A}$  \"aquivalent zu:
\begin{equation}\label{eq:1.12}
\partial_i\, ^{\varphi}F(\varphi(x))\cdot \partial_j\, ^{\varphi}F(\varphi(x))=\, ^{\varphi}g_{ij}(\varphi(x))
\end{equation}

f\"ur $x\in U$ und $1\leq i\leq j\leq n$ ist. Die Existenz einer freien Einbettung $F_0\in C^{\infty}(M,\mathbb{R}^q)$ mit der Eigenschaft, dass $g-F_0^{\ast}(g^{can})$ eine Riemannsche Metrik auf $M$ ist, wird Gegenstand von \autoref{chap:Kap2} sein. In \autoref{chap:Kap3} beginnt der Beweis von Hauptsatz \ref{haupt1}, indem zun\"achst gezeigt wird, dass endlich viele symmetrische Tensorfelder $h^{(1)},...,h^{(m)}\in \mathcal{T}^2(M)$ existieren, so dass $h:= g-F_0^{\ast}(g^{can})=\sum_{l=1}^m{h^{(l)}}$  gilt. Die hierbei konstruierten Tensorfelder $h^{(1)},...,h^{(m)}\in \mathcal{T}^2(M)$ haben die Eigenschaft, dass sie jeweils kompakt in einer Koordinatenumgebung liegen, bez\"uglich der sie eine gewisse besondere Koordinatendarstellung haben. Ist diese Zerlegungseigenschaft der Metrik $h$ bewiesen, dann reicht es zu zeigen, dass f\"ur eine gegebene freie Abbildung $F_l\in C^{\infty}(M,\mathbb{R}^q)$, f\"ur $l\in\{0,...,m-1\}$, eine freie Abbildung $F_{l+1}\in C^{\infty}(M,\mathbb{R}^q)$ existiert, so dass die Gleichungen:
\begin{align*}
&F_{l+1}^{\ast}(g^{can})=F_l^{\ast}(g^{can})+h^{(l)}\\
&F_{l+1}^{\ast}(g^{can})(x)=F_{l}^{\ast}(g^{can})(x)\hspace{0.5cm}\text{f\"ur alle }x\in M\backslash U_l
\end{align*}

und die Absch\"atzung:
\begin{align*}
&\max_{x\in M}{|F_{l+1}(x)-F_{l}(x)|_{\mathbb{R}^q}}\leq\delta
\end{align*}

wobei $\left\{(U_l,\varphi_l)\right\}_{1\leq l\leq m}\subseteq\mathcal{A}$, die in \autoref{chap:Kap3} konstruierten Karten beschreiben, erf\"ullt sind. Mithilfe dieses Resultates wird das Problem in eine lokale Fragestellung \"uberf\"uhrt. Da f\"ur die Karten $\varphi_l(U_l)=B_{1+\tau}(0)$, f\"ur $l\in\{1,...,m\}$ gilt, werden in den darauffolgenden Betrachtungen lokale Resultate auf $\gls{mathbb{B}}:= B_1(0)\subseteq\mathbb{R}^n$ hergeleitet. In \autoref{chap:Kap4} wird ein $l\in \{1,...,m\}$ fixiert, und es werden aus einer gegebenen freien Abbildung $F_0\in C^{\infty}(\overline{\mathbb{B}},\mathbb{R}^q)$ f\"ur $k\in\mathbb{N}\backslash\{0\}$, freie Abbildungen $F_{\epsilon,k}\in(\overline{\mathbb{B}},\mathbb{R}^q)$, mit $\epsilon\in(0,\epsilon_k]$ konstruiert, so dass f\"ur $k\geq 2$ die Eigenschaften:
\begin{align*}
&F_{\epsilon,k}^{\ast}(g^{can})=F_0^{\ast}(g^{can})+h^{(l)}+\epsilon^{k+1} f_{\epsilon,k}\\
&F_{\epsilon,k}(x)-F_0(x)\in C^{\infty}_0(\mathbb{B},\mathbb{R}^q)\\
&\max_{x\in \mathbb{B}}{|F_{\epsilon,k}(x)-F_{0}(x)|_{\mathbb{R}^q}}\leq C(k)\cdot\epsilon
\end{align*}

f\"ur eine, von $\epsilon$ unabh\"angige Konstante $C(k)\in\mathbb{R}_{>0}$, auf $\mathbb{B}$ erf\"ullt sind. Hierbei ist $f_{\epsilon,k}\in C^{\infty}_0(\mathbb{B})$ eine St\"orung, die in einer kompakten, von $\epsilon$  und $k$ unabh\"angigen, Menge $K\subseteq \mathbb{B}$ liegt, welche auch den Tr\"ager von $h^{(l)}$ beinhaltet. Diese St\"orung hat besondere Eigenschaften. In \autoref{chap:Kap6} werden feste offene Mengen $U_1,U_2\subseteq\mathbb{B}$ mit $K\subseteq U_1$, $\overline{U}_1\subseteq U_2$ sowie $\overline{U}_2\subseteq \mathbb{B}$ gew\"ahlt, und gezeigt, dass sich f\"ur eine hinreichend gro\ss{}e Wahl von $k\in\mathbb{N}\backslash\{0,1\}$, und eine hinreichend kleine Wahl von $\epsilon\in (0,\epsilon_k]$, die Funktion $F_{\epsilon,k}\in(\overline{\mathbb{B}},\mathbb{R}^q)$, durch Perturbation, in eine freie Abbildung $F\in C^{\infty}(\mathbb{B},\mathbb{R}^q)$ \"uberf\"uhren l\"asst, so dass der St\"orterm $f_{\epsilon,k}$ verschwindet. Dabei wird gew\"ahrleistet, dass $F_{\epsilon,k}(x)=F(x)$, f\"ur alle $x\in \mathbb{B}\backslash U_2$, gilt. Das hierf\"ur verwendete Resultat, zur L\"osung dieses lokalen Perturbationsproblems, wird in \autoref{chap:Kap5} hergeleitet. Damit ist Hauptsatz \ref{haupt1} gezeigt. Mit \autoref{chap:Kap2} folgt dann:
\begin{haupt}\label{haupt2}
Es sei $(M,g)$ eine geschlossene Riemannsche Mannigfaltigkeit der Dimension $n$, dann existiert f\"ur $q(n):=\max\{\frac{n}{2}(n+5),\frac{n}{2}(n+3)+5 \}$ eine freie Einbettung $F\in C^{\infty}(M,\mathbb{R}^q)$, so dass $F^{\ast}(g^{can})=g$ gilt.
\end{haupt}

In \autoref{chap:Kap7} wird, aufbauend auf Hauptsatz \ref{haupt2}, eine selbst entwickelte Methode beschrieben, die es erm\"oglicht, eine 1-parametrige Familie von Riemannschen Metriken in einen euklidischen Vektorraum isometrisch einzubetten. Von Interesse ist hierbei die glatte Abh\"angigkeit der Einbettungen vom Parameter.

\chapter{Existenzsätze für freie Einbettungen}
\label{chap:Kap2}
\thispagestyle{fancy}

In diesem Kapitel wird der Begriff der freien Abbildung eingef\"uhrt, aus dem sich der Begriff der freien Einbettung ergibt. F\"ur diese Abbildungen werden wesentliche Eigenschaften zusammengetragen, und anschlie\ss{}end Existenzs\"atze herausgearbeitet. Die Vorgehensweise wurde hierbei \cite[2.3.]{andrews2002notes} und \cite[2.5.3.]{gromov1970embeddings} entnommen.

\section{Begriffsklärung}
\label{sec:Ab2.1}

\begin{defi}
Sei $M$ eine glatte Mannigfaltigkeit der Dimension $n$, $F\in C^{\infty}(M,\mathbb{R}^q)$, $p\in M$, und $(U,\varphi)\in\mathcal{A}$ eine Karte, mit $p\in U$. Dann wird der Vektorraum:
\begin{equation*}
\gls{Raum2Karte}:= lin\left(\left\{\partial_i\,  ^\varphi F(\varphi(p))\right\}_{1\leq i\leq n}\cup \left\{\partial_i \partial_j\,  ^\varphi F(\varphi(p))\right\}_{1\leq i\leq j\leq n}\right)\subseteq \mathbb{R}^q
\end{equation*}

als \textbf{Raum aller Ableitungen von} $\bm{F}$ \textbf {bis zur zweiten Ordnung in} $\bm{p}$ \textbf{bez\"uglich} $\bm{(U,\varphi)}$\index{Raum aller Ableitungen bis zur zweiten Ordnung>$\sim$ bzgl. Karte} bezeichnet.

\end{defi}

Es wird gezeigt, dass dieser Begriff unabh\"angig von der Wahl der Karte $(U,\varphi)$ ist, hierbei wird auch auf \cite[1.2.3.]{gromov1970embeddings} verwiesen:

\begin{lem}
Sei $M$ eine glatte Mannigfaltigkeit der Dimension $n$,  $F\in C^{\infty}(M,\mathbb{R}^q)$, $p\in M$, und $(U,\varphi)\in\mathcal{A}$ , sowie $(V,\psi)\in\mathcal{A}$ Karten, mit $p\in U\cap V$. Dann gilt:
\begin{equation*}
^\varphi D^2_p(F)=\, ^\psi D^2_p(F)
\end{equation*}

\end{lem}

\begin{bew}
Es gilt f\"ur $j\in\{1,...,n\}$:
\begin{align*}
&\partial_j\,  ^\varphi F(\varphi(p))=\partial_j (F\circ  \varphi^{-1})(\varphi(p))=\partial_j (F\circ\psi^{-1}\circ \psi\circ \varphi^{-1})(\varphi(p))\\
=&\sum_{k=1}^{n}{\partial_k (F\circ\psi^{-1})(\psi(p))\cdot\partial_j(\psi^k\circ \varphi^{-1})(\varphi(p))}=\sum_{k=1}^{n}{\partial_k\,  ^\psi F(\psi(p))\cdot\partial_j(\psi^k\circ \varphi^{-1})(\varphi(p))}
\end{align*}

und f\"ur $i\in\{1,...,n\}$:
\begin{align*}
\partial_i\partial_j\,  ^\varphi F(\varphi(p))=&\sum_{k,l=1}^{n}{\partial_k\partial_l (F\circ\psi^{-1})(\psi(p))\cdot\partial_j(\psi^k\circ \varphi^{-1})(\varphi(p))\cdot \partial_i(\psi^l\circ \varphi^{-1})(\varphi(p))} \\
&+\sum_{k=1}^{n}{\partial_k (F\circ\psi^{-1})(\psi(p))\cdot\partial_i\partial_j(\psi^k\circ \varphi^{-1})(\varphi(p))}\\
=&\sum_{k,l=1}^{n}{\partial_k\partial_l\, ^\psi F(\psi(p))\cdot\partial_j(\psi^k\circ \varphi^{-1})(\varphi(p))\cdot\partial_i(\psi^l\circ \varphi^{-1})(\varphi(p))}\\
&+\sum_{k=1}^{n}{\partial_k\, ^\psi F(\psi(p))\cdot\partial_i\partial_j(\psi^k\circ \varphi^{-1})(\varphi(p))}
\end{align*}

Also gilt $^\varphi D^2_p(F)\subseteq\, ^\psi D^2_p(F)$. Vertauschung von $\varphi$ und $\psi$ ergibt $^\psi D^2_p(F) \subseteq\, ^\varphi D^2_p(F)$, womit die Behauptung gezeigt ist.

\end{bew}

Dieses Lemma rechtfertigt die folgende Festlegung:

\begin{defi}
Sei $M$ eine glatte Mannigfaltigkeit, $F\in C^{\infty}(M,\mathbb{R}^q)$, $p\in M$, und $(U,\varphi)\in\mathcal{A}$, eine Karte mit $p\in U$. Dann wird der Vektorraum:
\begin{equation*}
\gls{Raum2}:= \, ^\varphi D^2_p(F)
\end{equation*}

als \textbf{Raum aller Ableitungen von} $\bm{F}$ \textbf{bis zur zweiten Ordnung in} $\bm{p}$\index{Raum aller Ableitungen bis zur zweiten Ordnung} bezeichnet.

\end{defi}

Es gilt:
\begin{align*}
dim\, ^\varphi D^2_p(F)&\leq dim\left[ lin\left(\left\{\partial_i\,  ^\varphi F(\varphi(p))\right\}_{1\leq i\leq n}\right)\right]+dim\left[ lin \left(\left\{\partial_i \partial_j\,  ^\varphi F(\varphi(p))\right\}_{1\leq i\leq j\leq n}\right)\right]\\
&\leq n+\frac{n}{2}(n+1)=\frac{2n}{2}+\frac{n}{2}(n+1)=\frac{n}{2}(n+3)
\end{align*}

Daraus ergibt sich die Absch\"atzung:
\begin{equation}\label{eq:2.1}
0\leq dim\, D^2_p(F)\leq \min\left\{q, \frac{n}{2}(n+3)\right\}
\end{equation}

Eine besondere Bedeutung haben Abbildungen, bei denen, f\"ur $q\geq \frac{n}{2}(n+3)$, die Absch\"at\-zung \eqref{eq:2.1} mit Gleichheit angenommen wird:

\begin{defi}
Eine Abbildung $F\in C^{\infty}(M,\mathbb{R}^q)$, auf einer glatten Mannigfaltigkeit $M$, hei\ss{}t \textbf{freie Abbildung}\index{Abbildung>freie $\sim$}, falls:
\begin{equation*}
dim\, D^2_p(F)=\frac{n}{2}(n+3)
\end{equation*}

f\"ur alle $p\in M$ gilt. Falls $F$ eine glatte Einbettung ist, dann wird $F$ als \textbf{freie Einbettung}\index{Einbettung>freie $\sim$} bezeichnet.

\end{defi}


\section{Grundlegende Eigenschaften von freien Abbildungen}

Mit \eqref{eq:B6} ist jede freie Abbildung eine Immersion. Ist $M$ kompakt, dann ist mit \cite[Proposition 7.4.]{lee2003introduction}, jede injektive freie Abbildung bereits eine freie Einbettung. Weiterhin gilt:

\begin{lem}\label{lem:2.2}
Gegeben seien $n$-dimensionale glatte Mannigfaltigkeiten $(M,\mathcal{O}_M,\mathcal{A}_M)$ und $(N,\mathcal{O}_N,\mathcal{A}_N)$, sowie eine Immersion $G\in C^{\infty}(M,N)$, und eine freie Abbildung $F\in C^{\infty}(N,\mathbb{R}^q)$, dann ist die Verkn\"upfung $F\circ G$ ebenfalls eine freie Abbildung.
\end{lem}

\begin{bew}
Sei $p\in M$, mit \cite[Theorem 7.13]{lee2003introduction} existieren Karten $(U,\varphi)\in\mathcal{A}_M$ und $(V,\psi)\in\mathcal{A}_N$, mit $p\in U$ und $F(U)\subseteq V$, so dass f\"ur alle $i,j\in\{1,...,n\}$ und $x\in U$ die Gleichung:
\begin{equation}\label{eq:2.2}
\partial_j\, ^{\varphi}_\psi G^i(\varphi(x))=\delta^i_j
\end{equation}

gilt, dann ist f\"ur $j\in\{1,...,n\}$:
\begin{align}\label{eq:2.3}
\begin{split}
&\partial_j\, ^{\varphi}(F\circ G)(\varphi(p))=\partial_j\, (F\circ\psi^{-1} \circ \psi \circ G\circ \varphi^{-1})(\varphi(p))\\
=&\sum_{k=1}^n{\partial_k\, ^{\psi}F(\psi(G(p)))\cdot \partial_j \, ^{\varphi}_\psi G^k(\varphi(p))}\stackrel{\eqref{eq:2.2}}{=}\sum_{k=1}^n{\partial_k\, ^{\psi}F(\psi(G(p)))\cdot \delta_j^k}=\partial_j\, ^{\psi}F(\psi(G(p)))
\end{split}
\end{align}

und f\"ur $i\in\{1,...,n\}$, ausgehend von der Rechnung \eqref{eq:2.3}:
\begin{align*}
\partial_i\partial_j\, ^{\varphi}(F\circ G)(\varphi(p))=\partial_i \partial_j\, ^{\psi}F(\psi(G(p)))
\end{align*}

Damit gilt: $dim (D^2_p(F\circ G))=dim(D^2_{G(p)}(F))=\frac{n}{2}(n+3)$, womit gezeigt ist, dass $F\circ G$ eine freie Abbildung ist.

\end{bew}

In \autoref{sec:Ab2.5} wird das folgende Lemma von Bedeutung sein, welches direkt aus der Konstruktion des Atlas, einer eingebetteten Untermannigfaltigkeit folgt, siehe hierf\"ur auch Definition \ref{defi:B8}.

\begin{lem}\label{lem:2.5}

Gegeben sei eine eingebettete Untermannigfaltigkeit $S$, einer glatten Mannigfaltigkeit $M$, und eine freie Abbildung $F\in C^{\infty}(M,\mathbb{R}^q)$. Dann ist die Abbildung: $\left.F\right|_{S}: S\longrightarrow \mathbb{R}^q$ ebenfalls eine freie Abbildung.

\end{lem}


\section[Freie Einbettung einer offenen Menge \texorpdfstring{$\Omega\subseteq\mathbb{R}^n$}{Omega kleiner \mathbb{R}n}]{Freie Einbettung einer offenen Menge \texorpdfstring{$\bm{\Omega\subseteq\mathbb{R}^n}$}{Omega kleiner \mathbb{R}n}}

Eine freie Einbettung auf einer offenen Menge $\Omega\subseteq \mathbb{R}^n$ kann direkt angegeben werden, wie der Beweis des folgenden Satzes zeigt:

\begin{satz}\label{satz:2.2}
Sei $\Omega\subseteq\mathbb{R}^n$ eine offene Menge, dann existiert eine freie Einbettung $F_0\in C^{\infty}(\Omega,\mathbb{R}^{\frac{n}{2}(n+3)})$.
\end{satz}

\begin{bew}
Es sei $\{e_i\}_{1\leq i\leq n}\cup \{e_{ij}\}_{1\leq i\leq j\leq n}$ eine Orthonormalbasis des Raumes $\mathbb{R}^{\frac{n}{2}(n+3)}$. Damit wird die folgende Abbildung definiert:
\begin{align*}
F_0: \Omega&\longrightarrow\mathbb{R}^{\frac{n}{2}(n+3)}\\
(x^1,...,x^n)&\mapsto \sum_{i=1}^n{x^i e_i}+  \sum_{1\leq i\leq j\leq n}{x^i x^j e_{ij}}
\end{align*}

Zun\"achst wird gezeigt, dass die Abbildung $F_0$ injektiv ist. Seien dazu $x,y\in\Omega$ mit $F(x)=F(y)$, dann gilt insbesondere:
\begin{equation*}
\sum_{i=1}^n{x^i e_i}=\sum_{i=1}^n{y^i e_i}
\end{equation*}

also  sind $x$ und $y$ identisch. F\"ur $l\in\{1,...,n\}$ gilt:
\begin{align}\label{eq:2.4}
\begin{split}
\partial_l F_0(x)=&\sum_{i=1}^n{\delta_l^i e_i}+ \sum_{1\leq i\leq j\leq n}{\partial_l(x^i x^j) e_{ij}}=e_l+ \sum_{1\leq i\leq j\leq n}{\delta_l^i x^j e_{ij}}+ \sum_{1\leq i\leq j\leq n}{x^i \delta_l^j e_{ij}}\\
=& e_l+ \sum_{l\leq j\leq n}{ x^j e_{lj}}+ \sum_{1\leq i\leq l}{x^i  e_{il}}
\end{split}
\end{align}

F\"ur ein weiteres $k\in\{1,...,l\}$ ist dann:
\begin{equation*}
\partial_k\partial_l F_0(x)=
\begin{cases}
e_{kl} &\text{falls }k< l\\
2e_{kk} &\text{falls }k=l
\end{cases}
\end{equation*}

womit gezeigt ist, dass $F_0$ eine freie Abbildung ist. Dar\"uber hinaus ist die inverse Abbildung $F_0^{-1}: F_0(\Omega)\longrightarrow \Omega$ stetig, denn sind f\"ur ein gegebenes $\epsilon \in\mathbb{R}$, $x,y\in\Omega$ mit $|F_0(x)-F_0(y)|_{\mathbb{R}^{\frac{n}{2}(n+3)}}<\epsilon$, so folgt daraus direkt:
\begin{align*}
&|x-y|_{\mathbb{R}^{n}}=\left|\sum_{i=1}^n{x^i e_i}-\sum_{i=1}^n{y^i e_i}\right|_{\mathbb{R}^{\frac{n}{2}(n+3)}}=\left|\sum_{i=1}^n{(x^i-y^i) e_i}\right|_{\mathbb{R}^{\frac{n}{2}(n+3)}}\\
\leq & \left|\sum_{i=1}^n{(x^i-y^i) e_i}+ \sum_{1\leq i\leq j\leq n}{(x^i x^j-y^i y^j) e_{ij}}\right|_{\mathbb{R}^{\frac{n}{2}(n+3)}}=|F_0(x)-F_0(y)|_{\mathbb{R}^{\frac{n}{2}(n+3)}}<\epsilon
\end{align*}

\end{bew}


\section{Einbettungssatz f\"ur eine kompakte Mannigfaltigkeit}

Um einen Existenzsatz f\"ur eine freie Abbildung auf einer Mannigfaltigkeit zu zeigen, erweist es sich als n\"utzlich, die Mannigfaltigkeit zun\"achst in einen Vektorraum $\mathbb{R}^q$ einzubetten, wobei man vor\"ubergehend keine Kontrolle \"uber die Raumdimension $q\in\mathbb{N}$ haben m\"ochte. Dazu dient das folgende Lemma, vergleiche hierf\"ur auch \cite[Chapter 1, 3.4. Theorem]{hirsch1997differential}:

\begin{lem}\label{lem:2.1}
Gegeben sei eine kompakte $n$-dimensionale glatte Mannigfaltigkeit $M$, dann existiert f\"ur ein $q(n,M)\in \mathbb{N}$ eine glatte Einbettung $F\in C^{\infty}(M,\mathbb{R}^q)$.
\end{lem}

\begin{bew}
Mit \cite[Proposition 2.24]{lee2003introduction} und der Kompaktheit von $M$, existieren endlich viele Karten $(W_1,\varphi_1),...,(W_r,\varphi_m)\in\mathcal{A}$, mit $\varphi(W_i)=B_3(0)$ f\"ur alle $i\in\{1,...m\}$, so dass die Mengen $U_i:= \varphi_i^{-1}(\mathbb{B})$, f\"ur $i\in\{1,...,m\}$, in ihrer Gesamtheit die Mannigfaltigkeit $M$ \"uberdecken. Weiterhin wird mit \cite[Lemma 2.22]{lee2003introduction} eine Funktion $H\in C^{\infty}(\mathbb{R}^n,\mathbb{R})$ gew\"ahlt, welche die Eigenschaften $\left. H\right|_{\overline{\mathbb{B}}}\equiv 1$, $supp(H)=\overline{B_2(0)} $ und $H(\mathbb{R}^n)=[0,1]$ erf\"ullt. In der folgenden Definition werden f\"ur alle $i\in\{1,...,m\}$ die Funktionen $H\circ \varphi_i$ und $\varphi_i^T\ H\circ \varphi_i$ au\ss{}erhalb von $W_i$ durch $0$ konstant fortgesetzt, eine neue Bezeichnung soll nicht verwendet werden. Aus diesen Funktionen wird eine Funktion $F\in C^{\infty}(M,\mathbb{R}^{m(n+1)})$ wie folgt konstruiert:
\begin{equation*}
F:=\left((H\circ \varphi_i)_{1\leq i\leq m}, (\varphi_i^T\ H\circ \varphi_i)_{1\leq i\leq m} \right)^T
\end{equation*}

Um zu zeigen, dass $F$ eine glatte Einbettung ist, gen\"ugt es, aufgrund der Kompaktheit von $M$, mit \cite[Proposition 7.4]{lee2003introduction} zu zeigen, dass $F$ eine injektive Immersion ist. Zun\"achst die Injektivit\"at:
Sind $x,y\in M$ mit $F(x)=F(y)$, so existiert ein $i\in\{1,...,m\}$, mit $x\in U_i$. Dann gilt, wegen $F(x)=F(y)$, die Gleichung $(H\circ \varphi_i)(x)=(H\circ \varphi_i)(y)=1$ und es muss die Aussage $y\in U_i$ gelten. Aufgrund der Definition von $F$ gilt dann die Gleichung $\varphi_i(x)=\varphi_i(y)$. Als Kartenabbildung ist $\varphi_i$ insbesondere injektiv, und es folgt $x=y$. $F$ ist auch eine Immersion, denn f\"ur ein festes $x\in M$ existiert ein $i\in\{1,...,m\}$ mit $x\in U_i$. Die Jacobimatrix der Komponente $(\varphi_i^T\ H\circ \varphi_i)$ von $F$, bez\"uglich der Kartenabbildung $\varphi_i$, ist im Punkt $x$, nach Wahl von $H$, die $n$-dimensionale Einheitsmatrix. Die Ableitung von $F$ besitzt damit in $x$, unter Beachtung von \eqref{eq:B6}, den maximalen Rang $n$. Da $x\in M$ beliebig war, ist $F$ eine Immersion.

\end{bew}

Nun wird gezeigt, dass sich die in Lemma \ref{lem:2.1} hergeleitete Raumdimension $q$ immer, unabh\"angig von der konkreten Gestalt von $M$, auf $2n+1$ reduzieren l\"asst. Nebenbei ergibt sich noch ein Existenzsatz f\"ur Immersionen. Dieser Satz wird auch als \glqq Einbettungssatz von Whitney\grqq \,bezeichnet, vergleiche hierf\"ur \cite[Chapter 1, 3.5. Theorem]{hirsch1997differential}, oder mit einer etwas anderen Beweisidee \cite[Theorem 10.11]{lee2003introduction}. 

\begin{satz}\label{satz:2.1}
F\"ur jede kompakte n-dimensionale glatte Mannigfaltigkeit $M$ existiert eine glatte Einbettung in den Raum $\mathbb{R}^{2n+1}$, und eine Immersion in den Raum $\mathbb{R}^{2n}$.
\end{satz}

F\"ur den Beweis dieses Satzes werden einige ma\ss{}theoretische Begriffe und Zusammen\-h\"ange eingef\"uhrt:

\begin{defi}
Sei $M$ eine n-dimensionale glatte Mannigfaltigkeit, dann hei\ss{}t eine Menge $A\subseteq M$ \textbf{Nullmenge in} $\bm{M}$\index{Menge>Null-$\sim$ in $M$}, falls f\"ur alle $(U,\varphi)\in\mathcal{A}$ die Menge $\varphi(A\cap U)$ eine Lebesgue-Nullmenge in $\mathbb{R}^n$ ist.
\end{defi}

Die folgende wichtige Tatsache wird in \cite[Theorem 10.5]{lee2003introduction} gezeigt:

\begin{lem}\label{lem:2.3}
Seien $M$ und $N$ glatte Mannigfaltigkeiten mit $dim(M) < dim(N)$ und sei $F\in C^{\infty}(M,N)$. Dann ist die Menge $F(M)$ eine Nullmenge in N. Insbesondere ist $N\backslash F(M)$ eine dichte Teilmenge von $N$.
\end{lem}

Des Weiteren wird eine alternative Charakterisierung von eingebetteten Untermannigfaltigkeiten verwendet. 

\begin{lem}\label{lem:2.4}
Seien $M$ und $N$ glatte Mannigfaltigkeiten und $F\in C^{\infty}(M,N)$ eine glatte Einbettung, dann ist die Menge $F(M)$ eine eingebettete Untermannigfaltigkeit.
\end{lem}

Lemma \ref{lem:2.4} entspricht \cite[Theorem 8.3.]{lee2003introduction}. Nun zum Beweis von Satz \ref{satz:2.1}:

\begin{bew}
Mit Lemma \ref{lem:2.1}, Lemma \ref{lem:2.4} und Lemma \ref{lem:2.2} reicht es, die Aussagen f\"ur eine kompakte eingebettete Untermannigfaltigkeit eines Vektorraumes $\mathbb{R}^q$ zu zeigen. Es wird gezeigt, dass wenn $q>2n+1$, beziehungsweise $q>2n$ gilt, dann existiert ein $v\in\gls{Sq}$, so dass die Orthogonalprojektion:
\begin{align}\label{eq:2.7}
\begin{split}
\pi_{v}: M&\longrightarrow lin\{v \}^{\bot}\\
x&\mapsto x-(x\cdot v)\ v
\end{split}
\end{align}

eine glatte Einbettung, beziehungsweise eine Immersion ist. Die Behauptungen des Satzes folgen, wegen $dim(lin\{v \}^{\bot})=q-1$ , durch eine endliche Hintereinanderausf\"uhrung dieser Projektionen. F\"ur diese sukzessive Reduktion werden zwei Argumente ben\"otigt, die zur besseren \"Ubersicht zun\"achst getrennt voneinander behandelt werden:

\begin{lem}\label{argu:2.1}
Falls $q> 2n+1$ gilt, dann ist die Menge:
\begin{equation*}
V_1:=\left\{v\in\mathbb{S}^{q-1}: \pi_v \text{ ist injektiv }\right\}
\end{equation*}

eine dichte Teilmenge von $\mathbb{S}^{q-1}$. 

\end{lem}

\begin{bew}

Betrachte f\"ur $\Delta:= \left\{(x,y) \in M\times M: x=y \right\}$ die glatte Abbildung:
\begin{align*}
h: (M\times M)\backslash \Delta&\longrightarrow  \mathbb{S}^{q-1}\\
(x,y)&\mapsto \frac{x-y}{|x-y|}
\end{align*}

Hierbei ist die Menge $(M\times M)\backslash \Delta$ als offene Untermannigfaltigkeit der glatten Produktmannigfaltigkeit $M\times M$ aufzufassen, siehe hierf\"ur auch Definition \ref{defi:B9} und Definition \ref{defi:B10}.
Es wird gezeigt, dass die \"Aquivalenz:
\begin{equation}\label{eq:2.5}
v\notin h((M\times M)\backslash \Delta)\Longleftrightarrow \pi_v \text{ injektiv}
\end{equation}

gilt. Dazu werde angenommen, $\pi_v$ sei nicht injektiv. Dies ist gleichbedeutend damit, dass ein $(x,y)\in (M\times M)\backslash \Delta$ existiert, so dass $\pi_v(x)=\pi_v(y)$ gilt. Das hei\ss{}t:
\begin{align*}
  & &x-(x\cdot v)\ v&=y-(y\cdot v)\ v &\\
\Leftrightarrow & &x-y&=(x\cdot v)\ v-(y\cdot v)\ v &\\
& & &=\left[(x\cdot v)-(y\cdot v) \right]\ v\\
& & &=\left[(x-y)\cdot v \right]\ v
\end{align*}

Da $v\in\mathbb{S}^{q-1}$ gilt, ist diese Gleichheit a\"quivalent zu $v=\pm \frac{x-y}{|x-y|}$, beziehungsweise $v\in h(M\times M\backslash \Delta)$, womit \eqref{eq:2.5} gezeigt ist. Es folgt die Charakterisierung $V_1=\mathbb{S}^{q-1}\backslash h((M\times M)\backslash \Delta)$. Aufgrund der Absch\"atzung:
\begin{equation*}
\dim((M\times M)\backslash \Delta)=2n<q-1=\dim(\mathbb{S}^{q-1})
\end{equation*}

ist die Menge $V_1$, mit Lemma \ref{lem:2.3}, eine dichte Teilmenge von $\mathbb{S}^{q-1}$.

\end{bew}

\begin{lem}\label{argu:2.2}
Falls $q> 2n$ gilt, dann ist die Menge:
\begin{equation*}
V_2:=\left\{v\in\mathbb{S}^{q-1}: \pi_v \text{ ist eine Immersion }\right\}
\end{equation*}

eine offene dichte Teilmenge von $\mathbb{S}^{q-1}$.

\end{lem}

\begin{bew}

F\"ur die folgende Definition sei erw\"ahnt, dass f\"ur jedes $p\in M$ der Tangentialraum $T_p M$, \"uber die Ableitung der Inklusionsabbildung $i: M\longrightarrow \mathbb{R}^q$, welche mit \cite[Theorem 8.2]{lee2003introduction} glatt ist, als Unterraum von $T_p \mathbb{R}^q$ aufgefasst werden kann, siehe hierf\"ur auch \cite[Chapter 8, The Tangent Space to an Embedded Submanifold]{lee2003introduction}. Der Raum $T_p \mathbb{R}^q$ wird mit den Identifikationen in \cite[Chapter 3]{lee2003introduction} als $\mathbb{R}^q$ angesehen. Die Menge:
\begin{equation}
SM:= \coprod_{p\in M}{\left\{w\in T_p M: |w|_{\mathbb{R}^q}=1 \right\}}
\end{equation}

wird als glatte Mannigfaltigkeit der Dimension $2n-1$ aufgefasst, siehe hierf\"ur auch \cite[Chapter 1, Section 3]{hirsch1997differential}. Betrachte die glatte Abbildung:
\begin{align*}
\pi_1: SM&\longrightarrow \mathbb{S}^{q-1}\\
(w,p)&\mapsto  w
\end{align*}

Es wird gezeigt, dass die \"Aquivalenz:
\begin{equation}\label{eq:2.6}
v\notin \pi_1(SM)\Longleftrightarrow \pi_v\text{ ist eine Immersion}
\end{equation}

gilt. Betrachte dazu, f\"ur $v\in\mathbb{S}^{q-1}$, die Abbildung $\pi_v$ auf ganz $\mathbb{R}^q$, und bezeichne diese mit $\widehat{\pi}_v$. Das bedeutet konkret:
\begin{align*}
\widehat{\pi}_v: \mathbb{R}^q&\longrightarrow \mathbb{R}^q\\
x&\mapsto x-(x\cdot v)\ v
\end{align*}

F\"ur einen festen Punkt $p\in M$ gilt, wegen $v\in\mathbb{S}^{q-1}$, $Kern(\widehat{\pi}_{v{\ast}})=lin\{v\}\subseteq T_p \mathbb{R}^q$. Da wegen $\pi_v=\widehat{\pi}_v\circ i$, mit der Inklusion $i:M\longrightarrow\mathbb{R}^q$, die Gleichheit $Kern(\pi_{v{\ast}})=Kern(\widehat{\pi}_{v{\ast}})\cap T_p M=lin\{v\}\cap T_p M$ gilt, ist $\pi_v$ genau dann eine Immersion, wenn $v\notin T_p M$ f\"ur jedes $p\in M$, beziehungsweise $v\notin \pi_1(SM)$ gilt. Damit ist \eqref{eq:2.6} bewiesen, und es folgt $V_2=\mathbb{S}^{q-1}\backslash \pi_1(SM)$.\\
Aus der Kompaktheit der Mannigfaltigkeit $M$ folgt die Kompaktheit von $SM$. Wegen der Stetigkeit von $\pi_1$ ist die Menge $\pi_1(SM)$ kompakt, und damit insbesondere abgeschlossen in $\mathbb{S}^{q-1}$, woraus die Offenheit der Menge $V_2$ folgt. Ferner folgt aus der Absch\"atzung:
\begin{equation*}
dim(SM)=2n-1<q-1=dim(\mathbb{S}^{q-1})
\end{equation*}

mit Lemma \ref{lem:2.3} die Dichtheit der Menge $V_2$ in $\mathbb{S}^{q-1}$.
\end{bew}

Aus den Argumenten in Lemma \ref{argu:2.1} und Lemma \ref{argu:2.2} ergibt sich die Aussage:
\begin{equation*}
V_1\cap V_2=\left\{v\in\mathbb{S}^{q-1}: \pi_v\text{ ist eine injektive Immersion }\right\}\neq \emptyset
\end{equation*}

Wegen der Kompaktheit von $M$ ist, mit \cite[Proposition 7.4]{lee2003introduction}, jede injektive Immersion eine glatte Einbettung. Somit ist f\"ur $q>2n+1$:
\begin{equation*}
V_1\cap V_2=\left\{v\in\mathbb{S}^{q-1}: \pi_v\text{ ist eine glatte Einbettung }\right\}\neq \emptyset
\end{equation*}

Solange also die Raumdimension $q$ gr\"o\ss{}er als $2n+1$ ist, existiert ein $v\in\mathbb{S}^{q-1}$, so dass die Abbildung $\pi_v:M\longrightarrow lin\{v\}^{\bot}$ eine glatte Einbettung ist.
\end{bew}

In \cite{whit1944a} und \cite{whit1944b} wird, unter Verwendung anderer Beweismethoden, gezeigt, dass sich die Dimensionen in Satz \ref{satz:2.1}, f\"ur $n>1$, jeweils um $1$ reduzieren lassen, das hei\ss{}t es existiert eine glatte Einbettung in den Raum $\mathbb{R}^{2n}$ und eine Immersion in den Raum $\mathbb{R}^{2n-1}$.


\section[Freie Einbettung einer kompakten Riemannschen Mannigfaltigkeit]{Freie Einbettung einer kompakten Riemannschen Mannigfaltigkeit}
\label{sec:Ab2.5}

Die bisherigen Betrachtungen werden nun zum zentralen Existenzsatz, f\"ur eine freie Einbettung einer Riemannschen Mannigfaltigkeit  ausgebaut:

\begin{satz}\label{satz:2.3}
Gegeben sei eine kompakte n-dimensionale Riemannsche Mannigfaltigkeit\linebreak $(M,g)$, dann existiert eine freie Einbettung $F_0\in C^{\infty}(M,\mathbb{R}^{\frac{n}{2}(n+5)})$, so dass $g-F^{\ast}_0(g^{can})$ eine Riemannsche Metrik auf $M$ ist.
\end{satz}

\begin{bew}

Mit Satz \ref{satz:2.1} existiert eine glatte Einbettung $F_E\in C^{\infty}(M,\mathbb{R}^{2n+1})$, und mit Satz \ref{satz:2.2} l\"asst sich eine freie Einbettung $F_f\in C^{\infty}(\mathbb{R}^{2n+1},\mathbb{R}^{(2n+1)(n+2)})$ finden. Ferner kann mit Lemma \ref{lem:2.4} die Menge $F_E(M)$ als eingebettete Untermannigfaltigkeit des Raumes $\mathbb{R}^{2n+1}$ aufgefasst werden, und mit Lemma \ref{lem:2.5} ist die Einschr\"ankung $\left.F_f\right|_{F_E(M)}: F_E(M)\longrightarrow \mathbb{R}^{(2n+1)(n+2)}$ eine freie Abbildung auf $F_E(M)$. Daraus folgt, mit Lemma \ref{lem:2.2}, dass die Komposition $F_f\circ F_E: M\longrightarrow \mathbb{R}^{(2n+1)(n+2)}$
eine freie Einbettung auf $M$ ist. Die Raumdimension $(2n+1)(n+2)$ ist aber gr\"o\ss{}er als $\frac{n}{2}(n+5)$.\\
Aufbauend auf Satz \ref{satz:2.1} soll f\"ur eine gegebene freie Einbettung $F\in C^{\infty}(M,\mathbb{R}^q)$ ein geeignetes $v\in\mathbb{S}^{q-1}$ bestimmt werden, so dass wenn $q>\frac{n}{2}(n+5)$ erf\"ullt ist, die Abbildung $\pi_v\circ F:M\longrightarrow lin\{v\}^{\bot}$ eine freie Einbettung ist, hierbei ist $\pi_v$ die Orthogonalprojektion aus \eqref{eq:2.7}. Um das zu erreichen, wird, neben den beiden in Satz \ref{satz:2.1} verwendeten Reduktionsargumenten, ein weiteres ben\"otigt:

\begin{lem}
Falls $q> \frac{n}{2}(n+5)$ gilt, dann ist die Menge:
\begin{equation*}
V_3:=\left\{v\in \mathbb{S}^{q-1}:\pi_v\circ F\text{ ist eine freie Abbildung} \right\}
\end{equation*}

eine offene dichte Teilmenge in $\mathbb{S}^{q-1}$.

\end{lem}

Definiere:
\begin{equation}
SD^2(F):=\coprod_{p\in M}{\{w\in D^2_p(F): |w|_{\mathbb{R}^q}=1\}}
\end{equation}

Mit \cite[2.5]{gromov1970embeddings} kann die Menge $SD^2(F)$ als glatte Mannigfaltigkeit der Dimension $\frac{n}{2}(n+5)-1$ aufgefasst werden. Betrachte die glatte Abbildung:
\begin{align*}
\pi_2: SD^2(F)&\longrightarrow\mathbb{S}^{q-1}\\
(w,p)&\mapsto w
\end{align*}

Es wird die folgende Aussage gezeigt:
\begin{equation}\label{eq:2.8}
v\notin \pi_2(SD^2(F))\Longleftrightarrow \pi_v\circ F\text{ ist eine freie Abbildung}
\end{equation}

Angenommen $\pi_v\circ F$ ist keine freie Abbildung, was genau dann der Fall ist, wenn ein $p\in M$, ein $(U,\varphi)\in\mathcal{A}$ , und eine nichttriviale Auswahl von Koeffizienten $\left\{\lambda_{i}\right\}_{1\leq i\leq n}\cup \left\{\lambda_{ij}\right\}_{1\leq i\leq j\leq n}\subseteq\mathbb{R}$ existiert, so dass:
\begin{align*}
\sum_{i=1}^n{\lambda_i\, \partial_i\, ^{\varphi}(\pi_v\circ F)(\varphi(p))}+\sum_{1\leq i\leq j\leq n}{\lambda_{ij}\, \partial_i\partial_j\, ^{\varphi}(\pi_v\circ F)(\varphi(p))}=0
\end{align*}

Das ist gleichbedeutend mit:
\begin{align*}
\sum_{i=1}^n{\lambda_i\, \pi_v\left[ \partial_i\, ^{\varphi} F(\varphi(p))\right]}+\sum_{1\leq i\leq j\leq n}{\lambda_{ij}\, \pi_v \left[\partial_i\partial_j\, ^{\varphi} F(\varphi(p))\right]}=0
\end{align*}

beziehungsweise, nach eventueller Skalierung der Koeffizienten, und mit der Definition der Abbildung $\pi_v$ in \eqref{eq:2.7}:
\begin{align*}
D^2_p(F)\cap \mathbb{S}^{q-1}\ni&\sum_{i=1}^n{\lambda_i\,  \partial_i\, ^{\varphi} F(\varphi(p))}+\sum_{1\leq i\leq j\leq n}{\lambda_{ij}\,  \partial_i\partial_j\, ^{\varphi} F(\varphi(p))}\\
&=\sum_{i=1}^n{\lambda_i\, v\cdot  \partial_i\, ^{\varphi} F(\varphi(p))}\ v+\sum_{1\leq i\leq j\leq n}{\lambda_{ij}\, v\cdot  \partial_i\partial_j\, ^{\varphi} F(\varphi(p))\ v}\\
&=v\cdot\left[\sum_{i=1}^n{\lambda_i\,  \partial_i\, ^{\varphi} F(\varphi(p))}+\sum_{1\leq i\leq j\leq n}{\lambda_{ij}\,  \partial_i\partial_j\, ^{\varphi} F(\varphi(p))}\right]\ v
\end{align*}

Diese Gleichung ist \"aquivalent zu:
\begin{align*}
v=\pm \sum_{i=1}^n{\lambda_i\,  \partial_i\, ^{\varphi} F(\varphi(p))}+\sum_{1\leq i\leq j\leq n}{\lambda_{ij}\,  \partial_i\partial_j\, ^{\varphi} F(\varphi(p))}\in D^2_p(F)\cap \mathbb{S}^{q-1}
\end{align*}

beziehungsweise $v\in \pi_2(S D^2(F))$. Aussage \eqref{eq:2.8} ist damit bewiesen, und es folgt die Charakterisierung $V_3=\mathbb{S}^{q-1}\backslash \pi_2(SD^2(F))$.\\
Mit der Kompaktheit von $M$ l\"asst sich die Kompaktheit von $S D^2(F)$ zeigen, woraus, mit der Stetigkeit von $\pi_2$, die Offenheit von $V_3$ folgt. Die Dichtheit der Menge $V_3$ in $\mathbb{S}^{q-1}$ folgt mit Lemma \ref{lem:2.3} aus der Absch\"atzung:
\begin{equation*}
dim(SD^2(F))=\frac{n}{2}(n+5)-1< q-1=dim(\mathbb{S}^{q-1})
\end{equation*}

Ist $q>\frac{n}{2}(n+5)$ erf\"ullt, womit auch $q>2n+1$ erf\"ullt ist, so ergibt sich, unter Beachtung von Lemma \ref{argu:2.1} und Lemma \ref{argu:2.2}, die Aussage:
\begin{align*}
\left\{v\in\mathbb{S}^{q-1}: \pi_v\circ F \text{ ist eine freie Einbettung } \right\}\supseteq V_1\cap V_2\cap V_3\neq\emptyset
\end{align*}

Durch wiederholte Anwendung dieser Reduktionsargumente l\"asst sich eine freie Einbettung in den Raum $\mathbb{R}^{\frac{n}{2}(n+5)}$ konstruieren. Der Einfachheit halber werde diese Abbildung wieder mit $F$ bezeichnet.\\
Es bleibt zu untersuchen, wie aus dieser Abbildung eine freie Einbettung $F_0$ konstruiert werden kann, so dass $g-F_0^{\ast}(g^{can})$ eine Riemannsche Metrik auf $M$ ist. Das hei\ss{}t, dass in jedem Punkt $p\in M$ der Tensor $g(p)-F_0^{\ast}(g^{can})(p)\in T^2(T_p M)$ positiv definit sein soll. Dies l\"asst sich, unter Verwendung der Kompaktheit von $M$, durch eine Skalierung von $F$ erreichen:\\
Das \textbf{Einheitstangentialb\"undel}\index{Tangentialb\"undel>Einheits-$\sim$}:
\begin{equation*}
\gls{UM}:=\coprod_{p\in M}{\left\{X\in T_p M: |X|_g=1\right\}}
\end{equation*}

ist mit \cite[Theorem 8.8]{lee2003introduction} eine eingebettete Untermannigfaltigkeit des, in \cite[Lemma 4.1]{lee2003introduction} konstruierten, Tangentialb\"undels \index{Tangentialb\"undel} \gls{TM}. Ist $M$ kompakt, so ist auch $UM$ kompakt. F\"ur ein Tensorfeld $\sigma\in \mathcal{T}^2(M)$ ist die Abbildung:
\begin{equation*}
UM\ni (X,p)\mapsto \sigma(X,X)\in \mathbb{R}
\end{equation*}

stetig, und nimmt daher einen minimalen und einen maximalen Wert an. Sei nun $\epsilon\in\mathbb{R}_{>0}$, sowie $p\in M$ und $X\in T_p M$ mit $|X|_{g}=1$, dann ist:
\begin{equation}\label{eq:2.9}
\left[g-\epsilon F^{\ast}(g^{can})\right](X,X)=1-\epsilon F^{\ast}(g^{can})(X,X)\geq 1-\epsilon\max_{(Y,p)\in UM}{F^{\ast}(g^{can})(Y,Y)}
\end{equation}

Da $F^{\ast}(g^{can})$, als Riemannsche Metrik, insbesondere positiv definit ist, existiert ein $\epsilon_0\in\mathbb{R}_{>0}$ mit:
\begin{equation*}
\epsilon_0<\frac{1}{\max_{(Y,p)\in UM}{F^{\ast}(g^{can})(Y,Y)}}
\end{equation*}

Mit \eqref{eq:2.9} erf\"ullt die freie Einbettung $F_0:= \epsilon_0\cdot F\in C^{\infty}(M,\mathbb{R}^{\frac{n}{2}(n+5)})$ alle geforderten Eigenschaften, und der Satz ist bewiesen.

\end{bew}

Abschlie\ss{}end sei noch erw\"ahnt, dass die in Satz \ref{satz:2.3} beschriebene Raumdimension $\frac{n}{2}(n+5)$ nicht mehr reduziert werden kann, falls $n=2^k$ f\"ur ein $k\geq 2$ gilt. Hierf\"ur wird auf \cite{elias} verwiesen.

\chapter{Reduktion auf ein lokales Problem}
\label{chap:Kap3}
\thispagestyle{fancy}

Hauptgegenstand dieses Kapitels ist es zu zeigen, dass f\"ur eine gegebene Riemannsche Mannigfaltigkeit $(M,g)$, und die in \autoref{sec:Ab2.5} konstruierte freie Einbettung $F_0\in C^{\infty}(M,\mathbb{R}^{\frac{n}{2}(n+5)})$, die Riemannsche Metrik $h:= g-F_0^{\ast}(g^{can})$ als Summe von Tensorfeldern mit bestimmten w\"unschenswerten Eigenschaften dargestellt werden kann.

\section{Dekomposition der Metrik}

Der Inhalt des folgenden Satzes entspricht \cite[Satz 2.2.]{gunther1989einbettungssatz} f\"ur den kompakten Fall. Das, in diesem Satz verwendete, Lemma \ref{lem:3.4} wird in \autoref{sec:Ab3.2} nachtr\"aglich bewiesen. 

\begin{satz}\label{satz:3.1}

Es sei $(M,h)$ eine kompakte Riemannsche Mannigfaltigkeit, dann existieren glatte symmetrische kovariante $2-$Tensorfelder $h^{(1)},...,h^{(m)}\in\mathcal{T}^2(M)$ mit den Eigenschaften:
\begin{enumerate}[(i)]
\begin{item}
\begin{equation*}
h=\sum_{i=1}^m{h^{(i)}}
\end{equation*}
\end{item}

\begin{item}\label{satz:3.1.2}
F\"ur jedes $i\in\{1,...,m\}$ existiert eine Karte $(U_i,\varphi_i)\in\mathcal{A}$, mit $\varphi_i(U_i)=B_{1+\tau}(0)$, f\"ur ein $\tau>0$, und eine Abbildung $a_i\in C^{\infty}(U_i)$ mit $supp(a_i)\subseteq \varphi^{-1}(\mathbb{B})$, so dass:
\begin{align*}
h^{(i)}(x)=
\begin{cases}
a_i^4(x) \left.\, ^{\varphi_i}dx^1 \right|_x^2 &\text{f\"ur }x\in U_i\\
0 &\text{sonst}
\end{cases}
\end{align*}

\end{item}

\end{enumerate}

\end{satz}

Dem Beweis dieses Satzes wird eine Vor\"uberlegung vorangestellt: Betrachte die Menge von Abbildungen:
\begin{equation}\label{eq:3.2}
L(\mathbb{B}):= \left\{f\in C^{\infty}(\mathbb{B}): f(x)=\sum_{i=1}^n{c_i x^i}\text{ f\"ur ein }(c_1,...,c_n)\in\mathbb{R}^n \right\}
\end{equation}

Es sei $p\in M$, und $(U_p,\varphi_p)\in\mathcal{A}$ ein um $p$ zentrierter Koordinatenball mit dem Radius $1$, das hei\ss{}t: $p\in U$, $\varphi_p(p)=0$ und $\varphi_p(U_p)=\mathbb{B}$. Aus einer Funktion $f_p\in L(\mathbb{B})$ erh\"alt man, durch Verkn\"upfung  mit der Kartenabbildung $\varphi_p$, eine Abbildung $f_p\circ \varphi_p\in C^{\infty}(U_p)$. F\"ur diese Abbildung ist der Pullback $(f_p\circ \varphi_p)^{\ast}(g^{can})\in \mathcal{T}^2(U_p)$ wohldefiniert, und hat unter Beachtung von \eqref{eq:1.12} die Darstellung:
\begin{equation}\label{eq:3.1}
(f_p\circ \varphi_p)^{\ast}(g^{can})=\partial_i\> ^{\varphi_p}(f_p\circ \varphi_p)\cdot \partial_j \> ^{\varphi_p}(f_p\circ \varphi_p) \ ^{\varphi_p}dx^i\otimes\> ^{\varphi_p} dx^j=\partial_i f_p\cdot \partial_j f_p \ ^{\varphi_p}dx^i\otimes\> ^{\varphi_p} dx^j
\end{equation}

Sei dar\"uber hinaus eine Abbildung $\chi_p\in C^{\infty}_0(U_p)$ gegeben, dann ist $\chi_p^4\,\cdot\, (f_p\circ \varphi_p)^{\ast}(g^{can})\in \mathcal{T}^2(U_p)$ wohldefiniert, und hat mit \eqref{eq:3.1} die Darstellung:
\begin{equation*}
\chi_p^4\,\cdot\,(f_p\circ \varphi_p)^{\ast}(g^{can})=\chi_p^4\,\cdot\, \partial_i f_p\cdot \partial_j f_p\,\cdot\,^{\varphi_p}dx^i\otimes\> ^{\varphi_p} dx^j=\chi_p^4\,\cdot\,c_i^{(p)} c_j^{(p)} \,\cdot\, ^{\varphi_p}dx^i\otimes\> ^{\varphi_p} dx^j
\end{equation*}

Mittels trivialer Fortsetzung, wird $\chi_p^4\,\cdot\, (f_p\circ \varphi_p)^{\ast}(g^{can})$ absofort als glattes Tensorfeld in $\mathcal{T}^2(M)$ aufgefasst. Es wird gezeigt, dass f\"ur dieses Tensorfeld eine Karte $(U_p,\widehat{\varphi}_p)\in\mathcal{A}$, und eine Abbildung $a_p\in C^{\infty}_0(U_p)$ existiert, so dass:
\begin{align}\label{eq:3.3}
\begin{split}
\chi_p^4\,\cdot\, (f_p\circ \varphi_p)^{\ast}(g^{can})=
\begin{cases}
a_p^4\,\left.\cdot\, ^{\widehat{\varphi}_p}dx^1 \right|^2 &\text{in } U_p\\
0 &\text{sonst}
\end{cases}
\end{split}
\end{align}

Es sei $A_{(p)}\in GL(n,\mathbb{R})$ die Matrix, deren Inverse in der ersten Spalte den Vektor $c^{(p)}$ aus \eqref{eq:3.2} enth\"alt, und deren restlichen $n-1$ Spalten aus einer festen Basis des Untervektorraumes $lin\{c^{(p)}\}^{\bot}$ bestehen. Definiere nun eine neue Karte $(U_p,\widehat{\varphi}_p)\in \mathcal{A}$, wie folgt:
\begin{align*}
&\widehat{\varphi}_p: U_p\longrightarrow\mathbb{R}^n\\
&\widehat{\varphi}_p := A_{(p)}\cdot \varphi_p
\end{align*}

Dann gilt f\"ur $i\in\{1,...n\}$:
\begin{align*}
&\partial_i\> ^{\widehat{\varphi}_p}(f_p\circ \varphi_p)=\partial_i (f_p\circ \varphi_p\circ \widehat{\varphi}_p^{-1})=\partial_i (f_p\circ \varphi_p\circ (A_{(p)}\cdot \varphi_p)^{-1})\\
=&\partial_i (f_p\circ \varphi_p\circ \varphi_p^{-1}\circ A_{(p)}^{-1})=\partial_i (f_p\circ A_{(p)}^{-1})=\sum_{k=1}^n{\partial_k f_p\cdot (A_{(p)}^{-1})_{ki}}\\
=&\sum_{k=1}^n{c^{(p)}_k\cdot (A_{(p)}^{-1})_{ki}}=
\begin{cases}
|c^{(p)}|^2_{\mathbb{R}^n} &\text{falls }i=1\\
0 &\text{sonst}
\end{cases}
\end{align*}

Damit hat $\chi_p^4\,\cdot\, (f_p\circ \varphi_p)^{\ast}(g^{can})$ in $U_p$ die Darstellung:
\begin{align*}
&\chi_p^4\,\cdot\, (f_p\circ \varphi_p)^{\ast}(g^{can})=\chi_p^4 \ \partial_i\> ^{\widehat{\varphi}_p}(f_p\circ \varphi_p)\cdot \partial_j \> ^{\widehat{\varphi}_p}(f_p\circ \varphi_p) \> ^{\widehat{\varphi}_p}dx^i\otimes\> ^{\widehat{\varphi}_p} dx^j\\
=&(\chi_p\cdot |c^{(p)}|_{\mathbb{R}^n})^4\,\cdot\, \left. ^{\widehat{\varphi}_p}dx^1\right|^2
\end{align*}

womit $\eqref{eq:3.3}$ gezeigt ist. Nun zum Beweis von Satz \ref{satz:3.1}:

\begin{bew}

Es sei:
\begin{equation*}
F:= \left\{\left(\partial_i f\cdot \partial_j f \right)_{1\leq i\leq j\leq n}: f\in L(\mathbb{B}) \right\}\subseteq C^{\infty}(\overline{\mathbb{B}},\mathbb{R}^{\frac{n}{2}(n+1)})
\end{equation*}

definiert, wobei $L(\mathbb{B})$ bereits in \eqref{eq:3.2} festgelegt worden ist. F\"ur jedes $p\in M$ sei $(U_p,\varphi_p)$ ein, um $p$ zentrierter, Koordinatenball mit dem Radius 1. Aufgrund der stetigen Abh\"an\-gigkeit der Eigenwerte von den Matrixeintr\"agen \cite[Satz 1.2.]{werner1992numerische}, existiert f\"ur jedes $p\in M$ ein $\epsilon_p>0$, so dass die Elemente der Menge:
\begin{equation*}
 W_p:= \left\{z\in \mathbb{R}^{\frac{n}{2}(n+1)}: \max_{1\leq i\leq j\leq n}{\left|z_{ij}-\>^{\varphi_p}h_{ij}(0)\right|}<\epsilon_p \right\}
\end{equation*}

aufgefasst als symmetrische Matrizen in $\mathbb{R}^{n\times n}$, durch Spiegelung an der Hauptdiagonalen, genau wie $(\,^{\varphi_p}h_{ij}(0))_{1\leq i,j\leq n}\in\mathbb{R}^{n\times n}$ positiv definit sind. Mit \cite[Satz 1.11]{schwarz1997numerische} besitzt jede beliebige positiv definite Matrix $Z\in\mathbb{R}^{n\times n}$ eine Cholesky-Zerlegung, das hei\ss{}t, es existiert eine Matrix $C\in\mathbb{R}^{n\times n}$ mit $C_{ij}=0$, f\"ur $1\leq j < i\leq  n$, so dass $Z=C^{T} C$. Das bedeutet f\"ur alle $i,j\in\{1,...,n\}$:
\begin{equation}\label{eq:3.4}
z_{ij}=\sum_{k=1}^n{(C^{T})_{ik} C_{kj}}=\sum_{k=1}^n{C_{ki} C_{kj}}
\end{equation}

Definiere nun, f\"ur jedes $k\in\{1,...,n\}$, eine Funktion $f_k^z\in L(\mathbb{B})$, gem\"a\ss{} der Bildungsvorschrift:
\begin{align*}
f_k^z(x):=\sum_{m=1}^n{C_{km} x^m}
\end{align*}

Dann gilt mit \eqref{eq:3.4} f\"ur alle $i,j\in\{1,...,n\}$ die Gleichung $z_{ij}=\sum_{k=1}^n{\partial_i f^z_k(0)\cdot \partial_j f^z_k(0)}$. Damit sind alle Voraussetzungen von Lemma \ref{lem:3.4} auf Seite \pageref{lem:3.4} erf\"ullt und es existiert f\"ur jedes $p\in M$ ein $\delta_p>0$, und ein $l\in\mathbb{N}$, sowie Abbildungen $a_{p,1},...,a_{p,l}\in C^{\infty}(B_{\delta_p}(0))$, und $f_{p,1},...,f_{p,l}\in L(B_{\delta_p}(0))$, so dass:
\begin{equation*}
^{\varphi_p}h_{ij}(x)=\sum_{k=1}^{l}{a_{p,k}^4(x)\> \partial_i f_{p,k}(x)\cdot \partial_j f_{p,k}(x)}
\end{equation*}

f\"ur alle $x\in B_{\delta_p}(0)$ gilt. Nach Umskalierung der Kartenabbildungen wird nun angenommen, dass f\"ur alle $p\in M$ die Gleichheit $\delta_p=1+\tau$, f\"ur ein von $p$ unabh\"angiges $\tau\in(0,1)$ gilt, dabei sollen f\"ur die Kartenabbildungen keine neuen Bezeichnungen eingef\"uhrt werden. Nun werden, unter Verwendung der Kompaktheit von $M$, endlich viele $p_1,...,p_m\in M$ gew\"ahlt, so dass:
\begin{equation*}
M=\bigcup_{i=1}^m{\varphi_{p_i}^{-1}(B_{1-\tau}(0))}
\end{equation*}

Mit \cite[Theorem 2.25]{lee2003introduction} existieren Abbildungen $\psi_1,...,\psi_m\in C^{\infty}(M)$, mit $supp(\psi_i)\subseteq \varphi_{p_i}^{-1}(B_{1-\tau}(0))$ f\"ur alle $i\in\{1,...,m\}$, und $\sum_{i=1}^m{\psi_i^4}\equiv 1$ auf ganz $M$. Mit der Vor\"uberlegung am Anfang des Abschnittes, ist gezeigt, dass f\"ur jedes $i\in\{1,...,m\}$ das Tensorfeld $\psi_i^4\cdot h\in\mathcal{T}^2(M)$ eine Summe von Tensorfeldern mit der Eigenschaft \eqref{satz:3.1.2} ist. Wegen $\sum_{i=1}^m{\psi_i^4}\equiv 1$ ist der Satz bewiesen.
\end{bew}

Eine weitere M\"oglichkeit, Satz \ref{satz:3.1} zu beweisen, wird in \cite[Lemma 1.3.1]{hong2006isometric} gezeigt, dabei wird das lokale Einbettungsresultat aus \cite[2.8.1.]{gromov1970embeddings} verwendet. Ausgehend von Satz \ref{satz:3.1}, wird nun der folgende Satz formuliert, welcher in \autoref{chap:Kap6} bewiesen werden wird. Dieser \"uberf\"uhrt die globale Problemstellung in eine lokale. Der Beweis von Satz \ref{satz:3.2} erfordert eine Vielzahl von technischen Hilfsmitteln, die in den Kapiteln \ref{chap:Kap4} und \ref{chap:Kap5} behandelt werden.

\begin{satz}\label{satz:3.2}
Gegeben sei eine n-dimensionale glatte Mannigfaltigkeit $M$, und eine freie Einbettung $F_0\in C^{\infty}(M,\mathbb{R}^q)$, mit $q\geq\frac{n}{2}(n+3)+5$. Weiterhin sei ein glattes symmetrisches kovariantes $2-$Tensorfeld $h\in \mathcal{T}^2(M)$ gegeben, f\"ur das eine Karte $(U,\varphi)\in\mathcal{A}$, mit $\varphi(U)=B_{1+\tau}(0)$, f\"ur $\tau>0$, und eine Abbildung $a\in C^{\infty}_0(U)$ existiert, welche $supp(a)\subseteq \varphi^{-1}(\mathbb{B})$ erf\"ullt, so dass:
\begin{align*}
h(x)=
\begin{cases}
a^4(x) \left.\, ^{\varphi}dx^1 \right|_x^2 &\text{f\"ur }x\in U\\
0 &\text{sonst}
\end{cases}
\end{align*}

gilt. Dann existiert zu jedem $\epsilon\in\mathbb{R}_{>0}$ eine freie Einbettung $F\in C^{\infty}(M,\mathbb{R}^q)$, so dass:
\begin{align*}
&F^{\ast}(g^{can})=F_0^{\ast}(g^{can})+h \\
&supp(F-F_0)\subseteq U\\
&\max_{x\in M}{|F(x)-F_0(x)|_{\mathbb{R}^q}}\leq C\cdot \epsilon \\
\end{align*}

gilt, wobei die Konstante $C\in\mathbb{R}_{>0}$ nicht von $\epsilon$ abh\"angt.

\end{satz}

\section{Hilfsresultate}
\label{sec:Ab3.2}

Ziel dieses Abschnittes ist es, das in Satz \ref{satz:3.1} verwendete, Lemma \ref{lem:3.4} zu beweisen. Zun\"achst werden Begriffe aus der Konvexgeometrie eingef\"uhrt:

\begin{defi}
F\"ur eine Teilmenge $X\subseteq \mathbb{R}^q$ wird die Menge:
\begin{enumerate}[(i)]

\begin{item}
\begin{equation*}
\gls{ccone(X)}:=\left\{\sum_{i=1}^m{\lambda_i x_i}: m\in\mathbb{N},\ x_1,...,x_m\in X,\ \lambda_1,...,\lambda_m\in \mathbb{R}_{\geq 0} \right\}
\end{equation*}

als \textbf{konvex-konische H\"ulle von}\index{H\"ulle>konvex-konische $\sim$} $\bm{X}$ bezeichnet.

\end{item}

\begin{item}
\begin{equation*}
\gls{conv(X)}:=\left\{\sum_{i=1}^m{\lambda_i x_i}: m\in\mathbb{N},\ x_1,...,x_m\in X,\ \lambda_1,...,\lambda_m\in \mathbb{R}_{\geq 0}\text{ mit }\sum_{i=1}^m{\lambda_i}=1\right\}
\end{equation*}
\end{item}

als \textbf{konvexe H\"ulle von}\index{H\"ulle>konvexe $\sim$} $\bm{X}$ bezeichnet.

\end{enumerate}

\end{defi}

Die Inklusion $conv(X)\subseteq ccone(X)$ ist offensichtlich.

\begin{lem}\label{lem:3.1}

Gegeben seien linear unabh\"angige Vektoren $v_1,..,v_q\in\mathbb{R}^q$, dann gilt:
\begin{equation*}
z\in int(ccone\{v_1,..,v_q\})\Longleftrightarrow z=\sum_{i=1}^{q}{\lambda_i v_i} \hspace{0.5cm}\text{mit } \lambda_i\in\mathbb{R}_{>0}\ \forall i\in\{1,..,q\}
\end{equation*}

\end{lem}

\begin{bew}

Es sei $V\in\mathbb{R}^{q\times q}$ die Matrix, deren Spaltenvektoren $v_1,..,v_q$ sind. Wird $V$ als Abbildung $\mathbb{R}^{q}\longrightarrow\mathbb{R}^{q}$ aufgefasst, so ist diese Abbildung, aufgrund der Linearit\"at, stetig. Das gilt auch f\"ur die inverse Abbildung $V^{-1}: \mathbb{R}^{q}\longrightarrow\mathbb{R}^{q}$. Die Abbildung $V$ ist also ein Hom\"oomorphismus auf $\mathbb{R}^{q}$. F\"ur beliebige $z\in\mathbb{R}^q$ gilt die folgende Aussage:
\begin{equation*}
z\in ccone\{v_1,..,v_q\}\Leftrightarrow V^{-1}z\in\mathbb{R}_{\geq 0}^{q}
\end{equation*}

oder auch:
\begin{equation*}
ccone\{v_1,..,v_q\}=V(\mathbb{R}^q_{\geq 0})
\end{equation*}

Da $\partial \mathbb{R}_{\geq 0}^{q}=\left\{z\in \mathbb{R}_{\geq 0}^{q}: \exists i\in\{1,...,q\}: z_i=0 \right\}$, folgt mit:
\begin{equation*}
int(\mathbb{R}_{\geq 0})=\overline{\mathbb{R}^q}_{\geq 0}\backslash  \partial \mathbb{R}_{\geq 0}^{q}=\mathbb{R}^q_{\geq 0}\backslash  \partial \mathbb{R}_{\geq 0}^{q}
\end{equation*}

aus der Hom\"oomorphie von $V$, die Behauptung.

\end{bew}

\begin{lem}\label{lem:3.2}

Es sei $z\in \mathbb{R}^q\backslash \{0\}$, und $W$ eine Umgebung von $z$. Dann existiert eine linear unabh\"angige Teilmenge $\{v_1,..,v_q\}\subseteq W$, so dass $z\in int(ccone\{v_1,..v_q\})$.

\end{lem}

\begin{bew}

Mit \cite[Example 9.1 (e)]{lee2003introduction} und der Hom\"oomorphie-Eigenschaft von Drehungen und Streckungen, kann ohne Beschr\"ankung der Allgemeinheit angenommen werden, dass:
\begin{equation}\label{eq:3.5}
z=\sum_{i=1}^{q}{e_i}
\end{equation}

gilt. Hierbei sind $e_1,...,e_q$ die kanonischen Einheitsvektoren in $\mathbb{R}^q$. Unter Verwendung der Offenheit von $W$, existiert ein $\epsilon\in\mathbb{R}_{>0}$, so dass $B_\epsilon(z)\subseteq W$ gilt. Dann gilt f\"ur alle $i\in\{1,...,q\}$:
\begin{equation}\label{eq:3.6}
v_i:= z+\frac{\epsilon}{2}e_i\in B_\epsilon(z)
\end{equation}

Nun wird gezeigt, dass die Vektoren $v_1,...,v_q$ linear unabh\"angig sind. Angenommen, es existieren Koeffizienten ${\lambda_1,..,\lambda_q}\in\mathbb{R}$,
so dass:
\begin{equation*}
0=\sum_{i=1}^q{\lambda_i}{v_i}
\end{equation*}

gilt. Dann ist
\begin{equation}\label{eq:3.6a}
0=\sum_{i=1}^q{\lambda_i}{v_i}\stackrel{\eqref{eq:3.6}}{=}\sum_{i=1}^q{\lambda_i\left(z+\frac{\epsilon}{2}e_i\right)}=\sum_{i=1}^q{\lambda_i z}+\frac{\epsilon}{2}\sum_{i=1}^q{ \lambda_i e_i}
\end{equation} 

Ist nun $\sum_{i=1}^q{\lambda_i}=0$ gilt. Dann ist 
\begin{equation*}
0=\frac{\epsilon}{2}\sum_{i=1}^q{ \lambda_i e_i}
\end{equation*}

und es folgt $\lambda_i=0$ f\"ur alle $i\in\{1,...,q\}$. Ist $\sum_{i=1}^n{\lambda_i}\neq 0$ dann kann ohne Beschr\"ankung der Allgemeinheit angenommen werden, dass $\sum_{i=1}^q{\lambda_i}= 1$ gilt. In dem Fall folgt aus \eqref{eq:3.6a}
\begin{equation*}
\sum_{i=1}^{q}{e_i}\stackrel{\eqref{eq:3.5}}{=}z=-\frac{\epsilon}{2}\sum_{i=1}^q{\lambda_i e_i}
\end{equation*} 

und es gilt $\lambda_i=-\frac{2}{\epsilon}$ f\"ur alle $i\in\{1,...,q\}$. Dann ist
\begin{equation*}
\sum_{i=1}^q{\lambda_i}=-q\cdot \frac{2}{\epsilon}\neq 1
\end{equation*}

womit die lineare Unabh\"angigkeit der Vektoren $v_1,...,v_q$ bewiesen ist. Nun wird noch gezeigt, dass auch die Bedingung $z\in int\left(ccone\left\{v_1,..,v_q\right\}\right)$ erf\"ullt ist. Es ist:
\begin{equation*}
\sum_{i=1}^q{v_i}\stackrel{\eqref{eq:3.6}}{=}\sum_{i=1}^q{\left(z+\frac{\epsilon}{2}e_i\right)}=\sum_{i=1}^q{z}+\frac{\epsilon}{2}\sum_{i=1}^q{e_i}\stackrel{\eqref{eq:3.5}}{=}qz+\frac{\epsilon}{2}z=\left(q+\frac{\epsilon}{2}\right)z
\end{equation*}

beziehungsweise:
\begin{equation*}
z=\sum_{i=1}^q{\frac{1}{q+\frac{\epsilon}{2}}\cdot v_i}
\end{equation*}

Daraus folgt mit Lemma \ref{lem:3.1} die Behauptung.

\end{bew}

\begin{lem}\label{lem:3.3}

Es sei $U\subseteq\mathbb{R}^n$ eine offene Menge mit $0\in U$, und es seien $v_1,..,v_q\in C(U,\mathbb{R}^q)$ gegeben, so dass die Vektoren $v_1(0),..,v_q(0)\in\mathbb{R}^q$ linear unabh\"angig sind. Ferner sei $z\in int(ccone(v_1(0),..,v_q(0))$. Dann existieren $\delta,\rho\in\mathbb{R}_{>0}$, so dass f\"ur jedes $x\in B_{\rho}(0)$ die Inklusion:
\begin{equation*}
B_{\delta}(z)\subseteq int(ccone(v_1(x),..,v_q(x)))
\end{equation*}

gilt.

\end{lem}

\begin{bew}

Definiere eine Abbildung:
\begin{align*}
V: U&\longrightarrow \mathbb{R}^{q\times q}\\
x&\mapsto V(x)=[v_1(x),..,v_q(x)]
\end{align*}

Hierbei wird der Raum $\mathbb{R}^{q\times q}$ mit der, von der \textbf{euklidischen Norm induzierten Matrixnorm}\index{Norm>Matrix-$\sim$}:
\begin{align*}
\gls{Matrixnorm}: \mathbb{R}^{q\times q}&\longrightarrow\mathbb{R}\\
A&\mapsto \sup_{x\in\mathbb{S}^{q-1}}{|Ax|_{\mathbb{R}^q}}
\end{align*}

versehen. Nach Voraussetzung hat die Matrix $V(0)$ den Rang $q$. Aufgrund der Stetigkeit von $V$, und der Stetigkeit der Determinantenfunktion, kann angenommen werden, dass f\"ur alle $x\in U$ die Bedingung $V(x)\in GL(n,\mathbb{R})$ gilt. Ist dies nicht der Fall, so wird $U$ verkleinert. Demzufolge ist auch die Abbildung:
\begin{align*}
V^{-1}: U&\longrightarrow GL(q,\mathbb{R})\\
x&\mapsto V^{-1}(x)
\end{align*} 

wohldefiniert, und mit der Cramerschen Regel \cite[3.3.5.]{fischer2008lineare} auch stetig. Nach Voraussetzung gilt $V^{-1}(0)z \in \mathbb{R}_{>0}^q$. Dies motiviert die Definition der folgenden Abbildung, hierbei sei $R_1\in\mathbb{R}_{>0}$, so dass $\overline{B}_{R_1}(0)\subseteq U$ und $R_2\in\mathbb{R}_{>0}$ eine beliebige feste Zahl:
\begin{align*}
\Psi: B_{R_1}(0)\times B_{R_2}(z)&\longrightarrow\mathbb{R}^q\\
(x,b)&\mapsto V^{-1}(x)b
\end{align*}

Sei nun $x\in B_{R_1}(0)$ und $b\in B_{R_2}(z)$, dann gilt:
\begin{align}\label{eq:3.7}
\begin{split}
&\left|\Psi(x,b)-\Psi(0,z)\right|_{\mathbb{R}^{q}}=\left|V^{-1}(x)b-V^{-1}(0)z\right|_{\mathbb{R}^{q}}\\
\leq&\left|[V^{-1}(x)-V^{-1}(0)]b\right|_{\mathbb{R}^{q}}+\left|V^{-1}(0)[b-z] \right|_{\mathbb{R}^{q}}\\
\leq&\left\|V^{-1}(x)-V^{-1}(0)\right\|_{\mathbb{R}^{q\times q}}\left|b\right|_{\mathbb{R}^{q}}+\left\|V^{-1}(0)\right\|_{\mathbb{R}^{q\times q}}\left|b-z\right|_{\mathbb{R}^{q}}\\
\leq&C_1(z,R_2)\cdot\left\|V^{-1}(x)-V^{-1}(0)\right\|_{\mathbb{R}^{q\times q}}+C_2(V(0))\cdot \left|b-z\right|_{\mathbb{R}^{q}}\\
\end{split}
\end{align}

W\"ahle nun ein $\epsilon>0$, so dass $B_{\epsilon}(V^{-1}(0)z)\subseteq\mathbb{R}_{>0}^q$ gilt. Aufgrund der Stetigkeit der Abbildung $V^{-1}$, existiert ein $\rho(\epsilon)\in\mathbb{R}_{>0}$, so dass f\"ur $\left|x \right|_{\mathbb{R}^{n}}<\rho$ die Absch\"atzung:
\begin{equation*}
\left\|V^{-1}(x)-V^{-1}(0)\right\|_{\mathbb{R}^{q\times q}}<\frac{\epsilon}{2 C_1}
\end{equation*}

gilt. Sei nun $\left|b-z\right|_{\mathbb{R}^{q}}<\delta:=\frac{\epsilon}{2C_2}$, dann gilt, zusammen mit \eqref{eq:3.7}, die Absch\"atzung $\left|\Psi(x,b)-\Psi(0,z)\right|_{\mathbb{R}^{q}}<\epsilon$, das hei\ss{}t $\Psi(x,b)\in B_{\epsilon}(V^{-1}(0)z)\subseteq\mathbb{R}_{>0}$. Mit Lemma \ref{lem:3.1} folgt, aus der Wahl der Abbildung $\Psi$, die Behauptung.

\end{bew}

Zusammenfassend ergibt sich, aus Lemma \ref{lem:3.1} bis Lemma \ref{lem:3.3}, das in \cite[Lemma 2.3.]{gunther1989einbettungssatz} aufgef\"uhrte Resultat:

\begin{lem}\label{lem:3.4}
Es sei $U\subseteq\mathbb{R}^n$ eine offene Menge mit $0\in U$, eine Funktion $g\in C^{\infty}(U,\mathbb{R}^q)$ mit $g(0)\neq 0$, eine offene Menge $W\subseteq\mathbb{R}^q$ mit $g(0)\in W$, sowie eine Menge von Funktionen $F\subseteq C^{\infty}(U,\mathbb{R}^q)$ mit folgender Eigenschaft gegeben: F\"ur jedes $z\in W$ existiert eine Menge von Koeffizienten $\{\alpha^z_1,..,\alpha^z_{k}\}\subseteq\mathbb{R}_{>0}$ und eine Teilmenge von Funktionen $\{f_1^z,..,f_{k}^z\}\subseteq F$, so dass:
\begin{equation*}
z=\sum_{i=1}^{k}{\alpha_i^z f_i^z(0)}
\end{equation*}

Dann existiert ein $\rho>0$, ein $l(k,q)\in\mathbb{N}$, eine Menge von Koeffizientenfunktionen $\{\alpha_1,..,\alpha_l\}\subseteq C^{\infty}(B_{\rho}(0),\mathbb{R}_{>0})$, sowie eine Menge von Funktionen $\{f_1,..,f_l\}\subseteq F$, so dass:
\begin{equation*}
g(x)=\sum_{i=1}^{l}{\alpha_i(x)f_i(x)}
\end{equation*}

f\"ur alle $x\in B_{\rho}(0)$ gilt.

\end{lem}

\begin{bew}
Mit Lemma \ref{lem:3.2} existieren linear unabh\"angige Vektoren $\{v_1,..,v_q\}\subseteq W$, so dass $g(0)\in int(ccone\{v_1,..,v_q\})$. Schreibe, gem\"a\ss{} der Voraussetzung, f\"ur alle $j\in\{1,..,q \}$:
\begin{equation*}
v_j=\sum_{i=1}^{k}{\beta_{i}^{v_j}f_{i}^{v_j}(0)}
\end{equation*}

wobei die Koeffizienten $\beta_{1}^{v_j},...,\beta_{k}^{v_j}$ positiv sind, und die Abbildungen $f_{1}^{v_j},...,f_{k}^{v_j}$ in der Menge $F$ liegen. Daraus hervorgehend, wird f\"ur jedes $j\in\{1,...,q\}$ die Abbildung $v_j\in C^{\infty}(U,\mathbb{R}^q)$, mit der Bildungsvorschrift:
\begin{equation}\label{eq:3.8}
v_j(x):=\sum_{i=1}^{k}{\beta_{i}^{v_j}f_{i}^{v_j}(x)}\hspace{0.5cm}\text{f\"ur }x\in U
\end{equation}

definiert. Mit Lemma \ref{lem:3.3} existieren $\delta,\rho\in\mathbb{R}_{>0}$, so dass:
\begin{equation*}
B_{\delta}(g(0))\subseteq int(ccone(v_1(x),...,v_q(x)))
\end{equation*}

f\"ur alle $x\in B_{\rho}(0)$ gilt. Definiere nun:
\begin{align*}
V: B_{\rho}(0)&\longrightarrow GL(q,\mathbb{R})\\
x&\mapsto[v_1(x),...,v_q(x)]
\end{align*}

Wegen der Stetigkeit von $g$, kann angenommen werden, dass $g(B_{\rho}(0))\subseteq B_{\delta}(g(0))$ gilt. Dann ist $V(x)^{-1}(g(x))\in\mathbb{R}_{>0}^q$ f\"ur alle $x\in B_{\rho}(0)$. Definiere, daraus resultierend, eine glatte Abbildung:
\begin{align*}
\gamma: B_{\rho}(0)&\longrightarrow \mathbb{R}_{>0}^q\\
x&\mapsto V(x)^{-1}(g(x))
\end{align*}

und es gilt f\"ur jedes $x\in B_{\rho}(0)$:
\begin{align*}
g(x)=V(x)\cdot V(x)^{-1}(g(x))=\sum_{j=1}^q{\gamma_j(x)\,v_j(x)}\stackrel{\eqref{eq:3.8}}{=}\sum_{j=1}^{q}{\sum_{i=1}^{k}{\underbrace{\gamma_j(x)\beta_i^{v_j}}_{>0}f_i^{v_j}(x)}}
\end{align*}

woraus die Behauptung folgt.

\end{bew}

\chapter{Reduktion auf ein Perturbationsproblem}\label{chap:Kap4}
\thispagestyle{fancy}

In diesem Kapitel wird die lokale Problemstellung aus Satz \ref{satz:3.2} in ein Perturbationsproblem \"uberf\"uhrt, welches sich mit der Methode, die in \autoref{chap:Kap5} diskutiert werden wird, l\"osen l\"asst. Dazu wird der folgende Satz, welcher der Arbeit \cite[Kapitel 4]{gunther1989einbettungssatz}, beziehungsweise \cite[Theorem 1.3.9.]{hong2006isometric} entnommen ist, gezeigt:

\begin{satz}\label{satz:4.1}
Gegeben sei eine freie Abbildung $F_0\in C^{\infty}(\overline{\mathbb{B}},\mathbb{R}^q)$, f\"ur $q\geq \frac{n}{2}(n+3)+5$, und eine Abbildung $a\in C^{\infty}_0(\mathbb{B})$. Dann existiert eine kompakte Menge $K\subseteq \mathbb{B}$, mit $supp(a)\subseteq K$, so dass die folgende Aussage gilt: F\"ur jedes $k\in\mathbb{N}\backslash\{0,1\}$ existiert ein $\epsilon_k\in\mathbb{R}_{>0}$, so dass f\"ur jedes $\epsilon\in (0,\epsilon_k]$ eine freie Abbildung $F_{\epsilon,k}\in C^{\infty}(\overline{\mathbb{B}},\mathbb{R}^q)$ existiert, welche die Eigenschaften:
\begin{align}
\begin{split}\label{eq:4.1}
&\partial_1 F_{\epsilon,k}\cdot \partial_1 F_{\epsilon,k}=\partial_1 F_0\cdot \partial_1 F_0 +a^4+\epsilon^{k+1}f_{11}^{\epsilon,k}\\
&\partial_1 F_{\epsilon,k}\cdot \partial_i F_{\epsilon,k}=\partial_1 F_0\cdot \partial_i F_0 +\epsilon^{k+1}f_{1i}^{\epsilon,k}\hspace{0.5cm}\text{f\"ur }2\leq i\leq n\\
&\partial_i F_{\epsilon,k}\cdot \partial_j F_{\epsilon,k}=\partial_i F_0\cdot \partial_j F_0 +\epsilon^{k+1}f_{ij}^{\epsilon,k}\hspace{0.5cm}\text{f\"ur }2\leq i\leq j\leq n
\end{split}\\
\begin{split}\label{eq:4.2}
&supp(F_{\epsilon,k}-F_0)\in C^{\infty}_0(\mathbb{B},\mathbb{R}^q)
\end{split}\\
\begin{split}\label{eq:4.3}
&\max_{x\in\mathbb{B}}{|F_{\epsilon,k}(x)-F_0(x) |_{\mathbb{R}^q}}\leq C(n,k,F_0,a)\cdot \epsilon
\end{split}
\end{align}

erf\"ullt. Hierbei ist $f^{\epsilon,k}:=(f_{ij}^{\epsilon,k})_{1\leq i\leq j\leq n}\in C^{\infty}_0(\mathbb{B},\mathbb{R}^{\frac{n}{2}(n+1)})$ mit $supp(f^{\epsilon,k})\subseteq K$. Diese Abbildung erf\"ullt die Absch\"atzung:
\begin{equation}\label{eq:6.1}
\left\Vert f^{\epsilon,k} \right\Vert_{C^3(\mathbb{\overline{\mathbb{B}}} ,\mathbb{R}^{\frac{n}{2}(n+1)}   )}\leq C(n,k,F_0,a)\cdot \epsilon^{-3}
\end{equation}
\end{satz}

Um Satz \ref{satz:4.1} zu beweisen, wird eine Vielzahl von technischen \"Uberlegungen ben\"otigt, die im \autoref{sec:4.1} zusammengetragen und bewiesen werden, bevor die Abbildungen $F_{\epsilon,k}$ in \autoref{sec:4.2} induktiv hergeleitet werden.

\section{Hilfsresultate}
\label{sec:4.1}

Gegeben sei eine Menge von Funktionen $\{e_1,..,e_{n+5}\}\subseteq C^{\infty}(\overline{\mathbb{B}},\mathbb{R}^q)$, f\"ur $q\geq n+5$, so dass f\"ur jedes beliebige $x\in\overline{\mathbb{B}}$ die Vektoren $e_1(x),..,e_{n+5}(x)$ linear unabh\"angig sind. Mit \cite[5.4.9]{fischer2008lineare} existiert f\"ur jedes $z\in\mathbb{R}^q$ eine eindeutig bestimmte Orthogonalzerlegung, das hei\ss{}t, es existiert ein $z^{\top}\in lin\{e_1(x),..,e_{n+5}(x)\}$ und ein $z^{\bot}\in lin\{e_1(x),..,e_{n+5}(x)\}^{\bot}$, so dass $z=z^{\top}+z^{\bot}$. Mit:
\begin{align*}
P_x: \mathbb{R}^q\longrightarrow lin\{e_1(x),..,e_{n+5}(x)\}
\end{align*}
wird die Orthogonalprojektion auf den Unterraum $lin\{e_1(x),..,e_{n+5}(x)\}$ bezeichnet, das hei\ss{}t $P_x(z)=z^{\top}$. Das folgende Lemma entspricht \cite[Lemma 1.3.12]{hong2006isometric}, und stellt die Grundlage des gesamten Kapitels dar.

\begin{lem}\label{lem:4.2}
Es sei $\{e_1,..,e_{n+5}\}\subseteq C^{\infty}(\overline{\mathbb{B}},\mathbb{R}^q)$, f\"ur $q\geq n+5$, eine Menge von punktweise linear unabh\"angigen Funktionen, dann existieren Abbildungen $u,v\in C^{\infty}(\overline{\mathbb{B}},\mathbb{R}^q)$, mit der Eigenschaft:
\begin{equation*}
u(x),v(x)\in lin\{e_1(x),...,e_{n+5}(x) \}\hspace{0.5cm}\text{f\"ur alle }x\in\overline{\mathbb{B}}
\end{equation*}

so dass f\"ur alle $x\in\mathbb{\overline{B}}$ die $(n+2)$ Vektoren:
\begin{equation}\label{eq:4.4}
u(x),\, v(x),\, e_i(x)+sP_x\left(a\,\partial_i u(x) + b\,\partial_i v(x) \right) \hspace{0.5cm} i \in \{1,..,n\}
\end{equation}

f\"ur beliebige $s\in\mathbb{R}$, und $a,b\in\mathbb{R}$ mit $a^2+b^2=1$, linear unabh\"angig sind.

\end{lem}

\begin{bew}
Es sei $F\in C^{\infty}(\mathbb{R}^n,\mathbb{R}^{n+1})$ die Immersion aus Lemma \ref{lem:4.1}. F\"ur $\epsilon\in (0,1]$ werden Abbildungen $u, v\in C^{\infty}(\overline{\mathbb{B}},\mathbb{R}^q)$ wie folgt definiert:
\begin{align}\label{eq:4.5}
\begin{split}
&u(x):=\sum_{l=1}^{n+1}{\epsilon F_l(\epsilon^{-2}x)}\,e_{l+1}(x)+e_{n+4}(x)\\
&v(x):=\sum_{l=1}^{n+1}{\epsilon F_l(\epsilon^{-2}x)}\,e_{l+2}(x)+e_{n+5}(x)
\end{split}
\end{align}

Die Behauptung wird bewiesen, indem gezeigt wird, dass es ein $\epsilon\in(0,1]$ gibt, so dass die Abbildungen $u$ und $v$ die geforderten Eigenschaften \eqref{eq:4.4} erf\"ullen. Zun\"achst wird bemerkt, dass f\"ur jedes $x\in\overline{\mathbb{B}}$, und $a,b\in\mathbb{R}$ mit $a^2+b^2=1$, die glatten Abbildungen:
\begin{align}\label{eq:4.6}
w_l(x):= 
\begin{cases}
a\, e_{l+1}(x)+b\, e_{l+2}(x) &\text{f\"ur }l\in\{1,...,n+1\}\\
e_{l+2}(x) &\text{f\"ur }l\in\{n+2,...,n+3\}
\end{cases}
\end{align}

punktweise linear unabh\"angig sind. Seien dazu $x\in\overline{\mathbb{B}}$ und $c_1,...,c_{n+3}\in\mathbb{R}$, so dass:
\begin{equation*}
\sum_{l=1}^{n+3}{c_l\, w_l(x)}=0
\end{equation*}

 dann muss, wegen der punktweisen linearen Unabh\"angigkeit der Funktionen $e_1,...,e_{n+5}$, bereits $c_{n+2}=0=c_{n+3}$ gelten. Das bedeutet:
\begin{equation*}
c_1a\,e_2(x)+(c_1b+c_2a)\, e_3(x)+...+(c_nb+c_{n+1}a)\, e_{n+2}(x)+c_{n+1}b\, e_{n+3}(x)=0
\end{equation*}

beziehungsweise:
\begin{equation*}
c_1a=0,\hspace{0.5cm} c_1b+c_2a=0,\hspace{0.5cm}...,\hspace{0.5cm} c_nb+c_{n+1}a=0,\hspace{0.5cm} c_{n+1}b=0
\end{equation*}

Da nach Voraussetzung $a^2+b^2=1$ gilt, folgt $c_l=0$ f\"ur alle $l\in\{1,...,n+1\}$, womit die punktweise lineare Unabh\"angigkeit der Funktionen in \eqref{eq:4.6} gezeigt ist. Definiere nun Abbildungen $\widehat{w}_1,...,\widehat{w}_{n+2}\in C^{\infty}(\mathbb{R}^n\times\overline{\mathbb{B}},\mathbb{R}^q)$ wie folgt:
\begin{equation}\label{eq:4.7}
\widehat{w}_i(y,x):=
\begin{cases}
\sum_{l=1}^{n+1}{\partial_i F_l(y)\, w_l(x)} &\text{f\"ur }i\in\{1,...,n\}\\
w_{i+1}(x) &\text{f\"ur }i\in\{n+1,n+2\}
\end{cases}
\end{equation}

Es wird gezeigt, dass die Abbildungen $\widehat{w}_1,...,\widehat{w}_{n+2}$ punktweise linear unabh\"angig sind. F\"ur $y\in\mathbb{R}^n$ und ein $x\in\overline{\mathbb{B}}$ seien $\mu_1,...,\mu_{n+2}$, so dass:
\begin{equation*}
\sum_{i=1}^n{\mu_i\sum_{l=1}^{n+1}{\partial_i F_l(y)\, w_l(x)}}+\sum_{j=1}^2{\mu_{n+j}\, w_{n+j+1}(x)}=0
\end{equation*}

gilt, dann folgt, aus der linearen Unabh\"angigkeit der Vektoren $w_1(x),...,w_{n+3}(x)$, bereits $\mu_{n+1}=0=\mu_{n+2}$. Dies impliziert:
\begin{equation*}
\sum_{i=1}^{n}{\mu_i\, \partial_i F_l(y)}=0\hspace{0.5cm}\text{f\"ur alle }l\in\{1,...,n+1\}
\end{equation*}

beziehungsweise:
\begin{equation*}
\sum_{i=1}^{n}{\mu_i\, \partial_i F(y)}=0
\end{equation*}

Da $F$ eine Immersion ist, folgt $\mu_i=0$ f\"ur alle $i\in\{1,...n\}$. Es sei nun $\widehat{W}\in C^{\infty}(\mathbb{R}^n\times\overline{\mathbb{B}},\mathbb{R}^{n+2\times q})$ die matrixwertige Funktion, deren Zeilen die Funktionen $\widehat{w}_1^{T},...,\widehat{w}_{n+2}^{T}$ sind. Dann existiert, aufgrund der Kompaktheit der Menge $\overline{\mathbb{B}}$, und der Periodizit\"at der Abbildung $F$, ein $\theta\in\mathbb{R}_{>0}$, so dass:
\begin{equation}\label{eq:4.10}
det(\widehat{W}(y,x)\cdot \widehat{W}^T(y,x))\geq \theta
\end{equation}

f\"ur alle $y\in\mathbb{R}^n$ und $x\in\overline{\mathbb{B}}$ gilt. Im Folgenden seien $\epsilon\in(0,1]$ und $a,b\in\mathbb{R}$, mit $a^2+b^2=1$, sowie $s\in\mathbb{R}$, dann werden Funktionen $v_1,...,v_{n+2}\in C^{\infty}(\overline{\mathbb{B}},\mathbb{R}^{q})$ wie folgt definiert:
\begin{equation}\label{eq:4.13}
v_i(x):=
\begin{cases}
e_i(x)+sP_x(a\, \partial_i u(x)+b\,\partial_i v(x)) &\text{f\"ur }i\in\{1,...,n\}\\
u(x) &\text{f\"ur }i=n+1\\
v(x) &\text{f\"ur }i=n+2\\
\end{cases}
\end{equation}

Dann ist f\"ur $i\in\{1,...,n\}$:
\begin{align}\label{eq:4.8}
\begin{split}
v_i(x)\stackrel{\hphantom{\eqref{eq:4.5}}}{=}&e_i(x)+sP_x\left(a\,\partial_i u(x)+b\,\partial_i v(x) \right) \\
\stackrel{\eqref{eq:4.5}}{=}&e_i(x)+sP_x\left(\sum_{l=1}^{n+1}\frac{1}{\epsilon}\partial_i F_l(\epsilon^{-2}x)(a\, e_{l+1}(x)+b\, e_{l+2}(x)) \right)\\
\hphantom{\stackrel{\eqref{eq:4.5}}{=}}&+sP_x\left[\sum_{l=1}^{n+1}\epsilon F_l(\epsilon^{-2}x)\left(a \,\partial_i e_{l+1}(x)+b\, \partial_i e_{l+2}(x)\right) \right]\\
\hphantom{\stackrel{\eqref{eq:4.5}}{=}}&+sP_x\left[\partial_i e_{n+4}(x)+\partial_i e_{n+5}(x)\right]\\
\stackrel{\hphantom{\eqref{eq:4.5}}}{=}&e_i(x)+\frac{s}{\epsilon}\sum_{l=1}^{n+1}\partial_i F_l(x)(\epsilon^{-2}x)(a\, e_{l+1}(x)+b\, e_{l+2}(x))+s\mathcal{O}(1)\\
\stackrel{\eqref{eq:4.7}}{=}&e_i(x)+\frac{s}{\epsilon}\widehat{w}(\epsilon^{-2}x,x)+s\mathcal{O}(1)
\end{split}
\end{align}

ferner ist mit \eqref{eq:4.5}:
\begin{equation}\label{eq:4.9}
v_{n+1}=u(x)=e_{n+4}(x)+\mathcal{O}(\epsilon)\hspace{0.5cm}\text{ und }\hspace{0.5cm}v_{n+2}=v(x)=e_{n+5}(x)+\mathcal{O}(\epsilon)
\end{equation}

Zun\"achst wird angenommen, dass $|s|\geq \epsilon^{1-\frac{1}{2n}}$ gilt. Dann folgt aus \eqref{eq:4.8}, f\"ur $i\in\{1,...,n\}$:
\begin{equation}\label{eq:4.11}
v_i(x)=\frac{s}{\epsilon}\left[\epsilon\, e_i(x)+\widehat{w}(\epsilon^{-2}x,x)+\mathcal{O}(\epsilon)\right]=\frac{s}{\epsilon}\left[\widehat{w}_i(\epsilon^{-2}x,x)+\mathcal{O}(\epsilon^{\frac{1}{2n}})\right]
\end{equation}

und aus \eqref{eq:4.9}:
\begin{equation}\label{eq:4.12}
v_{n+1}=e_{n+4}(x)+\mathcal{O}(\epsilon^{\frac{1}{2n}})\hspace{0.5cm}\text{ und }\hspace{0.5cm}v_{n+2}=e_{n+5}(x)+\mathcal{O}(\epsilon^{\frac{1}{2n}})
\end{equation}

Nun sei $V\in C^{\infty}(\overline{\mathbb{B}},\mathbb{R}^{n+2\times q})$ die matrixwertige Funktion, deren Zeilen die Funktionen $v_1^{T},...,v_{n+2}^{T}$ sind. Aus \eqref{eq:4.11} und \eqref{eq:4.12} folgt, unter Beachtung von $|s|\geq \epsilon^{1-\frac{1}{2n}}$, f\"ur alle $x\in\mathbb{B}$:

\begin{align*}
det(V(x)\cdot V^T(x))=&\frac{s^{2n}}{\epsilon^{2n}}\left[det(\widehat{W}(\epsilon^{-2}x,x)\cdot \widehat{W}^T(\epsilon^{-2}x,x))+\mathcal{O}(\epsilon^{\frac{1}{2n}}) \right]\stackrel{\eqref{eq:4.10}}{\geq}\frac{s^{2n}}{\epsilon^{2n}}\left(\theta-C\epsilon^{\frac{1}{2n}} \right)\\
\geq& \frac{1}{\epsilon}\left(\theta-C\epsilon^{\frac{1}{2n}} \right)\geq \frac{\theta}{2\epsilon}
\end{align*}

falls $\epsilon\in (0,\epsilon_0]$, f\"ur ein hinreichend kleines $\epsilon_0\in\mathbb{R}_{>0}$. Damit ist gezeigt, dass die in \eqref{eq:4.13} definierten Funktionen, f\"ur $|s|\geq \epsilon^{1-\frac{1}{2n}}$ und alle $x\in\overline{\mathbb{B}}$ punktweise linear unabh\"angig sind, sofern $\epsilon\in\mathbb{R}_{>0}$ klein genug ist. Nun sei $|s|< \epsilon^{1-\frac{1}{2n}}$, und f\"ur ein $x\in\overline{\mathbb{B}}$, seien $c_1,...,c_{n+2}\in\mathbb{R}$, so dass $\sum_{i=1}^{n+2}{c_i\, v_i(x)}=0$ gilt. Setze: 
\begin{equation}\label{eq:4.17}
A:=\max_{1\leq i\leq n+2}{|c_i|}
\end{equation}

Es wird gezeigt, dass $A=0$ f\"ur hinreichend kleine $\epsilon\in\mathbb{R}_{>0}$ gelten muss. Aus $\eqref{eq:4.8}$ folgt f\"ur $i\in\{1,...,n\}$:
\begin{align}\label{eq:4.14}
\begin{split}
v_i(x)\stackrel{\hphantom{\eqref{eq:4.5}}}{=}&e_i(x)+\frac{s}{\epsilon}\widehat{w}_i(\epsilon^{-2}x,x)+\mathcal{O}(\epsilon^{1-\frac{1}{2n}})\\
\stackrel{\eqref{eq:4.5}}{=}&e_i(x)+\frac{s}{\epsilon}\sum_{l=1}^{n+1}\partial_i F_l(x)(\epsilon^{-2}x)(a\, e_{l+1}(x)+b\, e_{l+2}(x))+\mathcal{O}(\epsilon^{1-\frac{1}{2n}})
\end{split}
\end{align}

und aus $\eqref{eq:4.9}$:
\begin{equation}\label{eq:4.15}
v_{n+1}=e_{n+4}(x)+\mathcal{O}(\epsilon^{1-\frac{1}{2n}})\hspace{0.5cm}\text{ und }\hspace{0.5cm}v_{n+2}=e_{n+5}(x)+\mathcal{O}(\epsilon^{1-\frac{1}{2n}})\\
\end{equation}

Mit \eqref{eq:4.14} und \eqref{eq:4.15} folgt aus $\sum_{i=1}^{n+2}{c_i\, v_i(x)}=0$:
\begin{align*}
0\stackrel{\hphantom{\eqref{eq:4.16},\,\eqref{eq:4.17}}}{=}&\sum_{i=1}^n{c_i\left[e_i(x)+\frac{s}{\epsilon}\sum_{l=1}^{n+1}\partial_i F_l(\epsilon^{-2}x)\,(a\, e_{l+1}(x)+b\, e_{l+2}(x))+\mathcal{O}(\epsilon^{1-\frac{1}{2n}}) \right]}\\
&+c_{n+1}\left[e_{n+4}(x)+\mathcal{O}(\epsilon^{1-\frac{1}{2n}})\right]+c_{n+2}\left[e_{n+5}(x)+\mathcal{O}(\epsilon^{1-\frac{1}{2n}})\right]\\
\stackrel{\eqref{eq:4.16},\,\eqref{eq:4.17}}{=}&\sum_{i=1}^n{c_i\, e_i(x)}+\frac{s}{\epsilon}\sum_{l=1}^{n+1}{\sum_{i=1}^{min(l,n)}{c_i\left[\partial_i F_l(\epsilon^{-2}x)\,(a\, e_{l+1}(x)+b\, e_{l+2}(x))\right]}}\\
&+c_{n+1}\,e_{n+4}(x)+c_{n+2}\,e_{n+5}(x)+\mathcal{O}(A\epsilon^{1-\frac{1}{2n}})\\
\stackrel{\hphantom{\eqref{eq:4.16},\,\eqref{eq:4.17}}}{=}&\sum_{i=1}^n{c_i\, e_i(x)}+\frac{s}{\epsilon}\sum_{l=2}^{n+2}{\sum_{i=1}^{min(l-1,n)}{c_i\,\partial_i F_{l-1}(\epsilon^{-2}x)\,a\, e_{l}(x)}}\\
&+\frac{s}{\epsilon}\sum_{l=3}^{n+3}{\sum_{i=1}^{min(l-2,n)}{c_i\,\partial_i F_{l-2}(\epsilon^{-2}x)\,b\, e_{l}(x)}}\\
&+c_{n+1}\,e_{n+4}(x)+c_{n+2}\, e_{n+5}(x)+\mathcal{O}(A\epsilon^{1-\frac{1}{2n}})\\
\stackrel{\hphantom{\eqref{eq:4.16},\,\eqref{eq:4.17}}}{=}&\sum_{i=1}^n{c_i\, e_i(x)}+\frac{s}{\epsilon}\sum_{i=2}^{n+2}{\sum_{l=1}^{min(i-1,n)}{c_l\,\partial_l F_{i-1}(\epsilon^{-2}x)\,a\, e_{i}(x)}}\\
&+\frac{s}{\epsilon}\sum_{i=3}^{n+3}{\sum_{l=1}^{min(i-1,n)}{c_l\,\partial_l F_{i-2}(\epsilon^{-2}x)\,b\, e_{i}(x)}}\\
&+c_{n+1}\,e_{n+4}(x)+c_{n+2}\, e_{n+5}(x)+\mathcal{O}(A\epsilon^{1-\frac{1}{2n}})
\end{align*}

Aus der punktweisen linearen Unabh\"angigkeit der Abbildungen $e_1,....,e_{n+5}$ folgt:
\begin{align}\label{eq:4.18}
\begin{split}
&c_1+\mathcal{O}(A\epsilon^{1-\frac{1}{2n}})=0\\
&c_2+\frac{s}{\epsilon}c_1\, a\, \partial_1 F_1(\epsilon^{-2}x)+\mathcal{O}(A\epsilon^{1-\frac{1}{2n}})=0
\end{split}
\end{align}

sowie f\"ur $i\in\{3,...,n\}$:
\begin{equation}\label{eq:4.19}
c_i+\frac{s}{\epsilon}\sum_{l=1}^{i-1}{c_l\left[a\, \partial_l F_{i-1}(\epsilon^{-2}x)+b\, \partial_l F_{i-2}(\epsilon^{-2}x) \right]}=0
\end{equation}

und:
\begin{align}\label{eq:4.20}
\begin{split}
&c_{n+1}+\mathcal{O}(A\epsilon^{1-\frac{1}{2n}})=0\\
&c_{n+2}+\mathcal{O}(A\epsilon^{1-\frac{1}{2n}})=0\\
\end{split}
\end{align}

Aus \eqref{eq:4.18} folgt:
\begin{align*}
|c_1|\leq C\,A\epsilon^{1-\frac{1}{2n}}\\
|c_2|\leq C\,A\epsilon^{1-\frac{2}{2n}}
\end{align*}

woraus, mit \eqref{eq:4.19}, induktiv f\"ur $i\in\{3,...,n\}$ die Absch\"atzung:
\begin{align*}
|c_i|\leq C\,A\epsilon^{1-\frac{i}{2n}}
\end{align*}

hervorgeht. Ferner impliziert \eqref{eq:4.20}:
\begin{equation*}
|c_{n+1}|+|c_{n+2}|\leq C\,A\epsilon^{1-\frac{1}{2n}}\\
\end{equation*}

f\"ur eine bestimmte, von $x$ unabh\"angige Konstante, $C\in\mathbb{R}_{>0}$. Daraus folgt insgesamt $|c_i|\leq C\, A\sqrt{\epsilon}$ f\"ur alle $i\in\{1,...,n+2\}$, und schlie\ss{}lich, unter Beachtung von \eqref{eq:4.17}:
\begin{equation*}
A=\max_{1\leq i\leq n+2}{|c_i|}\leq C\,A\sqrt{\epsilon}\leq \frac{1}{2}A
\end{equation*}

falls $\epsilon\in(0,\epsilon_0]$, mit $C\sqrt{\epsilon_0}\leq \frac{1}{2}$. Daraus folgt die punktweise lineare Unabh\"angigkeit der Funktionen $v_1,...,v_{n+2}$, f\"ur $|s|< \epsilon^{1-\frac{1}{2n}}$, und hinreichend kleine $\epsilon>0$. Aus der Definition der Funktionen $v_1,...,v_{n+2}$ in \eqref{eq:4.13} folgt die Behauptung.

\end{bew}

\begin{bem}\label{bem:4.1}
Es kann angenommen werden, dass die in Lemma \ref{lem:4.2} konstruierten Funktionen $u$ und $v$ punktweise orthonormal zueinander sind.
\end{bem}

\begin{bew}
Es wird gezeigt, dass eine Orthonormalisierung der Funktionen $u$ und $v$, aus Lemma \ref{lem:4.2}, die lineare Unabh\"angigkeit der Vektoren in \eqref{eq:4.4}, f\"ur beliebige $x\in\overline{\mathbb{B}}$, nicht beeintr\"achtigt. Dazu werden punktweise orthonormalisierte Funktionen $\widehat{u},\widehat{v}\in C^{\infty}(\overline{\mathbb{B}})$, gem\"a\ss{} \cite[5.4.9]{fischer2008lineare}, wie folgt definiert: F\"ur $x\in\mathbb{B}$ ist, unter Beachtung der punktweisen linearen Unabh\"angigkeit der Funktionen $u$ und $v$:
\begin{align*}
\widehat{u}(x)&:=\frac{u(x)}{|u(x)|}=C_1(u(x))\cdot u(x) \\
\widehat{v}(x)&:=\frac{v(x)-\left(\frac{u(x)}{|u(x)|}\cdot v(x)\right)\frac{u(x)}{|u(x)|}}{\left|v(x)-\left(\frac{u(x)}{|u(x)|}\cdot v(x)\right)\frac{u(x)}{|u(x)|}\right|}=C_2(u(x),v(x))\cdot u(x)+C_3(u(x),v(x))\cdot v(x)
\end{align*}

Dann gilt f\"ur alle $i\in\{1,..,n\}$, $x\in\overline{\mathbb{B}}$, $s\in\mathbb{R}$ und $a,b\in\mathbb{R}$ mit $a^2+b^2=1$:
\begin{align*}
&e_i(x)+sP_x\left(a\,\partial_i \widehat{u}(x)+b\partial_i \widehat{v}(x) \right)\\
=&e_i(x)+sP_x\left[a\, C_1(u(x))\cdot\partial_i u(x)+b\, C_2(u(x),v(x))\cdot\partial_i u(x)+b \, C_3(u(x),v(x))\cdot\partial_i v(x) \right]\\
&+ sP_x\left[a\,C_4(u(x),i)\cdot u(x)+b\, C_5(u(x),v(x),i)\cdot u(x)+b \, C_6(u(x),v(x),i)\cdot v(x) \right]\\
=&e_i(x)+sP_x\left[C_7(u(x),v(x),a,b)\cdot\partial_i u(x)+C_8(u(x),v(x),b)\cdot\partial_i v(x) \right]\\
&+ C_9(u(x),v(x),i,s,a,b)\cdot u(x)+ C_{10}(u(x),v(x),i,s,b)\cdot v(x) \\
=&e_i(x)\\
&+C_{11}(u(x),v(x),a,b)\cdot P_x\left[C_{12}(u(x),v(x),a,b)\cdot\partial_i u(x)+C_{13}(u(x),v(x),a,b)\cdot\partial_i v(x) \right]\\
&+ C_9(u(x),v(x),i,s,a,b)\ u(x)+ \ C_{10}(u(x),v(x),i,s,b)\ v(x) 
\end{align*}

mit $C^2_{12}(u(x),v(x),a,b)+C^2_{13}(u(x),v(x),a,b)=1$. Mit linearen Unabh\"angigkeit der Vektoren in \eqref{eq:4.4} folgt die Behauptung.

\end{bew}

F\"ur weitere Konstruktionen ist das folgende Hilfsresultat wichtig:

\begin{lem}\label{lem:4.4}
Gegeben seien punktweise linear unabh\"angige Abbildungen $e_1,...,e_k\in\linebreak C^{\infty}(\overline{\mathbb{B}},\mathbb{R}^q)$, wobei $q\geq k+1$ gilt. Dann existieren Abbildungen $e_{k+1},...,e_q\in C^{\infty}(\overline{\mathbb{B}},\mathbb{R}^q)$, so dass die Abbildungen $e_1,...,e_q$ in $\overline{\mathbb{B}}$ punktweise linear unabh\"angig sind.
\end{lem}

\begin{bew}
Nach dem Fortsetzungssatz von Seeley \cite{seeley1964extension} existiert ein $\epsilon\in\mathbb{R}_{>0}$, so dass Abbildungen $\widehat{e}_1,...,\widehat{e}_k\in C^{\infty}(B_{1+\epsilon}(0),\mathbb{R}^q)$ mit der Eigenschaft existieren, dass $\left.\widehat{e}_i\right|_{\mathbb{B}}\equiv e_i$ f\"ur alle $i\in\{1,...,k\}$ gilt. Wegen der Stetigkeit der Determinante, kann angenommen werden, dass $\widehat{e}_1,...,\widehat{e}_k$ auf $B_{1+\epsilon}(0)$ punktweise linear unabh\"angig sind. Die Menge:
\begin{equation*}
N:=\coprod_{x\in B_{1+\epsilon}(0)}{lin\{\widehat{e}_1(x),...,\widehat{e}_k(x)\}^{\bot}}
\end{equation*}

kann, mit \cite[Lemma 8.41]{lee2003introduction}, als glattes Vektorb\"undel \"uber $B_{1+\epsilon}(0)$ aufgefasst werden. Mit \cite[Chapter 4, 2.5. Corollary]{hirsch1997differential} ist das Vektorb\"undel $N$ glatt trivialisierbar. Mit \cite[Corollary 5.11.]{lee2003introduction} folgt die Behauptung.

\end{bew}

Das folgende Lemma basiert auf \cite[Lemma 1.3.12.]{hong2006isometric}, beziehungsweise \cite[Satz 3.3.]{gunther1989einbettungssatz}, und bildet, aufbauend auf Lemma \ref{lem:4.2}, ebenfalls eine wichtige Grundlage f\"ur weitere Konstruktionen, hierbei ist $\gls{mathbb{S}}:= [-\pi,\pi]$.

\begin{lem}\label{lem:4.5}
Es sei $F_0\in C^{\infty}(\overline{\mathbb{B}},\mathbb{R}^q)$ eine freie Abbildung, wobei $q\geq \frac{n}{2}(n+3)+5$ gilt. Dann existiert eine Funktion $v_1\in C^{\infty}(\mathbb{S}\times\overline{\mathbb{B}},\mathbb{R}^q)$, so dass f\"ur ein beliebiges $a\in C^{\infty}_0(\mathbb{B})$ die Funktion $u_1:= a^2\,v_1 \in C^{\infty}(\mathbb{S}\times\overline{\mathbb{B}},\mathbb{R}^q)$ die folgenden Eigenschaften erf\"ullt:
\begin{enumerate}[(i)]
\begin{item}
Es gilt im punktweisen Sinne in $\mathbb{S}\times\overline{\mathbb{B}}$:
\begin{align}\label{eq:4.27}
\begin{split}
&\partial_i F_0\cdot u_1=0\hspace{0.5cm}\text{und}\hspace{0.5cm}\partial_i F_0\cdot \partial_t u_1=0\hspace{0.5cm}\text{f\"ur }1\leq i \leq n\\
&\partial_i\partial_j F_0\cdot u_1=0\hspace{0.5cm}\text{f\"ur }2\leq i\leq j\leq n
\end{split}
\end{align}
\end{item}
\begin{item}
Die Funktionen:
\begin{align}\label{eq:4.28}
\begin{array}{lll}
\partial_1 F_0+s\partial_t u_1 & \partial_i F_0\hspace{0.5cm}\text{f\"ur }2\leq i \leq n  & \partial_{ij}F_0\hspace{0.5cm}\text{f\"ur }2\leq i\leq j\leq n \\
\partial_1^2 F_0+2s\,\partial_t\partial_1 u_1  & \partial_1 \partial_i F_0+s\, \partial_t\partial_i u_1\hspace{0.5cm}\text{f\"ur }2\leq i \leq n &\\
\partial_t v_1 & \partial_t^2 v_1 &
\end{array}
\end{align}
sind im punktweisen Sinne in $\mathbb{S}\times\overline{\mathbb{B}}$, f\"ur alle $s\in\mathbb{R}$, linear unabh\"angig.
\end{item}
\begin{item}
Es gilt, f\"ur eine positive Funktion $\varrho\in C^{\infty}(\mathbb{S})$:
\begin{equation}\label{eq:4.29}
|\partial_t u_1(t,x)|_{\mathbb{R}^q}^2=a^4(x)\,\varrho^2(t)
\end{equation}
f\"ur alle $(t,x)\in \mathbb{S}\times\overline{\mathbb{B}}$.
\end{item}
\begin{item}
F\"ur jedes $x\in\overline{\mathbb{B}}$ gelten die Gleichungen:
\begin{align}\label{eq:4.35}
\begin{split}
&\int_{-\pi}^{\pi}{\partial_t u_1(t,x)\cdot \partial_i u_1(t,x)\> dt}=0\hspace{0.5cm}\text{f\"ur }1\leq i\leq n\\
&\int_{-\pi}^{\pi}{\varrho(t)\,\partial_i F_0(x)\cdot \partial_j u_1(t,x)\> dt}=0\hspace{0.5cm}\text{f\"ur }1\leq i, j\leq n
\end{split}
\end{align}
\end{item}
\end{enumerate}

\end{lem}

\begin{bew}
Da $F_0$ eine freie Abbildung ist, und $q\geq \frac{n}{2}(n+3)+5$ gilt, existieren mit Lemma \ref{lem:4.4} Abbildungen $f_1,...,f_5\in C^{\infty}(\overline{\mathbb{B}},\mathbb{R}^q)$, so dass die Abbildungen:
\begin{equation}\label{eq:4.25}
\partial_i F_0\hspace{0.5cm}\text{f\"ur }{1\leq i\leq n},\hspace{1cm}\partial_k\partial_l F_0\hspace{0.5cm}\text{f\"ur }{1\leq k\leq l\leq n},\hspace{1cm}f_1,...,f_5
\end{equation}

in $\overline{\mathbb{B}}$ punktweise linear unabh\"angig sind. F\"ur jedes $x\in\overline{\mathbb{B}}$ ist:
\begin{equation*}
L_x:= lin\left(\{\partial_i F_0 (x)\}_{1\leq i\leq n}\cup  \{\partial_i\partial_j F_0 (x)\}_{2\leq i\leq j\leq n}\right)\subseteq \mathbb{R}^q
\end{equation*}

und $P_x\in C^{\infty}(\mathbb{R}^q,L_x^{\bot})$ die Orthogonalprojektion auf den Raum $L_x^{\bot}$. Weiterhin seien $e_1,....,e_{n+5}\in C^{\infty}(\overline{\mathbb{B}},\mathbb{R}^q)$  wie folgt definiert:
\begin{align}\label{eq:4.31}
\begin{split}
e_i(x):=
\begin{cases}
\frac{1}{2}P_x \partial_1^2 F_0(x) &\text{f\"ur }i=1\\
P_x \partial_1 \partial_i F_0(x) &\text{f\"ur }i\in\{1,...,n\}\\
P_x f_{i-n}(x) &\text{f\"ur }i\in\{n+1,...,n+5\}
\end{cases}
\end{split}
\end{align}

Wegen der punktweisen linearen Unabh\"angigkeit der Abbildungen in \eqref{eq:4.25}, sind die Abbildungen $e_1,....,e_{n+5}$ ebenfalls punktweise linear unabh\"angig in $\overline{\mathbb{B}}$. Mit Lemma \ref{lem:4.2} existieren $u,v\in C^{\infty}(\overline{\mathbb{B}},\mathbb{R}^q)$, mit der Eigenschaft $u(x),v(x)\in lin\{e_1(x),...,e_{n+5}(x) \}$ f\"ur alle $x\in \overline{\mathbb{B}}$, so dass f\"ur alle $x\in\mathbb{\overline{B}}$ die $(n+2)$ Vektoren:
\begin{equation}\label{eq:4.32}
u(x),\, v(x),\, e_i(x)+sP_x\left(a\,\partial_i u(x) + b\,\partial_i v(x) \right) \hspace{0.5cm} i \in \{1,..,n\}
\end{equation}

f\"ur beliebige $s\in\mathbb{R}$, und $a,b\in\mathbb{R}$ mit $a^2+b^2=1$, linear unabh\"angig sind. Mit Bemerkung \ref{bem:4.1} kann angenommen werden, dass:
\begin{equation}\label{eq:4.26}
|u(x)|=1=|v(x)|\hspace{0.5cm}\text{und}\hspace{0.5cm}u(x)\cdot v(x)=0
\end{equation}

f\"ur alle $x\in\overline{\mathbb{B}}$ gilt. Definiere nun, mit den in Lemma \ref{lem:4.3} beschriebenen Funktionen $\alpha_1,\alpha_2\in C^{\infty}(\mathbb{S})$, die Abbildung $v_1\in C^{\infty}(\mathbb{S}\times\overline{\mathbb{B}},\mathbb{R}^q)$ mit:
\begin{equation}\label{eq:4.30}
v_1(t,x):= \alpha_1(t)\, u(x)+\alpha_2(t)\, v(x)
\end{equation}

Ist $a\in C^{\infty}_0(\mathbb{B})$ eine beliebige Abbildung, dann gilt f\"ur $u_1:= a^2\, v_1$, unter Beachtung von \eqref{eq:4.26}, die Gleichung:
\begin{equation*}
|\partial_t u_1(t,x)|^2=a^4\, (\alpha_1'(t)^2+\alpha_2'(t)^2)=a^4\, \varrho^2(t)
\end{equation*}

wobei $\varrho\in C^{\infty}(\mathbb{S})$, mit:
\begin{equation}\label{eq:4.34}
\varrho(t)=\sqrt{\alpha_1'(t)^2+\alpha_2'(t)^2}
\end{equation}

definiert ist. Dies entspricht \eqref{eq:4.29}. Da $u(x),v(x)\in lin\{e_1(x),...,e_{n+5}(x) \}\subseteq L_x^{\bot}$ f\"ur alle $x\in \overline{\mathbb{B}}$ gilt, ist \eqref{eq:4.27} bewiesen. Aus \eqref{eq:4.22} folgt, mit \eqref{eq:4.30}, die punktweise lineare Unabh\"angigkeit der Abbildungen $\partial_t v_1$ und $\partial_t^2 v_1$, sowie:
\begin{equation}\label{eq:4.33}
lin\{\partial_t v_1(t,x), \partial_t^2 v_1(t,x)\}=lin\{u(x), v(x)\}
\end{equation}

f\"ur alle $x\in\overline{\mathbb{B}}$. Ferner ist mit \eqref{eq:4.31} und \eqref{eq:4.30}, f\"ur alle $x\in\overline{\mathbb{B}}$ und $s\in\mathbb{R}$:
\begin{equation*}
P_x(\partial_1^2 F_0(x)+2s \partial_t \partial_1 v_1(t,x))=2e_i(x)+2s\, P_x(\alpha_1'(t)\, \partial_1 u(x)+\alpha_2'(t)\, \partial_1 v(x))
\end{equation*}

und f\"ur $i\in\{2,...,n\}$:
\begin{equation*}
P_x(\partial_1 \partial_i F_0(x)+s \partial_t \partial_i v_1(t,x))=e_i(x)+s\,  P_x(\alpha_1'(t)\, \partial_i u(x)+\alpha_2'(t)\, \partial_i v(x))
\end{equation*}

Aus der linearen Unabh\"angigkeit der Vektoren in \eqref{eq:4.32}, folgt mit \eqref{eq:4.33} die punktweise lineare Unabh\"angigkeit der Abbildungen in \eqref{eq:4.28}. Schlie\ss{}lich folgen aus \eqref{eq:4.24}, zusammen mit \eqref{eq:4.30}, die Identit\"aten in \eqref{eq:4.35}.

\end{bew}

Ist eine Matrix $A\in\mathbb{R}^{k\times q}$ gegeben, deren Zeilen linear unabh\"angig sind, dann l\"asst sich, f\"ur jedes $b\in \mathbb{R}^k$, das Gleichungssystem $Ax=b$ l\"osen. Falls $q>k$ gilt, so ist die L\"osung dieses Gleichungssystems nicht eindeutig bestimmt. In \autoref{sec:4.2} sollen Gleichungssysteme gel\"ost werden, bei denen die Koeffizienten der Matrix glatte Funktionen auf $\overline{\mathbb{B}}$ sind. Eine vektorwertige L\"osungsfunktion soll wieder eine glatte Funktion auf $\overline{\mathbb{B}}$ sein. Um das zu gew\"ahrleisten, wird gezeigt, dass im L\"osungsraum $L(A,b):=\{x\in\mathbb{R}^q: Ax=b \}$ genau ein Element minimaler euklidischer L\"ange existiert.

\begin{lem}\label{lem:4.6}
Gegeben sei eine Matrix $A\in\mathbb{R}^{k\times q}$ mit $rg(A)=k$, und ein Vektor $b\in \mathbb{R}^k$. Dann existiert genau ein $x_0\in L(A,b)$, so dass:
\begin{equation*}
\left|x_0 \right|_{\mathbb{R}^q}=\min\{\left|x\right|_{\mathbb{R}^q} : x\in L(A,b) \}
\end{equation*}

gilt. Diese L\"osung $x_0$ hat die Darstellung:
\begin{equation}\label{eq:4.36}
x_0= A^{T}(AA^{T})^{-1}b
\end{equation}

\end{lem}

\begin{bew}
Ein $x\in L(A,b)$ besitzt, mit \cite[5.4.9.]{fischer2008lineare}, eine eindeutige Orthogonalzerlegung $x=x^{\top}+x^{\bot}$, wobei $x^{\top}\in lin\{A^{T}\}$, und $x^{\bot}\in lin\{A^{T}\}^{\bot}$ erf\"ullt ist. Dann gilt:
\begin{equation}\label{eq:4.37}
\left|x \right|_{\mathbb{R}^q}= \left|x^{\top} \right|_{\mathbb{R}^q}+\left|x^{\bot} \right|_{\mathbb{R}^q}\geq \left|x^{\top} \right|_{\mathbb{R}^q}
\end{equation}

Angenommen, es existiert ein $\widehat{x}_0\in L(A,b)$, so dass $\left|\widehat{x}_0 \right|_{\mathbb{R}^q}=\min\{\left|x\right|_{\mathbb{R}^q} : x\in L(A,b) \}$ erf\"ullt ist, dann gilt, mit \eqref{eq:4.37}, die Gleichung $\widehat{x}_0=A^{T} y$, f\"ur ein $y\in\mathbb{R}^k$. Dann ist $b=A\widehat{x}_0=AA^{T} y$, woraus wegen $rg(AA^{T})=k$, die Gleichheit $y=(AA^{T})^{-1} b$ folgt. Also ist $\widehat{x}_0=x_0$, mit dem in \eqref{eq:4.36} definierten, $x_0$. Da $x_0\in L(A,b)$ gilt, folgt die Behauptung.

\end{bew}

Die Matrix $AA^T$ wird im Allgemeinen auch als \textbf{Gramsche Matrix}\index{Gramsche>$\sim$ Matrix} bezeichnet. In diesem Zusammenhang wird $det(AA^T)\in\mathbb{R}_{\geq 0}$ als \textbf{Gramsche Determinante}\index{Gramsche>$\sim$ Determinante} bezeichnet. Es gilt: $det(AA^T)=0$, genau dann wenn $rg(A)<k$ ist. Die Gramsche Matrix einer endlichen Menge von Zeilenvektoren gleicher Dimension, wird als die Gramsche Matrix der Matrix, in welcher die Zeilenvektoren, unter Beibehaltung der gegebenen Reihenfolge, stehen, definiert. Lemma \ref{lem:4.6} wird nun verwendet, um folgenden Sachverhalt zu zeigen:

\begin{fol}\label{fol:4.1}
Es sei $\Omega\subseteq\mathbb{R}^n$ eine beschr\"ankte offene Menge. Gegeben seien punktweise linear unabh\"angige Funktionen $e_1,...,e_k\in C^{\infty}(\overline{\Omega},\mathbb{R}^q)$, sowie ein $h\in C^{\infty}(\overline{\Omega},\mathbb{R}^k)$, dann existiert ein $v\in C^{\infty}(\overline{\Omega},\mathbb{R}^q)$, so dass f\"ur alle $i\in \{1,...,k\}$, und f\"ur alle $x\in \overline{\Omega}$ die Gleichung:
\begin{equation}\label{eq:4.38}
e_i(x)\cdot v(x)=h_i(x)
\end{equation}

erf\"ullt ist.
\end{fol}

\begin{bew}
Es sei $A\in C^{\infty}(\overline{\Omega},\mathbb{R}^{k\times q})$ die matrixwertige Funktion, deren Zeilenvektoren, unter Beibehaltung der Reihenfolge, die Funktionen $e_1^T,...,e_k^T$ sind. Daraus wird eine Abbildung $\Theta\in C^{\infty}(\overline{\Omega},\mathbb{R}^{q\times k})$ wie folgt konstruiert:
\begin{equation*}
\Theta(x):= A^{T}(x)\cdot(A(x)A^{T}(x))^{-1}
\end{equation*}

Mit Lemma \ref{lem:4.6} ist die Abbildung $v\in C^{\infty}(\overline{\Omega},\mathbb{R}^q)$ mit der Bildungsvorschrift:
\begin{equation}\label{eq:4.47}
v(x):= A^{T}(x)\cdot (A(x)A^{T}(x))^{-1}h(x)
\end{equation}

eine L\"osung des Gleichungssystems \eqref{eq:4.38}.
\end{bew}

\section{Konstruktion der lokalen Hilfsabbildungen \texorpdfstring{$F_{\epsilon,k}$}{F{eps,k}}}
\label{sec:4.2}

Es sei $\varrho\in C^{\infty}(\mathbb{S})$, die in \eqref{eq:4.29} eingef\"uhrte Funktion, welche in \eqref{eq:4.34} genau definiert worden ist. Diese Funktion wird, mittels periodischer Fortsetzung, als glatte Funktion auf ganz $\mathbb{R}$ aufgefasst. Mit $P\in C^{\infty}(\mathbb{R})$ wird eine beliebige, aber feste, Stammfunktion zu $\varrho$ bezeichnet. Da $\varrho(t)>0$ f\"ur alle $t\in\mathbb{R}$ gilt, ist diese Funktion streng monoton steigend, mit $\lim_{t\to \pm \infty}{P(t)}=\pm \infty$. Demzufolge exstiert eine Funktion $\beta\in C^{\infty}(\mathbb{R})$, so dass $P(\beta(t))=t$ f\"ur alle $t\in\mathbb{R}$ gilt. Mit der Kettenregel gilt dann:

\begin{equation}\label{eq:4.39}
\beta'(t)=\frac{1}{\varrho(\beta(t))}
\end{equation}

Die in diesem Abschnitt zu bestimmenden Hilfsabbildungen $F_{\epsilon,k}\in C^{\infty}(\overline{\mathbb{B}},\mathbb{R}^q)$ werden wie folgt konstruiert: Zun\"achst wird, f\"ur jedes $k\in\mathbb{N}\backslash\{0,1\}$, eine Abbildung $F_k\in C^{\infty}(I\times\mathbb{S}\times\overline{\mathbb{B}},\mathbb{R}^q)$ mit gewissen Eigenschaften konstruiert, wobei $\gls{I}:=[0,1]\subseteq\mathbb{R}$ ist. Anschlie\ss{}end wird, unter Beachtung von Lemma \ref{lem:C.3}, f\"ur $\epsilon\in [0,1]$:
\begin{equation}\label{eq:4.82}
F_{\epsilon,k}(x):= F_k(\epsilon,\beta(\epsilon^{-1}x_1),x)
\end{equation}

gesetzt. Dann gilt:
\begin{align}\label{eq:4.40}
\begin{split}
\partial_1 F_{\epsilon,k}(x)&\stackrel{\hphantom{\eqref{eq:4.39}}}{=}\partial_1 F_k(\epsilon,\beta(\epsilon^{-1}x_1),x)+\frac{\beta'(\epsilon^{-1}x_1)}{\epsilon}\,\partial_t F_k(\epsilon,\beta(\epsilon^{-1}x_1),x)\\
&\stackrel{\eqref{eq:4.39}}{=}\partial_1 F_k(\epsilon,\beta(\epsilon^{-1}x_1),x)+\frac{1}{\epsilon\varrho(\beta(\epsilon^{-1}x_1))}\,\partial_t F_k(\epsilon,\beta(\epsilon^{-1}x_1),x)
\end{split}
\end{align}

wobei $\partial_t$ die Ableitung nach der zweiten Komponente bezeichnet. Ferner gilt f\"ur $i\in \{2,...,n\}$:
\begin{equation}\label{eq:4.41}
\partial_i F_{\epsilon,k}(x)=\partial_i F_k(\epsilon,\beta(\epsilon^{-1}x_1),x)
\end{equation}

\subsection{Konstruktion der Hilfsabbildung \texorpdfstring{$F_2$}{F unten 2}}
\label{subsec:4.2.1}

Der folgende Satz stellt den Anfang, f\"ur die induktive Konstruktion der Abbildungen $F_k\in C^{\infty}(I\times\mathbb{S}\times\overline{\mathbb{B}},\mathbb{R}^q)$, f\"ur $k\geq 2$, dar.
\begin{satz}\label{satz:4.2}
Gegeben sei eine freie Abbildung $F_0\in C^{\infty}(\overline{\mathbb{B}},\mathbb{R}^q)$, mit $q\geq \frac{n}{2}(n+3)+5$, und eine Abbildung $a\in C^{\infty}_0(\mathbb{B})$. Dann existiert eine kompakte Menge $K\subseteq \mathbb{B}$ mit $supp(a)\subseteq K$, so dass eine Abbildung $F_2\in C^{\infty}(I\times\mathbb{S}\times\overline{\mathbb{B}})$ existiert, welche die Bedingung $supp(F_2(\epsilon,t,\cdot)-F_0)\subseteq K$, f\"ur alle $(\epsilon,t)\subseteq I\times\mathbb{S}$ erf\"ullt und welche die folgenden Gleichungen erf\"ullt:
\begin{enumerate}[(i)]
\begin{item}
F\"ur alle $(\epsilon,t,x)\in I\times\mathbb{S}\times\overline{\mathbb{B}}$  gilt:
\begin{equation}\label{eq:4.58}
\left|\partial_1 F_2(\epsilon,t,x)+\frac{1}{\epsilon \varrho(t)}\ \partial_t F_2(\epsilon,t,x) \right|^2_{\mathbb{R}^q}=|\partial_1 F_0(x)|^2_{\mathbb{R}^q}+a^4(x)+\epsilon^3 f_{11}^{(2)}(\epsilon,t,x)
\end{equation}
\end{item}
\begin{item}
F\"ur alle $(\epsilon,t,x)\in I\times\mathbb{S}\times\overline{\mathbb{B}}$ und $i\in\{2,...,n\}$ gilt:
\begin{equation}\label{eq:4.57}
\left(\partial_1 F_2(\epsilon,t,x)+\frac{1}{\epsilon \varrho(t)}\ \partial_t F_2(\epsilon,t,x)  \right)\cdot \partial_i F_2(\epsilon,t,x)=\partial_1 F_0\cdot\partial_i F_0+\epsilon^3 f_{1i}^{(2)}(\epsilon,t,x)
\end{equation}
\end{item}
\begin{item}
F\"ur alle $(\epsilon,t,x)\in I\times\mathbb{S}\times\overline{\mathbb{B}}$ und $i,j\in\{2,...,n\}$ mit $i\leq j$ gilt:
\begin{equation}\label{eq:4.56}
\partial_i F_2(\epsilon,t,x)\cdot \partial_j F_2(\epsilon,t,x)=\partial_i F_0\cdot\partial_j F_0+\epsilon^3f_{ij}^{(2)}(\epsilon,t,x)
\end{equation}
\end{item}

\end{enumerate}

Hierbei ist $f^{(2)}:=(f_{ij}^{(2)})_{1\leq i\leq j \leq n}\in C^{\infty}(I\times\mathbb{S}\times\overline{\mathbb{B}},\mathbb{R}^{\frac{n}{2}(n+1)})$ mit $supp(f^{(2)}(\epsilon,t,\cdot))\subseteq K$, f\"ur alle $(\epsilon,t)\subseteq I\times\mathbb{S}$. Diese Abbildung erf\"ullt die
Eigenschaft:
\begin{equation}\label{f2est}
\sup_{\epsilon\in I}{\left\Vert \partial^{\alpha} f^{(2)}(\epsilon,\cdot,\cdot) \right\Vert}_{\overline{\mathbb{B}},\mathbb{R}^{\frac{n}{2}(n+1)}}\leq C(n,|\alpha|,F_0,a)
\end{equation} 
f\"ur alle $\alpha=(\partial_0,...,\partial_n)\in \mathbb{N}^{n+1}$. Ferner gilt:
\begin{equation}\label{eq:4.65}
\max_{(\epsilon,t,x)\in I\times\mathbb{S}\times\overline{\mathbb{B}}}{|F_2(\epsilon,t,x)-F_0(x)|_{\mathbb{R}^q}}\leq C(n,F_0,a)\cdot \epsilon
\end{equation}

\end{satz}

\begin{bew}
Es seien $v_1$ wie in Lemma \ref{lem:4.5}, und $u_1:= a^2\, v_1\in C^{\infty}(\mathbb{S}\times\overline{\mathbb{B}},\mathbb{R}^q)$, sowie $\varrho\in C^{\infty}(\mathbb{S})$ wie in \eqref{eq:4.29}. Die zu bestimmende Abbildung $F_2$ soll die Form $F_0+\epsilon u_1+\epsilon^2 u_2$ haben, wobei $u_2$ in geeigneter Weise zu bestimmen ist. Zun\"achst wird bemerkt, dass f\"ur alle $i,j\in \{1,...,n\}$ die Gleichung:
\begin{equation}\label{eq:4.43}
\partial_i F_0\cdot \partial_j u_1=\partial_j(\partial_i F_0\cdot u_1)-\partial_i\partial_j F_0\cdot u_1\stackrel{\eqref{eq:4.27}}{=}-\partial_i\partial_j F_0\cdot u_1=-a^2  \partial_i\partial_j F_0\cdot v_1
\end{equation}

gilt. Mit \eqref{eq:4.35} existieren Funktionen $h_1,...,h_n\in C^{\infty}(\mathbb{S}\times\overline{\mathbb{B}})$, deren periodische Fortsetzungen, bez\"uglich der ersten Komponente, auf ganz $\mathbb{R}$ glatt sind, so dass die folgenden Gleichungen erf\"ullt sind:
\begin{equation}\label{eq:4.44}
\partial_t h_1\stackrel{\hphantom{\eqref{eq:4.43}}}{=}-\varrho\, \partial_1 F_0\cdot \partial_1 u_1-\partial_t u_1\cdot \partial_1 u_1\stackrel{\eqref{eq:4.43}}{=}a^2 \varrho  \,\partial_1^2F_0\cdot v_1-a^2 \partial_t v_1 \cdot \partial_1 u_1
\end{equation}
und f\"ur alle $i\in\{2,...,n\}$:
\begin{align}\label{eq:4.45}
\begin{split}
\partial_t h_i&\stackrel{\hphantom{\eqref{eq:4.43}}}{=}-\varrho\, \partial_1 F_0\cdot \partial_i u_1-\varrho\, \partial_i F_0\cdot \partial_1 u_1-\partial_t u_1\cdot \partial_i u_1\\
&\stackrel{\eqref{eq:4.43}}{=}a^2 \varrho\,  \partial_1\partial_i F_0\cdot v_1+a^2 \varrho\,  \partial_i\partial_1 F_0\cdot v_1-a^2 \partial_t v_1 \cdot \partial_i u_1
\end{split}
\end{align}

Die Funktionen $h_1,...,h_n$ k\"onnen so gew\"ahlt werden, dass sie jeweils ein Produkt der Funktion $a^2$, und einer Funktion aus dem Raum $C^{\infty}(\mathbb{S}\times\overline{\mathbb{B}})$ sind. Mit der punktweisen linearen Unabh\"angigkeit der Abbildungen in \eqref{eq:4.28} und Folgerung \ref{fol:4.1}, existiert dann ein $\widehat{u}_2\in C^{\infty}(\mathbb{S}\times\overline{\mathbb{B}},\mathbb{R}^q)$, so dass das folgende Gleichungssystem gilt:
\begin{align}\label{eq:4.46}
\notag\left(\partial_1 F_0+\frac{1}{\varrho}\partial_t u_1 \right)\cdot \widehat{u}_2&=h_1\\
\notag\partial_i F_0 \cdot \widehat{u}_2&=h_i\hspace{0.5cm}\text{f\"ur }2\leq i\leq n\\
\notag\left(\partial_1^2 F_0+\frac{2}{\varrho} \partial_t\partial_1 u_1 \right)\cdot \widehat{u}_2&=\frac{1}{2}(|\partial_1 u_1|_{\mathbb{R}^q}^2+2\partial_1 h_1)\\
\left(\partial_1\partial_i F_0+\frac{1}{\varrho} \partial_t\partial_i u_1 \right)\cdot \widehat{u}_2&=\frac{1}{2}(\partial_1 u_1\cdot \partial_i u_1+\partial_1 h_i+\partial_i h_1)\hspace{0.5cm}\text{f\"ur } 2\leq i\leq n\\
\notag\partial_i \partial_j F_0 \cdot \widehat{u}_2&=\frac{1}{2}(\partial_i u_1\cdot \partial_j u_1+\partial_i h_j+\partial_j h_i)\hspace{0.5cm}\text{f\"ur }2\leq i\leq j\leq n\\
\notag\partial_t v_1\cdot \widehat{u}_2&=0\hspace{2cm}\partial_t^2 v_1\cdot \widehat{u}_2=0
\end{align}

Dabei ist zu beachten, dass die rechten Seiten im Gleichungssystem \eqref{eq:4.46} jeweils ein Produkt der Abbildung $a$, mit einer Abbildung aus dem Raum $C^{\infty}(\mathbb{S}\times\overline{\mathbb{B}})$ sind, weshalb sich, unter Verwendung der Formel \eqref{eq:4.47}, auch $\widehat{u}_2$ als Produkt von $a$ und einer Funktion aus $C^{\infty}(\mathbb{S}\times\overline{\mathbb{B}},\mathbb{R}^q)$ darstellen l\"asst.  Nun wird ein $\widehat{v}_2\in C^{\infty}(\mathbb{S}\times\overline{\mathbb{B}},\mathbb{R}^q)$ bestimmt, so dass das Gleichungssystem:
\begin{align}\label{eq:4.48}
\begin{split}
\left(\partial_1 F_0+\frac{1}{\varrho}\partial_t u_1 \right)\cdot \widehat{v}_2&=-\partial_1 u_1\cdot \widehat{u}_2\\
\partial_i F_0\cdot \widehat{v}_2&=-\partial_i u_1\cdot \widehat{u}_2\hspace{0.5cm}\text{f\"ur }2\leq i\leq n\\
\partial_t v_1 \cdot \widehat{v}_2&=0\\
\partial_t^2 v_1 \cdot \widehat{v}_2&=\frac{1}{2 a^2}|\partial_t \widehat{u}_2|_{\mathbb{R}^q}^2
\end{split}
\end{align}

gilt. Dabei ist zu beachten, dass mit der obigen Bemerkung und \cite[Proposition 2.26]{lee2003introduction} eine Funktion $\overline{u}\in C^{\infty}(\mathbb{S}\times\overline{\mathbb{B}},\mathbb{R}^q)$ existiert, die in einer kompakten Teilmenge $K\subseteq \mathbb{B}$ mit $supp(a)\subseteq K$ enthalten ist, so dass $\partial_t \widehat{u}_2=a\,\overline{u}$ gilt. Unter Verwendung der Formel \eqref{eq:4.47} l\"asst sich dann $supp(\widehat{v}_2(t,\cdot))\subseteq K$, f\"ur alle $t\in\mathbb{S}$, erreichen. Definiere nun eine Abbildung $u_2\in C^{\infty}(I\times\mathbb{S}\times\overline{\mathbb{B}})$ mit:
\begin{equation}\label{eq:4.59}
u_2(\epsilon,t,x):= \widehat{u}_2(t,x)+\epsilon \widehat{v}_2(t,x)
\end{equation}

Mit \eqref{eq:4.46} und \eqref{eq:4.48} gelten dann die Gleichungen:
\begin{align}
\begin{split}\label{eq:4.49}
\left(\partial_1 F_0+\frac{1}{\varrho}\partial_t u_1+\epsilon \partial_1 u_1 \right)\cdot u_2&=h_1+\Lambda(\epsilon^2)
\end{split}\\
\begin{split}\label{eq:4.50}
(\partial_i F_0+\epsilon \partial_i u_1 )\cdot u_2&=h_i+\Lambda(\epsilon^2)\hspace{1cm}\text{f\"ur } 2\leq i\leq n
\end{split}\\
\begin{split}\label{eq:4.51}
\left(\partial_1^2 F_0+\frac{2}{\varrho}\partial_t \partial_1 u_1 \right)\cdot u_2&=\frac{1}{2}(|\partial_1 u_1|_{\mathbb{R}^q}^2+2\partial_1 h_1)+\Lambda(\epsilon)
\end{split}\\
\begin{split}\label{eq:4.52}
\left(\partial_1\partial_i F_0+\frac{1}{\varrho}\partial_t \partial_i u_1 \right)\cdot u_2&=\frac{1}{2}(\partial_1 u_1\cdot \partial_i u_1+\partial_1 h_i+\partial_i h_1)+\Lambda(\epsilon)\hspace{1cm}\text{f\"ur } 2\leq i\leq n
\end{split}\\
\begin{split}\label{eq:4.53}
\partial_i \partial_j F_0 \cdot u_2&=\frac{1}{2}(\partial_i u_1\cdot \partial_j u_1+\partial_i h_j+\partial_j h_i)+\Lambda(\epsilon)\hspace{0.5cm}\text{f\"ur }2\leq i\leq  j\leq n
\end{split}\\
\begin{split}\label{eq:4.54}
\partial_t u_1\cdot u_2&=0
\end{split}\\
\begin{split}\label{eq:4.55}
\partial_t^2 u_1\cdot u_2&=\frac{\epsilon}{2} |\partial_t u_2|_{\mathbb{R}^q}^2+\Lambda(\epsilon^2)
\end{split}
\end{align}

Hierbei sei vereinbart, dass \gls{Lambda}, mit $k\in\mathbb{N}$, f\"ur einen geeigneten Ausdruck der Form $\epsilon^k\, h$, wobei $h\in C^{\infty}(I\times\mathbb{S}\times\overline{\mathbb{B}})$, mit $h(\epsilon,t,\cdot)\subseteq K$ f\"ur alle $(\epsilon,t)\in I\times\mathbb{S}$,  steht. Nun wird f\"ur $F_2\in C^{\infty}(I\times\mathbb{S}\times\overline{\mathbb{B}})$ mit:
\begin{equation}\label{eq:4.60}
F_2(\epsilon,t,x):= F_0(x)+\epsilon u_1(t,x)+\epsilon^2 u_2(\epsilon,t,x)
\end{equation}

gezeigt, dass die Gleichungen \eqref{eq:4.58} bis \eqref{eq:4.56} gelten. Vorher sei erw\"ahnt, dass wegen $u_1=a^2\, v_1$ und \eqref{eq:4.59} beziehungsweise \eqref{eq:4.60}, die Bedingung $supp(F_2(\epsilon,t,\cdot)-F_0)\subseteq K$ f\"ur alle $(\epsilon,t)\subseteq I\times\mathbb{S}$, erf\"ullt ist. Zuerst wird \eqref{eq:4.58} gezeigt:
\begin{align}\label{eq:4.64}
\notag&\left|\partial_1 F_2+\frac{1}{\epsilon\varrho}\partial_t F_2 \right|^2_{\mathbb{R}^q}\stackrel{\eqref{eq:4.60}}{=}\left|\partial_1 F_0+\epsilon \partial_1 u_1+\epsilon^2\partial_1 u_2+\frac{1}{\epsilon\varrho}\underbrace{\partial_t F_0}_{=0}+ \frac{1}{\varrho}\partial_t u_1+\frac{\epsilon}{\varrho}\partial_t u_2\right|^2_{\mathbb{R}^q}\\
\notag=&\left|\partial_1 F_0\right|^2_{\mathbb{R}^q}+\frac{1}{\varrho^2}\underbrace{\left|\partial_t u_1 \right|^2_{\mathbb{R}^q}}_{\stackrel{\eqref{eq:4.29}}{=}a^4\varrho^2}+\frac{2}{\varrho}\underbrace{\partial_1 F_0\cdot\partial_t u_1}_{\stackrel{\eqref{eq:4.27}}{=}0}+2\epsilon\left(\partial_1 u_1+\frac{1}{\varrho}\partial_t u_2 \right)\cdot\left(\partial_1 F_0+\frac{1}{\varrho}\partial_t u_1 \right)\\
\notag&+\epsilon^2 \left|\partial_1 u_1+\frac{1}{\varrho} \partial_t u_2 \right|^2+2\epsilon^2\left(\partial_1 F_0+\frac{1}{\varrho} \partial_t u_1 \right)\cdot \partial_1 u_2+\Lambda(\epsilon^3)\\ 
=&\left|\partial_1 F_0\right|^2_{\mathbb{R}^q}+a^4+2\epsilon\left(\partial_1 u_1+\frac{1}{\varrho}\partial_t u_2 \right)\cdot\left(\partial_1 F_0+\frac{1}{\varrho}\partial_t u_1 \right)\\
\notag&+\epsilon^2 \left(|\partial_1 u_1|^2_{\mathbb{R}^q} +\frac{2}{\varrho}\partial_1 u_1\cdot \partial_t u_2+\frac{1}{\varrho^2} |\partial_t u_2|^2_{\mathbb{R}^q}\right)+2\epsilon^2\left(\partial_1 F_0+\frac{1}{\varrho} \partial_t u_1 \right)\cdot \partial_1 u_2+\Lambda(\epsilon^3)
\end{align}

Der dritte Summand auf der rechten Seite wird weiter umgeformt:
\begin{align*}
&\left(\partial_1 F_0+\frac{1}{\varrho}\partial_t u_1 \right)\cdot\left(\partial_1 u_1+\frac{1}{\varrho}\partial_t u_2 \right)\\
\stackrel{\hphantom{\eqref{eq:4.44}}}{=}&\left(\partial_1 F_0+\frac{1}{\varrho}\partial_t u_1 \right)\cdot \partial_1 u_1+\frac{1}{\varrho}\partial_t\biggl[\underbrace{\left( \partial_1 F_0+\frac{1}{\varrho}\partial_t u_1\right)\cdot u_2}_{\stackrel{\eqref{eq:4.49}}{=}h_1-\epsilon \partial_1 u_1\cdot u_2+\Lambda(\epsilon^2)\ }\biggr] -\frac{1}{\varrho}\underbrace{\partial_t \partial_1 F_0\cdot u_2}_{=0}\\
&-\frac{1}{\varrho^2}\underbrace{\partial_t^2 u_1\cdot u_2}_{\stackrel{\eqref{eq:4.55}}{=}\frac{\epsilon}{2}|\partial_t u_2|^2_{\mathbb{R}^q}+\Lambda(\epsilon^2)}-\frac{1}{\varrho}\partial_t\left(\frac{1}{\varrho} \right)\underbrace{\partial_t u_1\cdot u_2}_{\stackrel{\eqref{eq:4.54}}{=}0}\\
\stackrel{\hphantom{\eqref{eq:4.44}}}{=}&\left(\partial_1 F_0+\frac{1}{\varrho}\partial_t u_1 \right)\cdot \partial_1 u_1+\frac{1}{\varrho}\partial_t h_1-\frac{\epsilon}{\varrho}\partial_t(\partial_1 u_1\cdot u_2)-\frac{\epsilon}{2\varrho^2}|\partial_t u_2|^2_{\mathbb{R}^q}+\Lambda(\epsilon^2)\\
\stackrel{\eqref{eq:4.44}}{=}&-\frac{\epsilon}{\varrho}\partial_t(\partial_1 u_1\cdot u_2)-\frac{\epsilon}{2\varrho^2}|\partial_t u_2|^2_{\mathbb{R}^q}+\Lambda(\epsilon^2)\\
\end{align*}

Es folgt mit \eqref{eq:4.64}:
\begin{align*}
&\left|\partial_1 F_2+\frac{1}{\epsilon\varrho}\partial_t F_2 \right|^2_{\mathbb{R}^q}=\left|\partial_1 F_0\right|^2_{\mathbb{R}^q}+a^4\\
&+\epsilon^2\biggl[|\partial_1 u_1|^2_{\mathbb{R}^q} +\frac{2}{\varrho}\partial_1 u_1\cdot \partial_t u_2-\frac{2}{\varrho}\partial_t(\partial_1 u_1\cdot u_2)+2\left(\partial_1 F_0+\frac{1}{\varrho}\partial_t u_1 \right)\cdot \partial_1 u_2\biggr]+\Lambda(\epsilon^3)\\
\stackrel{\hphantom{\eqref{eq:4.49}}}{=}&\left|\partial_1 F_0\right|^2_{\mathbb{R}^q}+a^4+\epsilon^2\biggl[|\partial_1 u_1|^2_{\mathbb{R}^q} -\frac{2}{\varrho}\partial_t\partial_1 u_1\cdot u_2+2\left(\partial_1 F_0+\frac{1}{\varrho}\partial_t u_1 \right)\cdot \partial_1 u_2\biggr]+\Lambda(\epsilon^3)\\
\stackrel{\hphantom{\eqref{eq:4.49}}}{=}&\left|\partial_1 F_0\right|^2_{\mathbb{R}^q}+a^4\\
&+\epsilon^2\biggl[|\partial_1 u_1|^2_{\mathbb{R}^q}-\underbrace{2\left(\partial_1^2 F_0+\frac{2}{\varrho}\partial_t \partial_1 u_1\right)\cdot u_2}_{\stackrel{\eqref{eq:4.51}}{=}\left|\partial_1 u_1\right|^2_{\mathbb{R}^q}+2\partial_1 h_1+\Lambda(\epsilon)} +2\partial_1\biggl(\underbrace{\left(\partial_1 F_0+\frac{1}{\varrho}\partial_t u_1 \right)\cdot  u_2}_{\stackrel{\eqref{eq:4.49}}{=}h_1+\Lambda(\epsilon)}\biggr)\biggr]+\Lambda(\epsilon^3)\\
\stackrel{\eqref{eq:4.54}}{=}&\left|\partial_1 F_0\right|^2_{\mathbb{R}^q}+a^4+\epsilon^3 f_{11}^{(2)}
\end{align*}

Nun wird \eqref{eq:4.57} gezeigt, es sei $i\in\{2,...,n\}$, dann ist unter Beachtung von \eqref{eq:4.27}
\begin{align}\label{eq:4.61}
\begin{split}
&\left(\partial_1 F_2+\frac{1}{\epsilon \varrho}\partial_t F_2 \right)\cdot \partial_i F_2\\
\stackrel{\eqref{eq:4.60}}{=}& \biggl(\partial_1 F_0+\epsilon \partial_1 u_1 +\epsilon^2 \partial_1 u_2+\frac{1}{\epsilon\varrho}\underbrace{\partial_t F_0}_{=0}+\frac{1}{\varrho} \partial_t u_1 +\frac{\epsilon}{\varrho} \partial_t u_2\biggr)\cdot \left(\partial_i F_0+\epsilon \partial_i u_1 +\epsilon^2 \partial_i u_2 \right)\\
=&\partial_1 F_0\cdot \partial_i F_0+\epsilon\biggl[\left(\partial_1 u_1+\frac{1}{\varrho} \partial_t u_2\right)\cdot \partial_i F_0 +\left(\partial_1 F_0+ \frac{1}{\varrho} \partial_t u_1\right)\cdot\partial_i u_1\biggr]\\
+&\epsilon^2\biggl[\left(\partial_1 F_0+\frac{1}{\varrho}\partial_t u_1\right)\cdot \partial_i u_2+\left(\partial_1 u_1+\frac{1}{\varrho} \partial_t u_2 \right)\cdot \partial_i u_1 + \partial_1 u_2\cdot\partial_i F_0\biggr]+\Lambda(\epsilon^3)
\end{split}
\end{align}

Der Term in der ersten eckigen Klammer wird unter dem Aspekt untersucht, $\epsilon$ zu isolieren:
\begin{align}\label{eq:4.62}
\begin{split}
&\left(\partial_1 u_1+\frac{1}{\varrho} \partial_t u_2\right)\cdot \partial_i F_0 +\left(\partial_1 F_0+ \frac{1}{\varrho} \partial_t u_1\right)\cdot\partial_i u_1\\
\stackrel{\hphantom{\eqref{eq:4.50}}}{=}&\frac{1}{\varrho} \partial_t u_2\cdot \partial_i F_0+ \partial_1 F_0\cdot\partial_i u_1+\partial_i F_0\cdot\partial_1 u_1+\frac{1}{\varrho} \partial_t u_1\cdot \partial_i u_1\\
\stackrel{\eqref{eq:4.45}}{=}&\frac{1}{\varrho} \partial_t (u_2\cdot \partial_i F_0)-\frac{1}{\varrho}\partial_t h_i\stackrel{\eqref{eq:4.50}}{=}\frac{1}{\varrho}\partial_t h_i-\frac{\epsilon}{\varrho} \partial_t ( \partial_i u_1 \cdot u_2)-\frac{1}{\varrho}\partial_t h_i+\Lambda(\epsilon^2)\\
\stackrel{\hphantom{\eqref{eq:4.50}}}{=}&-\frac{\epsilon}{\varrho} \partial_t ( \partial_i u_1 \cdot u_2)+\Lambda(\epsilon^2)
\end{split}
\end{align}

Nun wird der Ausdruck in der zweiten eckigen Klammer in \eqref{eq:4.61} umgeformt:
\begin{align}\label{eq:4.63}
\begin{split}
&\left(\partial_1 F_0+\frac{1}{\varrho}\partial_t u_1\right)\cdot \partial_i u_2+\left(\partial_1 u_1+\frac{1}{\varrho} \partial_t u_2 \right)\cdot \partial_i u_1 + \partial_1 u_2\cdot\partial_i F_0\\
=&\partial_1 F_0\cdot \partial_i u_2+\frac{1}{\varrho}\partial_t u_1\cdot \partial_i u_2+\left(\partial_1 u_1+\frac{1}{\varrho} \partial_t u_2 \right)\cdot \partial_i u_1 + \partial_1 u_2\cdot\partial_i F_0\\
=&\underbrace{\partial_i(\partial_1 F_0\cdot u_2)}_{\stackrel{\eqref{eq:4.49}}{=}\partial_i h_1+\Lambda(\epsilon)}-\partial_i \partial_1 F_0\cdot u_2+\frac{1}{\varrho}\partial_i (\underbrace{\partial_t u_1\cdot  u_2}_{\stackrel{\eqref{eq:4.54}}{=}0})-\frac{1}{\varrho}\partial_t \partial_i u_1\cdot  u_2\\
&+\left(\partial_1 u_1+\frac{1}{\varrho} \partial_t u_2 \right)\cdot \partial_i u_1+\underbrace{\partial_1(u_2\cdot \partial_i F_0)}_{\stackrel{\eqref{eq:4.50}}{=}\partial_1 h_i +\Lambda(\epsilon)}- u_2\cdot \partial_1\partial_i F_0\\
=&\partial_1 h_i+\partial_i h_1-\frac{1}{\varrho} \partial_t \partial_i u_1 \cdot u_2-\underbrace{2\partial_i \partial_1 F_0\cdot u_2}_{\stackrel{\eqref{eq:4.52}}{=}\partial_1 u_1\cdot \partial_i u_1+\partial_1 h_i+\partial_i h_1-\frac{2}{\varrho}\partial_t\partial_i u_1\cdot u_2+\Lambda(\epsilon)}\\
&+\left(\partial_1 u_1+\frac{1}{\varrho} \partial_t u_2 \right)\cdot \partial_i u_1+\Lambda(\epsilon)\\
=&\frac{1}{\varrho} \partial_t \partial_i u_1 \cdot u_2+\frac{1}{\varrho} \partial_t u_2 \cdot \partial_i u_1 +\Lambda(\epsilon)=\frac{1}{\varrho} \partial_t (\partial_i u_1\cdot u_2) +\Lambda(\epsilon)
\end{split}
\end{align}

Aus \eqref{eq:4.61}, \eqref{eq:4.62} und \eqref{eq:4.63} folgt:
\begin{align*}
&\left(\partial_1 F_2+\frac{1}{\epsilon \varrho}\partial_t F_2 \right)\cdot \partial_i F_2\\
=&\partial_1 F_0\cdot \partial_i F_0+\epsilon^2\left(-\frac{1}{\varrho} \partial_t (\partial_i u_1\cdot u_2)+\frac{1}{\varrho} \partial_t (\partial_i u_1\cdot u_2) \right)+\Lambda(\epsilon^3)\\
=&\partial_1 F_0\cdot \partial_i F_0+\epsilon^3 f_{1i}^{(2)}
\end{align*}

Schlie\ss{}lich wird noch \eqref{eq:4.56} gezeigt. Es seien $i,j\in\{2,...,n\}$ mit $i\leq j$:
\begin{align*}
&\partial_i F_2 \cdot \partial_j F_2=\left(\partial_i F_0+\epsilon \partial_i u_1+\epsilon^2\partial_i u_2\right)\cdot \left(\partial_j F_0+\epsilon \partial_j u_1+\epsilon^2\partial_j u_2\right)\\
\stackrel{\eqref{eq:4.60}}{=}& \partial_i F_0\cdot \partial_j F_0+\epsilon\left(\partial_i F_0\cdot \partial_j u_1+\partial_i u_1\cdot \partial_j F_0\right)\\
&+\epsilon^2(\partial_i F_0\cdot\partial_j u_2 +\partial_i u_1\cdot \partial_j u_1+\partial_i u_2\cdot \partial_j F_0 )+\Lambda(\epsilon^3)\\
\stackrel{\eqref{eq:4.43}}{=}&  \partial_i F_0\cdot \partial_j F_0-\epsilon a^2(\underbrace{\partial_i \partial_j F_0\cdot  v_1}_{\stackrel{\eqref{eq:4.27}}{=}0}+ \underbrace{v_1\cdot \partial_j\partial_i F_0}_{\stackrel{\eqref{eq:4.27}}{=}0})\\
&+\epsilon^2\bigl[\underbrace{\partial_j(\partial_i F_0\cdot u_2)}_{\stackrel{\eqref{eq:4.50}}{=}\partial_j h_i+\Lambda(\epsilon)}-\partial_j\partial_i F_0\cdot u_2+\partial_i u_1\cdot \partial_j u_1+\underbrace{\partial_i( u_2\cdot \partial_j F_0)}_{\stackrel{\eqref{eq:4.50}}{=}\partial_i h_j+\Lambda(\epsilon)}- u_2\cdot \partial_i\partial_j F_0 \bigr]+\Lambda(\epsilon^3)\\
\stackrel{\eqref{eq:4.53}}{=}&  \partial_i F_0\cdot \partial_j F_0+\epsilon^3 f_{ij}^{(2)}
\end{align*}

Abschlie\ss{}end sei noch erw\"ahnt, dass mit \eqref{eq:4.60}, wegen $u_2\in C^{\infty}(I\times\mathbb{S}\times\overline{\mathbb{B}})$, die Absch\"at\-zung \eqref{eq:4.65} erf\"ullt ist.

\end{bew}

\subsection{Konstruktion der Hilfsabbildung \texorpdfstring{$F_k$}{Fk} für \texorpdfstring{$k\geq 3$}{k groesser 3}}
\label{subsec:4.2.2}

Satz \ref{satz:4.2} wird nun zu folgendem Satz verallgemeinert:
\begin{satz}\label{satz:4.3}
Gegeben sei eine freie Abbildung $F_0\in C^{\infty}(\overline{\mathbb{B}},\mathbb{R}^q)$, mit $q\geq \frac{n}{2}(n+3)+5$, und eine Abbildung $a\in C^{\infty}_0(\mathbb{B})$. Dann existiert eine kompakte Menge $K\subseteq \mathbb{B}$, mit $supp(a)\subseteq K$, so dass die folgende Aussage gilt: F\"ur jedes $k\in\mathbb{N}\backslash\{0,1\}$ existiert eine Abbildung $F_k\in C^{\infty}(I\times\mathbb{S}\times\overline{\mathbb{B}})$, welche die Bedingung $supp(F_k(\epsilon,t,\cdot)-F_0)\subseteq K$ f\"ur alle $(\epsilon,t)\subseteq I\times\mathbb{S}$ erf\"ullt, und welche die folgenden Gleichungen erf\"ullt:
\begin{enumerate}[(i)]
\begin{item}
F\"ur alle $(\epsilon,t,x)\in I\times\mathbb{S}\times\overline{\mathbb{B}}$  gilt:
\begin{equation}\label{eq:4.68}
\left|\partial_1 F_k(\epsilon,t,x)+\frac{1}{\epsilon \varrho(t)}\ \partial_t F_k(\epsilon,t,x) \right|^2_{\mathbb{R}^q}=|\partial_1 F_0(x)|^2_{\mathbb{R}^q}+a^4(x)+\epsilon^{k+1} f_{11}^{(k)}(\epsilon,t,x)
\end{equation}
\end{item}
\begin{item}
F\"ur alle $(\epsilon,t,x)\in I\times\mathbb{S}\times\overline{\mathbb{B}}$ und $i\in\{2,...,n\}$ gilt:
\begin{equation}\label{eq:4.67}
\left(\partial_1 F_k(\epsilon,t,x)+\frac{1}{\epsilon \varrho(t)}\ \partial_t F_k(\epsilon,t,x)  \right)\cdot \partial_i F_k(\epsilon,t,x)=\partial_1 F_0\cdot\partial_i F_0+\epsilon^{k+1} f_{1i}^{(k)}(\epsilon,t,x)
\end{equation}
\end{item}
\begin{item}
F\"ur alle $(\epsilon,t,x)\in I\times\mathbb{S}\times\overline{\mathbb{B}}$ und $i,j\in\{2,...,n\}$ mit $i\leq j$ gilt:
\begin{equation}\label{eq:4.66}
\partial_i F_k(\epsilon,t,x)\cdot \partial_j F_k(\epsilon,t,x)=\partial_i F_0\cdot\partial_j F_0+\epsilon^{k+1} f_{ij}^{(k)}(\epsilon,t,x)
\end{equation}
\end{item}
\end{enumerate}

Hierbei ist $f^{(k)}:=(f_{ij}^{(k)})_{1\leq i\leq j \leq n}\in C^{\infty}(I\times\mathbb{S}\times\overline{\mathbb{B}},\mathbb{R}^{\frac{n}{2}(n+1)})$ mit $supp(f^{(k)}(\epsilon,t,\cdot))\subseteq K$ f\"ur alle $(\epsilon,t)\subseteq I\times\mathbb{S}$. Diese Abbildung erf\"ullt die
Eigenschaft:
\begin{equation}\label{fkest}
\sup_{\epsilon\in I}{\left\Vert \partial^{\alpha} f^{(k)}(\epsilon,\cdot,\cdot) \right\Vert}_{\overline{\mathbb{B}},\mathbb{R}^{\frac{n}{2}(n+1)}}\leq C(n,k,|\alpha|,F_0,a)
\end{equation} 
f\"ur alle $\alpha=(\partial_0,...,\partial_n)\in \mathbb{N}^{n+1}$. Ferner gilt:
\begin{equation}\label{eq:4.69}
\max_{(\epsilon,t,x)\in I\times\mathbb{S}\times\overline{\mathbb{B}}}{|F_k(\epsilon,t,x)-F_0(x)|_{\mathbb{R}^q}}\leq C(n,k,F_0,a)\cdot \epsilon
\end{equation}

\end{satz}

\begin{bew}
Der Satz wird mittels Induktion \"uber $k$ bewiesen. Mit Satz \ref{satz:4.2} stimmt die Aussage f\"ur $k=2$. Nun sei die Aussage f\"ur ein $k\in\mathbb{N}\backslash\{0,1\}$ richtig. Es wird ein $u_{k+1}\in C^{\infty}(I\times\mathbb{S}\times\overline{\mathbb{B}},\mathbb{R}^q)$ bestimmt, so dass die Abbildung $F_{k+1}:= F_k+\epsilon^{k+1} u_{k+1}$ die gew\"unschten Eigenschaften besitzt. Da $F_k$ die Gleichungen \eqref{eq:4.68} bis \eqref{eq:4.66} erf\"ullt, existiert ein $h:= (h_{ij})_{1\leq i\leq j\leq n}\subseteq C^{\infty}(\mathbb{S}\times\overline{\mathbb{B}},\mathbb{R}^{\frac{n}{2}(n+1)})$, mit $supp(h(t,\cdot))\subseteq K$ f\"ur alle $t\in \mathbb{S}$, so dass:
\begin{align}
\begin{split}\label{eq:4.81}
\left|\partial_1 F_k+\frac{1}{\epsilon \varrho}\ \partial_t F_k \right|^2_{\mathbb{R}^q}&=\left|\partial_1 F_0\right|^2_{\mathbb{R}^q}+a^4+\epsilon^{k+1} h_{11}+\Lambda(\epsilon^{k+2})
\end{split}\\
\begin{split}\label{eq:4.80}
\left(\partial_1 F_k+\frac{1}{\epsilon \varrho}\ \partial_t F_k  \right)\cdot \partial_i F_k&=\partial_1 F_0\cdot \partial_i F_0+\epsilon^{k+1} h_{1i}+\Lambda(\epsilon^{k+2})\hspace{0.5cm}\text{f\"ur }2\leq i\leq n
\end{split}\\
\begin{split}\label{eq:4.79}
\partial_i F_k\cdot \partial_j F_k&=\partial_i F_0\cdot \partial_j F_0+\epsilon^{k+1} h_{ij}+\Lambda(\epsilon^{k+2})\hspace{0.5cm}\text{f\"ur }2\leq i\leq j\leq n
\end{split}
\end{align}

erf\"ullt ist. Dabei sei f\"ur den Induktionsschritt von $k=2$ auf $k+1=3$ auf das Gleichungssystem \eqref{eq:4.49} bis \eqref{eq:4.55} verwiesen. Analog zu \eqref{eq:4.59} soll die Abbildung $u_{k+1}$ von der Form $\widehat{u}_{k+1}+\epsilon \widehat{v}_{k+1}$ sein, wobei $\widehat{u}_{k+1},\widehat{v}_{k+1}\in C^{\infty}(\mathbb{S}\times\overline{\mathbb{B}},\mathbb{R}^q)$ mit $supp(\widehat{u}_{k+1}(t,\cdot)), supp(\widehat{v}_{k+1}(t,\cdot))\subseteq K$ f\"ur alle $t\in\mathbb{S}$. Dazu soll die Abbildung $u_{k+1}$ so gew\"ahlt werden, dass das Gleichungssystem:
\begin{align}
\begin{split}\label{eq:4.70}
(\partial_i F_0+\epsilon \partial_i u_1)\cdot u_{k+1}&=\Lambda(\epsilon^2)\hspace{0.5cm}\text{f\"ur }2\leq i\leq n
\end{split}\\
\begin{split}\label{eq:4.71}
\left[\partial_1 F_0+\frac{1}{\varrho} \partial_t u_1+\epsilon\left(\partial_1 u_1+\frac{1}{\varrho}\partial_t u_2\right) \right]\cdot u_{k+1}&=\Lambda(\epsilon^2)
\end{split}\\
\begin{split}\label{eq:4.72}
\partial_i \partial_j F_0\cdot u_{k+1}&=\frac{1}{2}h_{ij}+\Lambda(\epsilon)\hspace{0.5cm}\text{f\"ur }2\leq i\leq j\leq n
\end{split}\\
\begin{split}\label{eq:4.73}
\left(\partial_1 \partial_i F_0+\frac{1}{\varrho}\partial_t\partial_i u_1 \right)\cdot u_{k+1}&=\frac{1}{2}h_{1i}+\Lambda(\epsilon)\hspace{0.5cm}\text{f\"ur }2\leq i\leq n
\end{split}\\
\begin{split}\label{eq:4.74}
\left[\partial_1^2 F_0+\frac{2}{\varrho}\partial_t \partial_1 u_1 +\frac{1}{\epsilon\varrho}\partial_t\left(\frac{\partial_t u_1+\epsilon \partial_t u_2}{\varrho} \right)\right]\cdot u_{k+1}&=\frac{1}{2}h_{11}+\Lambda(\epsilon)
\end{split}
\end{align}

erf\"ullt ist. Nun wird beschrieben, wie $\widehat{u}_{k+1}$ und $\widehat{v}_{k+1}$ gew\"ahlt werden k\"onnen, damit die Abbildung $u_{k+1}=\widehat{u}_{k+1}+\epsilon \widehat{v}_{k+1}$ das Gleichungssystem \eqref{eq:4.70} bis \eqref{eq:4.74} erf\"ullt. In \eqref{eq:4.59} wurde die Funktion $u_2$ als Summe $u_2=\widehat{u}_2+\epsilon \widehat{v}_2$ eingef\"uhrt. Nach Konstruktion von $\widehat{u}_2$ in \eqref{eq:4.46} gilt die Absch\"atzung:
\begin{equation*}
\left| \partial_t \widehat{u}_2(t,x) \right|_{\mathbb{R}^q}+\left| \partial_t^2 \widehat{u}_2(t,x) \right|_{\mathbb{R}^q} \leq C\cdot |a(x)|
\end{equation*}

f\"ur alle $(t,x)\in \mathbb{S}\times \overline{\mathbb{B}}$. Hierbei ist $C\in\mathbb{R}_{>0}$ eine Konstante, die weder von $t$ noch von $x$ abh\"angt, aber dar\"uber hinaus nicht weiter spezifiziert wird. Wegen der stetigen Abh\"angigkeit der Determinante von  den Matrixeintr\"agen, existiert, mit der punktweisen linearen Unabh\"angigkeit der Funktionen in \eqref{eq:4.28}, ein $\delta\in\mathbb{R}_{>0}$, so dass die $\frac{n}{2}(n+3)+2$ Vektoren:
\begin{align}\label{eq:4.75}
\begin{array}{ll}
\partial_1 F_0(x)+\frac{1}{\varrho(t)}\partial_t u_1(t,x) & \partial_i F_0(x)\hspace{0.5cm}\text{f\"ur }2\leq i\leq n \\
\partial_1^2 F_0(x)+\frac{2}{\varrho(t)}\partial_t\partial_1 u_1(t,x)+\frac{1}{\varrho(t)}\partial_t\left(\frac{\partial_t \widehat{u}_2(t,x)}{\varrho(t)} \right)\\
\partial_1 \partial_i F_0(x)+\frac{1}{\varrho(t)}\partial_t\partial_i u_1(x,t)\hspace{0.5cm}\text{f\"ur }2\leq i\leq n & \partial_i\partial_j F_0(x)\hspace{0.5cm}\text{f\"ur }2\leq i\leq j\leq n  \\
\partial_t v_1(t,x) & \partial_t^2 v_1(t,x)
\end{array}
\end{align}

f\"ur alle $(t,x)\in \mathbb{S}\times\overline{\mathbb{B}}$, mit $|a(x)|<\delta$, linear unabh\"angig sind. Definiere nun offene Mengen $U_1,U_2\subseteq \mathbb{B}$ mit:
\begin{align}\label{eq:4.93}
&U_1:= \left\{x\in \mathbb{B} : |a(x)|<\delta\right\} & &U_2:= \left\{x\in \mathbb{B} : |a(x)|>\frac{\delta}{2}\right\}
\end{align}

und w\"ahle eine Zerlegung der Eins, bez\"uglich der offenen \"Uberdeckung $\left\{U_1,U_2 \right\}$, also $\varphi_1,\varphi_2\in C^{\infty}(\mathbb{B})$ mit $supp(\varphi_l)\subseteq U_l$, f\"ur $l\in\{1,2\}$, und $\varphi_1+\varphi_2\equiv 1$, und definiere f\"ur $l\in\{1,2\}$ jeweils $h^{(l)}\in C^{\infty}(\mathbb{S}\times \overline{\mathbb{B}},\mathbb{R}^{\frac{n}{2}(n+1)})$ mit:
\begin{equation*}
(h^{(l)}_{ij})_{1\leq i\leq j\leq n}:=(\varphi_l\cdot h_{ij})_{1\leq i\leq j\leq n}
\end{equation*}

Dann ist die Bedingung $supp(h^{(l)}(t,\cdot))\subseteq K$ f\"ur alle $t\in\mathbb{S}$ erf\"ullt. Nun werden Abbildungen $u_{k+1}^{(l)}\in C^{\infty}(\mathbb{S}\times \overline{\mathbb{B}},\mathbb{R}^q)$ konstruiert, so dass $u_{k+1}^{(l)}$ jeweils das Gleichungssystem \eqref{eq:4.70} bis \eqref{eq:4.74}, f\"ur $h=h^{(l)}$, l\"ost. Die Funktionen $u_{k+1}^{(l)}$ sollen dabei die Form $\widehat{u}_{k+1}^{(l)}+\epsilon\widehat{u}_{k+1}^{(l)}$ haben, wobei $\widehat{u}_{k+1}^{(l)},\widehat{v}_{k+1}^{(l)}\in C^{\infty}(\mathbb{S}\times \overline{\mathbb{B}},\mathbb{R}^q)$, mit $\widehat{u}_{k+1}^{(l)}(t,\cdot)\subseteq K$, und $\widehat{v}_{k+1}^{(l)}(t,\cdot)\subseteq K$ f\"ur alle $t\in\mathbb{S}$, erf\"ullt sein soll. Mit der linearen Unabh\"angigkeit der Vektoren in \eqref{eq:4.75}, existieren, unter Beachtung dass $u_1=a^2\, v_1$ gilt, Abbildungen $\widehat{u}_{k+1}^{(1)}\in C^{\infty}(\mathbb{S}\times \overline{\mathbb{B}},\mathbb{R}^q)$ und $\widehat{v}_{k+1}^{(1)}\in C^{\infty}(\mathbb{S}\times \overline{\mathbb{B}},\mathbb{R}^q)$, so dass:
\begin{align*}
\partial_i F_0\cdot \widehat{u}_{k+1}^{(1)}&=0\hspace{0.5cm}\text{f\"ur }2\leq i\leq n \\
\left(\partial_1 F_0+\frac{1}{\varrho} \partial_t u_1 \right)\cdot\widehat{u}_{k+1}^{(1)}&=0\\
\partial_i \partial_j F_0\cdot \widehat{u}_{k+1}^{(1)}&=\frac{h_{ij}^{(1)}}{2}\hspace{0.5cm}\text{f\"ur }2\leq i\leq j\leq n\\
\left(\partial_1 \partial_i F_0+\frac{1}{\varrho}\partial_t\partial_i u_1 \right)\cdot \widehat{u}_{k+1}^{(1)}&=\frac{h_{1i}^{(1)}}{2}\hspace{0.5cm}\text{f\"ur }2\leq i\leq n\\
\left[\partial_1^2 F_0+\frac{2}{\varrho}\partial_t \partial_1 u_1 +\frac{1}{\varrho}\partial_t\left(\frac{\partial_t \widehat{u}_2}{\varrho} \right)\right]\cdot \widehat{u}_{k+1}^{(1)}&=\frac{h_{11}^{(1)}}{2}\\
\partial_t v_1\cdot \widehat{u}_{k+1}^{(1)}&=0\\
\partial_t^2 v_1\cdot \widehat{u}_{k+1}^{(1)}&=0
\end{align*}

und:
\begin{align*}
\partial_i F_0\cdot \widehat{v}_{k+1}^{(1)}&=- \partial_i u_1\cdot \widehat{u}_{k+1}^{(1)}\hspace{0.5cm}\text{f\"ur }2\leq i\leq n\\
\left(\partial_1 F_0+\frac{1}{\varrho} \partial_t u_1 \right)\cdot\widehat{v}_{k+1}^{(1)}&=-\left(\partial_1 u_1+\frac{1}{\varrho}\partial_t u_2\right)\cdot \widehat{u}_{k+1}^{(1)}\\
\partial_t v_1\cdot \widehat{v}_{k+1}^{(1)}&=0\\
\partial_t^2 v_1\cdot \widehat{v}_{k+1}^{(1)}&=0
\end{align*}

woraus, f\"ur $u_{k+1}^{(1)}:=\widehat{u}_{k+1}^{(1)}+\epsilon \widehat{v}_{k+1}^{(1)}$, die G\"ultigkeit des Gleichungssystems \eqref{eq:4.70} bis \eqref{eq:4.74} f\"ur $h=h^{(1)}$ folgt. Um das System \eqref{eq:4.70} bis \eqref{eq:4.74} f\"ur den Fall $l=2$ zu l\"osen, wird folgende Vor\"uberlegung betrachtet: Die Gleichung \eqref{eq:4.74} wird mit $\epsilon$ multipliziert:
\begin{equation*}
\left[\epsilon\partial_1^2 F_0+\frac{2\epsilon}{\varrho}\partial_t \partial_1 u_1 +\frac{1}{\varrho}\partial_t\left(\frac{\partial_t u_1+\epsilon \partial_t u_2}{\varrho} \right)\right]\cdot u_{k+1}^{(2)}=\frac{\epsilon h_{11}^{(2)}}{2}+\Lambda(\epsilon^2)
\end{equation*}

Diese Gleichung ist mit \eqref{eq:4.59} aber erf\"ullt, falls die Gleichungen:
\begin{align*}
\left[\frac{1}{\varrho^2}\partial_t^2 u_1+\epsilon \partial_1^2 F_0+\frac{2\epsilon}{\varrho}\partial_t\partial_1 u_1+\frac{\epsilon}{\varrho}\partial_t\left(\frac{\partial_t \widehat{u}_2}{\varrho} \right)\right]\cdot u_{k+1}^{(2)}&=\frac{\epsilon h_{11}^{(2)}}{2}+\Lambda(\epsilon^2)\\
\partial_t v_1\cdot u_{k+1}^{(2)}&=0\\
\end{align*}

erf\"ullt sind. Wegen $|a(x)|>\frac{\delta}{2}$, f\"ur alle $x\in U_2$, sind die $\frac{n}{2}(n+3)+1$ Vektoren:
\begin{align}\label{eq:4.92}
\begin{array}{ll}
\partial F_0(x)+\frac{1}{\varrho(t)} \partial_t u_1(t,x) & \partial_i F_0(x) \hspace{0.5cm}\text{f\"ur }2\leq i\leq n \\
\partial_1 \partial_i F_0(x)+\frac{1}{\varrho(t)} \partial_t \partial_i u_1(t,x)  \hspace{0.5cm}\text{f\"ur }2\leq i\leq n & \partial_i \partial_j F_0(x)  \hspace{0.5cm}\text{f\"ur }2\leq i\leq j\leq n\\
\partial_t v_1(t,x) & \partial_t^2 u_1(t,x)
\end{array}
\end{align}

f\"ur alle $(t,x)\in \mathbb{S}\times U_2$ linear unabh\"angig. Um das Gleichungssystem \eqref{eq:4.70} bis \eqref{eq:4.74} f\"ur $h=h^{(2)}$ zu l\"osen, wird, unter Beachtung der linearen Unabh\"angigkeit der Vektoren in \eqref{eq:4.92}, das folgende Gleichungssystem gel\"ost:
\begin{align*}
\partial_i F_0\cdot \widehat{u}_{k+1}^{(2)}&=0\hspace{0.5cm}\text{f\"ur }2\leq i\leq n\\
\left(\partial_1 F_0+\frac{1}{\varrho} \partial_t u_1 \right)\cdot\widehat{u}_{k+1}^{(2)}&=0\\
\partial_i \partial_j F_0\cdot \widehat{u}_{k+1}^{(2)}&=\frac{h_{ij}^{(2)}}{2}\hspace{0.5cm}\text{f\"ur }2\leq i\leq j\leq n\\
\left(\partial_1 \partial_i F_0+\frac{1}{\varrho}\partial_t\partial_i u_1 \right)\cdot \widehat{u}_{k+1}^{(2)}&=\frac{h_{1i}^{(2)}}{2}\hspace{0.5cm}\text{f\"ur }2\leq i\leq n\\
\partial_t^2 u_1\cdot \widehat{u}_{k+1}^{(2)}&=0\\
\partial_t v_1\cdot \widehat{u}_{k+1}^{(2)}&=0
\end{align*}

sowie:
\begin{align*}
\partial_i F_0\cdot \widehat{v}_{k+1}^{(2)}&=- \partial_i u_1\cdot \widehat{u}_{k+1}^{(2)}\hspace{0.5cm}\text{f\"ur }2\leq i\leq n\\
\left(\partial_1 F_0+\frac{1}{\varrho} \partial_t u_1 \right)\cdot\widehat{v}_{k+1}^{(2)}&=-\left(\partial_1 u_1+\frac{1}{\varrho}\partial_t u_2\right)\cdot \widehat{u}_{k+1}^{(2)}\\
\frac{1}{\varrho^2}\partial_t^2 u_1\cdot \widehat{v}_{k+1}^{(2)}&=\frac{h_{11}^{(2)}}{2}-\left[ \partial_1^2 F_0+\frac{2}{\varrho}\partial_t\partial_1 u_1+\frac{1}{\varrho}\partial_t\left(\frac{\partial_t \widehat{u}_2}{\varrho} \right)\right]\cdot \widehat{u}_{k+1}^{(2)}\\
\partial_t v_1\cdot \widehat{v}_{k+1}^{(2)}&=0
\end{align*}

Damit erf\"ullt die Funktion $u_{n+1}\in C^{\infty}(I\times\mathbb{S}\times\overline{\mathbb{B}})$ mit:
\begin{equation*}
u_{k+1}:= u_{k+1}^{(1)}+u_{k+1}^{(2)}=\underbrace{\widehat{u}_{k+1}^{(1)}+\widehat{u}_{k+1}^{(2)}}_{=: \widehat{u}_{k+1}}+\epsilon \underbrace{\bigl(\widehat{v}_{k+1}^{(1)}+\widehat{v}_{k+1}^{(2)}\bigr)}_{=:\widehat{v}_{k+1}}
\end{equation*}

wegen $h=h^{(1)}+h^{(2)}$ das Gleichungssystem \eqref{eq:4.70} bis \eqref{eq:4.74}, und alle weiteren geforderten Bedingungen. Es gilt $\widehat{u}_{k+1},\widehat{v}_{k+1}\in C^{\infty}(\mathbb{S}\times\overline{\mathbb{B}},\mathbb{R}^q)$ mit $supp(\widehat{u}_{k+1}(t,\cdot))\subseteq K,\linebreak supp(\widehat{v}_{k+1}(t,\cdot))\subseteq K$ f\"ur alle $t\in\mathbb{S}$. Diese Eigenschaften \"ubertragen sich auf die Abbildung $u_{k+1}$. Nun wird $F_{k+1}\in C^{\infty}(I\times\mathbb{S}\times\overline{\mathbb{B}},\mathbb{R}^q)$ mit:
\begin{equation}\label{eq:4.78}
F_{k+1}(\epsilon,t,x):= F_k(\epsilon,t,x)+\epsilon^{k+1}u_{k+1}(\epsilon,t,x)
\end{equation}

definiert. Dann ist nach Konstruktion von $u_{k+1}$, zusammen mit der Induktionsvoraussetzung, die Bedingung $supp(F_{k+1}(\epsilon,t,\cdot)-F_0(\cdot))\subseteq K$ f\"ur alle $(\epsilon,t)\in I\times\mathbb{S}$ erf\"ullt. Ferner ist die Absch\"atzung \eqref{eq:4.69} erf\"ullt. Es wird nun gezeigt, dass die Gleichungen \eqref{eq:4.68} bis \eqref{eq:4.66} ebenfalls gelten. Dazu werden die Gleichungen \eqref{eq:4.70} bis \eqref{eq:4.74} verwendet. Vorher sei bemerkt, dass mit \eqref{eq:4.60} und \eqref{eq:4.78} die Gleichung:
\begin{equation}\label{eq:4.91}
F_{k}(\epsilon,t,x)=F_0(x)+\epsilon u_1(t,x)+\sum_{i=3}^k{\epsilon^i\,u_i(\epsilon,t,x)}
\end{equation}

f\"ur $(\epsilon,t,x)\in I\times\mathbb{S}\times \overline{\mathbb{B}}$ gilt. Zuerst wird \eqref{eq:4.68} gezeigt:
\begin{align*}
&\left|\partial_1 F_{k+1}+\frac{1}{\epsilon\varrho}\partial_t F_{k+1} \right|^2_{\mathbb{R}^q}\stackrel{\eqref{eq:4.78}}{=}\left|\partial_1 F_k+\frac{1}{\epsilon\varrho}\partial_t F_k+\epsilon^{k+1}\left(\partial_1 u_{k+1}+\frac{1}{\epsilon \varrho}\partial_t u_{k+1} \right) \right|^2\\
\stackrel{\hphantom{\eqref{eq:4.80}}}{=}&\left|\partial_1 F_k+\frac{1}{\epsilon\varrho}\partial_t F_{k}\right|^2_{\mathbb{R}^q}+2\epsilon^{k+1}\left(\partial_1 F_k+\frac{1}{\epsilon\varrho}\partial_t F_k\right)\cdot \left(\partial_1 u_{k+1}+\frac{1}{\epsilon \varrho}\partial_t u_{k+1} \right)+\Lambda(\epsilon^{2k+2})\\
\stackrel{\eqref{eq:4.81}}{=}&\left|\partial_1 F_0\right|^2_{\mathbb{R}^q}+a^4+\epsilon^{k+1}h_{11}\\
&+2\epsilon^{k+1}\left(\partial_1 F_k+\frac{1}{\epsilon\varrho}\partial_t F_k\right)\cdot \left(\partial_1 u_{k+1}+\frac{1}{\epsilon \varrho}\partial_t u_{k+1} \right)+\Lambda(\epsilon^{k+2})\\
\stackrel{\eqref{eq:4.78}}{=}&\left|\partial_1 F_0\right|^2_{\mathbb{R}^q}+a^4+\epsilon^{k+1}h_{11}\\
&+2\epsilon^{k+1}\left(\partial_1 F_0+\epsilon \partial_1 u_1+\frac{1}{\varrho}\partial_t u_1+\frac{\epsilon}{\varrho}\partial_t u_2\right)\cdot \left(\partial_1 u_{k+1}+\frac{1}{\epsilon \varrho}\partial_t u_{k+1} \right)+\Lambda(\epsilon^{k+2})\\
\stackrel{\hphantom{\eqref{eq:4.80}}}{=}&\left|\partial_1 F_0\right|^2_{\mathbb{R}^q}+a^4+\epsilon^{k+1}h_{11}+\frac{2}{\varrho}\epsilon^k\left[\left(\partial_1 F_0+\frac{1}{\varrho}\partial_t u_1 \right)+\epsilon\left(\partial_1 u_1+\frac{1}{\varrho}\partial_t u_2 \right) \right]\cdot \partial_t u_{k+1}\\
&+2\epsilon^{k+1}\left(\partial_1 F_0+\frac{1}{\varrho}\partial_t u_1 \right)\cdot \partial_1 u_{k+1}+\Lambda(\epsilon^{k+2})\\ \stackrel{\hphantom{\eqref{eq:4.80}}}{=}&\left|\partial_1 F_0\right|^2_{\mathbb{R}^q}+a^4+\epsilon^{k+1}h_{11}\\
&+\frac{2}{\varrho}\epsilon^k\underbrace{\partial_t\left[\left(\partial_1 F_0+\frac{1}{\varrho}\partial_t u_1 \right)\cdot  u_{k+1}\right]}_{\stackrel{\eqref{eq:4.71}}{=}-\epsilon\partial_t\left(\partial_1 u_1+\frac{1}{\varrho}\partial_t u_2 \right)\cdot u_{k+1}+\Lambda(\epsilon^2)}-\frac{2}{\varrho}\epsilon^k\left[\underbrace{\partial_t\partial_1 F_0}_{=0}+\partial_t\left(\frac{1}{\varrho}\partial_t u_1\right) \right]\cdot  u_{k+1}\\
&+\frac{2}{\varrho}\epsilon^{k+1}\partial_t\left[\left(\partial_1 u_1+\frac{1}{\varrho}\partial_t u_2 \right)\cdot u_{k+1}\right]-\frac{2}{\varrho}\epsilon^{k+1}\left[\partial_t\partial_1 u_1+\partial_t\left(\frac{1}{\varrho}\partial_t u_2\right) \right]\cdot u_{k+1} \\
&+2\epsilon^{k+1}\underbrace{\partial_1\left[\left(\partial_1 F_0+\frac{1}{\varrho}\partial_t u_1 \right)\cdot  u_{k+1}\right]}_{\stackrel{\eqref{eq:4.71}}{=}\Lambda(\epsilon)}-2\epsilon^{k+1}\left[\left(\partial^2_1 F_0+\frac{1}{\varrho}\partial_t\partial_1 u_1 \right)\cdot  u_{k+1}\right]+\Lambda(\epsilon^{k+2})\\
\stackrel{\hphantom{\eqref{eq:4.80}}}{=}&\left|\partial_1 F_0\right|^2_{\mathbb{R}^q}+a^4+\epsilon^{k+1}h_{11}\\
&-2\epsilon^{k+1}\left[\partial_1^2 F_0+\frac{2}{\varrho}\partial_t \partial_1 u_1+\frac{1}{\epsilon\varrho}\partial_t\left(\frac{\partial_t u_1+\epsilon \partial_t u_2}{\varrho} \right) \right]\cdot u_{k+1}+\Lambda(\epsilon^{k+2})\\
\stackrel{\eqref{eq:4.74}}{=}&\left|\partial_1 F_0\right|^2_{\mathbb{R}^q}+a^4+\epsilon^{k+2}f_{11}^{(k+1)}
\end{align*}

Nun wird \eqref{eq:4.67} gezeigt, es sei $i\in\{2,...,n\}$:
\begin{align*}
&\left(\partial_1 F_{k+1}+\frac{1}{\epsilon \varrho} \partial_t F_{k+1}\right)\cdot \partial_i F_{k+1}\\
\stackrel{\eqref{eq:4.78}}{=}&\left[\partial_1 F_k+\frac{1}{\epsilon \varrho} \partial_t F_k+\epsilon^{k+1}\left(\partial_1 u_{k+1}+\frac{1}{\epsilon \varrho} \partial_t u_{k+1}\right)\right]\cdot \left(\partial_i F_k+\epsilon^{k+1}\partial_i u_{k+1} \right)\\
\stackrel{\hphantom{\eqref{eq:4.78}}}{=}&\left(\partial_1 F_k+\frac{1}{\epsilon \varrho} \partial_t F_k \right)\cdot \partial_i F_k+\epsilon^{k+1}\left(\partial_1 u_{k+1}+\frac{1}{\epsilon \varrho} \partial_t u_{k+1}\right)\cdot \partial_i F_k\\
&+\epsilon^{k+1}\left(\partial_1 F_k+\frac{1}{\epsilon \varrho} \partial_t F_k \right)\cdot \partial_i u_{k+1}+\Lambda(\epsilon^{2k+2})\\
\stackrel{\eqref{eq:4.80}}{=}&\partial_1 F_0\cdot \partial_i F_0+\epsilon^{k+1} h_{1i}+\epsilon^{k+1}\left(\partial_1 u_{k+1}+\frac{1}{\epsilon \varrho} \partial_t u_{k+1}\right)\cdot \partial_i F_k\\
&+\epsilon^{k+1}\left(\partial_1 F_k+\frac{1}{\epsilon \varrho} \partial_t F_k \right)\cdot \partial_i u_{k+1}+\Lambda(\epsilon^{k+2})\\
\stackrel{\eqref{eq:4.78}}{=}&\partial_1 F_0\cdot \partial_i F_0+\epsilon^{k+1} h_{1i}+\epsilon^{k+1} \partial_i F_0\cdot \partial_1 u_{k+1}+\frac{\epsilon^k}{\varrho}\left( \partial_i F_0+\epsilon \partial_i u_1\right)\cdot \partial_t u_{k+1}\\
&+\epsilon^{k+1}\left(\partial_1 F_0+\frac{1}{\varrho} \partial_t u_1 \right)\cdot \partial_i u_{k+1}+\Lambda(\epsilon^{k+2})\\
\stackrel{\hphantom{\eqref{eq:4.78}}}{=}&\partial_1 F_0\cdot \partial_i F_0+\epsilon^{k+1} h_{1i}+\epsilon^{k+1} \biggl[\underbrace{\partial_1 \left(\partial_i F_0\cdot  u_{k+1}\right)}_{\stackrel{\eqref{eq:4.70}}{=}\Lambda(\epsilon)}-\partial_1\partial_i F_0\cdot  u_{k+1}\biggr]\\
&+\frac{\epsilon^k}{\varrho}\biggl[\underbrace{\partial_t\left( \partial_i F_0\cdot  u_{k+1}\right)}_{\stackrel{\eqref{eq:4.70}}{=}-\epsilon\partial_t (\partial_i u_1\cdot u_{k+1})+\Lambda(\epsilon^2)}-\underbrace{\partial_t \partial_i F_0}_{=0}\cdot  u_{k+1}\biggr]+\frac{\epsilon^{k+1}}{\varrho}\biggl[ \partial_t(\partial_i u_1\cdot u_{k+1})- \partial_t\partial_i u_1\cdot u_{k+1}\biggr]\\
&+\epsilon^{k+1}\partial_i\biggl[\underbrace{\left(\partial_1 F_0+\frac{1}{\varrho} \partial_t u_1 \right)\cdot  u_{k+1}}_{\stackrel{\eqref{eq:4.71}}{=}\Lambda(\epsilon)}\biggr]-\epsilon^{k+1}\underbrace{\left(\partial_i\partial_1 F_0+\frac{1}{\varrho} \partial_i\partial_t u_1 \right)\cdot  u_{k+1}}_{\stackrel{\eqref{eq:4.73}}{=}\frac{1}{2}h_{1i}+\Lambda(\epsilon)}+\Lambda(\epsilon^{k+2})\\
\stackrel{\hphantom{\eqref{eq:4.78}}}{=}&\partial_1 F_0\cdot \partial_i F_0+\frac{1}{2}\epsilon^{k+1} h_{1i}-\epsilon^{k+1}\underbrace{\left(\partial_1\partial_i F_0+\frac{1}{\varrho}\partial_t\partial_i u_1\right)\cdot  u_{k+1}}_{\stackrel{\eqref{eq:4.73}}{=}\frac{1}{2}h_{1i}+\Lambda(\epsilon)}+\Lambda(\epsilon^{k+2})\\
\stackrel{\hphantom{\eqref{eq:4.78}}}{=}&\partial_1 F_0\cdot \partial_i F_0+\epsilon^{k+2}f_{1i}^{(k+1)}\\
\end{align*}

Und schlie\ss{}lich \eqref{eq:4.66}, es gilt f\"ur $i,j\in\{2,...,n\}$ mit $i\leq j$:
\begin{align*}
&\partial_i F_{k+1}\cdot \partial_j F_{k+1}\stackrel{\eqref{eq:4.78}}{=}\partial_i (F_k+\epsilon^{k+1}u_{k+1})\cdot \partial_j (F_k+\epsilon^{k+1}u_{k+1})\\
\stackrel{\hphantom{\eqref{eq:4.79}}}{=}&\partial_i F_k\cdot \partial_j F_k+\epsilon^{k+1}\partial_i F_k\cdot\partial_j u_{k+1}+\epsilon^{k+1}\partial_j F_k\cdot\partial_i u_{k+1} +\Lambda(\epsilon^{2k+2})\\
\stackrel{\eqref{eq:4.79}}{=}&\partial_i F_0\cdot \partial_j F_0+\epsilon^{k+1}h_{ij}+\epsilon^{k+1}\partial_i F_k\cdot\partial_j u_{k+1}+\epsilon^{k+1}\partial_j F_k\cdot\partial_i u_{k+1} +\Lambda(\epsilon^{k+2})\\
\stackrel{\hphantom{\eqref{eq:4.79}}}{=}&\partial_i F_0\cdot \partial_j F_0+\epsilon^{k+1}h_{ij}+\epsilon^{k+1}\underbrace{\partial_j(\partial_i F_0\cdot u_{k+1})}_{\stackrel{\eqref{eq:4.70}}{=}\Lambda(\epsilon)}-\epsilon^{k+1}\underbrace{\partial_j\partial_i F_0\cdot u_{k+1}}_{\stackrel{\eqref{eq:4.72}}{=}\frac{h_{ij}}{2}+\Lambda(\epsilon)}\\
&+\epsilon^{k+1}\underbrace{\partial_i(\partial_j F_0\cdot u_{k+1})}_{\stackrel{\eqref{eq:4.70}}{=}\Lambda(\epsilon)}-\epsilon^{k+1}\underbrace{\partial_i\partial_j F_0\cdot u_{k+1}}_{\stackrel{\eqref{eq:4.72}}{=}\frac{h_{ij}}{2}+\Lambda(\epsilon)}+\Lambda(\epsilon^{k+2})\\
\stackrel{\hphantom{\eqref{eq:4.79}}}{=}&\partial_i F_0\cdot \partial_j F_0+\epsilon^{k+2}f_{ij}^{(k+1)}
\end{align*}

\end{bew}

Um nun den Beweis von Satz \ref{satz:4.1} abzuschließen wird aus der in \ref{satz:4.3} beschriebenen Funktion $F_k$ wird die Funktion $F_{\epsilon,k}$ gem\"a\ss{} \eqref{eq:4.82} konstruiert. Das bedeutet $F_{\epsilon,k}(x):= F_k(\epsilon,\beta(\epsilon^{-1}x_1),x)$, daraus 
ergibt sich die in \eqref{eq:4.1} bis \eqref{eq:4.3} beziehungsweise \eqref{eq:6.1} geforderte Abbildung $f_{ij}^{\epsilon,k}$ wie folgt: f\"ur $i,j\in\{1,...,n\}$ mit $i\leq j$ ist:
\begin{equation}\label{eq:4.102}
f_{ij}^{\epsilon,k}(x):=f_{ij}^{(k)}(\epsilon,\beta(\epsilon^{-1}x_1),x)
\end{equation}
Aus den Absch\"atzungen \eqref{f2est} und \eqref{fkest} ergibt sich \eqref{eq:6.1}. Zun\"achst folgt:
\begin{fol}\label{fol:4.2}
Gegeben sei eine freie Abbildung $F_0\in C^{\infty}(\overline{\mathbb{B}},\mathbb{R}^q)$, f\"ur $q\geq \frac{n}{2}(n+3)+5$, und eine Abbildung $a\in C^{\infty}_0(\mathbb{B})$. Dann existiert eine kompakte Menge $K\subseteq \mathbb{B}$, mit $supp(a)\subseteq K$, so dass die folgende Aussage gilt: F\"ur jedes $k\in\mathbb{N}\backslash\{0,1\}$, existiert f\"ur alle $\epsilon\in (0,1]$, eine Abbildung $F_{\epsilon,k}\in C^{\infty}(\overline{\mathbb{B}},\mathbb{R}^q)$, welche die Eigenschaften \eqref{eq:4.1} bis \eqref{eq:6.1} erf\"ullt.
\end{fol}

Es bleibt zu zeigen dass f\"ur hinreichend kleine $\epsilon\in\mathbb{R}_{>0}$, die Abbildung $F_{\epsilon,k}$, genau wie $F_0$, eine freie Abbildung ist, womit dann Satz \ref{satz:4.1} bewiesen ist:

\section{Perturbation einer freien Abbildung}
\label{sec:4.3}

\begin{lem}\label{lem:4.9}
Es sei $K\subseteq\mathbb{R}^n$ eine kompakte Menge, und es seien $e_1,...,e_k\subseteq C^{0}(K,\mathbb{R}^q)$ punktweise linear unabh\"angige Abbildungen, dann existiert ein $\epsilon_0\in\mathbb{R}_{>0}$, mit der folgenden Eigenschaft: sind Abbildungen $\delta_1,...,\delta_k\subseteq C^{0}(K,\mathbb{R}^q)$ gegeben mit:
\begin{equation*}
\max_{1\leq i \leq k}{\max_{x\in K}{|\delta_i(x)|_{\mathbb{R}^q}}\leq \epsilon_0}
\end{equation*}

dann ist, f\"ur $x\in K$, die Gramsche Determinante der Vektoren:
\begin{equation}\label{eq:4.90}
(w_1+\delta_1)^T(x),...,(w_k+\delta_k)^T(x)
\end{equation}

nach unten, gegen eine von $x$ unabh\"angige Konstante, beschr\"ankt. Insbesondere sind die Vektoren in \eqref{eq:4.90} linear unabh\"angig. 
\end{lem}

Lemma \ref{lem:4.9} folgt direkt aus der Stetigkeit der Determinante. Lemma \ref{lem:4.9} soll nun auf eine, in Folgerung \ref{fol:4.2} beschriebene, Funktion $F_{\epsilon,k}\in C^{\infty}(\overline{\mathbb{B}},\mathbb{R}^q)$ angewandt werden, um zu zeigen, dass $F_{\epsilon,k}$ f\"ur hinreichend kleine $\epsilon\in\mathbb{R}_{>0}$ eine freie Abbildung ist. Dann ist mit \eqref{eq:4.82}, \eqref{eq:4.40} und \eqref{eq:4.41} f\"ur $x\in\overline{\mathbb{B}}$ und $2\leq i\leq j\leq n$:
\begin{align}\label{eq:4.94}
\notag\partial_1 F_{\epsilon,k}(x)\stackrel{\hphantom{\eqref{eq:4.91}}}{=}&\partial_1 F_k(\epsilon,\beta(\epsilon^{-1} x_1),x)+\frac{1}{\epsilon \varrho(\beta(\epsilon^{-1}x_1))}\partial_t F_k(\epsilon,\beta(\epsilon^{-1} x_1),x)\\
\notag\stackrel{\eqref{eq:4.91}}{=}&\partial_1 F_0(x)+\frac{1}{\varrho(\beta(\epsilon^{-1}x_1))}\partial_t u_1(\beta(\epsilon^{-1} x_1),x)+\Lambda(\epsilon)\\
\notag\partial_i F_{\epsilon,k}(x)\stackrel{\hphantom{\eqref{eq:4.91}}}{=}&\partial_i F_k(\epsilon,\beta(\epsilon^{-1} x_1),x)\\
\notag\stackrel{\eqref{eq:4.91}}{=}&\partial_i F_0(x)+\Lambda(\epsilon) \\
\partial_1^2 F_{\epsilon,k}(x)\stackrel{\hphantom{\eqref{eq:4.91}}}{=}&\partial_1^2 F_k(\epsilon,\beta(\epsilon^{-1} x_1),x)+\frac{2}{\epsilon \varrho(\beta(\epsilon^{-1}x_1))}\partial_t \partial_1 F_k(\epsilon,\beta(\epsilon^{-1} x_1),x)\\
\notag&+\frac{1}{\epsilon^2 \varrho(\beta(\epsilon^{-1}x_1))}\left.\partial_t\left[\frac{1}{ \varrho(\cdot)}\partial_t F_k(\epsilon,\cdot,x)\right]\right|_{t=\beta(\epsilon^{-1} x_1)}\\
\notag\stackrel{\eqref{eq:4.91}}{=}&\partial_1^2 F_0(x)+\frac{2}{ \varrho(\beta(\epsilon^{-1}x_1))}\partial_t \partial_1 u_1(\beta(\epsilon^{-1} x_1),x)\\
\notag&+\frac{1}{\epsilon \varrho(\beta(\epsilon^{-1}x_1))}\left.\partial_t\left[\frac{\partial_t u_1(\cdot,x) +\epsilon \partial_t u_2(\epsilon,\cdot,x) }{\varrho(\cdot)}\right]\right|_{t=\beta(\epsilon^{-1} x_1)}+\Lambda(\epsilon)\\
\notag\partial_1 \partial_i F_{\epsilon,k}(x)\stackrel{\hphantom{\eqref{eq:4.91}}}{=}&\partial_1 \partial_i F_0(x)+\frac{1}{\varrho(\beta(\epsilon^{-1}x_1))}\partial_t \partial_i u_1(\beta(\epsilon^{-1} x_1),x)+\Lambda(\epsilon)\\
\notag\partial_i \partial_j F_{\epsilon,k}(x)\stackrel{\hphantom{\eqref{eq:4.91}}}{=}&\partial_i \partial_j F_0(x)+\Lambda(\epsilon)
\end{align}

Aus der linearen Unab\"angigkeit der Vektoren in \eqref{eq:4.75} und \eqref{eq:4.92} folgt, f\"ur hinreichend kleine $\epsilon\in\mathbb{R}_{>0}$ und alle $x\in\overline{\mathbb{B}}$, die lineare Unabh\"angigkeit der Vektoren in \eqref{eq:4.94}. Um dies einzusehen, kann eine Fallunterscheidung wie in \eqref{eq:4.93} durchgef\"uhrt werden.\\
Damit ist gezeigt, dass ein $\epsilon_k\in\mathbb{R}_{>0}$ existiert, so dass die Abbildung $F_{\epsilon,k}\in C^{\infty}(\overline{\mathbb{B}},\mathbb{R}^q)$ in Folgerung \ref{fol:4.2}, f\"ur $\epsilon\in (0,\epsilon_k]$, eine freie Abbildung ist, womit Satz \ref{satz:4.1} vollst\"andig bewiesen ist.\\
F\"ur den Beweis von Hauptsatz \ref{haupt1} in \autoref{chap:Kap6}, wird noch gezeigt, dass $F_{\epsilon,k}$ f\"ur hinreichend kleine $\epsilon\in\mathbb{R}_{>0}$ injektiv ist, sofern $F_0$ injektiv ist. Daf\"ur wird eine topologische Hilfsgr\"o\ss{}e, die sogenannte \textbf{Lebesgue-Zahl}\index{Lebesgue-Zahl}, eingef\"uhrt. Der folgende Sachverhalt wird in \cite[Lemma 27.5]{munkres2000topology} bewiesen, dabei wird, f\"ur eine Teilmenge $S$ eines metrischen Raumes $(M,d)$, der Ausdruck $\gls{durchmesser}:= \sup\{ d(x,y): x,y\in S \}\in [0,\infty ]$ als der \textbf{Durchmesser}\index{Durchmesser} von $S$ bezeichnet.

\begin{lem}\label{lem:4.7}
Es sei $\mathcal{X}$ eine offene \"Uberdeckung eines kompakten metrischen Raumes $(M,d)$. Dann existiert ein $\delta_L> 0$ mit der folgenden Eigenschaft: Ist $S\subseteq M$ eine Menge mit $diam(S)<\delta_L$, dann existiert ein $U\in\mathcal{X}$ so dass $S\subseteq U$.
\end{lem}

Um Missverst\"andnisse zu vermeiden, sei darauf hingewiesen, dass sich die Bezeichnung Lebesgue-Zahl, auf das $\delta_L$ bezieht. Diese Zahl ist nicht eindeutig bestimmt, und ist immer von der vorliegenden offenen \"Uberdeckung abh\"angig.

\begin{lem}\label{lem:4.10}
Ist die freie Abbildung $F_0\in C^{\infty}(\overline{\mathbb{B}},\mathbb{R}^q)$ f\"ur $q\geq \frac{n}{2}(n+3)+5$ injektiv, dann existiert ein $\epsilon_0 \in\mathbb{R}_{>0}$, so dass f\"ur alle $\epsilon\in [-\epsilon_0,\epsilon_0]$ und alle $(t_1,x),(t_2,y)\in \mathbb{S}\times\overline{\mathbb{B}}$, mit $x\neq y$ f\"ur die, in Lemma \ref{lem:4.5} konstruierte, Abbildung $u_1 \in C^{\infty}(\mathbb{S}\times\overline{\mathbb{B}},\mathbb{R}^q)$ die Bedingung:
\begin{equation*}
F_0(x)+\epsilon u_1(t_1,x)\neq F_0(y)+\epsilon u_1(t_2,y)
\end{equation*}
\end{lem}

gilt.

\begin{bew}
Es ist $u_1=a^2 v_1$, und mit \eqref{eq:4.30} ist $v_1(t,x)=\alpha_1(t)\,u(x)+\alpha_2(t)\,v(x)$, mit den, in Lemma \ref{lem:4.2} konstruierten, Abbildungen $u,v\in C^{\infty}(\overline{\mathbb{B}},\mathbb{R}^q)$. Wenn $u\equiv 0 \equiv v$ gilt, so ist nichts zu zeigen. Es sei also mindestens eine der beiden Funktionen $u$ und $v$ nicht identisch $0$. Zun\"achst wird gezeigt, dass f\"ur ein hinreichend kleines $\widehat{\epsilon}\in\mathbb {R}_{>0}$ die Abbildung $\Phi_{\widehat{\epsilon}}\in C^{\infty}([-\widehat{\epsilon},\widehat{\epsilon}]^2\times \overline{\mathbb{B}},\mathbb{R}^q)$ mit:
\begin{equation}\label{eq:4.97}
\Phi_{\widehat{\epsilon}}(\alpha,\beta,x):= F_0(x)+\alpha\, u(x)+\beta\, v(x)
\end{equation}

injektiv ist. Mit dem Fortsetzungssatz von Seeley \cite{seeley1964extension} werden Abbildungen $\widehat{F}_0, \widehat{u}, \widehat{v}\in C^{\infty}(B_{1+\epsilon}(0),\mathbb{R}^q)$ gew\"ahlt, welche die jeweiligen Abbildungen, \"uber den Rand der Menge $\overline{\mathbb{B}}$ hinaus, fortsetzen. Dann wird eine Abbildung $\Phi\in C^{\infty}(\mathbb{R}^2\times B_{1+\epsilon}(0),\mathbb{R}^q)$ wie folgt definiert:
\begin{equation}\label{eq:4.96}
\Phi(\alpha,\beta,x):= \widehat{F}_0(x)+\alpha\, \widehat{u}(x)+\beta\, \widehat{v}(x)
\end{equation}

Dann ist f\"ur jedes $x\in\overline{\mathbb{B}}$:
\begin{equation*}
D\Phi(0,0,x)=(u(x),v(x),DF_0(x))\in\mathbb{R}^{q\times n+2}
\end{equation*}

Nach Wahl von $u$ und $v$ im Beweis von Lemma \ref{lem:4.5}, folgt unmittelbar $rg(D\Phi(0,0,x))=n+2$ f\"ur alle $x\in\overline{\mathbb{B}}$. Dann existiert, mit \cite[Theorem 7.8]{lee2003introduction}, ein $\epsilon_x >0$ und eine, in $B_{1+\epsilon}(0)$ offene Umgebung $U_x$, von $x$, so dass die Abbildung $\Phi_x:= \left.\Phi \right|_{(-\epsilon_x,\epsilon_x)^2\times U_x}\in C^{\infty}([-\epsilon_x,\epsilon_x]^2\times \overline{U}_x,\mathbb{R}^q)$ injektiv ist. Unter Verwendung der Kompaktheit der Menge $\overline{\mathbb{B}}$, werden nun endlich viele, in $B_{1+\epsilon}(0)$ offene Mengen $U_1,...,U_m\subseteq B_{1+\epsilon}(0)$ und $\epsilon_1,...,\epsilon_m\in\mathbb{R}_{>0}$ gew\"ahlt, so dass:
\begin{equation*}
\overline{\mathbb{B}}\subseteq \bigcup_{i=1}^m{U_i}
\end{equation*}

und f\"ur alle $i\in \{1,...,m\}$ die Abbildung $\Phi_i:= \left.\Phi \right|_{(-\epsilon_i,\epsilon_i)^2\times U_i}\in C^{\infty}([-\epsilon_i,\epsilon_i]^2\times \overline{U}_i)$ injektiv ist. Nun wird, f\"ur jedes $i\in\{1,...,m\}$, die Menge $V_i:= U_i\cap \overline{\mathbb{B}}$ so definiert, dass das Mengensystem $\mathcal{X}:= \{V_i \}_{1\leq i\leq m}\subseteq \mathcal{P}(\overline{\mathbb{B}})$ eine offene \"Uberdeckung der Menge $\overline{\mathbb{B}}$ bildet. Es sei $\delta_L\in \mathbb{R}_{>0}$ eine Lebesgue-Zahl zur \"Uberdeckung $\mathcal{X}$, und definiere nun die in $\overline{\mathbb{B}}\times \overline{\mathbb{B}}$ offene Menge:
\begin{align}\label{eq:4.99}
\Delta^{\delta_L}:= \left\{(x,y)\in \overline{\mathbb{B}}\times \overline{\mathbb{B}}: |x-y|<\delta_L \right\}
\end{align}

Ist $(x,y)\in \Delta^{\delta_L}$, dann existiert ein $i\in\{1,...,m\}$, so dass $\{x,y\}\subseteq V_i\subseteq U_i$ gilt. Ist dann $x\neq y$, und sind $(\alpha_1,\beta_1),(\alpha_2,\beta_2)\in (-\epsilon_i,\epsilon_i)^2$, dann ist $\Phi(\alpha_1,\beta_1,x)\neq \Phi(\alpha_2,\beta_2,y)$. Nun sei $(x,y)\in \overline{\mathbb{B}}\times \overline{\mathbb{B}}$, so dass $(x,y)\notin \Delta^{\delta_L}$. Unter Verwendung der Injektivit\"at von $F_0$, und der Kompaktheit von $\overline{\mathbb{B}}\times \overline{\mathbb{B}}\backslash \Delta^{\delta_L}$, seien nun $(\alpha_1,\beta_1), (\alpha_2,\beta_2)\in\mathbb{R}^2$, so dass:
\begin{equation}\label{eq:4.95}
|\alpha_1|+|\beta_1|+|\alpha_2|+|\beta_2|<\frac{\min_{(x,y)\in \overline{\mathbb{B}}\times\overline{\mathbb{B}}\backslash \Delta^{\delta_L}}{\left|F_0(x)- F_0(y)\right|_{\mathbb{R}^q}}}{\max_{x\in\overline{\mathbb{B}}}{|u(x)|_{\mathbb{R}^q}}+\max_{x\in\overline{\mathbb{B}}}{|v(x)|_{\mathbb{R}^q}}}
\end{equation}

dann ist:
\begin{align*}
&\left|\Phi(\alpha_1,\beta_1,x)- \Phi(\alpha_2,\beta_2,y)\right|_{\mathbb{R}^q}\\
\stackrel{\eqref{eq:4.96}}{=}&\left|F_0(x)+\alpha_1\, u(x)+\beta_1\, v(x)- F_0(y)-\alpha_2\, u(y)-\beta_2\, v(y)\right|_{\mathbb{R}^q}\\
\stackrel{\hphantom{\eqref{eq:4.95}}}{\geq}&\left|F_0(x)- F_0(y)\right|_{\mathbb{R}^q}- \left|\alpha_1\, u(x)\right|_{\mathbb{R}^q}-\left|\beta_1\, v(x)\right|_{\mathbb{R}^q}-|\alpha_2\, u(y)|_{\mathbb{R}^q}-|\beta_2\, v(y)|_{\mathbb{R}^q}\\
\stackrel{\hphantom{\eqref{eq:4.95}}}{\geq}&  \min_{(x,y)\in \overline{\mathbb{B}}\times\overline{\mathbb{B}}\backslash \Delta^{\delta_L}}{\left|F_0(x)- F_0(y)\right|_{\mathbb{R}^q}}\\
&-\left(|\alpha_1|+|\alpha_2|\right)\cdot\max_{x\in\overline{\mathbb{B}}}{|u(x)|_{\mathbb{R}^q}}-\left(|\beta_1|+|\beta_2|\right)\cdot\max_{x\in\overline{\mathbb{B}}}{|v(x)|_{\mathbb{R}^q}}\\
\stackrel{\hphantom{\eqref{eq:4.95}}}{\geq}&  \min_{(x,y)\in \overline{\mathbb{B}}\times\overline{\mathbb{B}}\backslash \Delta^{\delta_L}}{\left|F_0(x)- F_0(y)\right|_{\mathbb{R}^q}}\\
&-\left(|\alpha_1|+|\alpha_2|+|\beta_1|+|\beta_2|\right)\cdot\left(\max_{x\in\overline{\mathbb{B}}}{|u(x)|_{\mathbb{R}^q}}+\max_{x\in\overline{\mathbb{B}}}{|v(x)|_{\mathbb{R}^q}}\right)\stackrel{\eqref{eq:4.95}}{>}0
\end{align*}

Die Abbildung $\Phi_{\widehat{\epsilon}}$ in \eqref{eq:4.97} ist demnach injektiv, falls:
\begin{equation*}
\widehat{\epsilon}<\frac{1}{2}\cdot\min\left\{\left\{\epsilon_i \right\}_{1\leq i \leq m}\cup\left\{\frac{1}{4}\frac{\min_{(x,y)\in \overline{\mathbb{B}}\times\overline{\mathbb{B}}\backslash \Delta^{\delta_L}}{\left|F_0(x)- F_0(y)\right|_{\mathbb{R}^q}}}{\max_{x\in\overline{\mathbb{B}}}{|u(x)|_{\mathbb{R}^q}}+\max_{x\in\overline{\mathbb{B}}}{|v(x)|_{\mathbb{R}^q}}} \right\} \right\} 
\end{equation*}

Es sei nun $a\not\equiv 0$, da sonst, wegen $u_1=a^2\, v_1$, nichts zu zeigen w\"are. Ist:
\begin{equation*}
\epsilon_0<\frac{\widehat{\epsilon}}{\max_{x\in\overline{\mathbb{B}}}{a^2(x)}}\cdot \min\left\{\frac{1}{\max_{t\in\mathbb{S}}{|\alpha_1(t)|}},\frac{1}{\max_{t\in\mathbb{S}}{|\alpha_2(t)|}}\right\}
\end{equation*}

erf\"ullt, dann ist f\"ur $\epsilon\in [-\epsilon_0,\epsilon_0]$ und $(t_1,x)\in \mathbb{S}\times \overline{\mathbb{B}}$ sowie $(t_2,y)\in \mathbb{S}\times \overline{\mathbb{B}}$ mit $x\neq y$:
\begin{align*}
F_0(x)+\epsilon u_1(t_1,x)&=F_0(x)+\underbrace{\epsilon a^2(x)\,\alpha_1(t_1)}_{\in (-\widehat{\epsilon},\widehat{\epsilon}) }\, u(x)+\underbrace{\epsilon a^2(x)\,\alpha_2(t_1)}_{\in (-\widehat{\epsilon},\widehat{\epsilon}) }\,v(x)
\intertext{und}
F_0(y)+\epsilon u_1(t_2,y)&=F_0(y)+\underbrace{\epsilon a^2(y)\,\alpha_1(t_2)}_{\in (-\widehat{\epsilon},\widehat{\epsilon}) }\, u(y)+\underbrace{\epsilon a^2(y)\,\alpha_2(t_2)}_{\in (-\widehat{\epsilon},\widehat{\epsilon}) }\, v(y)
\end{align*}

somit folgt aus der Injektivit\"at der Abbildung $\Phi_{\widehat{\epsilon}}$, dass:
\begin{equation*}
F_0(x)+\epsilon u_1(t_1,x)\neq F_0(y)+\epsilon u_1(t_2,y)
\end{equation*}

erf\"ullt ist.

\end{bew}

Die Ausf\"uhrungen in diesem Abschnitt bilden die Grundlage f\"ur den folgenden Satz. Dieser Satz entspricht \cite[Lemma 4.3.]{gunther1989einbettungssatz}.

\begin{satz}\label{satz:4.4}
Es sei $F_0\in C^{\infty}(\overline{\mathbb{B}},\mathbb{R}^q)$, f\"ur $q\geq \frac{n}{2}(n+3)+5$, eine freie Abbildung. Dann existiert f\"ur jedes $k\in\mathbb{N}\backslash\{0,1\}$ ein $\epsilon_k\in\mathbb{R}_{>0}$, so dass f\"ur jedes $u\in C^{\infty}(\overline{\mathbb{B}},\mathbb{R}^q)$ mit:
\begin{equation*}
\left\|u \right\|_{C^2(\overline{\mathbb{B}},\mathbb{R}^q)}<\epsilon_k
\end{equation*}

f\"ur alle $\epsilon\in (0,\epsilon_k]$ die folgenden Aussagen gelten:
\begin{enumerate}[(i)]
\begin{item}\label{satz:4.4.i}
Die Abbildung $F_{\epsilon,k}+u$ ist eine freie Abbildung, und f\"ur jedes $x\in\overline{\mathbb{B}}$ ist die Gramsche Determinante der Vektoren:
\begin{align*}
&\partial_i(F_{\epsilon,k}+u)(x)\hspace{0.5cm}1\leq i\leq n & &\partial_i\partial_j(F_{\epsilon,k}+u)(x)\hspace{0.5cm}1\leq i\leq j\leq n
\end{align*}

nach unten gegen eine, von $\epsilon, x$ und $u$ unabh\"angige, Konstante beschr\"ankt.

\end{item}

\begin{item}\label{satz:4.4.ii}
Ist die Abbildung $F_0$ injektiv, so ist die Abbildung $F_{\epsilon,k}+u$ auch injektiv.

\end{item}

\end{enumerate}

\end{satz}

\begin{bew}
Die Aussage \eqref{satz:4.4.i} folgt aus den Gleichungen \eqref{eq:4.94} und den sich daran anschlie\ss{}enden Bemerkungen, sowie Lemma \ref{lem:4.9}. Nun wird \eqref{satz:4.4.ii} gezeigt. Nach Konstruktion von $F_{\epsilon,k}$ existiert zu jedem $\gamma\in \mathbb{R}_{>0}$ ein $\eta(\gamma)\in \mathbb{R}_{>0}$, so dass f\"ur alle $\epsilon<\eta$ die Aussage:
\begin{equation*}
\left\Vert F_0(\cdot)+\epsilon u_1(\beta(\epsilon^{-1}pr_1(\cdot)),\cdot)-F_{\epsilon,k}(\cdot) \right\Vert_{C^1(\overline{\mathbb{B}},\mathbb{R}^q)}<\gamma
\end{equation*}

erf\"ullt ist. Zusammen mit Lemma \ref{lem:4.10} folgt die Behauptung.

\end{bew}

\chapter{Lösung des lokalen Perturbationsproblems}
\label{chap:Kap5}
\thispagestyle{fancy}

Gegeben seien offene Mengen $U_1, U_2\subseteq \mathbb{B}$ mit $\overline{U}_1\subseteq U_2$ und $\overline{U}_2 \subseteq \mathbb{B}$, eine freie Abbildung $F_0\in C^{\infty}(\overline{\mathbb{B}},\mathbb{R}^q)$, sowie ein $f\in C^{m,\alpha}(\overline{\mathbb{B}},\mathbb{R}^{\frac{n}{2}(n+1)})$ mit der Eigenschaft $supp(f)\subseteq U_1$. In diesem Kapitel wird die folgende Fragestellung untersucht: Unter welchen Voraussetzungen an $f$ existiert ein $u\in C^{m,\alpha}(\overline{\mathbb{B}},\mathbb{R}^q)$ mit $supp(u)\subseteq U_2$, so dass f\"ur alle $i,j\in\{1,...,n\}$ mit $i\leq j$ die Gleichung:
\begin{equation*}
\partial_i(F_0+u)\cdot \partial_j (F_0+u) =\partial_i F_0\cdot \partial_j F_0+f_{ij}
\end{equation*}

erf\"ullt ist? Es wird gezeigt, dass ein solches $u$ immer existiert, sofern $f$ im $C^{2,\alpha}$-Sinne klein genug ist. Die, in diesem Kapitel vorgestellte, Methode basiert auf \cite[Kapitel 5.]{gunther1989einbettungssatz}. Im gesamten Kapitel sei $\alpha\in (0,1)$ fest gew\"ahlt.


\section{Definition verschiedener Hilfsoperatoren}
\label{sec:5.1}

Zun\"achst wird, f\"ur ein $a\in C^{\infty}_0(\mathbb{B})$ und $i\in \{1,...,n\}$, der Operator:
\begin{align}\label{eq:5.2}
\begin{split}
N_i: C^{2,\alpha}(\overline{\mathbb{B}},\mathbb{R}^q)&\longrightarrow C^{0,\alpha}(\overline{\mathbb{B}})\\
N_i[a](v)&:= 2\partial_i a \Delta v\cdot v+a\Delta v\cdot\partial_i v
\end{split}
\end{align}
betrachtet. Mit $\Delta^{-1}N_i[a](v)$ wird die, nach \cite[Corollary 4.14]{gilbarg2001elliptic} eindeutig bestimmte, L\"osung $w_i\in C^{2,\alpha}(\overline{\mathbb{B}})$ der Poisson-Gleichung mit trivialer Dirichlet-Randbedingung:
\begin{align}\label{eq:5.28}
\begin{cases}
\Delta w_i=N_i[a](v)  &\text{auf }\mathbb{B}\\
\hphantom{\Delta} w_i=0  &\text{auf }\partial\mathbb{B}
\end{cases}
\end{align}

bezeichnet. F\"ur $a\in C^{\infty}_0(\mathbb{B})$ und $v\in C^{3,\alpha}(\overline{\mathbb{B}},\mathbb{R}^q) $ werden, f\"ur $1\leq i\leq j\leq n$, die Funktionen:
\begin{align}\label{eq:5.1}
\begin{split}
u^{(1)}_{ij}[a](v)&:= a\partial_i\Delta^{-1} N_j[a](v)+a\partial_j\Delta^{-1} N_i[a](v)+3\partial_i a \Delta^{-1} N_j[a](v)+3\partial_j a \Delta^{-1} N_i[a](v)\\
u^{(2)}_{ij}[a](v)&:= 4\partial_i a \partial_j a\,  v\cdot v+2a \partial_i a\, \partial_j v\cdot v+2a \partial_j a\, \partial_i v\cdot v+a^2\, \partial_i v\cdot\partial_j v
\end{split}
\end{align}
definiert. Dann gilt, f\"ur $l\in\{1,2\}$, $u_{ij}^{(l)}[a](v)\in C^{2,\alpha}(\overline{\mathbb{B}},\mathbb{R}^q)$. In den folgenden beiden Lemmata wird $\Delta u_{ij}^{(l)}[a](v)$ f\"ur $l\in\{1,2\}$ formal berechnet.

\begin{lem}\label{lem:5.1}
Es sei $a\in C^{\infty}_0(\mathbb{B})$ und $1\leq i\leq j\leq n$. Dann existiert ein Operator:

\begin{equation*}
L_{ij}[a] : C^{2,\alpha}(\overline{\mathbb{B}},\mathbb{R}^q)\longrightarrow C^{0,\alpha}(\overline{\mathbb{B}})
\end{equation*}

der Form:
\begin{align*}
L_{ij}[a](v)=&\sum_{\substack{s_1,s_2\in\mathbb{N}^n\\ |s_1|+|s_2|=3,\ |s_2|\leq 2}}{\lambda(s_1,s_2)\cdot\partial^{s_1} a\cdot\partial^{s_2}}(\Delta^{-1}N_i[a](v))\\
+&\sum_{\substack{s_1,s_2\in\mathbb{N}^n\\|s_1|+|s_2|=3,\ |s_2|\leq 2}}{\mu(s_1,s_2)\cdot\partial^{s_1} a\cdot\partial^{s_2}}(\Delta^{-1}N_j[a](v))
\end{align*}

mit $\lambda(s_1,s_2),\  \mu(s_1,s_2)\in\mathbb{R}\text{ f\"ur alle }s_1,s_2\in\mathbb{N}^n$, mit $|s_1|+|s_2|=3$ und $|s_2|\leq 2$, so dass f\"ur jedes $v\in C^{3,\alpha}(\overline{\mathbb{B}},\mathbb{R}^q)$ die Gleichung:
\begin{align*}
&\Delta u_{ij}^{(1)}[a](v)=a\,\partial_i N_j[a](v)+a \,\partial_j N_i[a](v)-L_{ij}[a](v)
\end{align*}

erf\"ullt ist.

\end{lem}

\begin{bew}

Es gilt f\"ur $i,j\in\{1,...,n\}$:
\begin{align*}
\Delta (a\partial_i\Delta^{-1} N_j[a](v))=&\Delta a\, \partial_i\Delta^{-1} N_j[a](v)+2\nabla a\cdot\nabla\partial_i\Delta^{-1} N_j[a](v)+a\,\partial_i N_j[a](v)\\
\Delta (\partial_i a \Delta^{-1} N_j[a](v))=&\Delta \partial_i a\, \Delta^{-1} N_j[a](v)+2\nabla \partial_i a\cdot\nabla \Delta^{-1} N_j[a](v)+\partial_i a\,  N_j[a](v)
\end{align*}

Aus \eqref{eq:5.1} folgt die Behauptung.

\end{bew}

\begin{lem}\label{lem:5.2}

Es sei $a\in C^{\infty}_0(\mathbb{B})$ und $1\leq i\leq j\leq n$. Dann existiert ein Operator:
\begin{equation*}
R_{ij}[a] : C^{2,\alpha}(\overline{\mathbb{B}},\mathbb{R}^q)\longrightarrow C^{0,\alpha}(\overline{\mathbb{B}})
\end{equation*}

der Form:
\begin{equation*}
R_{ij}[a](v)=\sum_{\substack{s_1,s_2,s_3,s_4\in\mathbb{N}^n\\|s_1|+|s_2|+|s_3|+|s_4|=4,\ |s_3|,|s_4|\leq 2}}{C(s_1,s_2,s_3,s_4)\cdot\partial^{s_1} a\ \partial^{s_2} a\ \partial^{s_3} v \cdot \partial^{s_4} v}
\end{equation*}

mit $C(s_1,s_2,s_3,s_4)\in\mathbb{R}\text{ f\"ur alle }s_1,s_2,s_3,s_4\in\mathbb{N}^n,\text{ mit }\sum_{l=1}^4{|s_l|}=4$ und $|s_3|,|s_4|\leq 2$, so dass f\"ur jedes $v\in C^{3,\alpha}(\overline{\mathbb{B}},\mathbb{R}^q)$ die Gleichung:
\begin{align*}
\Delta u^{(2)}_{ij}[a](v)=a\partial_i N_j[a](v)+a\partial_j N_i[a](v)+R_{ij}[a](v)
\end{align*}

erf\"ullt ist.

\end{lem} 

\begin{bew}

Es wird zun\"achst bemerkt, dass mit \eqref{eq:5.2} f\"ur alle $i,j\in\{1,..,n\}$ die Gleichung:
\begin{align}\label{eq:5.6}
\begin{split}
&a\partial_i N_j[a](v)=2a \partial_j a\, \partial_i\Delta v \cdot v+a^2\,\partial_i\Delta v\cdot \partial_j v\\
+&\sum_{\substack{s_1,s_2,s_3,s_4\in\mathbb{N}^n\\|s_1|+|s_2|+|s_3|+|s_4|=4,\ |s_3|,|s_4|\leq 2}}{C_1(s_1,s_2,s_3,s_4)\cdot\partial^{s_1} a\ \partial^{s_2} a\ \partial^{s_3} v \cdot \partial^{s_4} v}
\end{split}
\end{align}

gilt. Nun wird der Laplace-Operator auf die einzelnen Summanden von $u^{(2)}_{ij}$ angewandt, siehe hierf\"ur auch \eqref{eq:5.1}:
\begin{align}\label{eq:5.3}
\begin{split}
&\Delta (\partial_i a \partial_j a \ v\cdot v)=\Delta (\partial_i a \partial_j a )\, v\cdot v+2\nabla (\partial_i a \partial_j a )\cdot\nabla (v\cdot v)+\partial_i a \partial_j a \, \Delta(v\cdot v)\\
=&\sum_{\substack{s_1,s_2,s_3,s_4\in\mathbb{N}^n\\|s_1|+|s_2|+|s_3|+|s_4|=4,\ |s_3|,|s_4|\leq 2}}{C_2(s_1,s_2,s_3,s_4)\cdot\partial^{s_1} a\ \partial^{s_2} a\cdot \partial^{s_3} v \ \partial^{s_4} v}
\end{split}
\end{align}

\begin{align}\label{eq:5.4}
\begin{split}
&\Delta (2a \partial_i a \, \partial_j v\cdot v)=\Delta [a \partial_i a \, \partial_j (v\cdot v)]\\
=&\Delta (a \partial_i a) \, \partial_j (v\cdot v)+2 \nabla (a \partial_i a) \cdot \nabla\partial_j (v\cdot v)+ a \partial_i a \, \partial_j \Delta(v\cdot v)\\
=&2a \partial_i a \, \partial_j \Delta v\cdot v+\sum_{\substack{s_1,s_2,s_3,s_4\in\mathbb{N}^n\\|s_1|+|s_2|+|s_3|+|s_4|=4,\ |s_3|,|s_4|\leq 2}}{C_3(s_1,s_2,s_3,s_4)\cdot\partial^{s_1} a\ \partial^{s_2} a \, \partial^{s_3} v \cdot \partial^{s_4} v}
\end{split}
\end{align}

und schlie\ss{}lich:
\begin{align}\label{eq:5.5}
\begin{split}
&\Delta (a^2 \, \partial_i v\cdot \partial_j v)=\Delta (a^2) \, \partial_i v\cdot \partial_j v+2 \nabla (a^2) \cdot \nabla(\partial_i v\cdot \partial_j v)+ a^2 \, \Delta (\partial_i v\cdot \partial_j v)\\
=&a^2\, \partial_i\Delta v\cdot\partial_j v+a^2\, \partial_j\Delta v\cdot\partial_i v\\
&+\sum_{\substack{s_1,s_2,s_3,s_4\in\mathbb{N}^n\\|s_1|+|s_2|+|s_3|+|s_4|=4,\ |s_3|,|s_4|\leq 2}}{C_4(s_1,s_2,s_3,s_4)\cdot\partial^{s_1} a\, \partial^{s_2} a\, \partial^{s_3} v \cdot \partial^{s_4} v}
\end{split}
\end{align}

Aus \eqref{eq:5.3}, \eqref{eq:5.4} und \eqref{eq:5.5} folgt:
\begin{align*}
\Delta u^{(2)}_{ij}[a](v)&\stackrel{\eqref{eq:5.1}}{=}4\Delta(\partial_i a \partial_j a\,  v\cdot v)+\Delta(2a \partial_i a\, \partial_j v\cdot v)+\Delta(2a \partial_j a\, \partial_i v\cdot v)+\Delta(a^2\, \partial_i v\cdot\partial_j v)\\
\stackrel{\hphantom{\eqref{eq:5.1}}}{=}&2a \partial_i a \, \partial_j \Delta v\cdot v+2a \partial_j a \, \partial_i \Delta v\cdot v+a^2\, \partial_i\Delta v\cdot\partial_j v+a^2\, \partial_j\Delta v\cdot\partial_i v\\
&+\sum_{\substack{s_1,s_2,s_3,s_4\in\mathbb{N}^n\\|s_1|+|s_2|+|s_3|+|s_4|=4,\ |s_3|,|s_4|\leq 2}}{C(s_1,s_2,s_3,s_4)\cdot\partial^{s_1} a\, \partial^{s_2} a\, \partial^{s_3} v \cdot \partial^{s_4} v}\\
\stackrel{\eqref{eq:5.6}}{=}&a\partial_i N_j[a](v)+a\partial_j N_i[a](v)+R_{ij}[a](v)
\end{align*}

\end{bew}

F\"ur $1\leq i\leq j\leq n$ sei nun: 
\begin{equation}\label{eq:5.13}
u_{ij}[a](v):= u_{ij}^{(2)}[a](v)-u_{ij}^{(1)}[a](v)
\end{equation}

 dann ist mit Lemma \ref{lem:5.1} und Lemma \ref{lem:5.2}, f\"ur $a\in C^{\infty}_0(\mathbb{B})$ und $v\in C^{3,\alpha}(\overline{\mathbb{B}},\mathbb{R}^q)$, die Poisson-Gleichung mit trivialer Dirichlet-Randbedingung:
\begin{align}\label{eq:5.29}
\begin{cases}
\Delta u_{ij}[a](v)=M_{ij}[a](v) &\text{auf }\mathbb{B}\\
\hphantom{\Delta} u_{ij}[a](v)=0 &\text{auf }\partial\mathbb{B}
\end{cases}
\end{align}

f\"ur den Operator:
\begin{align}\label{eq:5.11}
\begin{split}
M_{ij} : C^{2,\alpha}(\overline{\mathbb{B}},\mathbb{R}^q)&\longrightarrow C^{0,\alpha}(\overline{\mathbb{B}})\\
M_{ij}[a](v)&:= L_{ij}[a](v)+R_{ij}[a](v)
\end{split}
\end{align}

erf\"ullt. F\"ur $1\leq i\leq j\leq n$ wird der Ausdruck:
\begin{equation}\label{eq:5.12}
\Delta^{-1}M_{ij}[a](v):= u_{ij}[a](v)
\end{equation}

definiert. Das folgende Lemma beinhaltet die zentrale Idee f\"ur die L\"osung des Perturbationsproblems, unter Verwendung der  Operatoren $N_{i}$ und $M_{ij}$:

\begin{lem}\label{lem:5.3}
Gegeben seien Funktionen $F_0 \in C^{2,\alpha}(\overline{\mathbb{B}},\mathbb{R}^q)$, $a\in C^{\infty}_0(\mathbb{B}), v \in C^{3,\alpha}(\overline{\mathbb{B}},\mathbb{R}^q)$ und $f \in C^{2,\alpha}(\overline{\mathbb{B}},\mathbb{R}^{\frac{n}{2}(n+1)})$, so dass in $\mathbb{B}$ die folgenden Gleichungen erf\"ullt sind:
\begin{align}
\label{eq:5.7}
\partial_i F_0\cdot v&=-a\Delta^{-1}N_i[a](v)  & &\text{f\"ur } 1\leq i\leq n\\
\label{eq:5.8}
\partial_i \partial_j F_0\cdot v&=-\frac{1}{2}f_{ij}+\frac{1}{2}\Delta^{-1}M_{ij}[a](v) & &\text{f\"ur } 1\leq i\leq  j\leq n
\end{align}

Dann gilt f\"ur die Funktion $F:= F_0+a^2 v$ in $\mathbb{B}$, f\"ur alle $i,j\in\{1,...,n\}$ mit $i\leq j$, die Gleichung:
\begin{equation}\label{eq:5.17}
\partial_i F\cdot \partial_j F=\partial_i F_0\cdot \partial_j F_0+a^2 f_{ij}
\end{equation}

\end{lem}

\begin{bew}

Mit \eqref{eq:5.8} ist f\"ur $1\leq i\leq j\leq n$:
\begin{equation}\label{eq:5.9}
2a^2\partial_i\partial_j F_0\cdot v+a^2 f_{ij}=a^2\Delta^{-1}M_{ij}[a](v)\stackrel{\eqref{eq:5.12}}{=}a^2 u_{ij}[a](v)\stackrel{\eqref{eq:5.13}}{=}a^2 [u^{(2)}_{ij}[a](v)-u^{(1)}_{ij}[a](v)]
\end{equation}

Es gilt:
\begin{align*}
&a^2 u^{(1)}_{ij}[a](v)\\
\stackrel{\eqref{eq:5.1}}{=}&a^2[a\partial_i\Delta^{-1} N_j[a](v)+a\partial_j\Delta^{-1} N_i[a](v)+3\partial_i a \Delta^{-1} N_j[a](v)+3\partial_j a \Delta^{-1} N_i[a](v)]\\
\stackrel{\hphantom{\eqref{eq:5.1}}}{=}&a^2\partial_i(a\Delta^{-1} N_j[a](v))+a^2\partial_j(a\Delta^{-1} N_i[a](v))\\
&\hspace{4cm}+2a^2\partial_i a \Delta^{-1} N_j[a](v)+2a^2\partial_j a \Delta^{-1} N_i[a](v)\\
\stackrel{\eqref{eq:5.7}}{=}&-a^2\partial_i(\partial_j F_0\cdot v)-a^2\partial_j(\partial_i F_0\cdot v)-2a\partial_i a\partial_j F_0\cdot v-2a\partial_j a\partial_i F_0\cdot v
\end{align*}

und:
\begin{align*}
a^2\,u^{(2)}_{ij}[a](v)&\stackrel{\eqref{eq:5.1}}{=}a^2[4\partial_i a \partial_j a \,v\cdot v+2a \partial_i a \, \partial_j v\cdot v+2a \partial_j a \, \partial_i v\cdot v+a^2\, \partial_i v\cdot\partial_j v]\\
&\stackrel{\hphantom{\eqref{eq:5.1}}}{=}\partial_i(a^2 v)\cdot \partial_j(a^2 v)
\end{align*}

Daraus ergibt sich mit \eqref{eq:5.9} die Gleichung:
\begin{align*}
2a^2\partial_i\partial_j F_0\cdot v+a^2 f_{ij}=&\partial_i(a^2 v)\cdot \partial_j(a^2 v)+a^2\partial_i(\partial_j F_0\cdot v)+a^2\partial_j(\partial_i F_0\cdot v)\\
&+2a\partial_i a\partial_j F_0\cdot v+2a\partial_j a\partial_i F_0\cdot v
\end{align*}

Anwendung der Produktregel auf den dritten und vierten Summanden auf der rechten Seite, und anschlie\ss{}ende Subtraktion des Terms $2a^2\partial_i\partial_j F_0\cdot v$ auf beiden Seiten ergibt:
\begin{align*}
a^2 f_{ij}=&\partial_i(a^2 v)\cdot \partial_j(a^2 v)+a^2\partial_j F_0\cdot \partial_i  v+a^2\partial_i F_0\cdot \partial_j v+2a\partial_i a\partial_j F_0\cdot v+2a\partial_j a\partial_i F_0\cdot v\\
=&\partial_i(F_0+a^2 v)\cdot \partial_j(F_0+a^2 v)-\partial_i F_0\cdot \partial_j F_0=\partial_i F \cdot \partial_j F-\partial_i F_0\cdot \partial_j F_0
\end{align*}

womit die Behauptung bewiesen ist.

\end{bew}


\section{Formulierung des Fixpunktproblems}
\label{sec:5.2}

Um das lokale Perturbationsproblem zu l\"osen, wird untersucht, ob sich das, in Lemma \ref{lem:5.3} aufgestellte Gleichungssystem, f\"ur eine gegebene freie Abbildung $F_0\in C^{\infty}(\overline{\mathbb{B}},\mathbb{R}^q)$, l\"osen l\"asst. Daf\"ur soll dieses Gleichungssystem durch die Einf\"uhrung geeigneter Operatoren als Fixpunktgleichung umgeschrieben werden:\\
Mit $A[F_0]\in C^{\infty}(\overline{\mathbb{B}},\mathbb{R}^{\frac{n}{2}(n+3)\times q})$ wird die matrixwertige Funktion bezeichnet, deren ersten $n$ Zeilen aus den Funktionen $\partial_i F_0^{\top}$ f\"ur $ i\in\{1,..,n\}$, und deren restlichen Zeilen, in lexikografischer Reihenfolge, aus den Funktionen $\partial_i\partial_j F_0^{\top}$, f\"ur $ i,j\in\{1,..,n\}$ mit $i\leq j$, bestehen. F\"ur geordnete Mengen von Funktionen $\{h_i \}_{1\leq i\leq n}\subseteq C(\overline{\mathbb{B}})$ und $\{f_{ij} \}_{1\leq i\leq j\leq n}\subseteq C(\overline{\mathbb{B}})$, werden mit $h\in C(\overline{\mathbb{B}},\mathbb{R}^n)$ und $f\in C(\overline{\mathbb{B}},\mathbb{R}^{\frac{n}{2}(n+1)})$ die spaltenvektorwertigen Funktionen mit den entsprechenden Komponentenfunktionen bezeichnet, auch hier wird f\"ur die Zweifach-Indizierung die lexikografische Reihenfolge verwendet. Mit der Cramerschen Regel ist die Abbildung:
\begin{align}\label{eq:5.41}
\begin{split}
\Theta[F_0] : \mathbb{B}&\longrightarrow \mathbb{R}^{q\times \frac{n}{2}(n+3)}\\
x&\mapsto A^{\top}[F_0](x)\cdot( A[F_0](x)  A^{\top}[F_0](x))^{-1}
\end{split}
\end{align}

genau wie $F_0$ glatt, wobei sich jede Ableitung beliebiger Ordnung stetig bis zum Rand fortsetzen lassen. Das bedeutet $\Theta[F_0]\in C^{\infty}(\overline{\mathbb{B}},\mathbb{R}^{q\times \frac{n}{2}(n+3)})$, und der lineare Operator:
\begin{align}\label{eq:5.19}
\begin{split}
E[F_0]:  C^{m,\alpha}(\overline{\mathbb{B}},\mathbb{R}^n) \times C^{m,\alpha}(\overline{\mathbb{B}},\mathbb{R}^{\frac{n}{2}(n+1)})&\longrightarrow C^{m,\alpha}(\overline{\mathbb{B}},\mathbb{R}^q) \\  
E[F_0](h,f)(x)&:=\Theta[F_0](x)\cdot\begin{pmatrix} h(x)\\ f(x) \end{pmatrix}
\end{split}
\end{align}

ist wohldefiniert. Aus der Definition des Operators $\Theta[F_0]$ folgt, f\"ur alle $x\in\mathbb{B}$, die Gleichheit $A[F_0](x)\cdot \Theta[F_0](x)=I$, wobei $I\in\mathbb{R}^{\frac{n}{2}(n+3)\times\frac{n}{2}(n+3)}$ die Einheitsmatrix bezeichnet, siehe hierf\"ur auch Lemma \ref{lem:4.6} und Folgerung \ref{fol:4.1}. Daraus folgen, f\"ur alle $x\in\mathbb{B}$, nach Definition von $E[F_0]$, die Gleichungen:
\begin{align}\label{eq:5.10}
\begin{split}
\partial_i F_0(x)\cdot E[F_0](h,f)(x)&=h_i(x)\hspace{0.5cm}\text{f\"ur }1\leq i\leq n\\
\partial_i\partial_j F_0(x)\cdot E[F_0](h,f)(x)&=f_{ij}(x)\hspace{0.5cm}\text{f\"ur }1\leq i\leq j\leq n
\end{split}
\end{align}

Werden nun die Gleichungen $\eqref{eq:5.7}$ und $\eqref{eq:5.8}$ erneut betrachtet, so liegt unter Beachtung von \eqref{eq:5.10}, f\"ur $a\in C^{\infty}_0(\mathbb{B})$ und $m\geq 2$, die Definition der folgenden Operatoren nahe:
\begin{align}
\begin{split}\label{eq:5.14}
 P: C^{m,\alpha}(\overline{\mathbb{B}},\mathbb{R}^q)&\longrightarrow C^{m,\alpha}(\overline{\mathbb{B}},\mathbb{R}^n)  \\
(P_i[a](v))_{1\leq i\leq n}&:=(a\Delta^{-1} N_i[a](v))_{1\leq i\leq n}
\end{split}
\intertext{und}
\begin{split}\label{eq:5.15}
 Q: C^{m,\alpha}(\overline{\mathbb{B}},\mathbb{R}^q)&\longrightarrow C^{m,\alpha}(\overline{\mathbb{B}},\mathbb{R}^{\frac{n}{2}(n+1)})\\
(Q_{ij}[a](v))_{1\leq i\leq j\leq n}&:=(\Delta^{-1} M_{ij}[a](v))_{1\leq i\leq j\leq n}
\end{split}
\end{align}

mit den, in \eqref{eq:5.2} und \eqref{eq:5.11} definierten, Operatoren $N_i$ und $M_{ij}$. Schlie\ss{}lich wird, f\"ur eine gegebene freie Abbildung $F_0\in C^{\infty}(\overline{\mathbb{B}},\mathbb{R}^q)$ und $a\in C^{\infty}_0(\mathbb{B})$, sowie $f\in C^{m,\alpha}(\overline{\mathbb{B}},\mathbb{R}^{\frac{n}{2}(n+1)})$, mit $m\geq 2$ der Operator:
\begin{align}\label{eq:5.18}
\begin{split}
\Phi[F_0,a,f]: C^{m,\alpha}(\overline{\mathbb{B}},\mathbb{R}^q)&\longrightarrow C^{m,\alpha}(\overline{\mathbb{B}},\mathbb{R}^q)\\
\Phi[F_0,a,f](v)&:=-E[F_0]\left(P[a](v), \frac{1}{2}f- \frac{1}{2} Q[a](v)\right)
\end{split}
\end{align}

eingef\"uhrt. Das Gleichungssystem in Lemma \ref{lem:5.3} ist, f\"ur $v\in C^{3,\alpha}(\overline{\mathbb{B}},\mathbb{R}^q)$, wegen \eqref{eq:5.10} erf\"ullt, falls die Fixpunktgleichung:
\begin{equation}\label{eq:5.16}
v=\Phi[F_0,a,f](v)
\end{equation}

erf\"ullt ist. Wie zu Beginn des Kapitels erw\"ahnt, wird $supp(f)\subseteq U_1$ angenommen. Dann kann, mit \cite[Proposition 2.26]{lee2003introduction}, ein $a\in C^{\infty}_0(\mathbb{B})$ so gew\"ahlt werden, dass $\left. a\right|_{\overline{U}_1}=1$ und $supp(a)\subseteq U_2$ erf\"ullt sind. Damit gilt $f\equiv a^2 f$, und mit \eqref{eq:5.17} ist das Problem gel\"ost. Es bleibt also die Frage zu kl\"aren, ob der, in \eqref{eq:5.18} eingef\"uhrte, Operator einen Fixpunkt besitzt. Daf\"ur werden in \autoref{sec:5.4} diverse $C^{m,\alpha}$-Absch\"atzungen f\"ur die bereits eingef\"uhrten Hilfsoperatoren, bewiesen. Anhand der Definition der Operatoren $P_i$ und $Q_{ij}$ wird deutlich, dass daf\"ur Absch\"atzungen von der L\"osung der Poissongleichung mit Dirichlet-Randbedingung, ben\"otigt werden. Diese Thematik ist Gegenstand des n\"achsten Abschnittes.


\section{Die Poissongleichung mit Dirichlet-Randbedingung}
\label{sec:5.3}

Gegeben seien Funktionen $f\in C^{0,\alpha}(\overline{\mathbb{B}})$ und $\varphi\in C^{2,\alpha}(\overline{\mathbb{B}})$. Nach \cite[Corollary 4.14]{gilbarg2001elliptic} existiert genau ein $u\in C^{2,\alpha}(\overline{\mathbb{B}})$, welches das Problem:
\begin{align}\label{eq:5.21}
\begin{cases}
\Delta u=f &\text{auf }\mathbb{B}\\
\hphantom{\Delta} u=\varphi &\text{auf }\partial \mathbb{B}
\end{cases}
\end{align}

l\"ost. Mit \cite[Theorem 6.6]{gilbarg2001elliptic} ergibt sich, zusammen mit \eqref{eq:A.3} und \eqref{eq:A.5}, die Absch\"at\-zung:
\begin{equation}\label{eq:5.20}
\left|u \right|_{C^{2,\alpha}(\overline{\mathbb{B}})}\leq C(n,\alpha)\cdot \left( \left|f \right|_{C^{0,\alpha}(\overline{\mathbb{B}})} + \left|\varphi \right|_{C^{2,\alpha}(\overline{\mathbb{B}})}\right)
\end{equation}

Wegen \eqref{eq:5.28} und \eqref{eq:5.29} gen\"ugt es, den Fall $\varphi\equiv 0$ zu untersuchen. Es wird gezeigt, dass sich die Regularit\"at der rechten Seite $f$ auf die Regularit\"at von $u$ \"ubertr\"agt. Daf\"ur wird eine \"aquivalente Charakterisierung des Sobolev-Raumes $W^{1,p}_0(\mathbb{B})$ f\"ur $p\in[1,\infty)$ verwendet. Zun\"achst wird ein sogenannter \textbf{Spuroperator}\index{Spuroperator} eingef\"uhrt, dabei wird die Randregularit\"at von $\mathbb{B}$ direkt ausgenutzt. Das folgende Lemma wird in \cite[5.5, Theorem 1]{evans1998partial} gezeigt. Der Funktionenraum $L^p(\partial\mathbb{B})$ wird dabei wie in \cite[A 6.5 (2)]{alt2012lineare} definiert.

\begin{lem}\label{lem:5.4}
Es sei $p\in [1,\infty)$, dann existiert ein beschr\"ankter linearer Operator:
\begin{equation*}
T: W^{1,p}(\mathbb{B})\longrightarrow L^p(\partial\mathbb{B})
\end{equation*}

mit den Eigenschaften:

\begin{enumerate}[(i)]

\begin{item}
$Tu=\left.u\right|_{\partial \mathbb{B}}\hspace{0.5cm}\text{f\"ur alle }u\in W^{1,p}(\mathbb{B})\cap C(\overline{\mathbb{B}})$
\end{item}

\begin{item}
Es existiert eine Konstante $C(n,p)\in\mathbb{R}$, so dass f\"ur alle $u\in W^{1,p}(\mathbb{B})$ die Absch\"at\-zung:
\begin{equation*}
\left\|Tu \right\|_{L^p(\partial(\mathbb{B}))}\leq C(n,p)\cdot\left\|u \right\|_{W^{1,p}(\mathbb{B})}
\end{equation*}

\end{item}
 gilt.
\end{enumerate}

\end{lem}

Der Sobolev-Raum $W^{1,p}_0(\mathbb{B})$ wird dann wie folgt charakterisiert:

\begin{lem}\label{lem:5.5}
Es sei $p\in [1,\infty)$, dann gilt f\"ur $u\in W^{1,p}(\mathbb{B})$ die folgende Aussage:

\begin{equation*}
u\in W^{1,p}_0(\mathbb{B}) \Longleftrightarrow Tu\equiv 0
\end{equation*}

\end{lem}
Dieses Lemma wird in \cite[5.5, Theorem 2]{evans1998partial} bewiesen.

\begin{lem}
F\"ur $m\in\mathbb{N}$  sei $f\in C^{m,\alpha}(\overline{\mathbb{B}})$ gegeben. Dann existiert genau ein $u\in {C^{m+2,\alpha}(\overline{\mathbb{B}})}$, welches das Randwertproblem \eqref{eq:5.21} f\"ur $\varphi\equiv 0$ l\"ost. Diese L\"osung erf\"ullt die Absch\"atzung:
\begin{equation}\label{eq:5.22}
\left|u \right|_{C^{m+2,\alpha}(\overline{\mathbb{B}})}\leq C(n,m,\alpha)\cdot \left|f \right|_{C^{m,\alpha}(\overline{\mathbb{B}})}
\end{equation}

\end{lem}

\begin{bew}
Zu Beginn des Abschnittes wurde der Sachverhalt f\"ur $m=0$ und beliebige $\varphi\in C^{2,\alpha}(\overline{\mathbb{B}})$ diskutiert. Mit \eqref{eq:A.6} gilt $f\in {C^{0,\alpha}(\overline{\mathbb{B}})}$. Es sei $u\in C^{2,\alpha}(\overline{\mathbb{B}})$ die eindeutig bestimmte L\"osung von \eqref{eq:5.21} f\"ur $\varphi\equiv 0$. Um die h\"ohere Regularit\"at zu zeigen, wird gezeigt, dass f\"ur $k\in\mathbb{N}$ mit $k\in[2,m+1]$, aus $u\in {C^{k,\alpha}}(\overline{\mathbb{B}})$ immer $u\in C^{k+1,\alpha}(\overline{\mathbb{B}})$ folgt. Es sei $s \in\mathbb{N}^n$ ein Multiindex der Ordnung $|s|=k-1$. Mit $\eqref{eq:5.21}$ gilt, wegen $f\in C^{m,\alpha}(\overline{\mathbb{B}})$, die Gleichung:
\begin{equation}\label{eq:5.23}
\partial^s\Delta u=\partial^s f
\end{equation}

Nun sei $\phi\in C^{\infty}_0(\mathbb{B})$ eine beliebige Testfunktion, dann folgt aus \eqref{eq:5.23} mittels partieller Integration::
\begin{align}\label{eq:5.24}
\begin{split}
&\int_{\mathbb{B}}{\partial^s f\cdot \phi \ dx}=\int_{\mathbb{B}}{\partial^s\Delta u\cdot \phi \ dx}=(-1)^{|s|}\int_{\mathbb{B}}{\Delta u\cdot \partial^s\phi \ dx}\\
=&(-1)^{|s|+1}\int_{\mathbb{B}}{\nabla u\cdot \nabla\partial^s\phi \ dx}=(-1)^{|s|+1}\int_{\mathbb{B}}{\nabla u\cdot \partial^s\nabla\phi \ dx}\\
=&-\int_{\mathbb{B}}{\partial^s\nabla u\cdot \nabla\phi \ dx}=-\int_{\mathbb{B}}{\nabla \partial^s u\cdot \nabla\phi \ dx}
\end{split}
\end{align}

Wegen $\partial^s f\in {C^{0}}(\overline{\mathbb{B}})$ gilt $\partial^s f\in L^2(\mathbb{B})$, und aus $\partial^s u\in {C^{1}}(\overline{\mathbb{B}})$ folgt $\partial^s u\in W^{1,2}(\mathbb{B})$. Somit ist $\partial^s u$ die, nach \cite[Theorem 8.3]{gilbarg2001elliptic} eindeutig bestimmte, schwache L\"osung des Problems:
\begin{align}\label{eq:5.25}
\begin{cases}
\Delta v=\partial^s f &\text{ auf }\mathbb{B}\\
\hphantom{\Delta} v=\partial^s u &\text{ auf }\partial \mathbb{B}
\end{cases}
\end{align}

Mit $u_s\in C^{2,\alpha}(\overline{\mathbb{B}})$ wird die, nach \cite[Theorem 4.3]{gilbarg2001elliptic} eindeutig bestimmte, klassische L\"osung des Problems \eqref{eq:5.25} bezeichnet. Dann gilt mit \eqref{eq:5.24} f\"ur alle $\phi\in C^{\infty}_{0}(\mathbb{B})$:
\begin{align*}
-\int_{\mathbb{B}}{\nabla \partial^s u\cdot \nabla\phi \ dx}&=\int_{\mathbb{B}}{\partial^s f\cdot \phi \ dx}=\int_{\mathbb{B}}{\Delta u_s\cdot \phi \ dx}=-\int_{\mathbb{B}}{\nabla u_s \cdot \nabla\phi \ dx}
\end{align*}

beziehungsweise:
\begin{equation*}
\int_{\mathbb{B}}{\nabla (\partial^s u-u_s)\cdot \nabla\phi \ dx}=0
\end{equation*}

Da $T\partial^s u=Tu_s$ f\"ur den, in Lemma \ref{lem:5.4} definierten, Spuroperator, ist, unter Beachtung von Lemma \ref{lem:5.5}, die Funktion $\partial^s u-u_s\in W^{1,2}_0(\mathbb{B})$ die schwache L\"osung des Problems:
\begin{align*}
\begin{cases}
\Delta v=0 &\text{ auf }\mathbb{B}\\
\hphantom{\Delta} v=0 &\text{ auf }\partial \mathbb{B}
\end{cases}
\end{align*}

Mit \cite[Corollary 8.2]{gilbarg2001elliptic} folgt $\partial^s u\equiv u_s$ auf $\mathbb{B}$. Dies impliziert $\partial^s u\in  C^{2,\alpha}(\overline{\mathbb{B}})$, beziehungsweise, da $s\in\mathbb{N}^n$ mit $|s|=k-1$ beliebig ist, $u\in C^{k+1,\alpha}(\overline{\mathbb{B}})$. Damit ist die Existenz und die Eindeutigkeit gezeigt. Die Absch\"atzung \eqref{eq:5.22} kann induktiv, unter Verwendung von \cite[Theorem 6.6]{gilbarg2001elliptic}, gezeigt werden, hierf\"ur sei auch auf den Beweis des folgenden Lemmas verwiesen.

\end{bew}

Falls $supp(f)\subseteq \mathbb{B}$ gilt, dann kann die Absch\"atzung \eqref{eq:5.22} versch\"arft werden:

\begin{lem}
F\"ur $m\in\mathbb{N}\backslash\{0\}$ sei $f\in C^{m,\alpha}(\overline{\mathbb{B}})$ mit $supp(f)\subseteq B_R(0)$, f\"ur ein $R\in (0,1)$ gegeben. Dann gilt f\"ur die eindeutig bestimmte L\"osung $u\in C^{m+2,\alpha}(\overline{\mathbb{B}})$ von \eqref{eq:5.21} mit $\varphi\equiv 0$ die Absch\"atzung:
\begin{equation}\label{eq:5.30}
\left|u \right|_{C^{m+2,\alpha}(\overline{\mathbb{B}})}\leq K(n,\alpha,R)\cdot \left|f \right|_{C^{m,\alpha}(\overline{\mathbb{B}})}+C(n,m,\alpha,R)\cdot \left|f \right|_{C^{m-1,\alpha}(\overline{\mathbb{B}})}
\end{equation} 

\end{lem}

\begin{bew}
F\"ur ein geeignetes $\epsilon\in (0,1-R)$ sei $U_0:= B_{R+\epsilon}(0)\subseteq \mathbb{B}$, sowie $U_1,...,U_k\subseteq \overline{\mathbb{B}}$ relativ offene Mengen, so dass $\partial \mathbb{B} \subseteq \bigcup_{j=1}^k{U_j}$, und $U_j\cap supp(f)=\emptyset$ sowie $U_j\cap \partial \mathbb{B}\neq\emptyset$ f\"ur alle $j\in\{1,...,k\}$ gilt. Mit \cite[Theorem 2.25]{lee2003introduction} existieren $\psi_0,...,\psi_k\subseteq C^{\infty}(\overline{\mathbb{B}})$, so dass $supp(\psi_j)\subseteq U_j$ f\"ur alle $j \in\{0,...,k\}$ und $\sum_{j=0}^k{\psi_j}\equiv 1$ gilt. Dann muss $\left. \psi_0 \right|_{supp(f)}\equiv 1$ gelten, woraus direkt $f\equiv \psi_0\, f$ folgt. Nun werde angenommen, f\"ur ein $l\in\{2,...,m+1\}$ gelte die Absch\"atzung:
\begin{equation}\label{eq:5.30a}
\left|u \right|_{C^{l,\alpha}(\overline{\mathbb{B}})}\leq K(n,\alpha,R)\cdot \left|f \right|_{C^{l-2,\alpha}(\overline{\mathbb{B}})}+C(n,l,\alpha,R)\cdot \left|f \right|_{C^{l-3,\alpha}(\overline{\mathbb{B}})}
\end{equation} 

Dann wird f\"ur alle $i\in\{1,...,n\}$ und $j\in\{0,...,k\}$ die Abbildung $v_i^j\in C^{l-1,\alpha}(\overline{\mathbb{B}})$ mit $v_i^j:= \partial_i (\psi_j\, u)$ definiert. Mit \eqref{eq:A.10} gilt:
\begin{align}\label{eq:5.202}
\notag&\left| u\right|_{C^{l+1,\alpha}(\overline{\mathbb{B}})}\leq C(n,\alpha)\cdot\left| u\right|_{C^{2,\alpha}(\overline{\mathbb{B}})}+\sum_{i=1}^n{\left| \partial_i u\right|_{C^{l,\alpha}(\overline{\mathbb{B}})}}\\
\notag\leq &C(n,\alpha)\cdot\left| u\right|_{C^{2,\alpha}(\overline{\mathbb{B}})}+\sum_{i=1}^n{\sum_{j=0}^k{\left| \partial_i (\psi_j\, u)\right|_{C^{l,\alpha}(\overline{\mathbb{B}})}}}\\
= &C(n,\alpha)\cdot\left| u\right|_{C^{2,\alpha}(\overline{\mathbb{B}})}+\sum_{i=1}^n{\sum_{j=0}^k{\left|v_i^j \right|_{C^{l,\alpha}(\overline{\mathbb{B}})}}}
\end{align}
 
F\"ur alle $i\in\{1,...,n\}$ und $j\in\{0,...,k\}$ gilt:
\begin{align}\label{eq:5.201}
\begin{split}
&\Delta v_i^j=\Delta (\partial_i (\psi_j\, u))=\Delta (\partial_i \psi_j\, u+ \psi_j\, \partial_i u)\\
\stackrel{\hphantom{\eqref{eq:5.21}}}{=}&\Delta \partial_i \psi_j\, u+2\nabla \partial_i \psi_j\cdot \nabla u+\partial_i \psi_j\, \Delta u+\Delta\psi_j\, \partial_i u+2\nabla\psi_j\cdot \nabla \partial_i u+\psi_j\, \Delta\partial_i u\\
\stackrel{\eqref{eq:5.21}}{=}&\Delta \partial_i \psi_j\, u+2\nabla \partial_i \psi_j\cdot \nabla u+\partial_i \psi_j\, f+\Delta\psi_j\, \partial_i u+2\nabla\psi_j\cdot \nabla \partial_i u+\psi_j\, \partial_i f
\end{split}
\end{align}

F\"ur $j=0$ folgt aus \eqref{eq:5.201}:
\begin{equation*}
\Delta v_i^0=\Delta \partial_i \psi_0\, u+2\nabla \partial_i \psi_0\cdot \nabla u+\partial_i \psi_0\, f+\Delta\psi_0\, \partial_i u+2\nabla\psi_0\cdot \nabla \partial_i u+ \partial_i f
\end{equation*}

Mit \eqref{eq:5.30a} ergibt daraus, unter Beachtung von \eqref{eq:A.8} und \eqref{eq:A.9}:
\begin{align*}
\left| v_i^0\right|_{C^{l,\alpha}(\overline{\mathbb{B}})}\leq & K_1(n,\alpha,R)\cdot \left| f\right|_{C^{l-1,\alpha}(\overline{\mathbb{B}})}+C_1(n,l,\alpha,R)\cdot\left| f\right|_{C^{l-2,\alpha}(\overline{\mathbb{B}})}\\
&+K_2(n,\alpha,R)\cdot\left| u\right|_{C^{l,\alpha}(\overline{\mathbb{B}})}
\end{align*}

Anwendung der Randabsch\"atzung \cite[Theorem 4.12]{gilbarg2001elliptic} auf $v_i^j$ f\"ur alle $i\in\{1,...,n\}$ und $j\in\{1,...,k\}$, ergibt mit \eqref{eq:5.202} die Absch\"atzung:
\begin{align*}
\left| u\right|_{C^{l+1,\alpha}(\overline{\mathbb{B}})}\leq & K_3(n,\alpha,R)\cdot \left| f\right|_{C^{l-1,\alpha}(\overline{\mathbb{B}})}+C_2(n,l,\alpha,R)\cdot\left| f\right|_{C^{l-2,\alpha}(\overline{\mathbb{B}})}\\
&+K_4(n,\alpha,R)\cdot\left| u\right|_{C^{l,\alpha}(\overline{\mathbb{B}})}
\end{align*}

Mit dem Maximumprinzip in \cite[Theorem 3.7]{gilbarg2001elliptic} und dem Ehrling-Lemma in \cite[U8.2]{alt2012lineare}, ergibt sich, aus der Induktionsvoraussetzung, die Absch\"atzung:
\begin{align*}
\left| u\right|_{C^{l+1,\alpha}(\overline{\mathbb{B}})}\leq & K(n,\alpha,R)\cdot \left| f\right|_{C^{l-1,\alpha}(\overline{\mathbb{B}})}+C(n,l,\alpha,R)\cdot\left| f\right|_{C^{l-2,\alpha}(\overline{\mathbb{B}})}
\end{align*}

womit die Behauptung gezeigt ist.

\end{bew}


\section{H\"olderabsch\"atzungen der Hilfsoperatoren}
\label{sec:5.4}

Um die, in \eqref{eq:5.14} und \eqref{eq:5.15} eingef\"uhrten, Operatoren $P$ und $Q$ abzusch\"atzen, werden wegen \eqref{eq:5.22}  und \eqref{eq:5.30} die Operatoren $N_i$ und $M_{ij}$ abgesch\"atzt. Da f\"ur $a\in C^{\infty}_0(\mathbb{B})$ und $v\in C^{2,\alpha}(\overline{\mathbb{B}},\mathbb{R}^q)$ stets $supp(N_i[a](v)),\, supp(M_{ij}[a](v)))\subseteq supp(a)\subseteq \mathbb{B}$ gilt, wird f\"ur Absch\"atzungen der Komponenten von $P$ und $Q$ vor allem die Ungleichung \eqref{eq:5.30} verwendet. Anschlie\ss{}end werden noch Eigenschaften des linearen Operators $E[F_0]$, wie zum Beispiel die Stetigkeit, bewiesen.

\begin{lem}\label{lem:5.b}
F\"ur jedes $i\in \{1,...,n\}$ gelten f\"ur den in \eqref{eq:5.2} definierten Operator $N_i$ die Absch\"atzungen:\\
F\"ur $a\in C_0^{\infty}(\mathbb{B})$ und $v_1,v_2\in C^{2,\alpha}(\overline{\mathbb{B}},\mathbb{R}^q)$:
\begin{equation}\label{eq:5.31}
\left|N_i[a](v_1) -N_i[a](v_2) \right|_{C^{0,\alpha}(\overline{\mathbb{B}})} \leq K(n,\alpha,a)\cdot \left(\left|v_1 \right|_{C^{2,\alpha}(\overline{\mathbb{B}},\mathbb{R}^q)}+\left|v_2 \right|_{C^{2,\alpha}(\overline{\mathbb{B}},\mathbb{R}^q)} \right)\left|v_1-v_2 \right|_{C^{2,\alpha}(\overline{\mathbb{B}},\mathbb{R}^q)}
\end{equation}

F\"ur $m\in\mathbb{N}$, $a\in C_0^{\infty}(\mathbb{B})$ und $v\in {C^{m+2,\alpha}(\overline{\mathbb{B}},\mathbb{R}^q)}$:
\begin{equation}\label{eq:5.32}
\left|N_i[a](v) \right|_{C^{m,\alpha}(\overline{\mathbb{B}})}\leq K(n,\alpha,a)\cdot\left|v \right|_{C^{m+2,\alpha}(\overline{\mathbb{B}},\mathbb{R}^q)}\left|v \right|_{C^{2,\alpha}(\overline{\mathbb{B}},\mathbb{R}^q)}+C(n,m,\alpha,a)\cdot\left|v \right|^2_{C^{m+1,\alpha}(\overline{\mathbb{B}},\mathbb{R}^q)}
\end{equation}

\end{lem}

\begin{bew}
Sind $a\in C^{\infty}(\mathbb{B})$, $v_1,v_2\in C^{2,\alpha}(\overline{\mathbb{B}},\mathbb{R}^q)$ und $i\in \{1,...,q\}$, dann ist nach Definition des Operators $N_i$ in \eqref{eq:5.2}:
\begin{align*}
&\left|N_i[a](v_1) -N_i[a](v_2) \right|_{C^{0,\alpha}(\overline{\mathbb{B}})}\\
\stackrel{\phantom{\eqref{eq:A.2}}}{=}&\left|2\partial_i a\, \Delta v_1\cdot v_1+a\,\Delta v_1\cdot\partial_i v_1 -2\partial_i a\, \Delta v_2\cdot v_2-a\,\Delta v_2\cdot\partial_i v_2\right|_{C^{0,\alpha}(\overline{\mathbb{B}})}\\
\stackrel{\phantom{\eqref{eq:A.2}}}{\leq}&\left|2\partial_i a\, (\Delta v_1\cdot v_1-\Delta v_2\cdot v_2)\right|_{C^{0,\alpha}(\overline{\mathbb{B}})}+\left|a\,(\Delta v_1\cdot\partial_i v_1 -\Delta v_2\cdot\partial_i v_2)\right|_{C^{0,\alpha}(\overline{\mathbb{B}})}\\
\stackrel{\eqref{eq:A.2}}{\leq} &K_1(\alpha,a)\cdot\left| \Delta v_1\cdot v_1-\Delta v_2\cdot v_2\right|_{C^{0,\alpha}(\overline{\mathbb{B}})}+K_1(\alpha,a)\cdot\left|\Delta v_1\cdot\partial_i v_1 -\Delta v_2\cdot\partial_i v_2\right|_{C^{0,\alpha}(\overline{\mathbb{B}})}\\
\stackrel{\phantom{\eqref{eq:A.2}}}{\leq} &K_1(\alpha,a)\cdot\left| \Delta (v_1-v_2)\cdot v_1+\Delta v_2\cdot (v_1-v_2)\right|_{C^{0,\alpha}(\overline{\mathbb{B}})}\\
&\hspace{3cm}+K_1(\alpha,a)\cdot\left|\Delta (v_1-v_2)\cdot\partial_i v_1 +\Delta v_2\cdot\partial_i(v_1- v_2)\right|_{C^{0,\alpha}(\overline{\mathbb{B}})}\\
\stackrel{\phantom{\eqref{eq:A.2}}}{\leq}&K_1(\alpha,a)\cdot\left| \Delta (v_1-v_2) \cdot v_1\right|_{C^{0,\alpha}(\overline{\mathbb{B}})}+K_1(\alpha,a)\cdot\left|\Delta v_2\cdot (v_1-v_2)\right|_{C^{0,\alpha}(\overline{\mathbb{B}})}\\
&+K_1(\alpha,a)\cdot\left|\Delta (v_1-v_2)\cdot\partial_i v_1\right|_{C^{0,\alpha}(\overline{\mathbb{B}})} +K_1(\alpha,a)\cdot\left|\Delta v_2\cdot\partial_i(v_1- v_2)\right|_{C^{0,\alpha}(\overline{\mathbb{B}})}\\
\stackrel{\eqref{eq:A.16}}{\leq}&K_1(\alpha,a)\cdot\left| \Delta (v_1-v_2)\right|_{C^{0,\alpha}(\overline{\mathbb{B}},\mathbb{R}^q)} \cdot \left|v_1\right|_{C^{0,\alpha}(\overline{\mathbb{B}},\mathbb{R}^q)}\\
&+K_1(\alpha,a)\cdot\left|\Delta v_2\right|_{C^{0,\alpha}(\overline{\mathbb{B}},\mathbb{R}^q)}\cdot \left|v_1-v_2\right|_{C^{0,\alpha}(\overline{\mathbb{B}},\mathbb{R}^q)}\\
&+K_1(\alpha,a)\cdot\left|\Delta (v_1-v_2)\right|_{C^{0,\alpha}(\overline{\mathbb{B}},\mathbb{R}^q)}\cdot\left|\partial_i v_1\right|_{C^{0,\alpha}(\overline{\mathbb{B}},\mathbb{R}^q)}\\
&+K_1(\alpha,a)\cdot\left|\Delta v_2\right|_{C^{0,\alpha}(\overline{\mathbb{B}},\mathbb{R}^q)}\cdot\left|\partial_i(v_1- v_2)\right|_{C^{0,\alpha}(\overline{\mathbb{B}},\mathbb{R}^q)}\\
\stackrel{\eqref{eq:A.12}}{\leq}& K(n,\alpha,a)\cdot \left(\left|v_1 \right|_{C^{2,\alpha}(\overline{\mathbb{B}},\mathbb{R}^q)}+\left|v_2 \right|_{C^{2,\alpha}(\overline{\mathbb{B}},\mathbb{R}^q)} \right)\left|v_1-v_2 \right|_{C^{2,\alpha}(\overline{\mathbb{B}},\mathbb{R}^q)}
\end{align*}

Damit ist \eqref{eq:5.31} gezeigt. Nun zu \eqref{eq:5.32}, es sei $a\in C_0^{\infty}(\mathbb{B})$ und $v\in {C^{m+2,\alpha}(\overline{\mathbb{B}},\mathbb{R}^q)}$, dann ist:
\begin{align*}
&\left|N_i[a](v)\right|_{C^{m,\alpha}(\overline{\mathbb{B}})}=\left|2\partial_i a\, \Delta v\cdot  v+a\Delta v\cdot \partial_i v\right|_{C^{m,\alpha}(\overline{\mathbb{B}})}\\
\stackrel{\hphantom{\eqref{eq:A.15}}}{\leq} & \left|2\partial_i a\, \Delta v \cdot v\right|_{C^{m,\alpha}(\overline{\mathbb{B}})}+\left|a\,\Delta v\cdot \partial_i v\right|_{C^{m,\alpha}(\overline{\mathbb{B}})}\\
\stackrel{\eqref{eq:A.15}}{\leq}&\left|2\partial_i a\, \Delta v\right|_{C^{0,\alpha}(\overline{\mathbb{B}},\mathbb{R}^q)} \left|v\right|_{C^{m,\alpha}(\overline{\mathbb{B}},\mathbb{R}^q)}+\left|2\partial_i a\, \Delta v\right|_{C^{m,\alpha}(\overline{\mathbb{B}},\mathbb{R}^q)}\left|v\right|_{C^{0,\alpha}(\overline{\mathbb{B}},\mathbb{R}^q)}\\
&+C_1(n,m,\alpha)\cdot \left|2\partial_i a\, \Delta v\right|_{C^{m-1,\alpha}(\overline{\mathbb{B}},\mathbb{R}^q))} \left|v\right|_{C^{m-1,\alpha}(\overline{\mathbb{B}},\mathbb{R}^q)}\\
&+\left|a\,\Delta v\right|_{C^{0,\alpha}(\overline{\mathbb{B}},\mathbb{R}^q)}\left|\partial_i v\right|_{C^{m,\alpha}(\overline{\mathbb{B}},\mathbb{R}^q))}+\left|a\,\Delta v\right|_{C^{m,\alpha}(\overline{\mathbb{B}},\mathbb{R}^q)}\left|\partial_i v\right|_{C^{0,\alpha}(\overline{\mathbb{B}},\mathbb{R}^q)}\\
&+C_1(n,m,\alpha)\cdot \left|a\,\Delta v\right|_{C^{m-1,\alpha}(\overline{\mathbb{B}},\mathbb{R}^q)} \left|\partial_i v\right|_{C^{m-1,\alpha}(\overline{\mathbb{B}},\mathbb{R}^q)}\\
\stackrel{\eqref{eq:A.13}}{\leq}&\left|2\partial_i a\right|_{C^{0,\alpha}(\overline{\mathbb{B}})} \left|\Delta v\right|_{C^{0,\alpha}(\overline{\mathbb{B}},\mathbb{R}^q)} \left|v\right|_{C^{m,\alpha}(\overline{\mathbb{B}},\mathbb{R}^q)}\\
&+\left|2\partial_i a\right|_{C^{0,\alpha}(\overline{\mathbb{B}})}  \left|\Delta v\right|_{C^{m,\alpha}(\overline{\mathbb{B}},\mathbb{R}^q)}  \left|v\right|_{C^{0,\alpha}(\overline{\mathbb{B}},\mathbb{R}^q)}+\left|2\partial_i a\right|_{C^{m,\alpha}(\overline{\mathbb{B}})}  \left|\Delta v\right|_{C^{0,\alpha}(\overline{\mathbb{B}},\mathbb{R}^q)}  \left|v\right|_{C^{0,\alpha}(\overline{\mathbb{B}},\mathbb{R}^q)}\\
&+C_2(n,m,\alpha)\cdot \left|2\partial_i a\right|_{C^{m-1,\alpha}(\overline{\mathbb{B}})} \left|\Delta v\right|_{C^{m-1,\alpha}(\overline{\mathbb{B}},\mathbb{R}^q)}  \left|v\right|_{C^{0,\alpha}(\overline{\mathbb{B}},\mathbb{R}^q)}\\
&+C_3(n,m,\alpha)\cdot \left|2\partial_i a\right|_{C^{m-1,\alpha}(\overline{\mathbb{B}})} \left|\Delta v\right|_{C^{m-1,\alpha}(\overline{\mathbb{B}},\mathbb{R}^q))}  \left|v\right|_{C^{m-1,\alpha}(\overline{\mathbb{B}},\mathbb{R}^q)}\\
&+\left|a\right|_{C^{0,\alpha}(\overline{\mathbb{B}})} \left|\Delta v\right|_{C^{0,\alpha}(\overline{\mathbb{B}},\mathbb{R}^q)}  \left|\partial_i v\right|_{C^{m,\alpha}(\overline{\mathbb{B}},\mathbb{R}^q))}\\
&+\left|a\right|_{C^{0,\alpha}(\overline{\mathbb{B}})} \left|\Delta v\right|_{C^{m,\alpha}(\overline{\mathbb{B}},\mathbb{R}^q)}  \left|\partial_i v\right|_{C^{0,\alpha}(\overline{\mathbb{B}},\mathbb{R}^q)}+\left|a\right|_{C^{m,\alpha}(\overline{\mathbb{B}})} \left|\Delta v\right|_{C^{0,\alpha}(\overline{\mathbb{B}},\mathbb{R}^q)}  \left|\partial_i v\right|_{C^{0,\alpha}(\overline{\mathbb{B}},\mathbb{R}^q)}\\
&+C_2(n,m,\alpha)\cdot \left|a\right|_{C^{m-1,\alpha}(\overline{\mathbb{B}},\mathbb{R}^q)}\left|\Delta v\right|_{C^{m-1,\alpha}(\overline{\mathbb{B}},\mathbb{R}^q)}  \left|\partial_i v\right|_{C^{0,\alpha}(\overline{\mathbb{B}},\mathbb{R}^q)}\\
&+C_3(n,m,\alpha)\cdot \left|a\right|_{C^{m-1,\alpha}(\overline{\mathbb{B}})}\left|\Delta v\right|_{C^{m-1,\alpha}(\overline{\mathbb{B}},\mathbb{R}^q)} \left|\partial_i v\right|_{C^{m-1,\alpha}(\overline{\mathbb{B}},\mathbb{R}^q)}\\
\stackrel{\hphantom{\eqref{eq:A.15}}}{\leq} &K(n,\alpha,a)\cdot \left|v \right|_{C^{m+2,\alpha}(\overline{\mathbb{B}},\mathbb{R}^q)}\left|v \right|_{C^{2,\alpha}(\overline{\mathbb{B}},\mathbb{R}^q)}+C(n,m,\alpha,a)\cdot \left|v \right|^2_{C^{m+1,\alpha}(\overline{\mathbb{B}},\mathbb{R}^q)}
\end{align*}

\end{bew} 

\begin{lem}\label{lem:5.c}
F\"ur $i,j\in \{1,...,n\}$ mit $i\leq j$ gelten f\"ur den in \eqref{eq:5.11} definierten Operator $M_{ij}$ die Absch\"atzungen:\\
F\"ur $a\in C_0^{\infty}(\mathbb{B})$ und $v_1,v_2\in C^{2,\alpha}(\overline{\mathbb{B}},\mathbb{R}^q)$:
\begin{align}\label{eq:5.33}
\begin{split}
&\left|M_{ij}[a](v_1) -M_{ij}[a](v_2) \right|_{C^{0,\alpha}(\overline{\mathbb{B}})}\\
 &\hspace{3cm}\leq K(n,\alpha,a)\cdot \left(\left|v_1 \right|_{C^{2,\alpha}(\overline{\mathbb{B}},\mathbb{R}^q)}+\left|v_2 \right|_{C^{2,\alpha}(\overline{\mathbb{B}},\mathbb{R}^q)} \right)\left|v_1-v_2 \right|_{C^{2,\alpha}(\overline{\mathbb{B}},\mathbb{R}^q)}
\end{split}
\end{align}

F\"ur $m\in\mathbb{N}$, $a\in C_0^{\infty}(\mathbb{B})$ und $v\in {C^{m+2,\alpha}(\overline{\mathbb{B}},\mathbb{R}^q)}$:
\begin{align}\label{eq:5.34}
\begin{split}
\left|M_{ij}[a](v) \right|_{C^{m,\alpha}(\overline{\mathbb{B}})}\leq &K(n,\alpha,a)\cdot\left|v \right|_{C^{m+2,\alpha}(\overline{\mathbb{B}},\mathbb{R}^q)}\left|v \right|_{C^{2,\alpha}(\overline{\mathbb{B}},\mathbb{R}^q)}\\ 
&\hspace{3cm}+C(n,m,\alpha,a)\cdot\left|v \right|^2_{C^{m+1,\alpha}(\overline{\mathbb{B}},\mathbb{R}^q)}
\end{split}
\end{align}

\end{lem}

\begin{bew}
Der Operator $M_{ij}[a]$ ist die Summe des in Lemma \ref{lem:5.1} definierten Operators $L_{ij}[a]$, und des in Lemma \ref{lem:5.2} definierten Operators $R_{ij}[a]$. Der Operator $L_{ij}[a]$ ist eine Summe von Termen der Form $\partial^{s_1} a \ \partial^{s_2} (\Delta^{-1}N_l)$ mit $|s_1|+|s_2|=3$, $|s_2|\leq 2$ und $l\in\{i,j\}$. Demzufolge wird, um \eqref{eq:5.33} zu zeigen, f\"ur $v_1,v_2\in C^{2,\alpha}(\overline{\mathbb{B}},\mathbb{R}^q)$ zun\"achst die Differenz:
\begin{equation*}
\partial^{s_1} a \ \partial^{s_2} (\Delta^{-1}N_l[a](v_1))-\partial^{s_1} a \ \partial^{s_2} (\Delta^{-1}N_l[a](v_2))
\end{equation*}

in der $\left|\cdot \right|_{C^{0,\alpha}(\overline{\mathbb{B}})}$-Norm abgesch\"atzt:
\begin{align*}
&\left|\partial^{s_1} a \ \partial^{s_2} (\Delta^{-1}N_l[a](v_1))-\partial^{s_1} a \ \partial^{s_2} (\Delta^{-1}N_l[a](v_2)) \right|_{C^{0,\alpha}(\overline{\mathbb{B}})}\\
=&\left|\partial^{s_1} a \ \partial^{s_2} [(\Delta^{-1}N_l[a](v_1))-(\Delta^{-1}N_l[a](v_2))]  \right|_{C^{0,\alpha}(\overline{\mathbb{B}})}\\
\stackrel{\eqref{eq:A.2}}{\leq}&\left|\partial^{s_1} a \right|_{C^{0,\alpha}(\overline{\mathbb{B}})} \cdot \left| \partial^{s_2} [(\Delta^{-1}N_l[a](v_1))-(\Delta^{-1}N_l[a](v_2))]  \right|_{C^{0,\alpha}(\overline{\mathbb{B}})}\\
\stackrel{\hphantom{\eqref{eq:A.2}}}{\leq}&K_1(\alpha,a)\cdot \left| \partial^{s_2} \Delta^{-1}[N_l[a](v_1)-N_l[a](v_2)]  \right|_{C^{0,\alpha}(\overline{\mathbb{B}})}\\
\stackrel{\hphantom{\eqref{eq:A.2}}}{\leq}&K_2(n,\alpha,a)\cdot \left|  \Delta^{-1}[N_l[a](v_1)-N_l[a](v_2)]  \right|_{C^{2,\alpha}(\overline{\mathbb{B}})}\\
\stackrel{\eqref{eq:5.20}}{\leq}&K_3(n,\alpha,a)\cdot \left|  N_l[a](v_1)-N_l[a](v_2)  \right|_{C^{0,\alpha}(\overline{\mathbb{B}})}\\
\stackrel{\eqref{eq:5.31}}{\leq}& K_4(n,\alpha,a)\cdot \left(\left|v_1 \right|_{C^{2,\alpha}(\overline{\mathbb{B}},\mathbb{R}^q)}+\left|v_2 \right|_{C^{2,\alpha}(\overline{\mathbb{B}},\mathbb{R}^q)} \right)\left|v_1-v_2 \right|_{C^{2,\alpha}(\overline{\mathbb{B}},\mathbb{R}^q)}
\end{align*}

Nun wird der Ausdruck $R_{ij}[a](v)$ untersucht. Der Operator ist Summe von Termen der Form:
\begin{equation*}
\partial^{s_1} a  \partial^{s_2} a \, \partial^{s_3}v\cdot \partial^{s_4}v
\end{equation*}

mit $\sum_{l=1}^4{|s_i|}=4$ und $|s_3|,|s_4|\leq 2$, also wird die Differenz:
\begin{equation*}
\partial^{s_1} a  \partial^{s_2} a \ \partial^{s_3}v_1\cdot \partial^{s_4}v_1-\partial^{s_1} a  \partial^{s_2} a \, \partial^{s_3}v_2\cdot \partial^{s_4}v_2
\end{equation*}

in der $\left|\cdot \right|_{C^{0,\alpha}(\overline{\mathbb{B}})}$-Norm abgesch\"atzt:
\begin{align*}
&\left|\partial^{s_1} a  \partial^{s_2} a \, \partial^{s_3}v_1\cdot \partial^{s_4}v_1-\partial^{s_1} a  \partial^{s_2} a \, \partial^{s_3}v_2\cdot \partial^{s_4}v_2 \right|_{C^{0,\alpha}(\overline{\mathbb{B}})}\\
\stackrel{\hphantom{\eqref{eq:A.9}}}{=}&\left|\partial^{s_1} a  \partial^{s_2} a \, [\partial^{s_3}v_1\cdot \partial^{s_4}v_1- \partial^{s_3}v_2\cdot \partial^{s_4}v_2] \right|_{C^{0,\alpha}(\overline{\mathbb{B}})}\\
\stackrel{\eqref{eq:A.9}}{\leq}&\left|\partial^{s_1} a \right|_{C^{0,\alpha}(\overline{\mathbb{B}})}\ \left| \partial^{s_2} a\right|_{C^{0,\alpha}(\overline{\mathbb{B}})} \ \left|\partial^{s_3}v_1\cdot \partial^{s_4}v_1- \partial^{s_3}v_2\cdot \partial^{s_4}v_2 \right|_{C^{0,\alpha}(\overline{\mathbb{B}})}\\
\stackrel{\hphantom{\eqref{eq:A.9}}}{\leq}& K_1(n,\alpha,a) \cdot \left|\partial^{s_3}v_1\cdot \partial^{s_4}v_1- \partial^{s_3}v_2\cdot \partial^{s_4}v_2 \right|_{C^{0,\alpha}(\overline{\mathbb{B}})}\\
\stackrel{\hphantom{\eqref{eq:A.9}}}{=}& K_1(n,\alpha,a) \cdot \left|\partial^{s_3}(v_1-v_2)\cdot \partial^{s_4}v_1+\partial^{s_3}v_2\cdot \partial^{s_4}(v_1-v_2) \right|_{C^{0,\alpha}(\overline{\mathbb{B}})}\\
\stackrel{\hphantom{\eqref{eq:A.9}}}{\leq} & K_1(n,\alpha,a) \cdot \left|\partial^{s_3}(v_1-v_2)\cdot \partial^{s_4}v_1\right|_{C^{0,\alpha}(\overline{\mathbb{B}})}+\left|\partial^{s_3}v_2\cdot \partial^{s_4}(v_1-v_2) \right|_{C^{0,\alpha}(\overline{\mathbb{B}})}\\
\stackrel{\eqref{eq:A.16}}{\leq}& K_2(n,\alpha,a) \cdot \left|\partial^{s_3}(v_1-v_2)\right|_{C^{0,\alpha}(\overline{\mathbb{B}},\mathbb{R}^q)}\left|\partial^{s_4}v_1\right|_{C^{0,\alpha}(\overline{\mathbb{B}},\mathbb{R}^q)}\\
&+\left|\partial^{s_3}v_2\right|_{C^{0,\alpha}(\overline{\mathbb{B}},\mathbb{R}^q)} \left|\partial^{s_4}(v_1-v_2) \right|_{C^{0,\alpha}(\overline{\mathbb{B}},\mathbb{R}^q)}\\
\stackrel{\eqref{eq:A.12}}{\leq}&K_3(n,\alpha,a)\cdot \left(\left|v_1\right|_{C^{2,\alpha}(\overline{\mathbb{B}},\mathbb{R}^q)}+\left|v_2\right|_{C^{2,\alpha}(\overline{\mathbb{B}},\mathbb{R}^q)}\right) \cdot \left|v_1-v_2\right|_{C^{2,\alpha}(\overline{\mathbb{B}},\mathbb{R}^q)}
\end{align*}

Damit ist \eqref{eq:5.33} gezeigt. Nun zu \eqref{eq:5.34}, es sei $a\in C_0^{\infty}(\mathbb{B})$ und $v\in {C^{m+2,\alpha}(\overline{\mathbb{B}},\mathbb{R}^q)}$. Die Summanden von $L_{ij}[a](v)$ und $R_{ij}[a](v)$ werden einzeln abgesch\"atzt. Zuerst wird ein Summand in $L_{ij}[a](v)$ untersucht: Es seien $s_1, s_2\in\mathbb{N}^n$ mit $|s_1|+|s_2|=3$, $|s_2|\leq 2$ sowie $l\in\{i,j\}$. Zun\"achst wird der Fall $m=0$ betrachtet:
\begin{align*}
&\left|\partial^{s_1} a \, \partial^{s_2} (\Delta^{-1}N_l[a](v))\right|_{C^{0,\alpha}(\overline{\mathbb{B}})}\stackrel{\eqref{eq:A.2}}{\leq} \left|\partial^{s_1} a\right|_{C^{0,\alpha}(\overline{\mathbb{B}})} \left| \partial^{s_2} (\Delta^{-1}N_l[a](v))\right|_{C^{0,\alpha}(\overline{\mathbb{B}})}\\
\stackrel{\hphantom{\eqref{eq:A.2}}}{\leq} & K_1(n,\alpha,a) \cdot \left| \Delta^{-1}N_l[a](v)\right|_{C^{|s_2|,\alpha}(\overline{\mathbb{B}})}\stackrel{\eqref{eq:A.6}}{\leq} K_2(n,\alpha,a) \cdot \left| \Delta^{-1}N_l[a](v)\right|_{C^{2,\alpha}(\overline{\mathbb{B}})}\\
\stackrel{\eqref{eq:5.20}}{\leq} &K_3(n,\alpha,a) \cdot \left| N_l[a](v)\right|_{C^{0,\alpha}(\overline{\mathbb{B}})}\\
\stackrel{\eqref{eq:5.32}}{\leq} & K_4(n,\alpha,a) \cdot \left| v\right|^2_{C^{2,\alpha}(\overline{\mathbb{B}},\mathbb{R}^q)}+C_1(n,m,\alpha,a)\cdot \left| v\right|_{C^{1,\alpha}(\overline{\mathbb{B}},\mathbb{R}^q)}^2
\end{align*}

Nun sei $m\geq 1$, dann ist:
\begin{align}\label{eq:5.35}
\notag&\left|\partial^{s_1} a \, \partial^{s_2} (\Delta^{-1}N_l[a](v))\right|_{C^{m,\alpha}(\overline{\mathbb{B}})}\\
\notag\stackrel{\eqref{eq:A.13}}{\leq}&\left|\partial^{s_1} a \right|_{C^{0,\alpha}(\overline{\mathbb{B}})}\ \left|\partial^{s_2} (\Delta^{-1}N_l[a](v))\right|_{C^{m,\alpha}(\overline{\mathbb{B}})}+\left|\partial^{s_1} a \right|_{C^{m,\alpha}(\overline{\mathbb{B}})}\ \left|\partial^{s_2} (\Delta^{-1}N_l[a](v))\right|_{C^{0,\alpha}(\overline{\mathbb{B}})}\\
\notag&+C_1(n,m,\alpha)\cdot\left|\partial^{s_1} a \right|_{C^{m-1,\alpha}(\overline{\mathbb{B}})}\ \left|\partial^{s_2} (\Delta^{-1}N_l[a](v))\right|_{C^{m-1,\alpha}(\overline{\mathbb{B}})}\\
\stackrel{\eqref{eq:A.12}}{\leq}&K_1(n,\alpha,a)\cdot \left|\Delta^{-1}N_l[a](v)\right|_{C^{m+|s_2|,\alpha}(\overline{\mathbb{B}})}+C_2(n,m,\alpha,a)\cdot \left|\Delta^{-1}N_l[a](v)\right|_{C^{2,\alpha}(\overline{\mathbb{B}})}\\
\notag&+C_3(n,m,\alpha,a)\cdot \left|\Delta^{-1}N_l[a](v)\right|_{C^{m+1,\alpha}(\overline{\mathbb{B}})}\\
\notag\stackrel{\eqref{eq:5.22}}{\leq}&K_1(n,\alpha,a)\cdot \left|\Delta^{-1}N_l[a](v)\right|_{C^{m+|s_2|,\alpha}(\overline{\mathbb{B}})}+C_4(n,m,\alpha,a)\cdot \left|N_l[a](v)\right|_{C^{0,\alpha}(\overline{\mathbb{B}})}\\
\notag&+C_5(n,m,\alpha,a)\cdot \left|N_l[a](v)\right|_{C^{m-1,\alpha}(\overline{\mathbb{B}})}\\
\notag\stackrel{\eqref{eq:5.32}}{\leq}&K_1(n,\alpha,a)\cdot\left|\Delta^{-1}N_l[a](v)\right|_{C^{m+|s_2|,\alpha}(\overline{\mathbb{B}})}+C_{6}(n,m,\alpha,a)\cdot \left| v\right|_{C^{m+1,\alpha}(\overline{\mathbb{B}},\mathbb{R}^q)}^2
\end{align}

Ist $|s_2|\leq 1$, dann folgt, nach erneuter Anwendung von \eqref{eq:5.22} und \eqref{eq:5.32}, aus \eqref{eq:5.35}, unter Beachtung von \eqref{eq:A.12}, die Absch\"atzung:
\begin{equation*}
\left|\partial^{s_1} a \, \partial^{s_2} (\Delta^{-1}N_l[a](v))\right|_{C^{m,\alpha}(\overline{\mathbb{B}})}\leq C_{7}(n,m,\alpha,a)\cdot \left| v\right|_{C^{m+1,\alpha}(\overline{\mathbb{B}},\mathbb{R}^q)}^2
\end{equation*}

und wenn $|s_2|=2$ gilt, dann ist:
\begin{align*}
&\left|\Delta^{-1}N_l[a](v)\right|_{C^{m+|s_2|,\alpha}(\overline{\mathbb{B}})}=\left|\Delta^{-1}N_l[a](v)\right|_{C^{m+2,\alpha}(\overline{\mathbb{B}})}\\
\stackrel{\eqref{eq:5.30}}{\leq}&K_2(n,\alpha,a)\cdot\left|N_l[a](v)\right|_{C^{m,\alpha}(\mathbb{B})}+C_8(n,m,\alpha,a)\cdot\left|N_l[a](v)\right|_{C^{m-1,\alpha}(\mathbb{B})}\\
\stackrel{\eqref{eq:5.32}}{\leq}&K_3(n,\alpha,a)\cdot\left|v\right|_{C^{m+2,\alpha}(\mathbb{B},\mathbb{R}^q)}\left|v\right|_{C^{2,\alpha}(\mathbb{B},\mathbb{R}^q)}+C_{9}(n,m,\alpha,a)\cdot \left| v\right|_{C^{m+1,\alpha}(\overline{\mathbb{B}},\mathbb{R}^q)}^2
\end{align*}

und es ergibt sich mit \eqref{eq:5.35} die Absch\"atzung:
\begin{align*}
\left|\partial^{s_1} a \, \partial^{s_2} (\Delta^{-1}N_l[a](v))\right|_{C^{m,\alpha}(\overline{\mathbb{B}})}\leq &K_4(n,\alpha,a)\cdot\left|v\right|_{C^{m+2,\alpha}(\mathbb{B},\mathbb{R}^q)}\left|v\right|_{C^{2,\alpha}(\mathbb{B},\mathbb{R}^q)}\\
&+C_{10}(n,m,\alpha,a)\cdot \left| v\right|_{C^{m+1,\alpha}(\overline{\mathbb{B}},\mathbb{R}^q)}^2
\end{align*}

Und schlie\ss{}lich zur Absch\"atzung eines Summanden in $R_{ij}[a](v)$. Es seien $s_1,s_2,s_3,s_4\in\mathbb{N}^n$ mit $\sum_{l=1}^4{|s_l|}=4$ und $|s_3|,|s_4|\leq 2$, wieder wird zuerst der Fall $m=0$ untersucht:
\begin{align*}
&\left|\partial^{s_1} a \partial^{s_2} a \, \partial^{s_3} v\cdot \partial^{s_4} v \right|_{C^{0,\alpha}(\overline{\mathbb{B}})}\\
\stackrel{\eqref{eq:A.14},\, \eqref{eq:A.16}}{\leq}&K_1(n,\alpha,a) \cdot\left|  v\right|_{C^{|s_3|,\alpha}(\overline{\mathbb{B}},\mathbb{R}^q)} \left| v \right|_{C^{|s_4|,\alpha}(\overline{\mathbb{B}},\mathbb{R}^q)}\stackrel{\eqref{eq:A.12}}{\leq}K_2(n,\alpha,a) \cdot\left|  v\right|_{C^{2,\alpha}(\overline{\mathbb{B}},\mathbb{R}^q)}^2
\end{align*}

Und f\"ur $m\geq 1$:
\begin{align}\label{eq:5.36}
\notag&\left|\partial^{s_1} a \partial^{s_2} a \, \partial^{s_3} v\cdot \partial^{s_4} v \right|_{C^{m,\alpha}(\overline{\mathbb{B}})}\\
\notag\stackrel{\eqref{eq:A.13}}{\leq} & \left|\partial^{s_1} a \partial^{s_2} a\, \right|_{C^{0,\alpha}(\overline{\mathbb{B}})}  \left|\partial^{s_3} v\cdot \partial^{s_4} v \right|_{C^{m,\alpha}(\overline{\mathbb{B}})}+\left|\partial^{s_1} a \partial^{s_2} a \right|_{C^{m,\alpha}(\overline{\mathbb{B}})} \left|\partial^{s_3} v\cdot \partial^{s_4} v \right|_{C^{0,\alpha}(\overline{\mathbb{B}})}\\
\notag&+C_1(n,m,\alpha)\cdot\left|\partial^{s_1} a\ \partial^{s_2} a \right|_{C^{m-1,\alpha}(\overline{\mathbb{B}})} \cdot \left|\partial^{s_3} v\cdot \partial^{s_4} v \right|_{C^{m-1,\alpha}(\overline{\mathbb{B}})}\\
\notag\stackrel{\eqref{eq:A.16}}{\leq} &K_1(n,\alpha,a)\cdot \left|\partial^{s_3} v\cdot \partial^{s_4} v \right|_{C^{m,\alpha}(\overline{\mathbb{B}})}+C_2(n,m,\alpha,a)\cdot \left|\partial^{s_3} v\right|_{C^{0,\alpha}(\overline{\mathbb{B}},\mathbb{R}^q)}\left|\partial^{s_4} v \right|_{C^{0,\alpha}(\overline{\mathbb{B}},\mathbb{R}^q)}\\
&+C_3(n,m,\alpha,a)\cdot\left|\partial^{s_3} v\right|_{C^{m-1,\alpha}(\overline{\mathbb{B}},\mathbb{R}^q)} \left|\partial^{s_4} v \right|_{C^{m-1,\alpha}(\overline{\mathbb{B}},\mathbb{R}^q)}\\
\notag\stackrel{\eqref{eq:A.12}}{\leq} &K_1(n,\alpha,a)\cdot\left|\partial^{s_3} v\cdot \partial^{s_4} v \right|_{C^{m,\alpha}(\overline{\mathbb{B}})}+C_4(n,m,\alpha,a)\cdot\left| v\right|^2_{C^{m+1,\alpha}(\overline{\mathbb{B}},\mathbb{R}^q)}\\
\notag\stackrel{\eqref{eq:A.15}}{\leq} &K_1(n,\alpha,a)\cdot\left(\left|\partial^{s_3} v\right|_{C^{0,\alpha}(\overline{\mathbb{B}},\mathbb{R}^q)}\left|\partial^{s_4} v \right|_{C^{m,\alpha}(\overline{\mathbb{B}},\mathbb{R}^q)}+\left|\partial^{s_3} v\right|_{C^{m,\alpha}(\overline{\mathbb{B}},\mathbb{R}^q)}\left|\partial^{s_4} v \right|_{C^{0,\alpha}(\overline{\mathbb{B}},\mathbb{R}^q)}\right)\\
\notag\hphantom{\stackrel{\eqref{eq:A.15}}{\leq}} &+C_5(n,m,\alpha,a)\cdot\left|\partial^{s_3} v\right|_{C^{m-1,\alpha}(\overline{\mathbb{B}},\mathbb{R}^q)}\left|\partial^{s_4} v \right|_{C^{m-1,\alpha}(\overline{\mathbb{B}},\mathbb{R}^q)}+C_4(n,m,\alpha,a)\cdot\left| v\right|^2_{C^{m+1,\alpha}(\overline{\mathbb{B}},\mathbb{R}^q)}\\
\notag\stackrel{\eqref{eq:A.12}}{\leq} &K_2(n,\alpha,a)\cdot\left[\left| v\right|_{C^{2,\alpha}(\overline{\mathbb{B}},\mathbb{R}^q)}\left(\left|\partial^{s_3} v \right|_{C^{m,\alpha}(\overline{\mathbb{B}},\mathbb{R}^q)}+\left|\partial^{s_4} v \right|_{C^{m,\alpha}(\overline{\mathbb{B}},\mathbb{R}^q)} \right)\right]\\
\notag\hphantom{\stackrel{\eqref{eq:A.15}}{\leq}} &+C_6(n,m,\alpha,a)\cdot\left| v\right|^2_{C^{m+1,\alpha}(\overline{\mathbb{B}},\mathbb{R}^q)}
\end{align}

Ist nun $|s_3|\leq 1$ und $|s_4|\leq 1$, so folgt, durch nochmalige Anwendung von \eqref{eq:A.12}, aus \eqref{eq:5.36} die Absch\"atzung:
\begin{equation*}
\left|\partial^{s_1} a \partial^{s_2} a \, \partial^{s_3} v\cdot \partial^{s_4} v \right|_{C^{m,\alpha}(\overline{\mathbb{B}})}\leq C_7(n,m,\alpha,a)\cdot\left| v\right|^2_{C^{m+1,\alpha}(\overline{\mathbb{B}},\mathbb{R}^q)}
\end{equation*}

und wenn $|s_3|=2$ oder $|s_4|=2$ gilt, dann folgt aus \eqref{eq:5.36} die Absch\"atzung:
\begin{align*}
\left|\partial^{s_1} a \partial^{s_2} a \, \partial^{s_3} v\cdot \partial^{s_4} v \right|_{C^{m,\alpha}(\overline{\mathbb{B}})}\leq &K_3(n,\alpha,a)\cdot\left|v \right|_{C^{m+2,\alpha}(\overline{\mathbb{B}},\mathbb{R}^q)}\left|v \right|_{C^{2,\alpha}(\overline{\mathbb{B}},\mathbb{R}^q)}\\ 
&+C_7(n,m,\alpha,a)\cdot\left|v \right|^2_{C^{m+1,\alpha}(\overline{\mathbb{B}},\mathbb{R}^q)}
\end{align*}

womit \eqref{eq:5.34} gezeigt ist.

\end{bew}

\begin{lem}\label{lem:5.d}
Die in \eqref{eq:5.14} und \eqref{eq:5.15} definierten Operatoren $P$ und $Q$ erf\"ullen die Absch\"at\-zungen:\\
F\"ur $a\in C_0^{\infty}(\mathbb{B})$ und $v_1,v_2\in C^{2,\alpha}(\overline{\mathbb{B}},\mathbb{R}^q)$:
\begin{align}\label{eq:5.37}
\begin{split}
&\left|P[a](v_1)-P[a](v_2) \right|_{C^{2,\alpha}(\overline{\mathbb{B}},\mathbb{R}^n)}+\left|Q[a](v_1)-Q[a](v_2) \right|_{C^{2,\alpha}(\overline{\mathbb{B}},\mathbb{R}^{\frac{n}{2}(n+1)})}\\
&\hspace{3cm}\leq K(n,\alpha,a)\cdot\left(\left|v_1 \right|_{C^{2,\alpha}(\overline{\mathbb{B}},\mathbb{R}^q)} +\left|v_2 \right|_{C^{2,\alpha}(\overline{\mathbb{B}},\mathbb{R}^q)}\right)\cdot \left|v_1-v_2 \right|_{C^{2,\alpha}(\overline{\mathbb{B}},\mathbb{R}^q)}
\end{split}
\end{align}

F\"ur $m\in\mathbb{N}\backslash\{0,1\}$, $a\in C_0^{\infty}(\mathbb{\mathbb{B}})$ und $v\in {C^{m,\alpha}(\overline{\mathbb{\mathbb{B}}},\mathbb{R}^q)}$:
\begin{align}\label{eq:5.38}
\begin{split}
&\left|P[a](v) \right|_{C^{m,\alpha}(\overline{\mathbb{B}},\mathbb{R}^n)}+\left|Q[a](v) \right|_{C^{m,\alpha}(\overline{\mathbb{B}},\mathbb{R}^{\frac{n}{2}(n+1)})}\\
&\hspace{3cm}\leq K(n,\alpha,a)\cdot\left|v \right|_{C^{m,\alpha}(\overline{\mathbb{B}},\mathbb{R}^q)}\cdot \left|v \right|_{C^{2,\alpha}(\overline{\mathbb{B}},\mathbb{R}^q)}+C(n,m,\alpha,a)\cdot \left|v \right|_{C^{m-1,\alpha}(\overline{\mathbb{B}},\mathbb{R}^q)}^2
\end{split}
\end{align}

\end{lem}

\begin{bew}

F\"ur $a\in C_0^{\infty}(\mathbb{B})$ und $v_1,v_2\in C^{2,\alpha}(\overline{\mathbb{B}},\mathbb{R}^q)$ gilt:
\begin{align*}
&\left|P[a](v_1)-P[a](v_2) \right|_{C^{2,\alpha}(\overline{\mathbb{B}},\mathbb{R}^n)}=\sum_{l=1}^n{\left|P_i[a](v_1)-P_i[a](v_2) \right|_{C^{2,\alpha}(\overline{\mathbb{B}})}}\\
\stackrel{\eqref{eq:5.14}}{=}&\sum_{i=1}^n{\left|a\Delta^{-1} N_i[a](v_1)-a\Delta^{-1} N_i[a](v_2) \right|_{C^{2,\alpha}(\overline{\mathbb{B}})}}\\
\stackrel{\hphantom{\eqref{eq:5.14}}}{=}&\sum_{i=1}^n{\left|a\Delta^{-1} \left[ N_i[a](v_1)- N_i[a](v_2)\right] \right|_{C^{2,\alpha}(\overline{\mathbb{B}})}}\\
\stackrel{\eqref{eq:A.9}}{\leq}&K_1(n,\alpha)\cdot\sum_{i=1}^n{\left|a\right|_{C^{2,\alpha}(\overline{\mathbb{B}})}\left|\Delta^{-1}\left[ N_i[a](v_1)- N_i[a](v_2)\right] \right|_{C^{2,\alpha}(\overline{\mathbb{B}})}}\\
\stackrel{\hphantom{\eqref{eq:A.9}}}{\leq}&K_2(n,\alpha,a)\cdot\sum_{i=1}^n{\left|\Delta^{-1}\left[ N_i[a](v_1)- N_i[a](v_2)\right] \right|_{C^{2,\alpha}(\overline{\mathbb{B}})}}\\
\stackrel{\eqref{eq:5.20}}{\leq}&K_3(n,\alpha,a)\cdot\left| N_i[a](v_1)- N_i[a](v_2)\right|_{C^{0,\alpha}(\overline{\mathbb{B}})}\\
\stackrel{\eqref{eq:5.31}}{\leq}&K_4(n,\alpha,a)\cdot \left(\left|v_1 \right|_{C^{2,\alpha}(\overline{\mathbb{B}},\mathbb{R}^q)} +\left|v_2 \right|_{C^{2,\alpha}(\overline{\mathbb{B}},\mathbb{R}^q)}\right) \left|v_1-v_2 \right|_{C^{2,\alpha}(\overline{\mathbb{B}},\mathbb{R}^q)}
\end{align*}

und:
\begin{align*}
&\left|Q[a](v_1)-Q[a](v_2) \right|_{C^{2,\alpha}(\overline{\mathbb{B}},\mathbb{R}^{\frac{n}{2}(n+1)})}\leq\sum_{1\leq i\leq j\leq n}{\left|Q_{ij}[a](v_1)-Q_{ij}[a](v_2) \right|_{C^{2,\alpha}(\overline{\mathbb{B}})}}\\
\stackrel{\eqref{eq:5.15}}{=}&\sum_{1\leq i\leq j\leq n}{\left|\Delta^{-1}M_{ij}[a](v_1)-\Delta^{-1}M_{ij}[a](v_2) \right|_{C^{2,\alpha}(\overline{\mathbb{B}})}}\\
\stackrel{\hphantom{\eqref{eq:5.15}}}{=}&\sum_{1\leq i\leq j\leq n}{\left|\Delta^{-1}\left[M_{ij}[a](v_1)-M_{ij}[a](v_2)\right] \right|_{C^{2,\alpha}(\overline{\mathbb{B}})}}\\
\stackrel{\eqref{eq:5.20}}{\leq}&\sum_{1\leq i\leq j\leq n}{K_5(n,\alpha,a)\cdot \left|M_{ij}[a](v_1)-M_{ij}[a](v_2)\right|_{C^{0,\alpha}(\overline{\mathbb{B}})}}\\
\stackrel{\eqref{eq:5.33}}{\leq}&K_6(n,\alpha,a)\cdot \left(\left|v_1 \right|_{C^{2,\alpha}(\overline{\mathbb{B}},\mathbb{R}^q)} +\left|v_2 \right|_{C^{2,\alpha}(\overline{\mathbb{B}},\mathbb{R}^q)}\right)\left|v_1-v_2 \right|_{C^{2,\alpha}(\overline{\mathbb{B}},\mathbb{R}^q)}
\end{align*}

Damit ist \eqref{eq:5.37} gezeigt. Nun sei $a\in C_0^{\infty}(\mathbb{\mathbb{B}})$ und $v\in {C^{m,\alpha}(\overline{\mathbb{\mathbb{B}}},\mathbb{R}^q)}$, f\"ur $m\in\mathbb{N}\backslash\{0\}$, dann gilt:
\begin{align}\label{eq:5.39}
\notag&\left|P[a](v) \right|_{C^{m,\alpha}(\overline{\mathbb{B}},\mathbb{R}^n)}=\sum_{i=1}^n{\left|P_i[a](v)\right|_{C^{m,\alpha}(\overline{\mathbb{B}})}}\stackrel{\eqref{eq:5.14}}{=}\sum_{i=1}^n{\left|a\Delta^{-1} N_i[a](v)\right|_{C^{m,\alpha}(\overline{\mathbb{B}})}}\\
\notag\stackrel{\eqref{eq:A.13}}{\leq}&\sum_{i=1}^n{\left|a\right|_{C^{0,\alpha}(\overline{\mathbb{B}})}\left|\Delta^{-1} N_i[a](v)\right|_{C^{m,\alpha}(\overline{\mathbb{B}})}}+\sum_{i=1}^n{\left|a\right|_{C^{m,\alpha}(\overline{\mathbb{B}})}\left|\Delta^{-1} N_i[a](v)\right|_{C^{0,\alpha}(\overline{\mathbb{B}})}}\\
&+C_1(n,m,\alpha) \cdot\sum_{i=1}^n{\left|a\right|_{C^{m-1,\alpha}(\overline{\mathbb{B}})}\left|\Delta^{-1} N_i[a](v)\right|_{C^{m-1,\alpha}(\overline{\mathbb{B}})}}\\
\notag\stackrel{\hphantom{\eqref{eq:A.13}}}{\leq}&K_1(\alpha,a)\cdot\sum_{i=1}^n{\left|\Delta^{-1} N_i[a](v)\right|_{C^{m,\alpha}(\overline{\mathbb{B}})}}+C_2(m,\alpha,a)\cdot\sum_{i=1}^n{\left|\Delta^{-1} N_i[a](v)\right|_{C^{0,\alpha}(\overline{\mathbb{B}})}}\\
\notag&+C_3(n,m,\alpha,a) \cdot\sum_{i=1}^n{\left|\Delta^{-1} N_i[a](v)\right|_{C^{m-1,\alpha}(\overline{\mathbb{B}})}}\\
\notag\stackrel{\eqref{eq:A.6}}{\leq}&K_1(\alpha,a)\cdot\sum_{i=1}^n{\left|\Delta^{-1} N_i[a](v)\right|_{C^{m,\alpha}(\overline{\mathbb{B}})}}+C_4(n,m,\alpha,a) \cdot\sum_{i=1}^n{\left|\Delta^{-1} N_i[a](v)\right|_{C^{m-1,\alpha}(\overline{\mathbb{B}})}}
\end{align}

Ist $m=2$, dann folgt mit \eqref{eq:A.6} aus \eqref{eq:5.39} die Absch\"atzung:
\begin{align*}
&\left|P[a](v) \right|_{C^{2,\alpha}(\overline{\mathbb{B}},\mathbb{R}^n)}\leq K_2(n,\alpha,a)\cdot\sum_{i=1}^n{\left|\Delta^{-1} N_i[a](v)\right|_{C^{2,\alpha}(\overline{\mathbb{B}})}}\\
\stackrel{\eqref{eq:5.20}}{\leq} &K_3(n,\alpha,a)\cdot\sum_{i=1}^n{\left| N_i[a](v)\right|_{C^{0,\alpha}(\overline{\mathbb{B}})}}\stackrel{\eqref{eq:5.32}}{\leq} K_4(n,\alpha,a)\cdot\left(\left|v \right|_{C^{2,\alpha}(\overline{\mathbb{B}},\mathbb{R}^q)}^2+\left|v \right|_{C^{1,\alpha}(\overline{\mathbb{B}},\mathbb{R}^q)}^2\right)
\end{align*}

und falls $m\geq 3$ gilt, dann folgt aus \eqref{eq:5.39} und \eqref{eq:5.30}, mit \eqref{eq:A.6}, die Absch\"atzung:
\begin{align*}
&\left|P[a](v) \right|_{C^{m,\alpha}(\overline{\mathbb{B}},\mathbb{R}^n)}\\
\stackrel{\hphantom{\eqref{eq:5.32}}}{\leq} &K_5(n,\alpha,a)\cdot\sum_{i=1}^n{\left|N_i[a](v)\right|_{C^{m-2,\alpha}(\overline{\mathbb{B}})}}+C_5(n,m,\alpha,a)\cdot \sum_{i=1}^n{\left|N_i[a](v)\right|_{C^{m-3,\alpha}(\overline{\mathbb{B}})}}\\
\stackrel{\eqref{eq:5.32}}{\leq}& K_6(n,\alpha,a)\cdot \left|v \right|_{C^{m,\alpha}(\overline{\mathbb{B}},\mathbb{R}^q)}\left|v \right|_{C^{2,\alpha}(\overline{\mathbb{B}},\mathbb{R}^q)}+C_6(n,m,\alpha,a)\cdot\left|v \right|_{C^{m-1,\alpha}(\overline{\mathbb{B}},\mathbb{R}^q)}^2
\end{align*}

und schlie\ss{}lich f\"ur $m\in \mathbb{N}\backslash\{0,1\}$:
\begin{equation}\label{eq:5.40}
\left|Q[a](v) \right|_{C^{m,\alpha}(\overline{\mathbb{B}},\mathbb{R}^{\frac{n}{2}(n+1)})}=\sum_{1\leq i\leq j\leq n}{\left|Q_{ij}[a](v)\right|_{C^{m,\alpha}(\overline{\mathbb{B}})}}=\sum_{1\leq i\leq j\leq n}{\left|\Delta^{-1}M_{ij}[a](v)\right|_{C^{m,\alpha}(\overline{\mathbb{B}})}}
\end{equation}

Falls $m=2$, dann folgt, unter Beachtung von \eqref{eq:A.6}, aus \eqref{eq:5.40}:
\begin{align*}
\left|Q[a](v) \right|_{C^{2,\alpha}(\overline{\mathbb{B}},\mathbb{R}^{\frac{n}{2}(n+1)})} \stackrel{\eqref{eq:5.20}}{\leq}& \sum_{1\leq i\leq j\leq n}{\left|M_{ij}[a](v)\right|_{C^{0,\alpha}(\overline{\mathbb{B}})}}\\
 \stackrel{\eqref{eq:5.34}}{\leq}&K_7(n,\alpha,a)\cdot\left(\left|v \right|_{C^{2,\alpha}(\overline{\mathbb{B}},\mathbb{R}^q)}^2+\left|v \right|_{C^{1,\alpha}(\overline{\mathbb{B}},\mathbb{R}^q)}^2\right)
\end{align*}

und wenn $m\geq 3$ gilt, dann folgt aus \eqref{eq:5.40} und \eqref{eq:5.30} mit \eqref{eq:A.6} die Absch\"atzung:
\begin{align*}
&\left|Q[a](v) \right|_{C^{m,\alpha}(\overline{\mathbb{B}},\mathbb{R}^{\frac{n}{2}(n+1)})}\\
\stackrel{\hphantom{\eqref{eq:5.34}}}{\leq} &K_8(n,\alpha,a)\cdot\sum_{1\leq i\leq j\leq n}{\left|M_{ij}[a](v)\right|_{C^{m-2,\alpha}(\overline{\mathbb{B}})}}+C_7(n,m,\alpha,a)\cdot\sum_{1\leq i\leq j\leq n}{\left|M_{ij}[a](v)\right|_{C^{m-3,\alpha}(\overline{\mathbb{B}})}}\\
 \stackrel{\eqref{eq:5.34}}{\leq}&K_9(n,\alpha,a)\cdot \left|v \right|_{C^{m,\alpha}(\overline{\mathbb{B}},\mathbb{R}^q)}\left|v \right|_{C^{2,\alpha}(\overline{\mathbb{B}},\mathbb{R}^q)}+C_8(n,m,\alpha,a)\cdot\left|v \right|_{C^{m-1,\alpha}(\overline{\mathbb{B}},\mathbb{R}^q)}^2
\end{align*}

womit \eqref{eq:5.38} gezeigt ist.

\end{bew}

Damit wurden die wesentlichen Absch\"atzungen f\"ur die Operatoren $P$ und $Q$ gezeigt. In den n\"achsten beiden Lemmata wird der in \eqref{eq:5.19} eingef\"uhrte lineare Operator $E[F_0]$ genauer untersucht. Zun\"achst wird die Stetigkeit gezeigt:

\begin{lem}\label{lem:5.a}
F\"ur eine beliebige freie Abbildung $F_0\in C^{\infty}(\overline{\mathbb{B}},\mathbb{R}^q)$ ist der lineare Operator $E[F_0]$ stetig , und es gilt f\"ur $h\in C^{m,\alpha}(\overline{\mathbb{B}},\mathbb{R}^n)$ und $f\in C^{m,\alpha}(\overline{\mathbb{B}},\mathbb{R}^{\frac{n}{2}(n+1)})$ die Absch\"atzung:
\begin{equation}\label{eq:5.46}
\left|E[F_0](h,f) \right|_{C^{m,\alpha}(\overline{\mathbb{B}},\mathbb{R}^q)}\leq C(n,q,m,\alpha,F_0)\cdot\left(\left|h \right|_{C^{m,\alpha}(\overline{\mathbb{B}},\mathbb{R}^n)}+\left|f \right|_{C^{m,\alpha}\left(\overline{\mathbb{B}},\mathbb{R}^{\frac{n}{2}(n+1)}\right)}\right)
\end{equation}

\end{lem}

\begin{bew}

Nach Definition gilt f\"ur die $l$-te Komponente von $E[F_0](h,f)$ die Gleichung:
\begin{equation}\label{eq:5.42}
E_l[F_0](h,f)=\sum_{i=1}^n{A_{l,i}h_i}+ \sum_{1\leq i\leq j\leq n}{B_{l,ij} f_{ij}}
\end{equation}

f\"ur Funktionen $\{A_{l,i}\}_{1\leq i\leq n}\subseteq C^{\infty}(\overline{\mathbb{B}})$ und $\{B_{l,ij}\}_{1\leq i\leq j\leq n}\subseteq C^{\infty}(\overline{\mathbb{B}})$, siehe hierf\"ur auch \eqref{eq:5.41}. Diese Funktionen h\"angen nur von Ableitungen der freien Abbildung $F_0$ ab. Sei nun $s\in\mathbb{N}^n$ ein Multiindex der Ordnung $k$ f\"ur ein $k\in\{0,...,m\}$, dann gilt:
\begin{align*}
&\left|\partial^s (E_l[F_0](h,f))\right|_{C^{0,\alpha}(\overline{\mathbb{B}})}=\left|\sum_{i=1}^n{\partial^s(A_{l,i}h_i)}+ \sum_{1\leq i\leq j\leq n}{\partial^s(B_{l,ij} f_{ij})}\right|_{C^{0,\alpha}(\overline{\mathbb{B}})}\\
\stackrel{\eqref{eq:A.7}}{=}&\left|\sum_{i=1}^n{\sum_{r\leq s}{\binom{s}{r}\partial^rA_{l,i}\partial^{s-r}h_i}}+ \sum_{1\leq i\leq j\leq n}{\sum_{r\leq s}{\binom{s}{r} \partial^r B_{l,ij} \partial^{s-r}f_{ij}}}\right|_{C^{0,\alpha}(\overline{\mathbb{B}})}\\
\stackrel{\hphantom{\eqref{eq:A.7}}}{\leq}&\sum_{i=1}^n{\sum_{r\leq s}{\binom{s}{r}\left|\partial^rA_{l,i}\partial^{s-r}h_i\right|_{C^{0,\alpha}(\overline{\mathbb{B}})}}}+ \sum_{1\leq i\leq j\leq n}{\sum_{r\leq s}{\binom{s}{r}\left|\partial^r B_{l,ij} \partial^{s-r}f_{ij}\right|_{C^{0,\alpha}(\overline{\mathbb{B}})}}}\\
\stackrel{\eqref{eq:A.2}}{\leq}&\sum_{i=1}^n{\sum_{r\leq s}{\binom{s}{r}\left|\partial^rA_{l,i}\right|_{C^{0,\alpha}(\overline{\mathbb{B}})}\left|\partial^{s-r}h_i\right|_{C^{0,\alpha}(\overline{\mathbb{B}})}}}\\
&\hspace{4cm}+ \sum_{1\leq i\leq j\leq n}{\sum_{r\leq s}{\binom{s}{r}\left|\partial^r B_{l,ij}\right|_{C^{0,\alpha}(\overline{\mathbb{B}})} \left|\partial^{s-r}f_{ij}\right|_{C^{0,\alpha}(\overline{\mathbb{B}})}}}\\
\stackrel{\hphantom{\eqref{eq:A.2}}}{\leq}&\sum_{i=1}^n{\sum_{r\leq s}{\binom{s}{r}\left|A_{l,i}\right|_{C^{|r|,\alpha}(\overline{\mathbb{B}})}\left|h_i\right|_{C^{|s-r|,\alpha}(\overline{\mathbb{B}})}}}\\
&\hspace{4cm}+ \sum_{1\leq i\leq j\leq n}{\sum_{r\leq s}{\binom{s}{r}\left|B_{l,ij}\right|_{C^{|r|,\alpha}(\overline{\mathbb{B}})} \left|f_{ij}\right|_{C^{|s-r|,\alpha}(\overline{\mathbb{B}})}}}\\
\stackrel{\eqref{eq:A.6}}{\leq}&C_1(n,k,\alpha,F_0)\cdot \left(\left|h \right|_{C^{k,\alpha}(\overline{\mathbb{B}},\mathbb{R}^n)}+\left|f \right|_{C^{k,\alpha}(\overline{\mathbb{B}},\mathbb{R}^{\frac{n}{2}(n+1)})}\right)
\end{align*}

Daraus ergibt sich die Absch\"atzung:
\begin{align*}
&\left|E_l[F_0](h,f)\right|_{C^{m,\alpha}(\overline{\mathbb{B}})}=\left|E_l[F_0](h,f)\right|_{C^{0,\alpha}(\overline{\mathbb{B}})}+\sum_{|s|=m}{\left|\partial^s (E_l[F_0](h,f)) \right|_{C^{0,\alpha}(\overline{\mathbb{B}})}}\\
\stackrel{\hphantom{\eqref{eq:A.12}}}{\leq}& C_2(n,\alpha,F_0)\cdot\left(\left|h \right|_{C^{0,\alpha}(\overline{\mathbb{B}},\mathbb{R}^n)}+\left|f \right|_{C^{0,\alpha}(\overline{\mathbb{B}},\mathbb{R}^{\frac{n}{2}(n+1)})}\right)\\
&+C_3(n,m,\alpha,F_0)\cdot\left(\left|h \right|_{C^{m,\alpha}(\overline{\mathbb{B}},\mathbb{R}^n)}+\left|f \right|_{C^{m,\alpha}(\overline{\mathbb{B}},\mathbb{R}^{\frac{n}{2}(n+1)})}\right)\\
\stackrel{\eqref{eq:A.12}}{\leq}&C_4(n,m,\alpha,F_0)\cdot\left(\left|h \right|_{C^{m,\alpha}(\overline{\mathbb{B}},\mathbb{R}^n)}+\left|f \right|_{C^{m,\alpha}\left(\overline{\mathbb{B}},\mathbb{R}^{\frac{n}{2}(n+1)}\right)}\right)
\end{align*}

und schlie\ss{}lich:
\begin{align*}
\left|E[F_0](h,f) \right|_{C^{m,\alpha}(\overline{\mathbb{B}},\mathbb{R}^q)}=&\sum_{l=1}^q{\left|E_l[F_0](h,f)\right|_{C^{m,\alpha}(\overline{\mathbb{B}})}}\\
\leq& C(n,q,m,\alpha,F_0)\cdot\left(\left|h \right|_{C^{m,\alpha}(\overline{\mathbb{B}},\mathbb{R}^n)}+\left|f \right|_{C^{m,\alpha}\left(\overline{\mathbb{B}},\mathbb{R}^{\frac{n}{2}(n+1)}\right)}\right)
\end{align*}

\end{bew}

F\"ur die nachfolgenden Betrachtungen wird, f\"ur $m\in\mathbb{N}$, die Operatornorm des stetigen linearen Operators $E[F_0]$ aus \eqref{eq:5.19} mit $\left\Vert E[F_0] \right\Vert_{m,\alpha}$ bezeichnet, siehe hierf\"ur auch \cite[3.2]{alt2012lineare}. Die Absch\"atzung \eqref{eq:5.46} kann noch versch\"arft werden:

\begin{lem}\label{lem:5.e}

Es sei $m\geq 3$ dann gilt f\"ur den Operator $E[F_0]$ die Absch\"atzung:
\begin{align}\label{eq:5.49}
\begin{split}
&\left|E[F_0](h,f) \right|_{C^{m,\alpha}(\overline{\mathbb{B}},\mathbb{R}^q)}\leq \left\Vert E[F_0] \right\Vert_{2,\alpha}\cdot \left(\left|h \right|_{C^{m,\alpha}(\overline{\mathbb{B}},\mathbb{R}^n)}+\left|f \right|_{C^{m,\alpha}\left(\overline{\mathbb{B}},\mathbb{R}^{\frac{n}{2}(n+1)}\right)}\right)\\
&\hspace{1cm}+C(n,q,m,\alpha,F_0)\cdot\left(\left|h \right|_{C^{m-1,\alpha}(\overline{\mathbb{B}},\mathbb{R}^n)}+\left|f \right|_{C^{m-1,\alpha}\left(\overline{\mathbb{B}},\mathbb{R}^{\frac{n}{2}(n+1)}\right)}\right)
\end{split}
\end{align}

f\"ur alle $h\in C^{m,\alpha}(\overline{\mathbb{B}},\mathbb{R}^n)$ und $f\in C^{m,\alpha}(\overline{\mathbb{B}},\mathbb{R}^{\frac{n}{2}(n+1)})$.

\end{lem}

\begin{bew}

Sei $l\in\{1,...,q\}$ und $s\in\mathbb{N}^n$ ein Multiindex der Ordnung $m-2$, dann gilt mit \eqref{eq:5.42}:
\begin{align*}
&\partial^s (E_l[F_0](h,f))=\sum_{l=1}^n{\partial^s (A_{l,i}h_i)}+\sum_{1\leq i\leq j\leq n}{\partial^s (B_{l,ij} f_{ij})}\\
=&\sum_{i=1}{\sum_{r\leq s}{\binom{s}{r}\partial^{r} A_{l,i} \partial^{s-r}h_i}}+\sum_{1\leq i\leq j\leq n}{\sum_{r\leq s}{\binom{s}{r}\partial^{r} B_{l,ij}\partial^{s-r} f_{ij}}}\\
=&\sum_{i=1}^n{A_{l,i}\partial^s h_i}+\sum_{i=1}^n{\sum_{0<r\leq s}{\binom{s}{r}\partial^{r} A_{l,i}\partial^{s-r}h_i}}\\
&+\sum_{1\leq i\leq j\leq n}{B_{l,ij}\partial^s f_{ij}}+\sum_{1\leq i\leq j\leq n}{\sum_{0<r\leq s}{\binom{s}{r}\partial^{r} B_{l,ij}\partial^{s-r} f_{ij}}}\\
=& E_l[F_0](\partial^s h,\partial^s f)+\sum_{0<r\leq s}{\binom{s}{r}\left[\partial^{r} A_{l,i} \partial^{s-r}h_i+\partial^{r} B_{l,ij}\partial^{s-r} f_{ij}\right]}
\end{align*}

Daraus folgt mit \eqref{eq:A.6}, wegen $|s|=m-2$:
\begin{align*}
&\left|\partial^s E_l[F_0](h,f)-E_l[F_0](\partial^s h,\partial^s f) \right|_{C^{2,\alpha}(\overline{\mathbb{B}})}\\
\leq &C_1(n,\alpha)\cdot\sum_{0<r\leq s}{\binom{s}{r}\left|\partial^{r} A_{l,i}\right|_{C^{2,\alpha}(\overline{\mathbb{B}})}\cdot \left|\partial^{s-r}h_i\right|_{C^{2,\alpha}(\overline{\mathbb{B}})}}\\
&+C_1(n,\alpha)\cdot\sum_{0<r\leq s}{\binom{s}{r}\left|\partial^{r} B_{l,ij}\right|_{C^{2,\alpha}(\overline{\mathbb{B}})}\cdot\left|\partial^{s-r} f_{ij}\right|_{C^{2,\alpha}(\overline{\mathbb{B}})}}\\
\leq &C_2(n,m,\alpha,F_0)\cdot \left(\left|h \right|_{C^{m-1,\alpha}(\overline{\mathbb{B}},\mathbb{R}^n)}+\left|f \right|_{C^{m-1,\alpha}\left(\overline{\mathbb{B}},\mathbb{R}^{\frac{n}{2}(n+1)}\right)}\right)
\end{align*}

Dies impliziert:
\begin{align}
\begin{split}
&\left|\partial^s E[F_0](h,f)-E[F_0](\partial^s h,\partial^s f) \right|_{C^{2,\alpha}(\overline{\mathbb{B}},\mathbb{R}^q)}\\
=&\sum_{l=1}^q{\left|\partial^s E_l[F_0](h,f)-E_l[F_0](\partial^s h,\partial^s f) \right|_{C^{2,\alpha}(\overline{\mathbb{B}})}}\\
\leq &C_3(n,q,m,\alpha,F_0)\cdot \left(\left|h \right|_{C^{m-1,\alpha}(\overline{\mathbb{B}},\mathbb{R}^n)}+\left|f \right|_{C^{m-1,\alpha}\left(\overline{\mathbb{B}},\mathbb{R}^{\frac{n}{2}(n+1)}\right)}\right)
\end{split}
\end{align}

Es folgt
\begin{align*}
&\left|E[F_0](h,f) \right|_{C^{m,\alpha}(\overline{\mathbb{B}},\mathbb{R}^q)}\leq \left|E[F_0](h,f) \right|_{C^{0,\alpha}(\overline{\mathbb{B}},\mathbb{R}^q)}+\sum_{|s|=m-2}{\left|\partial^s E[F_0](h,f) \right|_{C^{2,\alpha}(\overline{\mathbb{B}},\mathbb{R}^q)}}\\
\leq & \left\Vert E[F_0] \right\Vert_{0,\alpha}\cdot\left(\left|h \right|_{C^{0,\alpha}(\overline{\mathbb{B}},\mathbb{R}^n)}+\left|f \right|_{C^{0,\alpha}\left(\overline{\mathbb{B}},\mathbb{R}^{\frac{n}{2}(n+1)}\right)}\right)\\
&+\sum_{|s|=m-2}{\left|\partial^s E[F_0](h,f)-E[F_0](\partial^s h,\partial^s f) \right|_{C^{2,\alpha}(\overline{\mathbb{B}},\mathbb{R}^q)}}+\sum_{|s|=m-2}{\left|E[F_0](\partial^s h,\partial^s f) \right|_{C^{2,\alpha}(\overline{\mathbb{B}},\mathbb{R}^q)}}\\
\leq &\left\Vert E[F_0] \right\Vert_{0,\alpha}\cdot \left(\left|h \right|_{C^{0,\alpha}(\overline{\mathbb{B}},\mathbb{R}^n)}+\left|f \right|_{C^{0,\alpha}\left(\overline{\mathbb{B}},\mathbb{R}^{\frac{n}{2}(n+1)}\right)}\right)\\
&+C(n,q,m,\alpha,F_0)\cdot\left(\left|h \right|_{C^{m-1,\alpha}(\overline{\mathbb{B}},\mathbb{R}^n)}+\left|f \right|_{C^{m-1,\alpha}\left(\overline{\mathbb{B}},\mathbb{R}^{\frac{n}{2}(n+1)}\right)}\right)\\
&+\sum_{|s|=m-2}{\left\Vert E[F_0] \right\Vert_{2,\alpha}\cdot\left(\left|\partial^s h \right|_{C^{2,\alpha}(\overline{\mathbb{B}},\mathbb{R}^n)}+\left|\partial^s f \right|_{C^{2,\alpha}\left(\overline{\mathbb{B}},\mathbb{R}^{\frac{n}{2}(n+1)}\right)}\right)}
\end{align*}
woraus mit \eqref{eq:A.12} die gew\"unschte Absch\"atzung folgt.

\end{bew}


\section{L\"osung des Fixpunktproblems}
\label{sec:5.5}

Die in \autoref{sec:5.4} bewiesenen Absch\"atzungen werden nun verwendet, um zu zeigen, dass die Fixpunktgleichung \eqref{eq:5.16} eine L\"osung besitzt, sofern f\"ur die Abbildung $f\in C^{m,\alpha}(\overline{\mathbb{B}},\mathbb{R}^{\frac{n}{2}(n+1)})$, $m\geq 2 $, der Ausdruck $\left\Vert E[F_0]\right\Vert_{2,\alpha} \left|E[F_0](0,f) \right|_{C^{2,\alpha}(\overline{\mathbb{B}},\mathbb{R}^q)}$ klein genug ist. Dar\"uber hinaus wird gezeigt, dass sich die Regularit\"at von $f$ auf den konstruierten Fixpunkt \"ubertr\"agt.

\begin{lem}\label{lem:5.13}
Es sei $a\in C^{\infty}_0(\mathbb{B})$, dann existiert ein $\vartheta(n,\alpha,a)\in\mathbb{R}_{>0}$, mit der folgenden Eigenschaft: Ist eine freie Abbildung $F_0\in C^{\infty}(\overline{\mathbb{B}},\mathbb{R}^q)$, und ein $f\in C^{m,\alpha}(\overline{\mathbb{B}},\mathbb{R}^{\frac{n}{2}(n+1)})$, f\"ur $m\geq 2$, gegeben, so dass:
\begin{equation*}
\left\Vert E[F_0]\right\Vert_{2,\alpha}\left|E[F_0](0,f) \right|_{C^{2,\alpha}(\overline{\mathbb{B}},\mathbb{R}^q)}\leq \vartheta
\end{equation*}

erf\"ullt ist, dann besitzt der in \eqref{eq:5.18} definierte Operator $\Phi[F_0,a,f]$ einen Fixpunkt $v\in C^{m,\alpha}(\overline{\mathbb{B}},\mathbb{R}^q)$, welcher die Absch\"atzung:
\begin{equation}\label{eq:5.51}
\left|v \right|_{C^{2,\alpha}(\overline{\mathbb{B}},\mathbb{R}^q)}\leq \left|E[F_0](0,f) \right|_{C^{2,\alpha}(\overline{\mathbb{B}},\mathbb{R}^q)}
\end{equation}

erf\"ullt. Ferner gilt: $v\in C^{\infty}(\overline{\mathbb{B}},\mathbb{R}^q)$, falls $f\in C^{\infty}(\overline{\mathbb{B}},\mathbb{R}^{\frac{n}{2}(n+1)})$.

\end{lem}

\begin{bew}
Ist $f\equiv 0$, dann ist $v\equiv 0$ ein Fixpunkt des Operators $\Phi[F_0,a,f]$. Es sei also $f\not\equiv 0$. Definiere eine Folge $(v_k)_{k\in\mathbb{N}}\subseteq C^{m,\alpha}(\overline{\mathbb{B}},\mathbb{R}^q)$ wie folgt:
\begin{align}\label{eq:5.47}
v_k=
\begin{cases}
0 & \text{falls }  k=0\\
\Phi[F_0,a,f](v_{k-1})  &\text{falls } k\geq 1
\end{cases}
\end{align}

Zun\"achst wird gezeigt, dass diese Folge eine Cauchyfolge in der $\left|\cdot\right|_{C^{2,\alpha}(\overline{\mathbb{B}},\mathbb{R}^q)}$-Norm ist, falls der Ausdruck $\left\Vert E[F_0]\right\Vert_{2,\alpha} \left|E[F_0](0,f) \right|_{C^{2,\alpha}(\overline{\mathbb{B}},\mathbb{R}^q)}$ klein genug ist. Dann konvergiert diese Folge, aufgrund der Vollst\"andigkeit des Raumes $C^{2,\alpha}(\overline{\mathbb{B}},\mathbb{R}^q)$, gegen ein Grenzelement $v\in C^{2,\alpha}(\overline{\mathbb{B}},\mathbb{R}^q)$. Dazu wird zun\"achst gezeigt, dass die Folge beschr\"ankt ist, sofern der Ausdruck $\left\Vert E[F_0]\right\Vert_{2,\alpha} \left|E[F_0](0,f) \right|_{C^{2,\alpha}(\overline{\mathbb{B}},\mathbb{R}^q)}$ hinreichend klein ist. F\"ur $\mathbb{N}\backslash\{0\}$ gilt, unter Beachtung der Linearit\"at des Operators $E[F_0]$:
\begin{align*}
&\left|v_k \right|_{C^{2,\alpha}(\overline{\mathbb{B}},\mathbb{R}^q)}\stackrel{\eqref{eq:5.47}}{=}\left|\Phi[F_0,a,f](v_{k-1})\right|_{C^{2,\alpha}(\overline{\mathbb{B}},\mathbb{R}^q)}\\
\stackrel{\eqref{eq:5.18}}{=}&\left|E[F_0]\left(P[a](v_{k-1}),\frac{1}{2}f-\frac{1}{2}Q[a](v_{k-1}) \right)\right|_{C^{2,\alpha}(\overline{\mathbb{B}},\mathbb{R}^q)}\\
\stackrel{\hphantom{\eqref{eq:5.18}}}{\leq}&\left|E[F_0]\left(P[a](v_{k-1}),-\frac{1}{2}Q[a](v_{k-1}) \right)\right|_{C^{2,\alpha}(\overline{\mathbb{B}},\mathbb{R}^q)}+\left|E[F_0]\left(0,\frac{1}{2}f \right)\right|_{C^{2,\alpha}(\overline{\mathbb{B}},\mathbb{R}^q)}\\
\stackrel{\eqref{eq:5.46}}{\leq}& \left\Vert E[F_0] \right\Vert_{2,\alpha}\cdot\left[\left|P[a](v_{k-1}) \right|_{C^{2,\alpha}(\overline{\mathbb{B}},\mathbb{R}^n)} +\left|Q[a](v_{k-1}) \right|_{C^{2,\alpha}(\overline{\mathbb{B}},\mathbb{R}^{\frac{n}{2}(n+1)})}\right]\\
&+\frac{1}{2}\left|E[F_0]\left(0,f \right)\right|_{C^{2,\alpha}(\overline{\mathbb{B}},\mathbb{R}^q)}\\
\stackrel{\eqref{eq:5.38}}{\leq}&K_1(n,\alpha,a)\cdot \left\Vert E[F_0] \right\Vert_{2,\alpha}\cdot\left|v_{k-1} \right|_{C^{2,\alpha}(\overline{\mathbb{B}},\mathbb{R}^q)}^2+\frac{1}{2}\left|E[F_0]\left(0,f \right)\right|_{C^{2,\alpha}(\overline{\mathbb{B}},\mathbb{R}^q)}
\end{align*}

Falls nun:
\begin{equation*}
\left\Vert E[F_0]\right\Vert_{2,\alpha} \left|E[F_0](0,f) \right|_{C^{2,\alpha}(\overline{\mathbb{B}},\mathbb{R}^q)}\leq \frac{1}{2K_1(n,\alpha,a)}
\end{equation*}

gilt, so ist:
\begin{align*}
\left|v_k \right|_{C^{2,\alpha}(\overline{\mathbb{B}},\mathbb{R}^q)}\leq \frac{1}{2}\left(\frac{\left|v_{k-1} \right|_{C^{2,\alpha}(\overline{\mathbb{B}},\mathbb{R}^q)}^2}{ \left|E[F_0](0,f) \right|_{C^{2,\alpha}(\overline{\mathbb{B}},\mathbb{R}^q)}}+ \left|E[F_0]\left(0,f \right)\right|_{C^{2,\alpha}(\overline{\mathbb{B}},\mathbb{R}^q)}\right)
\end{align*}

woraus induktiv, wegen $v_0\equiv 0$, f\"ur alle $k\in\mathbb{N}$ die Absch\"atzung:
\begin{equation}\label{eq:5.48}
\left|v_k \right|_{C^{2,\alpha}(\overline{\mathbb{B}},\mathbb{R}^q)}\leq \left|E[F_0]\left(0,f \right)\right|_{C^{2,\alpha}(\overline{\mathbb{B}},\mathbb{R}^q)}
\end{equation}

folgt. Daraus resultiert f\"ur $k\in\mathbb{N}\backslash\{0\}$:
\begin{align*}
&\left|v_{k+1}-v_{k} \right|_{C^{2,\alpha}(\overline{\mathbb{B}},\mathbb{R}^q)}\stackrel{\eqref{eq:5.47}}{=}\left|\Phi[F_0,a,f](v_{k})-\Phi[F_0,a,f](v_{k-1}) \right|_{C^{2,\alpha}(\overline{\mathbb{B}},\mathbb{R}^q)}\\
\stackrel{\eqref{eq:5.18}}{=}&\left|E[F_0]\left(P[a](v_k)-P[a](v_{k-1}),-\frac{1}{2}(Q[a](v_k)-Q[a](v_{k-1}))\right) \right|_{C^{2,\alpha}(\overline{\mathbb{B}},\mathbb{R}^q)}\\
\stackrel{\eqref{eq:5.46}}{\leq}& \left\Vert E[F_0]\right\Vert_{2,\alpha}\cdot\Bigl(\left|P[a](v_k)-P[a](v_{k-1})\right|_{C^{2,\alpha}(\overline{\mathbb{B}},\mathbb{R}^n)}\\
&+\left| Q[a](v_k)-Q[a](v_{k-1}) \right|_{{C^{2,\alpha}(\overline{\mathbb{B}},\mathbb{R}^{\frac{n}{2}(n+1)})}}\Bigr)\\
\stackrel{\eqref{eq:5.37}}{\leq}&\widehat{K}_2(n,\alpha,a)\cdot\left\Vert E[F_0]\right\Vert_{2,\alpha}\cdot\left(\left|v_k \right|_{C^{2,\alpha}(\overline{\mathbb{B}},\mathbb{R}^q)} +\left|v_{k-1} \right|_{C^{2,\alpha}(\overline{\mathbb{B}},\mathbb{R}^q)}\right) \left|v_k-v_{k-1} \right|_{C^{2,\alpha}(\overline{\mathbb{B}},\mathbb{R}^q)}\\
\stackrel{\eqref{eq:5.48}}{\leq}& K_2(n,\alpha,a)\cdot \left\Vert E[F_0]\right\Vert_{2,\alpha}\left|E[F_0](0,f)\right|_{C^{2,\alpha}(\overline{\mathbb{B}},\mathbb{R}^q)} \left|v_k-v_{k-1}\right|_{C^{2,\alpha}(\overline{\mathbb{B}},\mathbb{R}^q)}\\
\end{align*}

Gilt dann die Absch\"atzung:
\begin{equation*}
\left\Vert E[F_0]\right\Vert_{2,\alpha} \left|E[F_0](0,f)\right|_{C^{2,\alpha}(\overline{\mathbb{B}},\mathbb{R}^q)} \leq \frac{1}{2K_2(n,\alpha,a)}
\end{equation*}

dann ist:
\begin{equation*}
\left|v_{k+1}-v_{k} \right|_{C^{2,\alpha}(\overline{\mathbb{B}},\mathbb{R}^q)}\leq \frac{1}{2} \left|v_k-v_{k-1}\right|_{C^{2,\alpha}(\overline{\mathbb{B}},\mathbb{R}^q)}
\end{equation*}

und die Folge $(v_k)_{k\in\mathbb{N}}\subseteq C^{m,\alpha}(\overline{\mathbb{B}},\mathbb{R}^q)$ ist eine Cauchy-Folge, die im $C^{2,\alpha}$-Sinne gegen ein $v\in C^{2,\alpha}(\overline{\mathbb{B}},\mathbb{R}^q)$ konvergiert, sofern die Absch\"atzung $\left\Vert E[F_0]\right\Vert_{2,\alpha} \left|E[F_0](0,f)\right|_{C^{2,\alpha}(\overline{\mathbb{B}},\mathbb{R}^q)} \leq \frac{1}{2}\min\left\{ K_1^{-1},K_2^{-1}\right\}$ erf\"ullt ist. Um zu zeigen, dass $v\in C^{m,\alpha}(\overline{\mathbb{B}},\mathbb{R}^q)$ gilt, wird gezeigt, dass die Folge $(v_k)_{k\in\mathbb{N}}$ in der $\left|\cdot \right|_{C^{m,\alpha}(\overline{\mathbb{B}},\mathbb{R}^q)}$-Norm beschr\"ankt ist. Dann existiert mit dem Satz von Arzel\`{a}-Ascoli \cite[2.12]{alt2012lineare} eine Teilfolge von $(v_k)_{k\in\mathbb{N}}$, die im $C^{m,\alpha}$-Sinne gegen $v$ konvergiert, womit dann $v\in C^{m,\alpha}(\overline{\mathbb{B}},\mathbb{R}^q)$ gezeigt ist. Mit \eqref{eq:5.48} ist die Behauptung f\"ur $m=2$ richtig. Es sei also $m\geq 3$, unter der Annahme, dass:
\begin{equation}\label{eq:5.50}
\left|v_{k} \right|_{C^{m-1,\alpha}(\overline{\mathbb{B}},\mathbb{R}^q)}\leq \eta
\end{equation}

f\"ur alle $k\in\mathbb{N}$ gilt, wobei $\eta$ nicht von $k$ abh\"angt. Sei nun $k\in\mathbb{N}\backslash\{0\}$ dann ist:
\begin{align*}
&\left|v_k \right|_{C^{m,\alpha}(\overline{\mathbb{B}},\mathbb{R}^q)}\stackrel{\eqref{eq:5.47}}{=}\left|\Phi[F_0,a,f](v_{k-1}) \right|_{C^{m,\alpha}(\overline{\mathbb{B}},\mathbb{R}^q)}\\
\stackrel{\eqref{eq:5.18}}{=}&\left|E[F_0]\left(P[a](v_{k-1}),\frac{1}{2}f-\frac{1}{2}Q[a](v_{k-1}) \right)\right|_{C^{m,\alpha}(\overline{\mathbb{B}},\mathbb{R}^q)}\\
\stackrel{\hphantom{\eqref{eq:5.18}}}{\leq}&\left|E[F_0]\left(P[a](v_{k-1}),-\frac{1}{2}Q[a](v_{k-1}) \right)\right|_{C^{m,\alpha}(\overline{\mathbb{B}},\mathbb{R}^q)}+\left|E[F_0]\left(0,\frac{1}{2}f \right)\right|_{C^{m,\alpha}(\overline{\mathbb{B}},\mathbb{R}^q)}\\
\stackrel{\eqref{eq:5.49}}{\leq}&\left\Vert E[F_0] \right\Vert_{2,\alpha}\cdot\left(\left|P[a](v_{k-1}) \right|_{C^{m,\alpha}(\overline{\mathbb{B}},\mathbb{R}^n)}+\left|Q[a](v_{k-1}) \right|_{C^{m,\alpha}(\overline{\mathbb{B}},\mathbb{R}^{\frac{n}{2}(n+1)})} \right)\\
&+C_1(n,q,m,\alpha,F_0)\cdot \left(\left|P[a](v_{k-1}) \right|_{C^{m-1,\alpha}(\overline{\mathbb{B}},\mathbb{R}^n)}+\left|Q[a](v_{k-1}) \right|_{C^{m-1,\alpha}(\overline{\mathbb{B}},\mathbb{R}^{\frac{n}{2}(n+1)})} \right)\\
&+\frac{1}{2}\left|E[F_0]\left(0,f \right)\right|_{C^{m,\alpha}(\overline{\mathbb{B}},\mathbb{R}^q)}\\
\stackrel{\eqref{eq:5.38}}{\leq}&K_3(n,\alpha,a) \cdot\left\Vert E[F_0] \right\Vert_{2,\alpha}\cdot\left|v_{k-1} \right|_{C^{m,\alpha}(\overline{\mathbb{B}},\mathbb{R}^q)}\left|v_{k-1} \right|_{C^{2,\alpha}(\overline{\mathbb{B}},\mathbb{R}^q)}+\frac{1}{2}\left|E[F_0]\left(0,f \right)\right|_{C^{m,\alpha}(\overline{\mathbb{B}},\mathbb{R}^q)}\\
&+C_2(n,q,m,\alpha,F_0,a)\cdot \left|v_{k-1} \right|_{C^{m-1,\alpha}(\overline{\mathbb{B}},\mathbb{R}^q)}^2\\
\stackrel{\eqref{eq:5.48}}{\leq}&K_3(n,\alpha,a) \cdot\left\Vert E[F_0] \right\Vert_{2,\alpha}\left|E[F_0](0,f) \right|_{C^{2,\alpha}(\overline{\mathbb{B}},\mathbb{R}^q)} \cdot\left|v_{k-1} \right|_{C^{m,\alpha}(\overline{\mathbb{B}},\mathbb{R}^q)}\\
&+\frac{1}{2}\left|E[F_0]\left(0,f \right)\right|_{C^{m,\alpha}(\overline{\mathbb{B}},\mathbb{R}^q)}+C_2(n,q,m,\alpha,F_0,a)\cdot \left|v_{k-1} \right|_{C^{m-1,\alpha}(\overline{\mathbb{B}},\mathbb{R}^q)}^2
\end{align*}

Falls nun:
\begin{equation*}
\left\Vert E[F_0] \right\Vert_{2,\alpha}\cdot\left|E[F_0](0,f) \right|_{C^{2,\alpha}(\overline{\mathbb{B}})} \leq \frac{1}{2K_3(n,\alpha,a)}
\end{equation*}

erf\"ullt ist, dann gilt:
\begin{align*}
\left|v_k \right|_{C^{m,\alpha}(\overline{\mathbb{B}},\mathbb{R}^q)}\stackrel{\hphantom{\eqref{eq:5.50}}}{\leq}& \frac{1}{2}\left(\left|v_{k-1} \right|_{C^{m,\alpha}(\overline{\mathbb{B}},\mathbb{R}^q)}+\left|E[F_0]\left(0,f \right)\right|_{C^{m,\alpha}(\overline{\mathbb{B}},\mathbb{R}^q)}\right)\\
&+C_2(n,q,m,\alpha,F_0,a)\cdot \left|v_{k-1} \right|_{C^{m-1,\alpha}(\overline{\mathbb{B}},\mathbb{R}^q)}^2\\
\stackrel{\eqref{eq:5.50}}{\leq}& \frac{1}{2}\left(\left|v_{k-1} \right|_{C^{m,\alpha}(\overline{\mathbb{B}},\mathbb{R}^q)}+\left|E[F_0]\left(0,f \right)\right|_{C^{m,\alpha}(\overline{\mathbb{B}},\mathbb{R}^q)}+C(n,q,m,\alpha,F_0,a,\eta)\right)
\end{align*}

und es folgt, wegen $v_0\equiv 0$, mit Lemma \ref{lem:7.8} induktiv f\"ur alle $k\in\mathbb{N}$ die Absch\"atzung:
\begin{equation}\label{eq:a}
\left|v_k \right|_{C^{m,\alpha}(\overline{\mathbb{B}},\mathbb{R}^q)}\leq \left|E[F_0]\left(0,f \right)\right|_{C^{m,\alpha}(\overline{\mathbb{B}},\mathbb{R}^q)}+C(n,q,m,\alpha,F_0,a,\eta)
\end{equation}

Die zu beweisenden Aussagen gelten also f\"ur die Wahl:
\begin{equation*}
\vartheta(n,\alpha,a):=\frac{1}{2}\min\left\{ K_1^{-1},K_2^{-1}, K_3^{-1}\right\}
\end{equation*}

\end{bew}

Es seien nun offene Mengen $U_1, U_2\subseteq \mathbb{B}$ mit $\overline{U}_1\subseteq U_2$ und $\overline{U}_2 \subseteq \mathbb{B}$, und es sei $f\in C^{m,\alpha}(\overline{\mathbb{B}},\mathbb{R}^{\frac{n}{2}(n+1)})$, f\"ur $m\geq 3$, mit der Eigenschaft $supp(f)\subseteq U_1$, wie am Anfang des Kapitels beschrieben, gegeben. Dann existiert mit \cite[Proposition. 2.26]{lee2003introduction} ein $a\in C^{\infty}_0(\mathbb{B})$, so dass $\left. a\right|_{\overline{U}_1}\equiv 1$ und $supp(a)\subseteq U_2$ gilt. Daraus folgt $f\equiv a^2 f$.  Es bezeichne $v\in C^{m,\alpha}(\overline{\mathbb{B}},\mathbb{R}^q)$ den, in Lemma \ref{lem:5.13} konstruierten, Fixpunkt des Operators $\Phi[F_0,a,f]$. Dann ist, unter Beachtung von Lemma \ref{lem:5.3}, f\"ur die Funktion $u:= a^2\, v\in C^{m,\alpha}(\overline{\mathbb{B}},\mathbb{R}^q)$ die Gleichung $\partial_i(F_0+u)\cdot \partial_j (F_0+u) =\partial_i F_0\cdot \partial_j F_0+f_{ij}$ f\"ur alle $i,j\in\{1,...,n\}$ mit $i\leq j$ erf\"ullt. Mit \eqref{eq:5.51} und \eqref{eq:A.14} gilt die Absch\"atzung $\left|u \right|_{C^{2,\alpha}(\overline{\mathbb{B}},\mathbb{R}^q)}\leq C(n,\alpha,a)\cdot \left|E[F_0]\left(0,f \right)\right|_{C^{2,\alpha}(\overline{\mathbb{B}},\mathbb{R}^q)}$. Daraus folgt der Existenzsatz:

\begin{satz}\label{satz:5.1}
Gegeben seien offene Mengen $U_1, U_2\subseteq \mathbb{B}$ mit $\overline{U}_1\subseteq U_2$ und $\overline{U}_2 \subseteq \mathbb{B}$. Dann existiert ein $\vartheta\in\mathbb{R}_{>0}$ mit der folgenden Eigenschaft: Ist eine freie Abbildung $F_0\in C^{\infty}(\overline{\mathbb{B}},\mathbb{R}^q)$ und ein $f\in C^{m,\alpha}(\overline{\mathbb{B}},\mathbb{R}^{\frac{n}{2}(n+1)})$, f\"ur $m\geq 3$, mit $supp(f)\subseteq U_1$ gegeben, so dass:
\begin{equation}\label{eq:5.52}
\left\Vert E[F_0]\right\Vert_{2,\alpha}\left|E[F_0](0,f) \right|_{C^{2,\alpha}(\overline{\mathbb{B}},\mathbb{R}^q)}\leq \vartheta
\end{equation}

erf\"ullt ist, dann existiert ein $u\in C^{m,\alpha}(\mathbb{B},\mathbb{R}^q)$ mit $supp(u)\subseteq U_2$, so dass f\"ur alle $i,j\in\{1,...,n\}$ mit $i\leq j$ die Gleichung:
\begin{equation*}
\partial_i(F_0+u)\cdot \partial_j (F_0+u) =\partial_i F_0\cdot \partial_j F_0+f_{ij}
\end{equation*}

erf\"ullt ist. Ferner gilt die Absch\"atzung:
\begin{equation}\label{eq:5.53}
\left|u \right|_{C^{2,\alpha}(\overline{\mathbb{B}},\mathbb{R}^q)}\leq C \cdot \left|E[F_0](0,f) \right|_{C^{2,\alpha}(\overline{\mathbb{B}},\mathbb{R}^q)}
\end{equation}

Die Konstanten $\vartheta$ und $C$ h\"angen dabei weder von $F_0$ noch von $m$ ab. Dar\"uber hinaus gilt: $u\in C^{\infty}(\overline{\mathbb{B}},\mathbb{R}^q)$ falls $f\in C^{\infty}(\overline{\mathbb{B}},\mathbb{R}^{\frac{n}{2}(n+1)})$. 

\end{satz}

\chapter{Konstruktion der isometrischen Einbettung}
\label{chap:Kap6}
\thispagestyle{fancy}

In \autoref{sec:6.1} werden Satz \ref{satz:4.1} und Satz \ref{satz:5.1} daf\"ur verwendet, um Satz \ref{satz:3.2} zu zeigen, womit das lokale Problem gel\"ost ist. Anschlie\ss{}end wird in \autoref{sec:6.2}, unter Verwendung von Satz \ref{satz:3.1}, der Hauptsatz \ref{haupt1} gezeigt. Mit dem Existenzsatz f\"ur eine freie Einbettung aus \autoref{chap:Kap2}, speziell Satz \ref{satz:2.3}, folgt die Existenz einer isometrischen Einbettung, f\"ur eine gegebene geschlossene Riemannsche Mannigfaltigkeit $(M,g)$. Dies entspricht Hauptsatz \ref{haupt2}.

\section{L\"osung des lokalen Problems}
\label{sec:6.1}

Zum Beweis von Satz \ref{satz:3.2}. Die Aussage des Satzes ist:

\begin{satzohnenum}
Gegeben sei eine n-dimensionale glatte Mannigfaltigkeit $M$, und eine freie Einbettung $F_0\in C^{\infty}(M,\mathbb{R}^q)$, mit $q\geq\frac{n}{2}(n+3)+5$. Weiterhin sei ein glattes symmetrisches kovariantes $2-$Tensorfeld $h\in \mathcal{T}^2(M)$ gegeben, f\"ur das eine Karte $(U,\varphi)\in\mathcal{A}$, mit $\varphi(U)=B_{1+\tau}(0)$, f\"ur $\tau>0$, und eine Abbildung $a\in C^{\infty}_0(U)$ existiert, welche $supp(a)\subseteq \varphi^{-1}(\mathbb{B})$ erf\"ullt, so dass:
\begin{equation}
h(x)=
\begin{cases}\label{eq:6.2}
a^4(x) \left.\, ^{\varphi}dx^1 \right|_x^2 &\text{f\"ur }x\in U\\
0 &\text{sonst}
\end{cases}
\end{equation}

gilt. Dann existiert zu jedem $\epsilon\in\mathbb{R}_{>0}$ eine freie Einbettung $F\in C^{\infty}(M,\mathbb{R}^q)$, so dass:
\begin{align}
\begin{split}\label{eq:6.8}
&F^{\ast}(g^{can})=F_0^{\ast}(g^{can})+h
\end{split}\\
\begin{split}\label{eq:6.9}
&supp(F-F_0)\subseteq U
\end{split}\\
\begin{split}\label{eq:6.10}
&\max_{x\in M}{|F(x)-F_0(x)|_{\mathbb{R}^q}}\leq C\cdot \epsilon
\end{split}
\end{align}

gilt, wobei die Konstante $C\in\mathbb{R}_{>0}$ nicht von $\epsilon$ abh\"angt.

\end{satzohnenum}

\begin{bew}
Die Abbildung $\left.\, ^{\varphi}F_0\right|_{\mathbb{B}}\in C^{\infty}(\overline{\mathbb{B}},\mathbb{R}^q)$ wird wieder mit $F_0$ bezeichnet. Mit Satz \ref{satz:4.1} existiert eine kompakte Menge $K\subseteq \mathbb{B}$, mit $supp(a)\subseteq K$, so dass die folgende Aussage gilt: F\"ur jedes $k\in\mathbb{N}\backslash\{0,1\}$ existiert ein $\epsilon_k>0$, so dass f\"ur alle $\epsilon\in (0,\epsilon_k]$ eine freie Abbildung $F_{\epsilon,k}\in C^{\infty}(\overline{\mathbb{B}},\mathbb{R}^q)$ existiert, welche die Eigenschaften:
\begin{enumerate}[(i)]
\begin{item}\label{6.1i}
\begin{align*}
&\partial_1 F_{\epsilon,k}\cdot \partial_1 F_{\epsilon,k}=\partial_1 F_0\cdot \partial_1 F_0 +a^4+\epsilon^{k+1}f_{11}^{\epsilon,k}\\
&\partial_1 F_{\epsilon,k}\cdot \partial_i F_{\epsilon,k}=\partial_1 F_0\cdot \partial_i F_0 +\epsilon^{k+1}f_{1i}^{\epsilon,k}\hspace{0.5cm}\text{f\"ur }2\leq i\leq n\\
&\partial_i F_{\epsilon,k}\cdot \partial_j F_{\epsilon,k}=\partial_i F_0\cdot \partial_j F_0 +\epsilon^{k+1}f_{ij}^{\epsilon,k}\hspace{0.5cm}\text{f\"ur }2\leq i\leq j\leq n
\end{align*}
\end{item}
\begin{item}\label{6.1ii}
\begin{equation*}
supp(F_{\epsilon,k}-F_0)\in C^{\infty}_0(\mathbb{B},\mathbb{R}^q)
\end{equation*}
\end{item}
\begin{item}\label{6.1iii}
\begin{equation*}
\max_{x\in\mathbb{B}}{|F_{\epsilon,k}(x)-F_0(x) |_{\mathbb{R}^q}}\leq C(n,k,F_0,a)\cdot \epsilon
\end{equation*}
\end{item}

erf\"ullt. Hierbei ist $f^{\epsilon,k}:=(f_{ij}^{\epsilon,k})_{1\leq i\leq j\leq n}\in C^{\infty}_0(\mathbb{B},\mathbb{R}^{\frac{n}{2}(n+1)})$ mit $supp(f^{\epsilon,k})\subseteq K$. Diese Abbildung erf\"ullt die Absch\"atzung:
\begin{item}\label{6.1iv}
\begin{equation*}
\left\Vert f^{\epsilon,k} \right\Vert_{C^3(\mathbb{\overline{\mathbb{B}}} ,\mathbb{R}^{\frac{n}{2}(n+1)}   )}\leq C(n,k,F_0,a)\cdot \epsilon^{-3}
\end{equation*}
\end{item}
\end{enumerate}

Nun sei $U_1\subseteq \mathbb{B}$ eine offene Menge mit $K\subseteq U_1$ und $\overline{U}_1\subseteq\overline{\mathbb{B}}$. Mit \eqref{eq:6.2}, \eqref{6.1i} und \eqref{6.1iv} gilt f\"ur $i,j\in\{1,...,n\}$ mit $i\leq j$ die Absch\"at\-zung:
\begin{align}\label{eq:6.3}
\begin{split}
&\left|\partial_i F_{\epsilon,k}\cdot \partial_j F_{\epsilon,k}-\partial_i F_0\cdot \partial_j F_0-\, ^{\varphi}h_{ij}\right|_{C^{2,\alpha}(\overline{\mathbb{B}})}=\epsilon^{k+1}\left|f_{ij}^{\epsilon,k}\right|_{C^{2,\alpha}(\overline{\mathbb{B}})}\\
\leq&C(n,\alpha)\cdot\epsilon^{k+1}\left\Vert f_{ij}^{\epsilon,k}\right\Vert_{C^{3}(\overline{\mathbb{B}})}\leq C(n,k,F_0,a)\cdot\epsilon^{k-2}
\end{split}
\end{align}

Hierbei sei $\alpha\in(0,1)$ auch f\"ur die folgenden Betrachtungen fest gew\"ahlt. Da $F_{\epsilon,k}$ f\"ur $\epsilon\in(0,\epsilon_k]$ eine freie Abbildung ist, so ist, f\"ur $m\in\mathbb{N}$, der in \eqref{eq:5.19} eingef\"uhrte Operator $E[F_{\epsilon,k}]$:
\begin{align*}
E[F_{\epsilon,k}]:  C^{m,\alpha}(\overline{\mathbb{B}},\mathbb{R}^n) \times C^{m,\alpha}(\overline{\mathbb{B}},\mathbb{R}^{\frac{n}{2}(n+1)})&\longrightarrow C^{m,\alpha}(\overline{\mathbb{B}},\mathbb{R}^q) \\  
E[F_{\epsilon,k}](\widehat{h},f)(x)&:=\Theta[F_{\epsilon,k}](x)\cdot\begin{pmatrix} \widehat{h}(x)\\ f(x) \end{pmatrix}
\end{align*}

wobei $\Theta[F_{\epsilon,k}]\in C^{\infty}(\overline{\mathbb{B}},\mathbb{R}^{q\times \frac{n}{2}(n+3)})$ mit:
\begin{equation*}
\Theta[F_{\epsilon,k}](x):= A^{\top}[F_{\epsilon,k}](x)\cdot( A[F_{\epsilon,k}](x)  A^{\top}[F_{\epsilon,k}](x))^{-1}
\end{equation*}

ist, wohldefiniert. Mit $A[F_{\epsilon,k}](x)\in \mathbb{R}^{\frac{n}{2}(n+3)\times q}$ wird hierbei die Matrix bezeichnet, deren Zeilenvektoren die Ableitungsvektoren, erster und zweiter Ordnung, der Abbildung $F_{\epsilon,k}$ sind. Nun soll die Operatornorm $\left\| E[F_{\epsilon,k}] \right\|_{2,\alpha}$ abgesch\"atzt werden.  Daf\"ur wird verwendet, dass f\"ur $x\in\mathbb{B}$ die Komponenten der Matrix 
\begin{equation*}
\Theta[F_{\epsilon,k}](x)=A^{\top}[F_{\epsilon,k}](x)\cdot( A[F_{\epsilon,k}](x)  A^{\top}[F_{\epsilon,k}](x))^{-1}
\end{equation*}

unter Beachtung der Cramerschen Regel, Polynome von Ableitungen von $F_{\epsilon,k}$, in erster und zweiter Ordnung in $x$, dividiert durch $det\left( (A[F_{\epsilon,k}](x)  A^{\top}[F_{\epsilon,k}](x))^{-1} \right)$, sind. Ferner gilt, wegen $F_{\epsilon,k}:= F_k(\epsilon,\beta(\epsilon^{-1}x_1),x)$, f\"ur $m\in\mathbb{N}$ die Absch\"atzung:
\begin{align*}
\left|F_{\epsilon,k} \right|_{C^{m,\alpha}(\overline{\mathbb{B}},\mathbb{R}^q)}&=\sum_{i=1}^q{\left|\left[ F_{\epsilon,k}\right]_i \right|_{C^{m,\alpha}(\overline{\mathbb{B}})}}\leq C\cdot\sum_{i=1}^q{\left|\left[ F_{\epsilon,k}\right]_i \right|_{C^{m+1}(\overline{\mathbb{B}})}}\leq C(k)\cdot\epsilon^{-(m+1)}
\end{align*}

Daraus folgt, f\"ur beliebige $\widehat{h}\in C^{2,\alpha}(\overline{\mathbb{B}},\mathbb{R}^n)$ und $f\in C^{2,\alpha}(\overline{\mathbb{B}},\mathbb{R}^{\frac{n}{2}(n+1)})$, die Absch\"atzung:
\begin{align*}
\left|E[F_{\epsilon,k}](\widehat{h},f) \right|_{C^{2,\alpha}(\overline{\mathbb{B}},\mathbb{R}^q)}\leq C(k)\cdot \epsilon^{-k_0}\cdot\left(|\widehat{h} |_{C^{2,\alpha}(\overline{\mathbb{B}},\mathbb{R}^n)}+|f |_{C^{2,\alpha}(\overline{\mathbb{B}},\mathbb{R}^{\frac{n}{2}(n+1)})} \right)
\end{align*}

f\"ur ein $k_0\in\mathbb{N}$. Dies impliziert die Ungleichung:
\begin{equation}\label{eq:6.4}
\left\| E[F_{\epsilon,k}] \right\|_{2,\alpha}\leq C(k)\cdot \epsilon^{-k_0}
\end{equation}

Unter der Annahme, dass $k\geq 2k_0+3$ gilt, ergibt sich, mit \eqref{eq:6.3} und \eqref{eq:6.4}, die Absch\"at\-zung:
\begin{align}\label{eq:6.5}
\begin{split}
&\left\| E[F_{\epsilon,k}] \right\|_{2,\alpha}\cdot \left| E[F_{\epsilon,k}](0,(\partial_i F_{\epsilon,k}\cdot \partial_j F_{\epsilon,k}-\partial_i F_0\cdot \partial_j F_0-\, ^{\varphi}h_{ij})_{1\leq i\leq j\leq n}) \right|_{C^{2,\alpha}(\overline{\mathbb{B}},\mathbb{R}^q)}\\
\leq &\left\| E[F_{\epsilon,k}] \right\|_{2,\alpha}^2\cdot | (\partial_i F_{\epsilon,k}\cdot \partial_j F_{\epsilon,k}-\partial_i F_0\cdot \partial_j F_0-\, ^{\varphi}h_{ij})_{1\leq i\leq j\leq n})    |_{C^{2,\alpha}(\overline{\mathbb{B}},\mathbb{R}^{\frac{n}{2}(n+1)})}\\
\leq & C(k)\cdot \epsilon^{-2 k_0+k-2}\leq C(k)\cdot \epsilon
\end{split}
\end{align}

Nun wird eine Menge $U_2\subseteq \mathbb{B}$, mit $\overline{U}_1\subseteq U_2$ und $\overline{U}_2\subseteq\mathbb{B}$, gew\"ahlt. Mit der Absch\"atzung \eqref{eq:6.5} existiert, unter der Annahme dass $\epsilon\in (0,\epsilon_k]$ klein genug ist, mit \eqref{6.1i} und Satz \ref{satz:5.1}, ein $u\in C^{\infty}(\overline{\mathbb{B}},\mathbb{R}^q)$ mit den Eigenschaften:
\begin{equation}\label{eq:6.6}
\partial_i(F_{\epsilon,k}+u)\cdot \partial_j(F_{\epsilon,k}+u)=\partial_i F_0\cdot \partial_j F_0+\,^{\varphi}h_{ij}\hspace{0.5cm}\text{f\"ur }1\leq i\leq j\leq n
\end{equation}
und
\begin{equation}\label{eq:6.7}
supp(u)\subseteq U_2
\end{equation}

Die Abbildung $u$ erf\"ullt mit \eqref{eq:5.53} die Absch\"atzung:
\begin{align}
\begin{split}\label{eq:6.11}
&\left|u \right|_{C^{2,\alpha}(\overline{\mathbb{B}},\mathbb{R}^q)}\\
\stackrel{\hphantom{\eqref{eq:6.3},\, \eqref{eq:6.4}}}{\leq} & C\cdot \left|E[F_{\epsilon,k}](0,(\partial_i F_{\epsilon,k}\cdot\partial_j F_{\epsilon,k}-\partial_i F_0\cdot \partial_j F_0-\,^{\varphi}h_{ij})_{1\leq i\leq j\leq n} )\right|_{C^{2,\alpha}(\overline{\mathbb{B}},\mathbb{R}^q)}\\
\stackrel{\hphantom{\eqref{eq:6.3},\, \eqref{eq:6.4}}}{\leq} &C\cdot \left\| E[F_{\epsilon,k}] \right\|_{2,\alpha}\cdot \left|(\partial_i F_{\epsilon,k}\cdot\partial_j F_{\epsilon,k}-\partial_i F_0\cdot \partial_j F_0-\,^{\varphi}h_{ij})_{1\leq i\leq j\leq n} \right|_{C^{2,\alpha}(\overline{\mathbb{B}},\mathbb{R}^{\frac{n}{2}(n+1)})}\\
\stackrel{\eqref{eq:6.3},\, \eqref{eq:6.4}}{\leq} &C(k)\cdot \epsilon^{- k_0+k-2}\stackrel{k\geq 2k_0+3}{\leq} C(k)\cdot \epsilon^{k_0+1}\stackrel{k_0\geq 0}{\leq}  C(k)\cdot \epsilon
\end{split}
\end{align}

Schlie\ss{}lich wird, unter Beachtung von \eqref{eq:6.7} und \eqref{6.1ii}, die Funktion $F\in C^{\infty}(M,\mathbb{R}^q)$ mit:
\begin{align*}
F(x):=
\begin{cases}
F_{\epsilon,k}(\varphi(x))+u(\varphi(x)) &\text{f\"ur }x\in \varphi^{-1}(\mathbb{B})\\
F_0(x) &\text{sonst}
\end{cases}
\end{align*}

definiert. Mit \eqref{eq:6.6} und \eqref{eq:6.2} ist  \eqref{eq:6.8} erf\"ullt. Aus \eqref{6.1iii} und \eqref{eq:6.11} folgt die Absch\"atzung \eqref{eq:6.10}. Es bleibt noch zu zeigen, dass $F\in C^{\infty}(M,\mathbb{R}^q)$, genau wie $F_0\in C^{\infty}(M,\mathbb{R}^q)$, eine freie Einbettung ist, sofern $\epsilon\in (0,\epsilon_k]$ klein genug ist. Mit Satz \ref{satz:4.4} \eqref{satz:4.4.i} ist $F\in C^{\infty}(M,\mathbb{R}^q)$, f\"ur kleine $\epsilon\in(0,\epsilon_k]$, eine freie Abbildung, und damit insbesondere eine Immersion. Mit \cite[Proposition 7.4.]{lee2003introduction} ist nur noch zu zeigen, dass $F\in C^{\infty}(M,\mathbb{R}^q)$ f\"ur kleine $\epsilon\in(0,\epsilon_k]$ auch injektiv ist. Angenommen, es existieren $x_1,x_2\in M$ mit $x_1\neq x_2$, so dass $F(x_1)=F(x_2)$ gilt. Mit Satz \ref{satz:4.4} \eqref{satz:4.4.ii} kann ausgeschlossen werden, dass $x_1,x_2\in U$ gilt. Weiterhin kann, wegen \eqref{eq:6.7}, die M\"oglichkeit $x_1,x_2\in M\backslash \varphi^{-1}(U_2)$ ebenfalls ausgeschlossen werden. Es bleibt noch der Fall zu untersuchen, dass $x_1\in M\backslash U$ und $x_2\in \varphi^{-1}(U_2)$ gilt. Da $F_0\in C^{\infty}(M,\mathbb{R}^q)$ eine freie Einbettung ist, und $\overline{U}_2\subseteq \mathbb{B}$ gilt, ist:
\begin{equation*}
\inf_{x\in M\backslash U,\, y\in \varphi^{-1}(U_2)}{|F_0(x)-F_0(y)|_{\mathbb{R}^q}}>0
\end{equation*}

Ferner folgt mit \eqref{eq:6.10}, unter Verwendung der Dreiecksungleichung, f\"ur alle $x,y\in M$:
\begin{equation*}
|F(x)-F(y)|_{\mathbb{R}^q}\geq |F_0(x)-F_0(y)|_{\mathbb{R}^q}-2\cdot C(k)\cdot \epsilon
\end{equation*}

Daraus ergibt sich, f\"ur hinreichend kleine $\epsilon\in (0,\epsilon_k]$, die Absch\"atzung, :
\begin{equation*}
\inf_{x\in M\backslash U,\, y\in \varphi^{-1}(U_2)}{|F(x)-F(y)|_{\mathbb{R}^q}}>0
\end{equation*}

womit die Injektivit\"at der Abbildung $F$ auf ganz $M$ gezeigt ist.

\end{bew}

Satz \ref{satz:3.2} wird nun daf\"ur verwendet, um Hauptsatz \ref{haupt1} zu zeigen.

\section{Beweis der Haups\"atze}
\label{sec:6.2}

Die Aussage von Hauptsatz \ref{haupt1} ist:

\begin{hauptohnenum}
Gegeben sei eine $n$-dimensionale geschlossene Riemannsche Mannigfaltigkeit $(M,g)$, und eine freie Einbettung $F_0\in C^{\infty}(M,\mathbb{R}^q)$, mit $q\geq \frac{n}{2}(n+3)+5$, so dass $g-F_0^{\ast}(g^{can})\in\mathcal{T}^2(M)$ eine Riemannsche Metrik ist. Ferner sei ein $\delta\in\mathbb{R}_{>0}$ gegeben, dann existiert eine freie Einbettung $F\in C^{\infty}(M,\mathbb{R}^q)$, so dass $F^{\ast}(g^{can})=g$ und:
\begin{equation}\label{eq:6.17}
\max_{x\in M}{|F(x)-F_0(x)|_{\mathbb{R}^q}}\leq \delta
\end{equation}

gilt.

\end{hauptohnenum}

\begin{bew}
Mit Satz \ref{satz:3.1} existieren symmetrische kovariante $2$-Tensorfelder $h_1,...,h_m\in\mathcal{T}^2(M)$, so dass:
\begin{equation}\label{eq:6.13}
h:= g-F_0^{\ast}(g^{can})=\sum_{i=1}^m{h^{(i)}}
\end{equation}

und f\"ur jedes $i\in\{1,...,m\}$ eine Karte $(U_i,\varphi_i)$, mit $\varphi_i(U_i)=B_{1+\tau}(0)$, f\"ur $\tau>0$, und eine Abbildung $a_i\in C^{\infty}(U_i)$, mit $supp(a_i)\subseteq \varphi^{-1}(\mathbb{B})$ existiert, so dass:
\begin{align*}
h^{(i)}(x)=
\begin{cases}
a_i^4(x) \left.\, ^{\varphi_i}dx^1 \right|_x^2 &\text{f\"ur }x\in U_i\\
0 &\text{sonst}
\end{cases}
\end{align*}

gilt. Mit dem, in \autoref{sec:6.1} bewiesenen, Satz \ref{satz:2.2}, k\"onnen sukzessiv freie Einbettungen $F_1,...,F_m\in C^{\infty}(M,\mathbb{R}^q)$ konstruiert werden, so dass f\"ur jedes $i\in\{1,...,m\}$ sowohl:
\begin{equation}\label{eq:6.14}
F^{\ast}_i(g^{can})=F^{\ast}_{i-1}(g^{can})+h^{(i)}
\end{equation}

als auch:
\begin{align}
\begin{split}\label{eq:6.15}
&\max_{x\in M}{|F_i(x)-F_{i-1}(x)|_{\mathbb{R}^q}}\leq \frac{\delta}{m}\\
&F_i(x)=F_{i-1}(x)\hspace{0.5cm}\text{f\"ur alle }x\in M\backslash{U_i}
\end{split}
\end{align}

gilt. Aus \eqref{eq:6.13} und \eqref{eq:6.14} folgt:
\begin{equation*}
F^{\ast}_m(g^{can})=F^{\ast}_0(g^{can})+\sum_{i=1}^m{h^{(i)}}\stackrel{\eqref{eq:6.13}}{=}g
\end{equation*}

Die Abbildung $F:= F_m\in C^{\infty}(M,\mathbb{R}^q)$ ist also eine isometrische Einbettung der Riemannschen Mannigfaltigkeit $(M,g)$. Abschlie\ss{}end wird noch \eqref{eq:6.17} gezeigt. Mit \eqref{eq:6.15} gilt f\"ur jedes $x\in M$ die Absch\"atzung:
\begin{equation*}
|F(x)-F_0(x)|_{\mathbb{R}^q}\leq\sum_{i=1}^m{|F_i(x)-F_{i-1}(x)|_{\mathbb{R}^q}}\leq\delta
\end{equation*}

Damit ist der Satz bewiesen.

\end{bew}

Zusammen mit dem Satz \ref{satz:2.3} folgt Hauptsatz \ref{haupt2}.

\chapter{Isometrische Einbettung einer Familie von Riemannschen Mannigfaltigkeiten}
\label{chap:Kap7}
\thispagestyle{fancy}

In diesem, abschlie\ss{}enden Kapitel wird eine M\"oglichkeit beschrieben, wie Hauptsatz \ref{haupt2} und die in \autoref{chap:Kap5} aufgef\"uhrten Resultate, daf\"ur verwendet werden k\"onnen, eine Familie von Riemannschen Mannigfaltigkeiten isometrisch in einen Vektorraum $\mathbb{R}^q$ einzubetten. Daf\"ur wird zun\"achst ein lokaler Einbettungssatz gezeigt, mit dessen Hilfe dann die gesamte Mannigfaltigkeit f\"ur eine kurze Zeit isometrisch eingebettet wird. Im gesamten Kapitel sei $\alpha\in (0,1)$ fest.

\section{Konstruktion einer lokaler L\"osung}
\label{sec:7.1}

Das folgende Lemma basiert direkt auf Lemma \ref{lem:5.13}. Es besagt, dass $\vartheta\in\mathbb{R}_{>0}$ so gew\"ahlt werden kann, dass zwei verschiedene Fixpunkte des Operators $\Phi$, zu jeweils gegebenem $f\in C^{\infty}(\overline{\mathbb{B}},\mathbb{R}^{\frac{n}{2}(n+1)})$ und $g\in C^{\infty}(\overline{\mathbb{B}},\mathbb{R}^{\frac{n}{2}(n+1)})$, eine gewisse Stetigkeitsabsch\"atzung erf\"ullen, falls die \glqq Energien\grqq\, klein genug sind.
 
\begin{lem}\label{lem:7.1}
Es sei $a\in C^{\infty}_0(\mathbb{B})$, dann existiert ein $\vartheta(n,\alpha,a)\in\mathbb{R}_{>0}$ mit der folgenden Eigenschaft: Ist eine freie Abbildung $F_0\in C^{\infty}(\overline{\mathbb{B}},\mathbb{R}^q)$ gegeben, und sind $f,g\in C^{\infty}(\overline{\mathbb{B}},\mathbb{R}^{\frac{n}{2}(n+1)})$ mit:
\begin{align}\label{eq:7.3}
\begin{split}
\left\Vert E[F_0]\right\Vert_{2,\alpha}\left|E[F_0](0,f) \right|_{C^{2,\alpha}(\overline{\mathbb{B}},\mathbb{R}^q)}&\leq \vartheta\\
\left\Vert E[F_0]\right\Vert_{2,\alpha}\left|E[F_0](0,g) \right|_{C^{2,\alpha}(\overline{\mathbb{B}},\mathbb{R}^q)}&\leq \vartheta
\end{split}
\end{align}

dann erf\"ullen die, in Lemma \ref{lem:5.13} konstruierten, Abbildungen $v^{(1)},v^{(2)}\in C^{\infty}(\overline{\mathbb{B}},\mathbb{R}^q)$, welche jeweils die Fixpunktgleichung:
\begin{align*}
&v^{(1)}=\Phi[F_0,a,f](v^{(1)}) & & v^{(2)}=\Phi[F_0,a,g](v^{(2)})
\end{align*}

erf\"ullen, die Absch\"atzung:
\begin{equation}\label{eq:7.5}
\left|v^{(1)}-v^{(2)}\right|_{C^{2,\alpha}(\overline{\mathbb{B}},\mathbb{R}^q)}\leq \left|E[F_0](0,f-g) \right|_{C^{2,\alpha}(\overline{\mathbb{B}},\mathbb{R}^q)}
\end{equation}

\end{lem}

\begin{bew}
Es seien $(v^{(1)}_k)_{k\in\mathbb{N}}, (v^{(2)}_k)_{k\in\mathbb{N}}\subseteq C^{\infty}(\overline{\mathbb{B}},\mathbb{R}^q)$, die jeweiligen, in \eqref{eq:5.47} konstruierten, Approximierenden. Das bedeutet konkret:
\begin{align*}
v^{(1)}_k=
\begin{cases}
0 & \text{falls }  k=0\\
\Phi[F_0,a,f](v^{(1)}_{k-1})  &\text{falls } k\geq 1
\end{cases}
\end{align*}

und:
\begin{align*}
v^{(2)}_k=
\begin{cases}
0 & \text{falls }  k=0\\
\Phi[F_0,a,g](v^{(2)}_{k-1})  &\text{falls } k\geq 1
\end{cases}
\end{align*}

Dann gilt f\"ur $k\in\mathbb{N}\backslash\{0\}$:
\begin{align}\label{eq:7.4}
\begin{split}
&\left|v^{(1)}_k-v^{(2)}_k\right|_{C^{2,\alpha}(\overline{\mathbb{B}},\mathbb{R}^q)}=\left|\Phi[F_0,a,f](v^{(1)}_{k-1})-\Phi[F_0,a,g](v^{(2)}_{k-1})\right|_{C^{2,\alpha}(\overline{\mathbb{B}},\mathbb{R}^q)}\\
\stackrel{\eqref{eq:5.18}}{=} & \left|E[F_0]\left(P[a](v^{(1)}_{k-1})-P[a](v^{(2)}_{k-1}), \frac{1}{2}\left[(f-g)-Q[a](v^{(1)}_{k-1})+ Q[a](v^{(2)}_{k-1})\right]\right)\right|_{C^{2,\alpha}(\overline{\mathbb{B}},\mathbb{R}^q)}\\
\stackrel{\hphantom{\eqref{eq:5.18}}}{\leq} & \left\Vert E[F_0]\right\Vert_{2,\alpha}\cdot  \left|P[a](v^{(1)}_{k-1})-P[a](v^{(1)}_{k-1}) \right|_{C^{2,\alpha}(\overline{\mathbb{B}},\mathbb{R}^n)}\\
&\hspace{3cm}+\left\Vert E[F_0]\right\Vert_{2,\alpha}\cdot \left|Q[a](v^{(1)}_{k-1})-Q[a](v^{(2)}_{k-1}) \right|_{C^{2,\alpha}(\overline{\mathbb{B}},\mathbb{R}^{\frac{n}{2}(n+1)})}\\
&\hspace{3cm}+\frac{1}{2}  \left|E[F_0](0,f-g)\right|_{C^{2,\alpha}(\overline{\mathbb{B}},\mathbb{R}^q)}\\
\stackrel{\eqref{eq:5.37}}{\leq} & K(n,\alpha,a)\cdot\left\Vert E[F_0]\right\Vert_{2,\alpha}\cdot \left(|v^{(1)}_{k-1} |_{C^{2,\alpha}(\overline{\mathbb{B}},\mathbb{R}^q)} +|v^{(2)}_{k-1} |_{C^{2,\alpha}(\overline{\mathbb{B}},\mathbb{R}^q)}\right)\cdot|v^{(1)}_{k-1}-v^{(2)}_{k-1} |_{C^{2,\alpha}(\overline{\mathbb{B}},\mathbb{R}^q)}\\
&\hspace{3cm}+\frac{1}{2}  \left|E[F_0](0,f-g)\right|_{C^{2,\alpha}(\overline{\mathbb{B}},\mathbb{R}^q)}\\
\stackrel{\eqref{eq:5.48}}{\leq} & K(n,\alpha,a)\cdot\left\Vert E[F_0]\right\Vert_{2,\alpha} \left(\left|E[F_0](0,f) \right|_{C^{2,\alpha}(\overline{\mathbb{B}},\mathbb{R}^q)}+\left|E[F_0](0,g) \right|_{C^{2,\alpha}(\overline{\mathbb{B}},\mathbb{R}^q)}\right)\\
&\hspace{2cm}\cdot|v^{(1)}_{k-1}-v^{(2)}_{k-1}|_{C^{2,\alpha}(\overline{\mathbb{B}},\mathbb{R}^q)}+\frac{1}{2}  \left|E[F_0](0,f-g)\right|_{C^{2,\alpha}(\overline{\mathbb{B}},\mathbb{R}^q)}
\end{split}
\end{align}

Ist nun:
\begin{equation*}
\left\Vert E[F_0]\right\Vert_{2,\alpha} \left(\left|E[F_0](0,f) \right|_{C^{2,\alpha}(\overline{\mathbb{B}},\mathbb{R}^q)}+\left|E[F_0](0,g) \right|_{C^{2,\alpha}(\overline{\mathbb{B}},\mathbb{R}^q)}\right)\leq \frac{1}{2 K(n,\alpha,a)}
\end{equation*}

dann folgt aus \eqref{eq:7.4} die Absch\"atzung:
\begin{align*}
&|v^{(1)}_k-v^{(2)}_k|_{C^{2,\alpha}(\overline{\mathbb{B}},\mathbb{R}^q)}\leq \frac{1}{2}\left(|v^{(1)}_{k-1}-v^{(2)}_{k-1} |_{C^{2,\alpha}(\overline{\mathbb{B}},\mathbb{R}^q)}+ \left|E[F_0](0,f-g)\right|_{C^{2,\alpha}(\overline{\mathbb{B}},\mathbb{R}^q)} \right)
\end{align*}

woraus induktiv, wegen $v^{(1)}_0=0=v^{(2)}_0$, die Absch\"atzung:
\begin{align*}
&|v^{(1)}_k-v^{(2)}_k|_{C^{2,\alpha}(\overline{\mathbb{B}},\mathbb{R}^q)}\leq  \left|E[F_0](0,f-g)\right|_{C^{2,\alpha}(\overline{\mathbb{B}},\mathbb{R}^q)}
\end{align*}

f\"ur alle $k\in\mathbb{N}$ folgt. Die Behauptung ist damit bewiesen.

\end{bew}

\begin{satz}\label{satz:7.1}
Gegeben sei ein kompakte $n$-dimensionale glatte Mannigfaltigkeit $M$, ferner sei eine Familie von Riemannschen Metriken $(g(\cdot,t))_{t\in[0,T]}\subseteq \mathcal{T}^2(M)$ gegeben, so dass f\"ur jedes $(U,\varphi)\in\mathcal{A}$ die Bedingung $(\, ^{\varphi}g_{ij}(\cdot,\cdot))_{1\leq i\leq j\leq n}\in C^{\infty}(\varphi(U)\times[0,T],\mathbb{R}^{\frac{n}{2}(n+1)})$ gilt. Dann existiert, f\"ur $q(n):=\max\left\{\frac{n}{2}(n+5),\frac{n}{2}(n+3)+5 \right\}$, f\"ur jedes $x\in M$ eine Umgebung $U_x\subseteq M$ von $x$, ein $T_x\in (0,T]$ und eine Familie von freien Einbettungen $(F(\cdot,t))_{t\in [0,T_x]}\subseteq C^{\infty}(\overline{U}_x,\mathbb{R}^q)$, so dass:
\begin{equation}\label{eq:7.1}
F(\cdot,t)^{\ast}(g^{can})=g(\cdot,t)
\end{equation}

auf $U_x$ f\"ur alle $t\in[0,T_x]$ und:
\begin{equation}\label{eq:7.2}
F\in C^0(\overline{U}_x \times[0,T_x],\mathbb{R}^q)
\end{equation}

gilt.
\end{satz}

\begin{bew}
Zun\"achst sei $F_0\in C^{\infty}(M,\mathbb{R}^q)$ eine freie isometrische Einbettung der Riemannschen Mannigfaltigkeit $(M,g(\cdot,0))$, gem\"a\ss{} Hauptsatz \ref{haupt2}. W\"ahle ein $(U,\varphi)\in \mathcal{A}$ mit $\varphi(U)=B_{1+\tau}(0)$ und $\varphi(x)=0$, sowie eine Funktion $\psi\in C^{\infty}(\mathbb{B})$ so dass $\left. \psi\right|_{\overline{V_x}}\equiv 1$ f\"ur eine Umgebung $V_x\subseteq\mathbb{B}$ von $0$, und $supp(\psi)\subseteq \mathbb{B}$ gilt. Die Abbildung $\left.\,^{\varphi_2}F_0 \right|_{\mathbb{B}}\in C^{\infty}(\overline{\mathbb{B}},\mathbb{R}^q)$ wird im Folgenden kurz mit $F_0$ bezeichnet. Nun sei $\widehat{g}\in C^{\infty}(\overline{\mathbb{B}},\mathbb{R}^{\frac{n}{2}(n+1)})$ mit:
\begin{equation}\label{eq:7.17}
\widehat{g}(x,t):=(\widehat{g}_{ij}(x,t))_{1\leq i\leq j\leq n}=\psi(x)\cdot ( \,^{\varphi}g_{ij}(x,t)-\,^{\varphi}g_{ij}(x,0))_{1\leq i\leq j\leq n}
\end{equation}

W\"ahle dann ein $a\in C^{\infty}_0(\mathbb{B})$ mit $\left. a\right|_{supp(\psi)}\equiv 1$, so dass $a^2\, \widehat{g}\equiv \widehat{g} $ gilt, dann sei $\widehat{T}_x\in (0,T]$, so dass f\"ur alle $t\in [0,T_x]$ die Absch\"atzung:
\begin{equation}
\left\Vert E[F_0]\right\Vert_{2,\alpha}\left|E[F_0](0,\widehat{g}(\cdot,t)) \right|_{C^{2,\alpha}(\overline{\mathbb{B}},\mathbb{R}^q)}\leq \vartheta(n,\alpha,a)
\end{equation}

mit dem in Satz \ref{satz:5.1}, beziehungsweise Lemma \ref{lem:7.1}, definierten $\vartheta\in\mathbb{R}_{>0}$ erf\"ullt ist. Dann existiert f\"ur jedes $t\in [0,T_x]$ ein $v(\cdot,t)\in C^{\infty}(\overline{\mathbb{B}},\mathbb{R}^q)$, so dass f\"ur $u(\cdot,t):= a^2\, v(\cdot,t)\in C^{\infty}_0(\overline{\mathbb{B}},\mathbb{R}^q) $, f\"ur alle $(x,t)\in \mathbb{B}\times [0,T_x]$, und f\"ur jede Auswahl $i,j\in\{1,...,n\}$ mit $i\leq j$, die Gleichung:
\begin{align}\label{eq:7.16}
\begin{split}
&\partial_i (\,^{\varphi}F_0(x)+u(x,t))\cdot \partial_j (\,^{\varphi}F_0(x)+u(x,t))\\
\stackrel{\hphantom{\eqref{eq:7.17}}}{=}&\partial_i \,^{\varphi}F_0(x)\cdot \partial_j \,^{\varphi}F_0(x)+\widehat{g}_{ij}(x,t)\\
\stackrel{\eqref{eq:7.17}}{=}&\,^{\varphi}g_{ij}(x,0)+\psi(x)\cdot (\,^{\varphi}g_{ij}(x,t)-\,^{\varphi}g_{ij}(x,0))
\end{split}
\end{align}

erf\"ullt ist. Ferner sei f\"ur sp\"atere Zwecke erw\"ahnt, dass f\"ur jedes $t\in[0,T_x]$ die Abbildung $v(\cdot,t)\in C^{\infty}(\overline{\mathbb{B}},\mathbb{R}^q)$, gem\"a\ss{} \eqref{eq:5.47}, durch die Folge $(v_k(\cdot,t))_{k\in\mathbb{N}}\subseteq C^{\infty}(\overline{\mathbb{B}},\mathbb{R}^q)$ mit:
\begin{align}\label{eq:7.15}
v_k(\cdot,t)=
\begin{cases}
0 & \text{falls }  k=0\\
\Phi[F_0,a,\widehat{g}(\cdot,t)](v_{k-1}(\cdot,t))  &\text{falls } k\geq 1
\end{cases}
\end{align}

im $C^{2,\alpha}$-Sinne approximiert wird. Mit \eqref{eq:5.48} gilt die Absch\"atzung:
\begin{equation}\label{eq:7.18}
\max_{t\in[0,T_x]}\left|v_k(\cdot,t) \right|_{C^{2,\alpha}(\overline{\mathbb{B}},\mathbb{R}^q)}\leq \max_{t\in[0,T_x]}\left|E[F_0]\left(0,\widehat{g}(\cdot,t)\right)\right|_{C^{2,\alpha}(\overline{\mathbb{B}},\mathbb{R}^q)}
\end{equation}

 Da $\left. \psi\right|_{\overline{V_x}}\equiv 1$ gilt, folgt f\"ur $U_x:= \varphi^{-1}(V_x)$ und $(F(\cdot,t))_{t\in[0,T]}\subseteq C^{\infty}(U_x,\mathbb{R}^q)$ mit:
\begin{equation}\label{7.26}
F(x,t)=F_0(x)+u(\varphi(x),t) 
\end{equation}

aus \eqref{eq:7.16}, unter Beachtung von \eqref{eq:7.5}, die Behauptung.

\end{bew}


\section{Regularit\"at der lokalen L\"osung}
\label{sec:7.2}

Um die Regularit\"at der, in \autoref{sec:7.1} konstruierten, lokalen L\"osung zu untersuchen, werden Eigenschaften der in \eqref{eq:7.15} eingef\"uhrten Approximierenden untersucht. Dabei sind vor allem die Ableitungen in Zeitrichtung, sofern existent, von Interesse. Hierf\"ur werden die in \autoref{sec:5.4} gezeigten Absch\"atzungen f\"ur die Betrachtung von Zeitableitungen aufgearbeitet. Im Folgenden sei $T\in\mathbb{R}_{>0}$ fest gew\"ahlt.

\begin{lem}\label{lem:7.2}
Es sei $a\in C_0^{\infty}(\mathbb{B})$ dann gilt f\"ur $v_1,v_2\in C^{\infty}(\overline{{\mathbb{B}}}\times [0,T],\mathbb{R}^q)$:
\begin{align}\label{eq:7.6}
\begin{split}
&\max_{t\in[0,T]}\left| N_i[a](v_1(\cdot,t) )-N_i[a](v_2(\cdot,t)) \right|_{C^{0,\alpha}(\overline{\mathbb{B}})} \\
\leq &K(n,\alpha,a)\cdot \left(\max_{t\in[0,T]}\left|v_1(\cdot,t) \right|_{C^{2,\alpha}(\overline{\mathbb{B}},\mathbb{R}^q)}+\max_{t\in[0,T]}\left|v_2(\cdot,t) \right|_{C^{2,\alpha}(\overline{\mathbb{B}},\mathbb{R}^q)} \right)\\
&\cdot\max_{t\in[0,T]}\left|v_1(\cdot,t)-v_2(\cdot,t) \right|_{C^{2,\alpha}(\overline{\mathbb{B}},\mathbb{R}^q)}
\end{split}
\end{align}

F\"ur $v\in C^{\infty}(\overline{{\mathbb{B}}}\times [0,T],\mathbb{R}^q)$ und $m,r\in\mathbb{N}$ gilt:
\begin{align}\label{eq:7.7}
\begin{split}
&\max_{t\in[0,T]}\left|\partial_t^r N_i[a](v(\cdot,t)) \right|_{C^{m,\alpha}(\overline{\mathbb{B}})}\\
\leq &K(n,\alpha,a)\cdot\sum_{s=0}^r{\binom{r}{s}\max_{t\in[0,T]}\left|\partial_t^s  v(\cdot,t) \right|_{C^{m+2,\alpha}(\overline{\mathbb{B}},\mathbb{R}^q)}\max_{t\in[0,T]}\left|\partial_t^{r-s} v(\cdot,t) \right|_{C^{2,\alpha}(\overline{\mathbb{B}},\mathbb{R}^q)}}\\
&+C(n,m,\alpha,a)\cdot\sum_{s=0}^r{\binom{r}{s}\max_{t\in[0,T]}\left|\partial_t^s v(\cdot,t)\right|_{C^{m+1,\alpha}(\overline{\mathbb{B}},\mathbb{R}^q)}\max_{t\in[0,T]}\left|\partial_t^{r-s} v(\cdot,t) \right|_{C^{m+1,\alpha}(\overline{\mathbb{B}},\mathbb{R}^q)}}
\end{split}
\end{align}

\end{lem}

\begin{bew}
Um \eqref{eq:7.6} zu zeigen, kann die Ungleichung \eqref{eq:5.31} in Lemma \ref{lem:5.b} auf feste Zeiten angewandt werden. Nun wird \eqref{eq:7.7} gezeigt:
\begin{align*}
&\max_{t\in[0,T]}\left|\partial_t^r N_i[a](v(\cdot,t))\right|_{C^{m,\alpha}(\overline{\mathbb{B}})}\\
\stackrel{\eqref{eq:5.2}}{=}&\max_{t\in[0,T]}\left|2\partial_i a\, \partial_t^r[ \Delta v(\cdot,t)\cdot  v(\cdot,t)]+a \partial_t^r[\Delta v(\cdot,t)\cdot \partial_i v(\cdot,t)]\right|_{C^{m,\alpha}(\overline{\mathbb{B}})}\\
\stackrel{\hphantom{\eqref{eq:5.2}}}{=}&\max_{t\in[0,T]}\left|2\partial_i a\, \sum_{s=0}^r \binom{r}{s} \Delta \partial^s_t v(\cdot,t)\cdot  \partial^{r-s}_t v(\cdot,t)+a \,\sum_{s=0}\binom{r}{s}\Delta \partial^{s}_t v(\cdot,t)\cdot \partial_i \partial^{r-s}_t v(\cdot,t)\right|_{C^{m,\alpha}(\overline{\mathbb{B}})}\\
\stackrel{\hphantom{\eqref{eq:A.15}}}{\leq} &\sum_{s=0}^r \binom{r}{s} \max_{t\in[0,T]}\left|2\partial_i a\, \Delta \partial^{s}_t v(\cdot,t) \cdot \partial^{r-s}_t v(\cdot,t)\right|_{C^{m,\alpha}(\overline{\mathbb{B}})}\\
&\hspace{4,5cm}+\sum_{s=0}^r \binom{r}{s}\max_{t\in[0,T]}\left|a\,\Delta \partial^{s}_t v(\cdot,t) \cdot \partial_i \partial^{r-s}_t v(\cdot,t)\right|_{C^{m,\alpha}(\overline{\mathbb{B}})}\\
\stackrel{\eqref{eq:A.15}}{\leq}&\sum_{s=0}^r \binom{r}{s}\left|2\partial_i a\, \Delta \partial^{s}_t v(\cdot,t)\right|_{C^{0,\alpha}(\overline{\mathbb{B}},\mathbb{R}^q)} \left|\partial^{r-s}_t v(\cdot,t)\right|_{C^{m,\alpha}(\overline{\mathbb{B}},\mathbb{R}^q)}\\
&+\sum_{s=0}^r \binom{r}{s}\max_{t\in[0,T]}\left|2\partial_i a\, \Delta \partial^{s}_t v(\cdot,t)\right|_{C^{m,\alpha}(\overline{\mathbb{B}},\mathbb{R}^q)}\max_{t\in[0,T]}\left|\partial^{r-s}_t v(\cdot,t)\right|_{C^{0,\alpha}(\overline{\mathbb{B}},\mathbb{R}^q)}\\
+&C_1(n,m,\alpha)\sum_{s=0}^r \binom{r}{s}\max_{t\in[0,T]}\left|2\partial_i a\, \Delta \partial^{s}_t v(\cdot,t)\right|_{C^{m-1,\alpha}(\overline{\mathbb{B}},\mathbb{R}^q))} \max_{t\in[0,T]}\left|\partial^{r-s}_t v(\cdot,t)\right|_{C^{m-1,\alpha}(\overline{\mathbb{B}},\mathbb{R}^q)}\\
&+\sum_{s=0}^r \binom{r}{s} \max_{t\in[0,T]}\left|a\,\Delta \partial^{s}_t v(\cdot,t)\right|_{C^{0,\alpha}(\overline{\mathbb{B}},\mathbb{R}^q)}\max_{t\in[0,T]}\left|\partial_i \partial^{r-s}_t  v(\cdot,t)\right|_{C^{m,\alpha}(\overline{\mathbb{B}},\mathbb{R}^q))}\\
&+\sum_{s=0}^r \binom{r}{s}\max_{t\in[0,T]}\left|a\,\Delta \partial^{s}_t v(\cdot,t)\right|_{C^{m,\alpha}(\overline{\mathbb{B}},\mathbb{R}^q)}\max_{t\in[0,T]}\left|\partial_i \partial^{r-s}_t v(\cdot,t)\right|_{C^{0,\alpha}(\overline{\mathbb{B}},\mathbb{R}^q)}\\
&+C_1(n,m,\alpha) \sum_{s=0}^r \binom{r}{s} \max_{t\in[0,T]}\left|a\,\Delta \partial^{s}_t v(\cdot,t)\right|_{C^{m-1,\alpha}(\overline{\mathbb{B}},\mathbb{R}^q)} \max_{t\in[0,T]}\left|\partial_i \partial^{r-s}_t v(\cdot,t)\right|_{C^{m-1,\alpha}(\overline{\mathbb{B}},\mathbb{R}^q)}\\
\stackrel{\eqref{eq:A.13}}{\leq} &K(n,\alpha,a)\cdot\sum_{s=0}^r{\binom{r}{s}\max_{t\in[0,T]}\left|\partial_t^s  v(\cdot,t) \right|_{C^{m+2,\alpha}(\overline{\mathbb{B}},\mathbb{R}^q)}\max_{t\in[0,T]}\left|\partial_t^{r-s} v(\cdot,t) \right|_{C^{2,\alpha}(\overline{\mathbb{B}},\mathbb{R}^q)}}\\
&+C(n,m,\alpha,a)\cdot\sum_{s=0}^r{\binom{r}{s}\max_{t\in[0,T]}\left|\partial_t^s v(\cdot,t)\right|_{C^{m+1,\alpha}(\overline{\mathbb{B}},\mathbb{R}^q)}\max_{t\in[0,T]}\left|\partial_t^{r-s} v(\cdot,t) \right|_{C^{m+1,\alpha}(\overline{\mathbb{B}},\mathbb{R}^q)}}
\end{align*}
\end{bew}

Aus dem folgenden Lemma wird sich die Glattheit der Operatoren $M_{ij}[a]$ aus \eqref{eq:5.11}, unter Ber\"ucksichtigung der neu eingef\"uhrten Zeitrichtung, ergeben.

\begin{lem}\label{lem:7.3}
F\"ur ein $f\in C^{\infty}(\overline{{\mathbb{B}}}\times [0,T])$ sei 
u: $\mathbb{B}\times [0,T]\longrightarrow \mathbb{R}$, mit $u(\cdot,t)\in C^{\infty}(\overline{\mathbb{B}})$  f\"ur alle $t\in[0,T]$, die eindeutig bestimmte L\"osung der Gleichung:
\begin{align}\label{eq:7.8}
\begin{cases}
\Delta u(\cdot,t)=f(\cdot,t) &\text{auf }\mathbb{B}\\
\hphantom{\Delta} u(\cdot,t)=0 &\text{auf }\partial \mathbb{B}
\end{cases}
\end{align}

auf $\mathbb{B}$, dann gilt: $u\in C^{\infty}(\overline{{\mathbb{B}}}\times [0,T])$.
\end{lem}

\begin{bew}
Es seien $\varrho\in (0,1)$ sowie $T_1,T_2\in (0,T)$ mit $T_1< T_2$. F\"ur $h\in (0,\max\{T_1, T-T_2\})$ sei $\gls{differ} \in C^{\infty}(\overline{\mathbb{B}}\times [T_1,T_2])$ mit:
\begin{equation*}
(\partial_{t,h}u)(x,t):= \frac{u(x,t+h)-u(x,t)}{h}
\end{equation*}

der \textbf{Differenzenquotient von} $\bm{u}$ \textbf{in} $\bm{t}$-\textbf{Richtung zur Schrittweite} $\bm{h}$\index{Differenzenquotient}. F\"ur $h_1,...,h_r\in\mathbb{R}_{>0}$ mit $\sum_{i=1}^r{h_i}< \max\{T_1, T-T_2\}$ folgt dann induktiv aus \eqref{eq:7.8} f\"ur alle $[T_1,T_2]$:
\begin{align}\label{eq:7.a}
\begin{cases}
\Delta\, \partial_{t,h_1}\partial_{t,h_2}...\partial_{t,h_r} u(\cdot,t)= \partial_{t,h_1}\partial_{t,h_2}...\partial_{t,h_r}f(\cdot,t) &\text{auf }\mathbb{B}\\
\hphantom{\Delta} \partial_{t,h_1}\partial_{t,h_2}...\partial_{t,h_r} u(\cdot,t)=0 &\text{auf }\partial \mathbb{B}
\end{cases}
\end{align}

F\"ur einen beliebigen Multiindex $s\in\mathbb{N}^n$, mit $|s|=m\in\mathbb{N}\backslash\{0,1\}$, folgt dann aus \eqref{eq:7.a}:
\begin{align*}
&\int_{T_1}^{T_2}{\int_{B_{\varrho}(0)}{\left|\partial_{t,h_1}\partial_{t,h_2}...\partial_{t,h_r} \partial^s u(x,t) \right|^2 dx}\,dt}\\
\stackrel{\hphantom{\eqref{eq:5.22}}}{=}&\int_{T_1}^{T_2}{\int_{B_{\varrho}(0)}{\left|\partial^s \partial_{t,h_1}\partial_{t,h_2}...\partial_{t,h_r}  u(x,t) \right|^2 dx}\,dt}\\
\stackrel{\hphantom{\eqref{eq:5.22}}}{\leq}&\int_{T_1}^{T_2}{\int_{\mathbb{B}}{\left| \partial_{t,h_1}\partial_{t,h_2}...\partial_{t,h_r}  u(\cdot,t) \right|^2_{C^{m,\alpha}(\overline{\mathbb{B}})} dx}\,dt}\\
\stackrel{\hphantom{\eqref{eq:5.22}}}{=}&C(n)\cdot\int_{T_1}^{T_2}{\left| \partial_{t,h_1}\partial_{t,h_2}...\partial_{t,h_r}  u(\cdot,t) \right|^2_{C^{m,\alpha}(\overline{\mathbb{B}})}\, dt}\\
\stackrel{\hphantom{\eqref{eq:5.22}}}{\leq}&C(n)\cdot(T_2-T_1)\cdot\max_{t\in [T_1,T_2]}\left| \partial_{t,h_1}\partial_{t,h_2}...\partial_{t,h_r}  u(\cdot,t) \right|^2_{C^{m,\alpha}(\overline{\mathbb{B}})}\\
\stackrel{\eqref{eq:5.22}}{\leq} &C(n,m,\alpha,T)\cdot \max_{t\in [T_1,T_2]}\left| \partial_{t,h_1}\partial_{t,h_2}...\partial_{t,h_r}  f(\cdot,t) \right|^2_{C^{m-2,\alpha}(\overline{\mathbb{B}})}
\end{align*}

Mit Lemma 7.23 und Lemma 7.24 in \cite{gilbarg2001elliptic} folgt aus der Glattheit von $f$ die Existenz einer schwachen Ableitung von $u$ beliebiger Ordnung. Mit dem Einbettungssatz von Sobolev \cite[Theorem 7.10]{gilbarg2001elliptic} folgt die Behauptung.

\end{bew}

\begin{lem}\label{lem:7.4}
Es sei $a\in C_0^{\infty}(\mathbb{B})$ dann gilt f\"ur $v_1,v_2\in C^{\infty}(\overline{{\mathbb{B}}}\times [0,T],\mathbb{R}^q)$:
\begin{align}\label{eq:7.9}
\begin{split}
&\max_{t\in[0,T]}\left| M_{ij}[a](v_1(\cdot,t) )-M_{ij}[a](v_2(\cdot,t)) \right|_{C^{0,\alpha}(\overline{\mathbb{B}})} \\
\leq &K(n,\alpha,a)\cdot \left(\max_{t\in[0,T]}\left|v_1(\cdot,t) \right|_{C^{2,\alpha}(\overline{\mathbb{B}},\mathbb{R}^q)}+\max_{t\in[0,T]}\left|v_2(\cdot,t) \right|_{C^{2,\alpha}(\overline{\mathbb{B}},\mathbb{R}^q)} \right)\\
&\cdot\max_{t\in[0,T]}\left|v_1(\cdot,t)-v_2(\cdot,t) \right|_{C^{2,\alpha}(\overline{\mathbb{B}},\mathbb{R}^q)}
\end{split}
\end{align}

F\"ur $v\in C^{\infty}(\overline{{\mathbb{B}}}\times [0,T],\mathbb{R}^q)$ und $m,r\in\mathbb{N}$ gilt:
\begin{align}\label{eq:7.10}
\begin{split}
&\max_{t\in[0,T]}\left|\partial_t^r M_{ij}[a](v(\cdot,t)) \right|_{C^{m,\alpha}(\overline{\mathbb{B}})}\\
\leq &K(n,\alpha,a)\cdot\sum_{s=0}^r{\binom{r}{s}\max_{t\in[0,T]}\left|\partial_t^s  v(\cdot,t) \right|_{C^{m+2,\alpha}(\overline{\mathbb{B}},\mathbb{R}^q)}\max_{t\in[0,T]}\left|\partial_t^{r-s} v(\cdot,t) \right|_{C^{2,\alpha}(\overline{\mathbb{B}},\mathbb{R}^q)}}\\
&+C(n,m,\alpha,a)\cdot\sum_{s=0}^r{\binom{r}{s}\max_{t\in[0,T]}\left|\partial_t^s v(\cdot,t)\right|_{C^{m+1,\alpha}(\overline{\mathbb{B}},\mathbb{R}^q)}\max_{t\in[0,T]}\left|\partial_t^{r-s} v(\cdot,t) \right|_{C^{m+1,\alpha}(\overline{\mathbb{B}},\mathbb{R}^q)}}
\end{split}
\end{align}

\begin{bew}
Um \eqref{eq:7.9} zu zeigen, kann, analog zu Lemma \ref{lem:7.2}, die Ungleichung \eqref{eq:5.33} in Lemma \ref{lem:5.c} auf feste Zeiten angewandt werden. Unter Beachtung der Definition von $M_{ij}[a]$ in \eqref{eq:5.11}, als Summe der Operatoren $L_{ij}[a]$ und $R_{ij}[a]$ aus Lemma \ref{lem:5.1} und Lemma \ref{lem:5.2}, ergibt sich die zu zeigende Absch\"atzung  \eqref{eq:7.10} einerseits aus den Absch\"atzungen \eqref{eq:5.22} und \eqref{eq:5.30} f\"ur den $L_{ij}[a]$-Anteil und andererseits analog zum Beweis von \eqref{eq:7.7} f\"ur den $R_{ij}[a]$-Anteil.
\end{bew}

Aus Lemma \ref{lem:7.1} und Lemma \ref{lem:7.3} folgt, mit den Definitionen in \eqref{eq:5.14} und \eqref{eq:5.15}, analog zum Beweis von Lemma \ref{lem:5.d}:

\end{lem}

\begin{lem}\label{lem:7.5}
Es sei $a\in C_0^{\infty}(\mathbb{B})$ dann gilt f\"ur $v_1,v_2\in C^{\infty}(\overline{{\mathbb{B}}}\times [0,T],\mathbb{R}^q)$, dann gilt:
\begin{align}\label{eq:7.11}
\begin{split}
&\max_{t\in[0,T]}\biggl(\left| P[a](v_1(\cdot,t) )-P[a](v_2(\cdot,t)) \right|_{C^{2,\alpha}(\overline{\mathbb{B}},\mathbb{R}^n)}\\
&\hspace{3cm}+\left| Q[a](v_1(\cdot,t) )-Q[a](v_2(\cdot,t)) \right|_{C^{2,\alpha}(\overline{\mathbb{B}},\mathbb{R}^{\frac{n}{2}(n+1)})} \biggr)\\
\leq &K(n,\alpha,a)\cdot \left(\max_{t\in[0,T]}\left|v_1(\cdot,t) \right|_{C^{2,\alpha}(\overline{\mathbb{B}},\mathbb{R}^q)}+\max_{t\in[0,T]}\left|v_2(\cdot,t) \right|_{C^{2,\alpha}(\overline{\mathbb{B}},\mathbb{R}^q)} \right)\\
&\cdot\max_{t\in[0,T]}\left|v_1(\cdot,t)-v_2(\cdot,t) \right|_{C^{2,\alpha}(\overline{\mathbb{B}},\mathbb{R}^q)}
\end{split}
\end{align}

F\"ur $v\in C^{\infty}(\overline{{\mathbb{B}}}\times [0,T],\mathbb{R}^q)$ und $m,r\in\mathbb{N}$ gilt:
\begin{align}\label{eq:7.12}
\begin{split}
&\max_{t\in[0,T]}\left(\left|\partial_t^r P[a](v(\cdot,t)) \right|_{C^{m,\alpha}(\overline{\mathbb{B}},\mathbb{R}^n)}+\left|\partial_t^r Q[a](v(\cdot,t)) \right|_{C^{m,\alpha}(\overline{\mathbb{B}},\mathbb{R}^{\frac{n}{2}(n+1)})}\right)\\
\leq &K(n,\alpha,a)\cdot\sum_{s=0}^r{\binom{r}{s}\max_{t\in[0,T]}\left|\partial_t^s  v(\cdot,t) \right|_{C^{m,\alpha}(\overline{\mathbb{B}},\mathbb{R}^q)}\max_{t\in[0,T]}\left|\partial_t^{r-s} v(\cdot,t) \right|_{C^{2,\alpha}(\overline{\mathbb{B}},\mathbb{R}^q)}}\\
&+C(n,m,\alpha,a)\cdot\sum_{s=0}^r{\binom{r}{s}\max_{t\in[0,T]}\left|\partial_t^s v(\cdot,t)\right|_{C^{m-1,\alpha}(\overline{\mathbb{B}},\mathbb{R}^q)}\max_{t\in[0,T]}\left|\partial_t^{r-s} v(\cdot,t) \right|_{C^{m-1,\alpha}(\overline{\mathbb{B}},\mathbb{R}^q)}}
\end{split}
\end{align}

\end{lem}

Die folgenden beiden Lemmata k\"onnen genau wie Lemma \ref{lem:5.a} und \ref{lem:5.e} gezeigt werden.

\begin{lem}\label{lem:7.6}
Es sei $F_0\in C^{\infty}(\overline{\mathbb{B}},\mathbb{R}^q)$ eine beliebige freie Abbildung und $r\in\mathbb{N}$. Dann gilt, f\"ur beliebige $h\in C^{\infty}(\overline{\mathbb{B}}\times[0,T],\mathbb{R}^n)$, $f\in C^{\infty}(\overline{\mathbb{B}}\times[0,T],\mathbb{R}^{\frac{n}{2}(n+1)})$, f\"ur den Operator $E[F_0]$  die Absch\"atzung:
\begin{align}\label{eq:7.13}
\begin{split}
&\max_{t\in[0,T]}\left|\partial_t^r E[F_0](h(\cdot,t),f(\cdot,t)) \right|_{C^{m,\alpha}(\overline{\mathbb{B}},\mathbb{R}^q)}\\
\leq&C(n,q,m,\alpha,F_0)\cdot\max_{t\in[0,T]}\left(\left|\partial_t^r  h(\cdot,t) \right|_{C^{m,\alpha}(\overline{\mathbb{B}},\mathbb{R}^n)}+\left|\partial_t^r  f(\cdot,t) \right|_{C^{m,\alpha}\left(\overline{\mathbb{B}},\mathbb{R}^{\frac{n}{2}(n+1)}\right)}\right)
\end{split}
\end{align}

\end{lem}

\begin{lem}\label{lem:7.7}
Es sei $F_0\in C^{\infty}(\overline{\mathbb{B}},\mathbb{R}^q)$ eine beliebige freie Abbildung und $r\in\mathbb{N}$. Dann gilt f\"ur den Operator $E[F_0]$ f\"ur $m\geq 3$ die Absch\"atzung:
\begin{align}\label{eq:7.21}
\begin{split}
&\max_{t\in[0,T]}\left|\partial_t^r E[F_0](h(\cdot,t),f(\cdot,t)) \right|_{C^{m,\alpha}(\overline{\mathbb{B}},\mathbb{R}^q)}\\
\leq &\left\Vert E[F_0] \right\Vert_{2,\alpha}\cdot \max_{t\in[0,T]}\left( \left|\partial_t^r h(\cdot,t)\right|_{C^{m,\alpha}(\overline{\mathbb{B}},\mathbb{R}^n)}+\left| \partial_t^r f(\cdot,t) \right|_{C^{m,\alpha}\left(\overline{\mathbb{B}},\mathbb{R}^{\frac{n}{2}(n+1)}\right)}\right)\\
&+C(n,q,m,\alpha,F_0)\cdot\max_{t\in[0,T]}\left(\left|\partial_t^r h(\cdot,t) \right|_{C^{m-1,\alpha}(\overline{\mathbb{B}},\mathbb{R}^n)}+\left|\partial_t^r f(\cdot,t)\right|_{C^{m-1,\alpha}\left(\overline{\mathbb{B}},\mathbb{R}^{\frac{n}{2}(n+1)}\right)}\right)
\end{split}
\end{align}

f\"ur beliebige $h\in C^{\infty}(\overline{\mathbb{B}}\times[0,T],\mathbb{R}^n)$, $f\in C^{\infty}(\overline{\mathbb{B}}\times[0,T],\mathbb{R}^{\frac{n}{2}(n+1)})$.

\end{lem}

\begin{satz}\label{satz:7.2}
In der Situation von Satz \ref{satz:7.1} existiert ein $\widehat{T}_x\in (0,T_x]$, so dass f\"ur die, in diesem Satz konstruierte, L\"osung $(F(\cdot,t))_{t\in [0,\widehat{T}_x]}\subseteq C^{\infty}(\overline{U}_x,\mathbb{R}^q)$ die Eigenschaft:
\begin{equation*}
F\in C^{\infty}(\overline{U}_x\times[0,\widehat{T}_x],\mathbb{R}^q)
\end{equation*}

gilt.

\end{satz}

\begin{bew}
Mit Lemma \ref{lem:7.3} ist die, in \eqref{eq:7.15} definierte, Folgen von Approximationen $(v_k)_{k\in\mathbb{N}}$ glatt, das hei\ss{}t $(v_k)_{k\in\mathbb{N}}\subseteq C^{\infty}(\overline{\mathbb{B}}\times[0,T_x],\mathbb{R}^q)$. Hierbei sei erneut, unter Beachtung der Glattheit der Metrik $g$, auf die Definition des Operators $\Phi[F_0,a,\widehat{g}]$ in \eqref{eq:5.18}, und die Definition der Operatoren $P$ und $Q$ in \eqref{eq:5.14} und \eqref{eq:5.15} verwiesen. Um zu zeigen, dass  $v\in C^{\infty}(\overline{\mathbb{B}}\times[0,\widehat{T}_x],\mathbb{R}^q)$ gilt, wird gezeigt, dass nach eventueller Verkleinerung von $\widehat{T}_x\in (0,T_x]$ f\"ur alle $r\in\mathbb{N}$ und $s\in\mathbb{N}^n$ eine Konstante $C(r,s)\in\mathbb{R}_{>0}$ existiert, so dass f\"ur alle $k\in\mathbb{N}$ die Absch\"atzung:
\begin{align}\label{eq:7.26}
\left\Vert \partial_t^r \partial^s v_k \right\Vert_{C^0(\overline{\mathbb{B}}\times [0,\widehat{T}_x],\mathbb{R}^q)}\leq C(r,s)
\end{align}

gilt. Die Konstante $C(r,s)$ soll dabei insbesondere nicht von $k$ abh\"angen. Es sei $r\in\mathbb{N}$ dann folgt aus \eqref{eq:7.15}, \eqref{eq:5.18}, \eqref{eq:7.13}, \eqref{eq:7.12} und \eqref{eq:7.18} f\"ur alle $k\in\mathbb{N}\backslash\{0\}$ die Absch\"atzung:
\begin{align}\label{eq:7.20}
\notag& \max_{t\in[0,\widehat{T}_x]}\left|\partial_t^r v_k(\cdot,t) \right|_{C^{2,\alpha}(\overline{\mathbb{B}},\mathbb{R}^q)}\\
\notag\leq &2\cdot K_1(n,\alpha,a)\cdot   \left\Vert E[F_0] \right\Vert_{2,\alpha}\cdot\max_{t\in[0,\widehat{T}_x]}\left|E[F_0]\left(0,\widehat{g}(\cdot,t)\right)\right|_{C^{2,\alpha}(\overline{\mathbb{B}},\mathbb{R}^q)}\max_{t\in[0,\widehat{T}_x]}\left|\partial_t^r v_{k-1}(\cdot,t) \right|_{C^{2,\alpha}(\overline{\mathbb{B}},\mathbb{R}^q)}\\
&+K_1(n,\alpha,a)\cdot   \left\Vert E[F_0] \right\Vert_{2,\alpha}\\
\notag&\hspace{2cm}\cdot\sum_{s=1}^{r-1}{\binom{r}{s}\max_{t\in[0,\widehat{T}_x]}\left|\partial_t^s  v_{k-1}(\cdot,t) \right|_{C^{2,\alpha}(\overline{\mathbb{B}},\mathbb{R}^q)}\max_{t\in[0,\widehat{T}_x]}\left|\partial_t^{r-s} v_{k-1}(\cdot,t) \right|_{C^{2,\alpha}(\overline{\mathbb{B}},\mathbb{R}^q)}}\\
\notag&+\frac{1}{2}\max_{t\in[0,\widehat{T}_x]}\left|E[F_0]\left(0,\partial_t^r \widehat{g}(\cdot,t) \right)\right|_{C^{2,\alpha}(\overline{\mathbb{B}},\mathbb{R}^q)}
\end{align}
 
Unter der Annahme, dass $\widehat{T}_x\in(0,T_x]$ so klein gew\"ahlt ist, dass:
\begin{equation*}
\left\Vert E[F_0] \right\Vert_{2,\alpha}\cdot\max_{t\in[0,\widehat{T}_x]}\left|E[F_0]\left(0,\widehat{g}(\cdot,t)\right)\right|_{C^{2,\alpha}(\overline{\mathbb{B}},\mathbb{R}^q)}\leq \frac{1}{4K_1(n,\alpha,a)}
\end{equation*}

gilt, folgt aus \eqref{eq:7.20} die Absch\"atzung:
\begin{align*}
&\max_{t\in[0,\widehat{T}_x]}\left|\partial_t^r v_k(\cdot,t) \right|_{C^{2,\alpha}(\overline{\mathbb{B}},\mathbb{R}^q)}\\
\leq &\frac{1}{2}\left(\max_{t\in[0,\widehat{T}_x]}\left|\partial_t^r v_{k-1}(\cdot,t) \right|_{C^{2,\alpha}(\overline{\mathbb{B}},\mathbb{R}^q)}+\max_{t\in[0,\widehat{T}_x]}\left|E[F_0]\left(0,\partial_t^r \widehat{g}(\cdot,t) \right)\right|_{C^{2,\alpha}(\overline{\mathbb{B}},\mathbb{R}^q)} \right)\\
&+K_1(n,\alpha,a)\cdot   \left\Vert E[F_0] \right\Vert_{2,\alpha}\\
&\hspace{2cm}\cdot\sum_{s=1}^{r-1}{\binom{r}{s}\max_{t\in[0,\widehat{T}_x]}\left|\partial_t^s  v_{k-1}(\cdot,t) \right|_{C^{2,\alpha}(\overline{\mathbb{B}},\mathbb{R}^q)}\max_{t\in[0,\widehat{T}_x]}\left|\partial_t^{r-s} v_{k-1}(\cdot,t) \right|_{C^{2,\alpha}(\overline{\mathbb{B}},\mathbb{R}^q)}}
\end{align*}

und es folgt per Induktion \"uber $r\in\mathbb{N}$ mit Lemma \ref{lem:7.8} die Absch\"atzung \eqref{eq:7.26} f\"ur alle $|s|\leq 2$ und $r\in\mathbb{N}$. Insbesondere existiert mit dem Satz von Arzel\`a-Ascoli eine Teilfolge von $\left(v_k \right)_{k\in\mathbb{N}}$, die auf $\overline{\mathbb{B}}\times [0,T_x]$ gleichm\"a\ss{}ig gegen ein $\widehat{v}\in C^{0}(\overline{\mathbb{B}}\times [0,T],\mathbb{R})$ konvergiert. Da f\"ur feste Zeiten $t\in[0,T_x]$ die Folge $\left(v_k(\cdot,t) \right)_{k\in\mathbb{N}}$ im $C^{2,\alpha}$-Sinne gegen $v(\cdot,t)$ konvergiert, muss aber $\widehat{v}\equiv v$ gelten.\\
Nun sei $m\geq 3$, unter der Annahme, dass die Absch\"atzung \eqref{eq:7.26} bereits f\"ur $|s|\leq m-1$ und alle $r\in\mathbb{N}$ gilt. Dann folgt aus \eqref{eq:7.15}, \eqref{eq:5.18}, \eqref{eq:7.21} und \eqref{eq:7.12} f\"ur alle $k\in\mathbb{N}\backslash\{0\}$ und $r\in\mathbb{N}$:
\begin{align}\label{eq:7.23}
\begin{split}
&\max_{t\in[0,\widehat{T}_x]}\left|\partial_t^r  v_k(\cdot,t) \right|_{C^{m,\alpha}(\overline{\mathbb{B}},\mathbb{R}^q)}\\
\leq &K_2(n,\alpha,a) \cdot\left\Vert E[F_0] \right\Vert_{2,\alpha}\max_{t\in[0,\widehat{T}_x]}\left|\partial_t^r v_{k-1}(\cdot,t) \right|_{C^{m,\alpha}(\overline{\mathbb{B}},\mathbb{R}^q)}\max_{t\in[0,\widehat{T}_x]}\left| v_{k-1}(\cdot,t) \right|_{C^{2,\alpha}(\overline{\mathbb{B}},\mathbb{R}^q)}\\
+&K_2(n,\alpha,a) \cdot\left\Vert E[F_0] \right\Vert_{2,\alpha}\cdot\sum_{s=0}^{r-1}{\max_{t\in[0,\widehat{T}_x]}\left|\partial_t^s v_{k-1}(\cdot,t) \right|_{C^{m,\alpha}(\overline{\mathbb{B}},\mathbb{R}^q)}\max_{t\in[0,\widehat{T}_x]}\left|\partial_t^{r-s} v_{k-1}(\cdot,t) \right|_{C^{2,\alpha}(\overline{\mathbb{B}},\mathbb{R}^q)}}\\
+&C_1(n,m,\alpha,F_0)\cdot \sum_{s=0}^r{\max_{t\in[0,\widehat{T}_x]}\left|\partial^s v_{k-1}(\cdot,t) \right|_{C^{m-1,\alpha}(\overline{\mathbb{B}},\mathbb{R}^q)}\max_{t\in[0,\widehat{T}_x]}\left|\partial^{r-s} v_{k-1}(\cdot,t) \right|_{C^{m-1,\alpha}(\overline{\mathbb{B}},\mathbb{R}^q)}}\\
&\hspace{6cm}+\frac{1}{2}\max_{t\in[0,\widehat{T}_x]}\left|E[F_0]\left(0,\partial_t^r\widehat{g}(\cdot,t) \right)\right|_{C^{m,\alpha}(\overline{\mathbb{B}},\mathbb{R}^q)}
\end{split}
\end{align}

Mit \eqref{eq:7.18} und \eqref{eq:7.26} folgt aus \eqref{eq:7.23}, zusammen mit der Induktionsvoraussetzung, die Absch\"atzung:
\begin{align}\label{eq:7.25}
\notag&\max_{t\in[0,\widehat{T}_x]}\left|\partial_t^r  v_k(\cdot,t) \right|_{C^{m,\alpha}(\overline{\mathbb{B}},\mathbb{R}^q)}\\
\notag\leq &K_2(n,\alpha,a) \cdot\left\Vert E[F_0] \right\Vert_{2,\alpha}\max_{t\in[0,T_x]}\left|E[F_0]\left(0,\widehat{g}(\cdot,t)\right)\right|_{C^{2,\alpha}(\overline{\mathbb{B}},\mathbb{R}^q)}\max_{t\in[0,\widehat{T}_x]}\left|\partial_t^r v_{k-1}(\cdot,t) \right|_{C^{m,\alpha}(\overline{\mathbb{B}},\mathbb{R}^q)}\\
&+C_2\left(n,\alpha,a,r,m,F_0,\widehat{T}_x,\widehat{g}\right)\cdot\sum_{s=0}^{r-1}{\max_{t\in[0,\widehat{T}_x]}\left|\partial_t^s v_{k-1}(\cdot,t) \right|_{C^{m,\alpha}(\overline{\mathbb{B}},\mathbb{R}^q)}}\\
\notag\notag&+C_3\left(n,\alpha,a,r,m,F_0,\widehat{T}_x,\widehat{g}\right)+\frac{1}{2}\max_{t\in[0,\widehat{T}_x]}\left|E[F_0]\left(0,\partial_t^r\widehat{g}(\cdot,t) \right)\right|_{C^{m,\alpha}(\overline{\mathbb{B}},\mathbb{R}^q)}
\end{align}

Falls nun $\widehat{T}_x\in (0,T_x]$ so klein gew\"ahlt ist, dass f\"ur alle $t\in [0,\widehat{T}_x]$ die Absch\"atzung:
\begin{equation*}
\left\Vert E[F_0] \right\Vert_{2,\alpha}\max_{t\in[0,T_x]}\left|E[F_0]\left(0,\widehat{g}(\cdot,t)\right)\right|_{C^{2,\alpha}(\overline{\mathbb{B}},\mathbb{R}^q)}\leq \frac{1}{2K_2(n,\alpha,a)}
\end{equation*}

erf\"ullt ist, dann folgt aus \eqref{eq:7.25}, per Induktion \"uber $r\in\mathbb{N}$, die Absch\"atzung \eqref{eq:7.26} f\"ur $|s|=m$ und alle $r\in\mathbb{N}$.

\end{bew}


\section{Konstruktion der globalen Einbettung}
\label{sec:7.3}

Die lokalen Resultate aus den ersten beiden Abschnitten dieses Kapitels werden nun zu einem globalen Satz ausgebaut:

\begin{satz}
Gegeben sei eine kompakte $n$-dimensionale glatte Mannigfaltigkeit $M$, ferner sei eine Familie von Riemannschen Metriken $(g(\cdot,t))_{t\in[0,T]}\subseteq \mathcal{T}^2(M)$ gegeben, so dass f\"ur jedes $(U,\varphi)\in\mathcal{A}$ die Bedingung $(\, ^{\varphi}g_{ij}(\cdot,\cdot))_{1\leq i\leq j\leq n}\in C^{\infty}(\varphi(U)\times[0,T],\mathbb{R}^{\frac{n}{2}(n+1)})$ gilt. Dann existiert, f\"ur $q(n):=\max\left\{\frac{n}{2}(n+5),\frac{n}{2}(n+3)+5 \right\}$ ein $\widehat{T}\in (0,T]$, und eine Familie von freien Einbettungen $(F(\cdot,t))_{t\in[0,\widehat{T}]}\subseteq C^{\infty}(M,\mathbb{R}^q)$, mit $F\in C^{\infty}(M\times[0,T],\mathbb{R}^q)$, so dass f\"ur alle $t\in [0,\widehat{T}]$ die Gleichung:
\begin{equation*}
F(\cdot,t)^{\ast}(g^{can})=g(\cdot,t)\hspace{0.5cm}
\end{equation*}

auf $M$ gilt.
\end{satz}

\begin{bew}
Es seien $\left\{(U_i,\varphi_i)\right\}_{1\leq i\leq m}\subseteq \mathcal{A}$ Karten, mit $\varphi_i(U_i)=B_{1+\tau}(0)$ f\"ur alle $i\in\{1,...,m\}$ f\"ur ein festes $\tau\in\mathbb{R}_{>0}$, so dass $M=\bigcup_{i=1}^m{\varphi_i^{-1}(\mathbb{B})}$ gilt. Ferner seien $\left\{\psi_i\right\}_{1\leq i\leq m}\subseteq C^{\infty}(M)$ mit $supp(\psi_i)\subseteq \varphi_i^{-1}(\mathbb{B})$, f\"ur alle $i\in\{1,...,m\}$ und $\sum_{i=1}^m{\psi_i}\equiv 1$, wie in \cite[Theorem 2.25]{lee2003introduction} gew\"ahlt. Dann gilt f\"ur jedes $(x,t)\in M\times [0,T]$ die Gleichung:
\begin{equation}\label{eq:7.41}
g(x,t)=g(x,0)+g(x,t)-g(x,0)=g(x,0)+\sum_{i=1}^m{\underbrace{\psi_i [g(x,t)-g(x,0)]}_{=: \widehat{g}^{(i)}(x,t)}}
\end{equation}

Gem\"a\ss{} Hauptsatz \ref{haupt2} wird zun\"achst, genau wie in Satz \ref{satz:7.1}, eine freie isometrische Einbettung $F_0\in C^{\infty}(M,\mathbb{R}^q)$ der Riemannschen Mannigfaltigkeit $(M,g(\cdot,0))$ gew\"ahlt. Mit der Methode aus Satz \ref{satz:7.1} kann, unter Beachtung der Regularit\"atsbetrachtungen aus \autoref{sec:7.2}, f\"ur ein $T_1\in (0,T]$ eine Familie von freien Einbettungen $F_1\in C^{\infty}(M\times[0,T_1])$ konstruiert werden, so dass f\"ur jedes $t\in [0,T_1]$ auf ganz $M$ die Gleichung:
\begin{equation}\label{eq:7.27}
F_1(\cdot,t)^{\ast}(g^{can})=g(\cdot,0)+\widehat{g}^{(1)}(\cdot,t)
\end{equation}

erf\"ullt ist. Diese Abbildung wird, wie in \eqref{eq:7.26}, \"uber die Bildungsvorschrift:
\begin{equation}
F_1(x,t):=
\begin{cases}
F_0(x)+u^{(1)}(\varphi_1(x),t) &\text{falls }x\in U_1 \\
F_0(x) &\text{sonst }
\end{cases} 
\end{equation}

wobei $u^{(1)}\in C^{\infty}(\overline{\mathbb{B}}\times [0,T_1],\mathbb{R}^q)$ eine Abbildung mit $supp(u^{(1)}(\cdot,t))\subseteq \mathbb{B}$ f\"ur alle $t\in [0,T_1]$ ist, bestimmt. Im Folgenden wird, f\"ur feste $t\in [0,T_1]$, die Abbildung $\left.\,^{\varphi_2}F_1(\cdot,t) \right|_{\mathbb{B}}\in C^{\infty}(\overline{\mathbb{B}},\mathbb{R}^q)$ ebenfalls mit $F_1(\cdot,t)$ bezeichnet. Es sei $a_2\in C^{\infty}_0(\mathbb{B})$ eine Abbildung, so dass $\left. a_2 \right|_{supp(\psi_2\circ \varphi_2^{-1})}\equiv 1$ erf\"ullt ist, dann kann $T_2\in (0,T_1]$ so gew\"ahlt werden, dass f\"ur alle $t\in[0,T_2]$ die Absch\"atzung:
\begin{equation*}
\max_{t\in[0,T_2]}\left\Vert E[F_1(\cdot,t)]\right\Vert_{2,\alpha}\max_{t\in[0,T_2]}\left|E[F_1(\cdot,t)](0,\,^{\varphi_2}\widehat{g}^{(2)}(\cdot,t)) \right|_{C^{2,\alpha}(\overline{\mathbb{B}},\mathbb{R}^q)}\leq \vartheta(n,\alpha,a_2)
\end{equation*}

mit dem in Satz \ref{satz:5.1} definierten $\vartheta\in\mathbb{R}_{>0}$, erf\"ullt ist. Dann existiert f\"ur jedes $t\in [0,T_2]$ ein $v^{(2)}(\cdot,t)\in C^{\infty}(\overline{\mathbb{B}},\mathbb{R}^q)$, so dass f\"ur $u^{(2)}(\cdot,t):= a_2^2\, v^{(2)}(\cdot,t)\in C^{\infty}_0(\mathbb{B},\mathbb{R}^q) $ f\"ur alle $x\in \mathbb{B}$, und f\"ur jedes $i,j\in\{1,...,n\}$ mit $i\leq j$, die Gleichung:
\begin{align}\label{eq:7.29}
\begin{split}
&\partial_i (F_1(x,t)+u^{(2)}(x,t))\cdot \partial_j (F_1(x,t)+u^{(2)}(x,t))\\
\stackrel{\hphantom{\eqref{eq:7.17}}}{=}&\partial_i F_1(x,t)\cdot \partial_j F_1(x,t)+\,^{\varphi_2}\widehat{g}^{(2)}_{ij}(x,t)\\
\stackrel{\eqref{eq:7.27}}{=}& ^{\varphi_2}g_{ij}(x,0)+\,^{\varphi_2}\widehat{g}^{(1)}_{ij}(x,t)+\,^{\varphi_2}\widehat{g}^{(2)}_{ij}(x,t)
\end{split}
\end{align}

erf\"ullt ist. F\"ur jedes $t\in [0,T_2]$ wird die Abbildung $v^{(2)}(\cdot,t)\in C^{\infty}(\overline{\mathbb{B}},\mathbb{R}^q)$, gem\"a\ss{} \eqref{eq:5.47}, durch die Folge: $(v^{(2)}_k(\cdot,t))_{k\in\mathbb{N}}\subseteq C^{\infty}(\overline{\mathbb{B}},\mathbb{R}^q)$ mit:
\begin{align}\label{eq:7.28}
v^{(2)}_k(\cdot,t)=
\begin{cases}
0 & \text{falls }  k=0\\
\Phi[F_1(\cdot,t),a_2,\,^{\varphi_2}\widehat{g}^{(2)}(\cdot,t)](v^{(2)}_{k-1}(\cdot,t))  &\text{falls } k\geq 1
\end{cases}
\end{align}

im $C^{2,\alpha}$-Sinne approximiert wird. Mit Lemma \ref{lem:7.3} gilt, unter Beachtung der Glattheit der Metrik $g$, der Defintion des Operators $\Phi[F_1(\cdot,t),a_2,\,^{\varphi_2}\widehat{g}^{(2)}(\cdot,t)]$ in \eqref{eq:5.18}, sowie der Definition der Operatoren $P$ und $Q$ in \eqref{eq:5.14} und \eqref{eq:5.15}, die Aussage $(v^{(2)}_k)_{k\in\mathbb{N}}\subseteq C^{\infty}(\overline{\mathbb{B}}\times[0,T_2],\mathbb{R}^q)$. Mit \eqref{eq:5.48} gilt die Absch\"atzung:
\begin{equation}\label{eq:7.30}
\max_{t\in[0,T_2]}\left|v^{(2)}_k(\cdot,t) \right|_{C^{2,\alpha}(\overline{\mathbb{B}},\mathbb{R}^q)}\leq \max_{t\in[0,T_2]}\left|E[F_1(\cdot,t)]\left(0,\,^{\varphi_2}\widehat{g}^{(2)}(\cdot,t)\right)\right|_{C^{2,\alpha}(\overline{\mathbb{B}},\mathbb{R}^q)}
\end{equation}

Definiere nun $(F_2(\cdot,t))_{t\in[0,T_2]}\subseteq C^{\infty}(M,\mathbb{R}^q)$ wie folgt:
\begin{equation}\label{eq:7.40}
F_2(x,t):=
\begin{cases}
F_1(x,t)+u^{(2)}(\varphi_2(x),t) &\text{falls }x\in U_2 \\
F_1(x,t) &\text{sonst}
\end{cases} 
\end{equation}

Dann ist mit \eqref{eq:7.29} f\"ur jedes $t\in [0,T_2]$ die Gleichung:
\begin{equation*}
F_2(\cdot,t)^{\ast}(g^{can})=g(\cdot,0)+\widehat{g}^{(1)}(\cdot,t)+\widehat{g}^{(2)}(\cdot,t)
\end{equation*}

erf\"ullt. Es wird gezeigt, dass, nach eventueller Verkleinerung von $T_2\in (0,T_1]$, die Aussage $F_2\in C^{\infty}(M\times [0,T_2],\mathbb{R}^q)$ gilt. Dazu wird wieder gezeigt, dass f\"ur jedes $r\in\mathbb{N}$ und $s\in\mathbb{N}^n$ eine Konstante $C(r,s)\in\mathbb{R}_{>0}$ existiert, so dass f\"ur alle $k\in\mathbb{N}$ die Absch\"atzung:
\begin{align}\label{eq:7.34}
\left\Vert \partial_t^r \partial^s v^{(2)}_k \right\Vert_{C^0(\overline{\mathbb{B}}\times [0,T_2],\mathbb{R}^q)}\leq C(r,s)
\end{align}

gilt, wobei die Konstante $C$ nicht von $k$ abh\"angt. Dabei ist nun die Zeitabh\"angigkeit des Operators $E[F_1(\cdot,t)]$ zu beachten. Im Folgenden ist f\"ur $t\in[0,T_2]$:
\begin{align*}
\partial_t^r E[F_1(\cdot,t)]: C^{\infty}(\overline{\mathbb{B}},\mathbb{R}^n) \times C^{\infty}(\overline{\mathbb{B}},\mathbb{R}^{\frac{n}{2}(n+1)})&\longrightarrow C^{\infty}(\overline{\mathbb{B}},\mathbb{R}^q)\\
\partial_t^r E[F_1(\cdot,t)](h,f)(x)&:= \partial_t^r\Theta[F_1(\cdot,t)](x)\cdot\begin{pmatrix} h(x)\\ f(x) \end{pmatrix}
\end{align*}

wobei $\Theta[F_1(\cdot,t)]$ in \eqref{eq:5.41} definiert worden ist. Die Absch\"atzung:
\begin{align}\label{eq:7.37}
\begin{split}
&\max_{t\in[0,T_2]}\left|\partial_t^r E[F_1(\cdot,t)](h,f)\right|_{C^{m,\alpha}(\overline{\mathbb{B}},\mathbb{R}^q)}\\
&\hspace{2cm}\leq  C(n,m,\alpha,r,F_1,T_2)\cdot\left(\left|h \right|_{C^{m,\alpha}(\overline{\mathbb{B}},\mathbb{R}^n)}+\left|f \right|_{C^{m,\alpha}\left(\overline{\mathbb{B}},\mathbb{R}^{\frac{n}{2}(n+1)}\right)}\right)
\end{split}
\end{align}

ist analog zu Lemma \ref{lem:5.a} zu beweisen. F\"ur die folgenden Absch\"atzungen sei erw\"ahnt, dass f\"ur $h\in C^{\infty}\left(\overline{\mathbb{B}}\times [0,T_2],\mathbb{R}^n\right)$ und $f\in C^{\infty}\left(\overline{\mathbb{B}}\times [0,T_2],\mathbb{R}^{\frac{n}{2}(n+1)}\right)$ f\"ur alle $t\in [0,T_2]$ die Gleichung:
\begin{equation}\label{eq:7.36}
\partial_t^r( E[F_1(\cdot,t)](h(\cdot,t),f(\cdot,t)))=\sum_{s=0}^r{\binom{r}{s}\partial_t^s E[F_1(\cdot,t)](\partial_t^{r-s} h(\cdot,t),\partial_t^{r-s} f(\cdot,t))}
\end{equation}

gilt. Dann ist f\"ur $r\in\mathbb{N}$ und $k\in\mathbb{N}\backslash\{0\}$, unter Beachtung von \eqref{eq:7.28}, \eqref{eq:5.18}, \eqref{eq:7.36}, \eqref{eq:7.37}, \eqref{eq:7.12} und \eqref{eq:7.30}:
\begin{align*}
\begin{split}
& \max_{t\in[0,T_2]}\left|\partial_t^r v^{(2)}_k(\cdot,t) \right|_{C^{2,\alpha}(\overline{\mathbb{B}},\mathbb{R}^q)}\\
\leq&2\cdot K_1(n,\alpha,a_2)\cdot\left\Vert  E[F_1(\cdot,t)]\right\Vert_{2,\alpha}\cdot\max_{t\in[0,T_2]}\left|E[F_1(\cdot,t)]\left(0,\,^{\varphi_2}\widehat{g}^{(2)}(\cdot,t)\right)\right|_{C^{2,\alpha}(\overline{\mathbb{B}},\mathbb{R}^q)}\\
&\hspace{2cm}\cdot\max_{t\in[0,T_2]}\left|\partial_t^{r} v^{(2)}_{k-1}(\cdot,t) \right|_{C^{2,\alpha}(\overline{\mathbb{B}},\mathbb{R}^q)}\\
&+K_1(n,\alpha,a_2)\cdot\max_{t\in[0,T_2]}\left\Vert  E[F_1(\cdot,t)] \right\Vert_{2,\alpha}\\
&\hspace{1cm}\cdot\sum_{s=1}^{r-1}{\binom{r}{s}\max_{t\in[0,T_2]}\left|\partial_t^{s}  v^{(2)}_{k-1}(\cdot,t) \right|_{C^{2,\alpha}(\overline{\mathbb{B}},\mathbb{R}^q)}\max_{t\in[0,T_2]}\left|\partial_t^{r-s} v^{(2)}_{k-1}(\cdot,t) \right|_{C^{2,\alpha}(\overline{\mathbb{B}},\mathbb{R}^q)}}\\
&+K_1(n,\alpha,a_2)\cdot \sum_{s=1}^r  \max_{t\in[0,T_2]} \left\Vert \partial_t^s E[F_1(\cdot,t)] \right\Vert_{2,\alpha}\\
&\hspace{1cm}\cdot\sum_{\widehat{s}=0}^{r-s}{\binom{r-s}{\widehat{s}}\max_{t\in[0,T_2]}\left|\partial_t^{\widehat{s}}  v^{(2)}_{k-1}(\cdot,t) \right|_{C^{2,\alpha}(\overline{\mathbb{B}},\mathbb{R}^q)}\max_{t\in[0,T_2]}\left|\partial_t^{r-s-\widehat{s}} v^{(2)}_{k-1}(\cdot,t) \right|_{C^{2,\alpha}(\overline{\mathbb{B}},\mathbb{R}^q)}}\\
&+C\left(n,\alpha,a_2,r,F_1,T_2,\,^{\varphi_2}\widehat{g}^{(2)}\right)
\end{split}
\end{align*}

Unter der Annahme, dass $T_2\in (0,T_1]$ so klein ist, dass:
\begin{equation*}
\max_{t\in[0,T_2]}\left\Vert  E[F_1(\cdot,t)]\right\Vert_{2,\alpha}\max_{t\in[0,T_2]}\left|E[F_1(\cdot,t)]\left(0,\,^{\varphi_2}\widehat{g}^{(2)}(\cdot,t)\right)\right|_{C^{2,\alpha}(\overline{\mathbb{B}},\mathbb{R}^q)}\leq \frac{1}{4 K_1(n,\alpha,a_2)}
\end{equation*}

gilt, folgt die Absch\"atzung mit Lemma \ref{lem:7.8}, per Induktion \"uber $r\in\mathbb{N}$, die Absch\"atzung \eqref{eq:7.34} f\"ur $|s|\leq 2$ und beliebige $r\in\mathbb{N}$. Nun sei $m\geq 3$, unter der Annahme, dass die Absch\"atzung \eqref{eq:7.34} bereits f\"ur $|s|\leq m-1$ und alle $r\in\mathbb{N}$ gezeigt ist. Dann gilt f\"ur $r\in\mathbb{N}$ und $k\in\mathbb{N}\backslash\{0\}$ unter Beachtung von \eqref{eq:7.28}, \eqref{eq:5.18}, \eqref{eq:7.36}, \eqref{eq:7.13}, \eqref{eq:7.12} und \eqref{eq:7.37}:
\begin{align*}
\begin{split}
& \max_{t\in[0,T_2]}\left|\partial_t^r v^{(2)}_k(\cdot,t) \right|_{C^{m,\alpha}(\overline{\mathbb{B}},\mathbb{R}^q)}\\
\leq&K_2(n,\alpha,a_2)\cdot\max_{t\in[0,T_2]}\left\Vert  E[F_1(\cdot,t)] \right\Vert_{2,\alpha}\max_{t\in[0,T]} \cdot\left|E[F_1(\cdot,t)]\left(0,\,^{\varphi_2}\widehat{g}^{(2)}(\cdot,t)\right)\right|_{C^{2,\alpha}(\overline{\mathbb{B}},\mathbb{R}^q)}\\
&\hspace{0.5cm}\cdot \max_{t\in[0,T_2]}\left|\partial_t^r  v^{(2)}_{k-1}(\cdot,t) \right|_{C^{m,\alpha}(\overline{\mathbb{B}},\mathbb{R}^q)}\\
&+C\left(n,m,\alpha,a_2,r,F_1,T_2,\,^{\varphi_2}\widehat{g}^{(2)}\right)\\
&\hspace{2cm}\cdot\sum_{s=0}^{r-1}{\max_{t\in[0,T_2]}\left|\partial_t^s  v^{(2)}_{k-1}(\cdot,t) \right|_{C^{m,\alpha}(\overline{\mathbb{B}},\mathbb{R}^q)}\max_{t\in[0,T_2]}\left|\partial_t^{r-s} v^{(2)}_{k-1}(\cdot,t) \right|_{C^{2,\alpha}(\overline{\mathbb{B}},\mathbb{R}^q)}}\\
&+C\left(n,m,\alpha,a_2,r,F_1,T_2,\,^{\varphi_2}\widehat{g}^{(2)}\right)\\
&\hspace{2cm}\cdot\sum_{s=0}^r{\max_{t\in[0,T_2]}\left|\partial_t^s v^{(2)}_{k-1}(\cdot,t)\right|_{C^{m-1,\alpha}(\overline{\mathbb{B}},\mathbb{R}^q)}\max_{t\in[0,T_2]}\left|\partial_t^{r-s} v^{(2)}_{k-1}(\cdot,t) \right|_{C^{m-1,\alpha}(\overline{\mathbb{B}},\mathbb{R}^q)}}\\
&+C\left(n,m,\alpha,a_2,r,F_1,T_2,\,^{\varphi_2}\widehat{g}^{(2)}\right)
\end{split}
\end{align*}

Sei nun $T_2 \in (0,T_1]$  klein genug, so dass:
\begin{align*}
\max_{t\in[0,T_2]}\left\Vert  E[F_1(\cdot,t)] \right\Vert_{2,\alpha}\max_{t\in[0,T]} \left|E[F_1(\cdot,t)]\left(0,\,^{\varphi_2}\widehat{g}^{(2)}(\cdot,t)\right)\right|_{C^{2,\alpha}(\overline{\mathbb{B}},\mathbb{R}^q)}\leq \frac{1}{2 K_3(n,\alpha,a_2)}
\end{align*}

gilt, dann folgt mit Lemma \ref{lem:7.8}, zusammen mit der Induktionsvoraussetzung, die Ab\-sch\"atzung \eqref{eq:7.34} f\"ur $|s|=m$ und alle $r\in \mathbb{N}$. Ist $F_2(\cdot,t)$ f\"ur ein $t\in[0,T_2]$ keine freie Einbettung, so kann dies aber mit Satz \ref{satz:4.4}, unter Beachtung der Absch\"atzung \eqref{eq:5.53}, durch eine Verkleinerung von $T_2$ vermieden werden. Dann kann die soeben beschriebene Konstruktion \eqref{eq:7.40} iterativ auf die anderen Karten $\left\{(U_i,\varphi_i)\right\}_{3\leq i\leq m}\subseteq \mathcal{A}$ angewandt werden, und es ergibt sich nach endlich vielen Schritten, unter Beachtung der Gleichung \eqref{eq:7.41}, die Behauptung.

\end{bew}

\setcounter{secnumdepth}{1}

\appendix

\addcontentsline{toc}{chapter}{Anhang} 

\chapter{Funktionenr\"aume}
\thispagestyle{fancy}

In diesem Kapitel werden die, in der Arbeit verwendeten, Funktionenr\"aume definiert und einige Absch\"atzungen bewiesen, welche vowiegend in \autoref{chap:Kap5} und \autoref{chap:Kap7} verwendet werden.

\begin{defi}\label{defi:A.1a}
Es sei $\Omega\subseteq\mathbb{R}^n$ eine offene Menge und $m\in\mathbb{N}$. Dann wird der folgende Vektorraum definiert:
\begin{align*}
\gls{C^{m}(Omega,Rq)}&:=\{u:\Omega\longrightarrow \mathbb{R}^q: \text{u ist m-mal stetig differenzierbar in }\Omega \}
\end{align*}

Wenn $\Omega$ zus\"atzlich beschr\"ankt ist, so wird ein weiterer Vektorraum definiert:
\begin{align*}
\gls{C^{m}(overline{Omega},Rq)}:=\{&u:\Omega\longrightarrow \mathbb{R}^q: \text{u ist m-mal stetig differenzierbar in }\Omega,\\
&\text{und f\"ur alle }|s|\leq m \text{ ist }  \partial^{s}u \text{ auf }\overline{\Omega}\text{ stetig fortsetzbar} \}
\end{align*}

Auf diesem Vektorraum wird die folgende Norm definiert:
\begin{align*}
\gls{NormaufCm}: C^{m}(\overline{\Omega},\mathbb{R}^q)&\longrightarrow \mathbb{R}\\
u&\mapsto\sum_{|s|\leq m}{\left\Vert\partial^{s} u\right\Vert}_{C^0(\overline{\Omega},\mathbb{R}^q)}
\end{align*}

hierbei ist $\left\Vert u\right\Vert_{C^{0}(\overline{\Omega},\mathbb{R}^q)}\:=\sup_{x\in\overline{\Omega}}{|u(x)|_{\mathbb{R}^q}}$. Ferner wird noch der Funktionenraum aller \textbf{glatten Funktionen auf einer offenen Menge} $\bm{\Omega}$ definiert:
\begin{equation*}
\gls{RaumGlatt}:= \bigcap_{m\in\mathbb{N}}{C^{m}(\Omega,\mathbb{R}^q)}
\end{equation*}

und falls $\Omega$ beschr\"ankt ist:
\begin{equation*}
\gls{RaumGlattFort}:= \bigcap_{m\in\mathbb{N}}{C^{m}(\overline{\Omega},\mathbb{R}^q)}
\end{equation*}

\end{defi}

Mit der Norm $\left\Vert\cdot\right\Vert_{C^{m}(\overline{\Omega},\mathbb{R}^q)}$ ist der Vektorraum $C^{m}(\overline{\Omega},\mathbb{R}^q)$ vollst\"andig \cite[1.6]{alt2012lineare}. Im Fall $q=1$  wird $\mathbb{R}^q$ in der Notation von $C^{m}(\overline{\Omega},\mathbb{R}^q)$ weggelassen.

\begin{defi}
Sei $S\subseteq\mathbb{R}^n$ und $\alpha\in (0,1)$. F\"ur eine Abbildung $u:S\longrightarrow\mathbb{R}$ hei\ss{}t:
\begin{equation*}
\gls{Holderkonstante}:=\sup\left\{\frac{|u(x)-u(y)|}{|x-y|^{\alpha}}: x,y\in S, x\neq y \right\}\in[0,\infty]
\end{equation*}

\textbf{H\"olderkonstante von} ${\bm{u}}$  \textbf{auf} ${\bm{S}}$ \textbf{zum Exponenten} $\bm{\alpha}$\index{H\"olderkonstante}.

\end{defi}

F\"ur zwei Funktionen $u,v: S\longrightarrow\mathbb{R}$ gilt die folgende Absch\"atzung:
\begin{equation}\label{eq:A.1}
[uv]_{\alpha,S}\leq \sup_{x\in S}{|u(x)|}\cdot[v]_{\alpha,S}+[u]_{\alpha,S}\cdot\sup_{x\in S}{|v(x)|}
\end{equation}

Die im Folgenden definierten H\"olderr\"aume werden, der Einfachheit halber, auf $\mathbb{B}$ definiert. Dies gew\"ahrleistet, neben der Konvexit\"at, eine hohe Randregularit\"at, siehe auch \cite[6.2]{gilbarg2001elliptic}. Da die analytischen Betrachtungen von lokaler Art sind, erweisen sich diese Definitionen als ausreichend.

\begin{defi}
F\"ur $\alpha\in (0,1)$ wird mit:
\begin{equation*}
\gls{C0alpha}:=\left\{u\in C^{0}(\overline{\mathbb{B}}): [u]_{\alpha,\overline{\mathbb{B}}}<\infty \right\}
\end{equation*} 

der Vektorraum aller $\bm{\alpha}$\textbf{-h\"olderstetigen reellen Funktionen auf} $\bm{\mathbb{B}}$\index{Stetigkeit>H\"older-$\sim$} definiert. Des Weiteren wird auf diesem Raum die Norm:
\begin{align*}
\gls{NormaufC0alpha}: C^{0,\alpha}(\overline{\mathbb{B}})&\longrightarrow\mathbb{R}\\
u&\mapsto \left\Vert u\right\Vert_{C^0(\overline{\mathbb{B}})}+[u]_{\alpha,\overline{\mathbb{B}}}
\end{align*}

definiert.

\end{defi}

In \cite[1.7]{alt2012lineare} wird gezeigt, dass der Raum $C^{0,\alpha}(\overline{\mathbb{B}})$, versehen mit der Norm $\left\Vert\cdot\right\Vert_{C^{0,\alpha}(\overline{\mathbb{B}})}$, ein Banachraum ist.

F\"ur zwei Funktionen $u,v\in C^{0,\alpha}(\overline{\mathbb{B}})$ ergibt sich mit \eqref{eq:A.1} die Absch\"atzung:
\begin{equation}\label{eq:A.2}
\left\Vert uv\right\Vert_{C^{0,\alpha}(\overline{\mathbb{B}})}\leq \left\Vert u\right\Vert_{C^{0,\alpha}(\overline{\mathbb{B}})}\left\Vert v\right\Vert_{C^{0,\alpha}(\overline{\mathbb{B}})}
\end{equation}

\begin{defi}\label{defi:A.1}
Auf $\mathbb{B}$ werden, f\"ur $m\in\mathbb{N}$ und $\alpha\in(0,1)$, die folgenden Vektorr\"aume definiert:
\begin{align*}
\gls{CmalphaB}:=\{u\in C^{m,\alpha}(\overline{\mathbb{B}}): \left\Vert\partial^s u\right\Vert_{C^{0,\alpha}(\overline{\mathbb{B}})}<\infty\hspace{0.25cm}\text{f\"ur alle }s\in\mathbb{N}^n\text{ mit }|s|=m\}
\end{align*}

Auf dieser Menge werden die Normen $\left|\cdot\right|_{C^{m,\alpha}(\overline{\mathbb{B}})},\left\Vert\cdot\right\Vert_{C^{m,\alpha}(\overline{\mathbb{B}})}: C^{m,\alpha}(\overline{\mathbb{B}})\longrightarrow\mathbb{R}$ wie folgt definiert:
\begin{align*}
\gls{Norm1}&:=
\begin{cases}
\left\Vert u \right\Vert_{C^{0,\alpha}(\overline{\mathbb{B}})} &\text{f\"ur }m=0\\
\left\Vert u \right\Vert_{C^{0,\alpha}(\overline{\mathbb{B}})}+\sum_{|s|=m}{\left\Vert\partial^s u\right\Vert_{C^{0,\alpha}(\overline{\mathbb{B}})}} &\text{f\"ur }m\geq 1
\end{cases}\\
\gls{Norm2}& :=\sum_{j=0}^m{\max_{|s|=j}{\left\Vert \partial^{s} u\right\Vert}_{C^{0}(\overline{\mathbb{B}})}}+\max_{|s|=m}{[\partial^{s} u]}_{\alpha,\mathbb{B}}
\end{align*}

definiert.
\end{defi}

Die Norm $\left\Vert \cdot \right\Vert_{C^{m,\alpha}(\overline{\mathbb{B}})}$ wird in \cite{gilbarg2001elliptic} verwendet. Existenz- und Eindeutigkeitss\"atze der Poisson-Gleichung mit Dirichlet-Randbedingung, sowie diverse Absch\"atzungen werden aus diesem Buch \"ubernommen und an den entsprechenden Stellen direkt angegeben. Aus praktischen Gr\"unden wird in dieser Arbeit die Norm $\left|\cdot \right|_{C^{m,\alpha}(\overline{\mathbb{B}})}$ verwendet. Dazu wird gezeigt, dass diese beiden definierten Normen \"aquivalent sind. F\"ur $m\geq 1$ gilt f\"ur jedes $u\in C^{m,\alpha}(\overline{\mathbb{B}})$:
\begin{align}\label{eq:A.3}
\begin{split}
\left|u \right|_{C^{m,\alpha}(\overline{\mathbb{B}})}:=\left\Vert u \right\Vert_{C^{0,\alpha}(\overline{\mathbb{B}})}+\sum_{|s|=m}{\left\Vert\partial^s u\right\Vert_{C^{0,\alpha}(\overline{\mathbb{B}})}}\leq& C(n,m)\cdot\left(\left\Vert u \right\Vert_{C^{0,\alpha}(\overline{\mathbb{B}})}+\max_{|s|=m}{\left\Vert\partial^s u\right\Vert_{C^{0,\alpha}(\overline{\mathbb{B}})}} \right)\\
\leq&C(n,m)\cdot \left\Vert u \right\Vert_{C^{m,\alpha}(\overline{\mathbb{B}})}
\end{split}
\end{align}

andererseits folgt, f\"ur $k\in\mathbb{N}$ und $\beta\in (0,1)$ mit $k+\beta<m+\alpha$, aus den Absch\"atzungen \cite[(6.82), (6.89)]{gilbarg2001elliptic}, zusammen mit \cite[(4.17)', (4.17)'']{gilbarg2001elliptic}, mit der Randregularit\"at der Menge $\mathbb{B}$, die Absch\"atzung:
\begin{equation}\label{eq:A.4}
\left\Vert u \right\Vert_{C^{k,\beta}(\overline{\mathbb{B}})}\leq C(n,m,k,\alpha,\beta)\cdot\left(\left\Vert u \right\Vert_{C^{0}(\overline{\mathbb{B}})}+\max_{|s|=m}{[\partial^{s}u]_{\alpha,\mathbb{B}}} \right)\leq C(n,m,k,\alpha,\beta)\cdot\left\Vert u \right\Vert_{C^{m,\alpha}(\overline{\mathbb{B}})}
\end{equation}

Aus der ersten Ungleichung in \eqref{eq:A.4} folgt:
\begin{align}\label{eq:A.5}
\begin{split}
\left\Vert u \right\Vert_{C^{m,\alpha}(\overline{\mathbb{B}})}\leq&\left\Vert u \right\Vert_{C^{m-1,\alpha}(\overline{\mathbb{B}})}+\max_{|s|=m}{\left\Vert \partial^{s} u\right\Vert}_{C^{0}(\overline{\mathbb{B}})}+\max_{|s|=m}{[\partial^{s} u]}_{\alpha,\mathbb{B}}\\
\leq& C(n,m,\alpha)\cdot\left(\left\Vert u \right\Vert_{C^{0}(\overline{\mathbb{B}})}+\max_{|s|=m}{[\partial^{s}u]_{\alpha,\mathbb{B}}} \right)+\max_{|s|=m}{\left\Vert \partial^{s} u\right\Vert}_{C^{0}(\overline{\mathbb{B}})}+\max_{|s|=m}{[\partial^{s} u]}_{\alpha,\mathbb{B}}\\
\leq& \widehat{C}(n,m,\alpha)\cdot \left|u\right|_{C^{m,\alpha}(\overline{\mathbb{B}})}
\end{split}
\end{align}

Die \"Aquivalenz der Normen $\left|\cdot\right|_{C^{m,\alpha}(\overline{\mathbb{B}})}$ und $\left\Vert \cdot \right\Vert_{C^{m,\alpha}(\overline{\mathbb{B}})}$ ist damit gezeigt. Es sei noch erw\"ahnt, dass ausgehend von \eqref{eq:A.4} aus \eqref{eq:A.3} und \eqref{eq:A.5}, die Absch\"atzung:
\begin{equation}\label{eq:A.6}
\left| u \right|_{C^{k,\beta}(\overline{\mathbb{B}})}\leq C(n,m,k,\alpha,\beta)\cdot\left| u \right|_{C^{m,\alpha}(\overline{\mathbb{B}})}
\end{equation}

f\"ur $k+\beta<m+\alpha$ folgt. F\"ur die n\"achste Absch\"atzung sei erw\"ahnt, dass f\"ur $u,v\in C^{m}(\Omega)$, wobei $\Omega\subseteq\mathbb{R}^n$ eine beliebige offene Menge ist, f\"ur $|s|\leq m$ stets die Gleichung:
\begin{equation}\label{eq:A.7}
\partial^s(uv)=\sum_{\beta\leq s}{\binom{s}{\beta}}\partial^{\beta}u\ \partial^{s-\beta}v
\end{equation}

gilt. Nun seien $u,v\in C^{m,\alpha}(\overline{\mathbb{B}})$, dann ist mit \eqref{eq:A.2} und \eqref{eq:A.6}:
\begin{align*}
&\left|\partial^s(uv) \right|_{C^{0,\alpha}(\overline{\mathbb{B}})}\leq \sum_{\beta\leq s}{\binom{s}{\beta}}\left|\partial^{\beta}u \right|_{C^{0,\alpha}(\overline{\mathbb{B}})}\left|\partial^{s-\beta}v \right|_{C^{0,\alpha}(\overline{\mathbb{B}})}\\
=&\left| u \right|_{C^{0,\alpha}(\overline{\mathbb{B}})}\left|\partial^{s}v \right|_{C^{0,\alpha}(\overline{\mathbb{B}})}+\left|\partial^{s}u \right|_{C^{0,\alpha}(\overline{\mathbb{B}})}\left| v \right|_{C^{0,\alpha}(\overline{\mathbb{B}})}+\sum_{0<\beta< s}{\binom{s}{\beta}}\left|\partial^{\beta}u \right|_{C^{0,\alpha}(\overline{\mathbb{B}})}\left|\partial^{s-\beta}v \right|_{C^{0,\alpha}(\overline{\mathbb{B}})}\\
\leq&\left| u \right|_{C^{0,\alpha}(\overline{\mathbb{B}})}\left|\partial^{s}v \right|_{C^{0,\alpha}(\overline{\mathbb{B}})}+\left|\partial^{s}u \right|_{C^{0,\alpha}(\overline{\mathbb{B}})}\left| v \right|_{C^{0,\alpha}(\overline{\mathbb{B}})}+\sum_{0<\beta< s}{\binom{s}{\beta}}\left| u \right|_{C^{|\beta|,\alpha}(\overline{\mathbb{B}})}\left| v \right|_{C^{|s-\beta|,\alpha}(\overline{\mathbb{B}})}\\
\leq& \left| u \right|_{C^{0,\alpha}(\overline{\mathbb{B}})}\left|\partial^{s}v \right|_{C^{0,\alpha}(\overline{\mathbb{B}})}+\left|\partial^{s}u \right|_{C^{0,\alpha}(\overline{\mathbb{B}})}\left| v \right|_{C^{0,\alpha}(\overline{\mathbb{B}})}+C(n,m,\alpha)\cdot\left|u\right|_{C^{m-1,\alpha}(\overline{\mathbb{B}})}\left|v\right|_{C^{m-1,\alpha}(\overline{\mathbb{B}})}
\end{align*}

Daraus folgt, zusammen mit \eqref{eq:A.2}:
\begin{align}\label{eq:A.8}
\begin{split}
&\left|uv \right|_{C^{m,\alpha}(\overline{\mathbb{B}})}=\left|uv \right|_{C^{0,\alpha}(\overline{\mathbb{B}})}+\sum_{|s|=m}{\left|\partial^s(uv) \right|_{C^{0,\alpha}(\overline{\mathbb{B}})}} \\
\leq&\left|u\right|_{C^{0,\alpha}(\overline{\mathbb{B}})}\left|v\right|_{C^{m,\alpha}(\overline{\mathbb{B}})}+\left|u\right|_{C^{m,\alpha}(\overline{\mathbb{B}})}\left|v\right|_{C^{0,\alpha}(\overline{\mathbb{B}})}+C(n,m,\alpha)\cdot\left|u\right|_{C^{m-1,\alpha}(\overline{\mathbb{B}})}\left|v\right|_{C^{m-1,\alpha}(\overline{\mathbb{B}})}
\end{split}
\end{align}

Durch nochmalige Anwendung von $\eqref{eq:A.6}$ folgt die etwas allgemeinere Absch\"atzung:
\begin{equation}\label{eq:A.9}
\left|uv \right|_{C^{m,\alpha}(\overline{\mathbb{B}})}\leq C(n,m,\alpha)\cdot \left|u\right|_{C^{m,\alpha}(\overline{\mathbb{B}})}\left|v\right|_{C^{m,\alpha}(\overline{\mathbb{B}})}
\end{equation}

Ferner gilt f\"ur $k\in \{1,...,m\}$:
\begin{align}\label{eq:A.10}
\begin{split}
\left|u \right|_{C^{m,\alpha}(\overline{\mathbb{B}})}\stackrel{\hphantom{\eqref{eq:A.6}}}{=}&\left\Vert u \right\Vert_{C^{0,\alpha}(\overline{\mathbb{B}})}+\sum_{|s|=m}{\left\Vert\partial^s u \right\Vert_{C^{0,\alpha}(\overline{\mathbb{B}})}}\leq\left| u \right|_{C^{0,\alpha}(\overline{\mathbb{B}})}+\sum_{|s|=k}{\left|\partial^s u \right|_{C^{m-k,\alpha}(\overline{\mathbb{B}})}}\\
\stackrel{\eqref{eq:A.6}}{\leq}&C(n,\alpha)\cdot\left| u \right|_{C^{2,\alpha}(\overline{\mathbb{B}})}+\sum_{|s|=k}{\left|\partial^s u \right|_{C^{m-k,\alpha}(\overline{\mathbb{B}})}} 
\end{split}
\end{align}

Nun soll ein geeigneter Funktionenraum f\"ur vektorwertige Funktionen, deren Komponentenfunktionen h\"olderstetig sind, eingef\"uhrt werden:

\begin{defi}
Auf $\mathbb{B}$ wird der folgende Vektorraum definiert:
\begin{align*}
\gls{CmalphaBRq}:=\bigl\{&(u_1,...,u_q)^{\top}: \mathbb{B}\longrightarrow \mathbb{R}^q \text{ mit }u_j\in C^{m,\alpha}(\overline{\mathbb{B}})\ \forall j\in\{1,..., q\}\bigr\}
\end{align*}

Auf dieser Menge wird die Norm:
\begin{align*}
\gls{NormaufCmalphaBRq}:= \sum _{j=1}^q{\left|u_j \right|_{C^{m,\alpha}(\overline{\mathbb{B}})}}
\end{align*}

definiert.

\end{defi}

Die Vollst\"andigkeit des Raumes $C^{m,\alpha}(\overline{\mathbb{B}},\mathbb{R}^q)$ folgt aus der Vollst\"andigkeit des Raumes $C^{m,\alpha}(\overline{\mathbb{B}})$. Mit \eqref{eq:A.6} gilt f\"ur $k\in\mathbb{N}$ und $\beta\in (0,1)$ mit $k+\beta<m+\alpha$ die Absch\"atzung:

\begin{equation}\label{eq:A.12}
\left|u\right|_{C^{k,\beta}(\overline{\mathbb{B}},\mathbb{R}^q)}\leq C(n,m,k,\alpha,\beta)\cdot\left|u\right|_{C^{m,\alpha}(\overline{\mathbb{B}},\mathbb{R}^q)}
\end{equation}

Im Folgenden sei $m\geq 1$. Mit \eqref{eq:A.8} folgt f\"ur $a\in C^{m,\alpha}(\overline{\mathbb{B}})$  und $u\in C^{m,\alpha}(\overline{\mathbb{B}},\mathbb{R}^q)$:
\begin{align}\label{eq:A.13}
\begin{split}
&\left|au\right|_{C^{m,\alpha}(\overline{\mathbb{B}},\mathbb{R}^q)}=\sum_{i=1}^q{\left|au_i\right|_{C^{m,\alpha}(\overline{\mathbb{B}})}}\leq\sum_{i=1}^q{\left| a \right|_{C^{0,\alpha}(\overline{\mathbb{B}})}\left|u_i \right|_{C^{m,\alpha}(\overline{\mathbb{B}})}}+\sum_{i=1}^q{\left|a\right|_{C^{m,\alpha}(\overline{\mathbb{B}})}\left| u_i \right|_{C^{0,\alpha}(\overline{\mathbb{B}})}}\\
&\hspace{6cm}+\sum_{i=1}^q{C(n,m,\alpha)\cdot\left|a\right|_{C^{m-1,\alpha}(\overline{\mathbb{B}})}\left|u_i\right|_{C^{m-1,\alpha}(\overline{\mathbb{B}})}}\\
=& \left| a \right|_{C^{0,\alpha}(\overline{\mathbb{B}})}\left| u \right|_{C^{m,\alpha}(\overline{\mathbb{B}},\mathbb{R}^q)}+\left| a \right|_{C^{m,\alpha}(\overline{\mathbb{B}})}\left| u \right|_{C^{0,\alpha}(\overline{\mathbb{B}},\mathbb{R}^q)}+C(n,m,\alpha)\cdot \left| a \right|_{C^{m-1,\alpha}(\overline{\mathbb{B}})}\left| u \right|_{C^{m-1,\alpha}(\overline{\mathbb{B}},\mathbb{R}^q)}
\end{split}
\end{align}

Zusammen mit \eqref{eq:A.6} und \eqref{eq:A.12} ist:
\begin{equation}\label{eq:A.14}
\left|au\right|_{C^{m,\alpha}(\overline{\mathbb{B}},\mathbb{R}^q)}\leq C(n,m,\alpha)\cdot \left| a \right|_{C^{m,\alpha}(\overline{\mathbb{B}})}\left| u \right|_{C^{m,\alpha}(\overline{\mathbb{B}},\mathbb{R}^q)}
\end{equation}

Die Absch\"atzung \eqref{eq:A.14} gilt auch f\"ur $m=0$ mit $C=1$. Sei nun ein weiteres $v\in C^{m,\alpha}(\overline{\mathbb{B}},\mathbb{R}^q)$ gegeben, dann ist mit \eqref{eq:A.8}:
\begin{align}\label{eq:A.15}
\notag&\left|u\cdot v\right|_{C^{m,\alpha}(\overline{\mathbb{B}})}\leq\sum_{i=1}^q{\left|u_i v_i\right|_{C^{m,\alpha}(\overline{\mathbb{B}})}}
\notag\leq \sum_{i=1}^q{\left| u_i \right|_{C^{0,\alpha}(\overline{\mathbb{B}})}\left|v_i \right|_{C^{m,\alpha}(\overline{\mathbb{B}})}}+\sum_{i=1}^q{\left|u_i\right|_{C^{m,\alpha}(\overline{\mathbb{B}})}\left| v_i \right|_{C^{0,\alpha}(\overline{\mathbb{B}})}}\\
\notag&\hspace{6cm}+\sum_{i=1}^q{C(n,m,\alpha)\cdot\left|u_i\right|_{C^{m-1,\alpha}(\overline{\mathbb{B}})}\left|v_i\right|_{C^{m-1,\alpha}(\overline{\mathbb{B}})}}\\
\leq &\left|u \right|_{C^{0,\alpha}(\overline{\mathbb{B}},\mathbb{R}^q)}\left|v \right|_{C^{m,\alpha}(\overline{\mathbb{B}},\mathbb{R}^q)}+\left|u \right|_{C^{m,\alpha}(\overline{\mathbb{B}},\mathbb{R}^q)}\left|v \right|_{C^{0,\alpha}(\overline{\mathbb{B}},\mathbb{R}^q)}\\
\notag&\hspace{6cm}+C(n,m,\alpha)\cdot \left|u \right|_{C^{m-1,\alpha}(\overline{\mathbb{B}},\mathbb{R}^q)}\left|v \right|_{C^{m-1,\alpha}(\overline{\mathbb{B}},\mathbb{R}^q)}
\end{align}

und es folgt mit \eqref{eq:A.12}:
\begin{equation}\label{eq:A.16}
\left|u\cdot v\right|_{C^{m,\alpha}(\overline{\mathbb{B}})}\leq C(n,m,\alpha)\cdot \left|u \right|_{C^{m,\alpha}(\overline{\mathbb{B}},\mathbb{R}^q)}\left|v \right|_{C^{m,\alpha}(\overline{\mathbb{B}},\mathbb{R}^q)}
\end{equation}

Auch die Absch\"atzung \eqref{eq:A.16} gilt f\"ur $m=0$ mit $C=1$. Schlie\ss{}lich gilt mit \eqref{eq:A.10} f\"ur $k\in \{1,...,m\}$:
\begin{align}\label{eq:A.11}
\begin{split}
\left|u \right|_{C^{m,\alpha}(\overline{\mathbb{B}},\mathbb{R}^q)}=\sum _{j=1}^q{\left|u_j \right|_{C^{m,\alpha}(\overline{\mathbb{B}})}}\leq &\sum _{j=1}^q{\left| u_j \right|_{C^{0,\alpha}(\overline{\mathbb{B}})}}+\sum _{j=1}^q{\sum_{|s|=k}{\left|\partial^s u_j \right|_{C^{m-k,\alpha}(\overline{\mathbb{B}})}}}\\
=&\left|u \right|_{C^{0,\alpha}(\overline{\mathbb{B}},\mathbb{R}^q)}+\sum_{|s|=k}{\left|\partial^s u \right|_{C^{m-k,\alpha}(\overline{\mathbb{B}},\mathbb{R}^q)}}
\end{split}
\end{align}

Abschlie\ss{}end werden noch die Sobolev-R\"aume eingef\"uhrt. Hierbei wird f\"ur $p\in[1,\infty]$ der \gls{LRaum}-Raum\index{Lebesgue-Raum} wie in \cite[D.1]{evans1998partial} definiert.

\begin{defi}
Es sei $\Omega\subseteq \mathbb{R}^n$ eine offene und beschr\"ankte Menge. Der \textbf{Sobolev-Raum}\index{Sobolev-Raum} der Ordnung $m\in \mathbb{N}$, mit Exponenten $p\in [1,\infty]$, wird wie folgt definiert:
\begin{align*}
\gls{sobo}:= \bigl\{f\in L^p(\Omega) :\ &\text{F\"ur alle } |s|\leq m \text{ existiert ein } f^{(s)}\in L^p(\Omega) \text{ mit }f^{(0)}=f \text{ und }\\
&\int_{\Omega}{\partial^s \zeta\, f^{(0)}\, dx}=(-1)^{|s|} \int_{\Omega}{\zeta\, f^{(s)} dx} \text{ f\"ur alle }\zeta\in C^{\infty}_0(\Omega)\bigr\}
\end{align*}

Auf diesem Raum wird die Norm $\gls{sobonorm}: W^{m,p}(\Omega)\longrightarrow\mathbb{R}$ wie folgt definiert:
\begin{equation*}
\left\Vert f\right\Vert_{W^{m,p}(\Omega)}:= \sum_{|s|\leq m}{\left\Vert f^{(s)}\right\Vert_{L^p(\Omega)}}
\end{equation*}

Ferner wird der Unterraum $W^{m,p}_0(\Omega)$ wie folgt definiert:
\begin{align*}
\gls{sobonull}:= \bigl\{f\in W^{m,p}(\Omega): \exists (f_k)_{k\in \mathbb{N}}\subseteq C^{\infty}_0(\Omega) \text{ mit } \lim_{k\to\infty}{\left\Vert f-f_k\right\Vert_{W^{m,p}(\Omega)}}=0 \bigr\}
\end{align*}

\end{defi}


\pagebreak
\chapter{Begriffe}
\thispagestyle{fancy}

In diesem Teil des Anhangs werden die grundlegenden topologischen und differentialgeometrischen Begriffe und Zusammenh\"ange, die in  der Arbeit verwendet werden, eingef\"uhrt. Diese wurden im Wesentlichen \cite{munkres2000topology} und \cite{lee2003introduction} entnommen.
\begin{defi}

Eine \textbf{Topologie} \index{Topologie|see{topologischer Raum}} auf einer Menge $X$ ist ein Teilsystem der Potenzmenge $\mathcal{T}\subseteq\mathcal{P}(X)$ mit den folgenden Eigenschaften: 

\begin{enumerate}[(i)]

\begin{item}
$X,\emptyset\in \mathcal{T}$
\end{item}

\begin{item}
Falls $\left\{U_i \right\}_{i\in\Lambda}\subseteq \mathcal{T}$ f\"ur eine beliebige Indexmenge $\Lambda$, dann ist auch $\bigcup_{i\in\Lambda}\in\mathcal{T}$
\end{item}

\begin{item}
Falls $\left\{U_i \right\}_{i\in\Lambda}\subseteq \mathcal{T}$ f\"ur eine endliche Indexmenge $\Lambda$, dann ist auch $\bigcap_{i\in\Lambda}\in\mathcal{T}$
\end{item}

\end{enumerate}

Das Tupel \gls{(X,mathcal{T})} hei\ss{}t \textbf{topologischer Raum}. \index{topologischer Raum} Die Elemente in $\mathcal{T}$ werden \textbf{offene Mengen} \index{Menge>offene $\sim$} genannt. Eine Menge $S\subseteq X$ hei\ss{}t \textbf{abgeschlossen}, \index{Menge>abgeschlossene $\sim$}wenn das Komplement $X\backslash S$ eine offene Menge ist. F\"ur einen Punkt $p \in X$ hei\ss{}t eine Menge $U\in\mathcal{T}$ mit $p\in U$ \textbf{Umgebung des Punktes} $\bm{p}$.
\end{defi}

Wenn aus dem Zusammenhang heraus klar ist, welche Topologie auf $X$ definiert ist, oder Verwechselungen ausgeschlossen sind, dann wird anstatt $(X,\mathcal{T})$ nur $X$ geschrieben.

\begin{defi}
Sei $(X,\mathcal{T})$ ein topologischer Raum und $S\subseteq X$. Die Menge:
\begin{equation*}
\gls{int(S)}:=\bigcup_{U\in\mathcal{T},\  U\subseteq S}{U}
\end{equation*}

wird als das \textbf{Innere von S} \index{Menge>Inneres einer $\sim$} und die Menge:
\begin{equation*}
\gls{overline{S}}:=\bigcap_{U\in\mathcal{T},\ X\backslash U\supseteq S}{X\backslash U}
\end{equation*}

als der \textbf{Abschluss von S} \index{Menge>Abschluss einer $\sim$} bezeichnet. Die Menge:
\begin{equation*}
\gls{randS}:=\overline{S}\backslash int(S)
\end{equation*}

hei\ss{}t \textbf{Rand von S}. \index{Menge>Rand einer $\sim$} Falls $\overline{S}=X$ gilt, dann wird die Menge $S$ als \textbf{dichte Teilmenge} \index{dichte Teilmenge} in $X$ bezeichnet. 

\end{defi}

\begin{defi}

Sei $(X,\mathcal{T})$ ein topologischer Raum, $V$ ein Vektorraum, sowie eine Abbildung $f: X\longrightarrow V$ gegeben. Dann wird die Menge:
\begin{equation*}
\gls{supp(f)}:=\overline{\{p\in M: f(p)\neq 0 \}} 
\end{equation*}

als der \textbf{Tr\"ager von f} bezeichnet.

\end{defi}

\begin{defi}

Sei $(X,\mathcal{T})$ ein topologischer Raum, und sei $\mathcal{U}\subseteq\mathcal{P}(X)$ ein Mengensystem, dann hei\ss{}t $\mathcal{U}$ \textbf{\"Uberdeckung von} $\bm{X}$, \index{\"Uberdeckung} wenn:
\begin{equation*}
\bigcup_{U\in\mathcal{U}}{U}=X
\end{equation*}

gilt. Falls $\mathcal{U}\subseteq\mathcal{T}$, dann wird  $\mathcal{U}$ \textbf{offene \"Uberdeckung von X} \index{\"Uberdeckung>offene $\sim$} genannt.

\end{defi}

\begin{defi}

Ein topologischer Raum $(X,\mathcal{T})$ hei\ss{}t \textbf{kompakt}, \index{topologischer Raum>kompakter $\sim$} falls jede offene \"Uber\-deckung $\mathcal{U}$ von $X$ eine endliche Teil\"uberdeckung $\left\{U_{1},...,U_{N}\right\}\subseteq\mathcal{U}$ besitzt. Letzteres bedeutet, dass:
\begin{equation*}
\bigcup_{j=1}^N{U_{j}}=X
\end{equation*}

gilt.

\begin{defi}

Gegeben sei eine Menge $M$, dann wird eine Abbildung $\gls{Metrik}: M\times M\longrightarrow\mathbb{R}$ als \textbf{Metrik} \index{Metrik|see{metrischer Raum}}bezeichnet, falls die folgenden Eigenschaften erf\"ullt sind:

\begin{enumerate}[(i)]

\begin{item}
$\forall x,y\in M$  ist $d(x,y)\geq 0$ und  $d(x,y)=0$ genau dann wenn $x=y$ gilt
\end{item}

\begin{item}
$\forall x,y\in M$  ist $d(x,y)=d(y,x)$.
\end{item}

\begin{item}
$\forall x,y,z\in M$  ist $d(x,z)\leq d(x,y)+d(y,z)$.
\end{item}

\end{enumerate}

Das Tupel $\gls{metrischer Raum}$ wird als \textbf{metrischer Raum} \index{metrischer Raum} bezeichnet. 

\end{defi}

Ist aus dem Zusammenhang klar, welche Metrik auf $M$ definiert ist, oder sind Verwechselungen auszuschlie\ss{}en, so wird f\"ur $(M,d)$ nur $M$ geschrieben. In einem metrischen Raum wird, f\"ur ein $x\in M$ und einen Radius $R\in\mathbb{R}_{>0}$, die Menge:
\begin{equation*}
\gls{offener Ball}:=\left\{y\in X: d(x,y)< R \right\}
\end{equation*}

als \textbf{offener Ball} \index{offener Ball} \textbf{um} $\bm{x}$ \textbf{mit dem Radius} $\bm{R}$ bezeichnet. Auf dem Vektorraum $\mathbb{R}^n$ wird immer die Metrik verwenden, die von der \textbf{euklidischen Norm} $\gls{euk}:=\sqrt{\sum_{i=1}^n{(x_i)^2}}$ induziert wird, das hei\ss{}t: $d(x,y)=|x-y|_{\mathbb{R}^n}$.

\begin{defi}\label{defi:B11}

Sei $X$ eine Menge. Eine \textbf{Basis} \index{Basis>$\sim$ f\"ur eine Topologie} \textbf{f\"ur eine Topologie auf} $\bm{X}$ ist ein Mengensystem $\mathcal{B}\subseteq\mathcal{P}(X)$ mit folgenden Eigenschaften:
\begin{enumerate}[(i)]
\begin{item}
\begin{equation*}
X=\bigcup_{B\in\mathcal{B}}{B}
\end{equation*}
\end{item}
\begin{item}
Ist $B_1, B_2\in\mathcal{B}$ und $x\in B_1\cap B_2$, dann existiert ein $B_{3}\in\mathcal{B}$, so dass $x\in B_3\subseteq B_1\cap B_2$ gilt.
\end{item}
\end{enumerate}

\end{defi}
Wie es die Wortwahl bereits nahelegt, kann f\"ur eine Menge $X$ eine, in Definition \ref{defi:B11} beschriebene, Basis $\mathcal{B}$ daf\"ur verwendet werden, eine Topologie auf $X$ zu erzeugen. Welche Gestalt diese Topologie haben soll, wird in der folgenden Definition festgelegt:

\begin{defi}

Sei $X$ eine Menge und $\mathcal{B}\subseteq\mathcal{P}(X)$ eine Basis f\"ur eine Topologie auf $X$, dann wird ein Mengensystem $\mathcal{T}_{\mathcal{B}}$ wie folgt definiert:
\begin{equation*}
\text{F\"ur }U\subseteq X\text{ gilt }U\in\mathcal{T}_{\mathcal{B}} \Longleftrightarrow \forall x\in U \ \exists B_x\in\mathcal{B}\text{ mit }x\in B_x\subseteq U
\end{equation*}

Das Mengensystem $\mathcal{T}_{\mathcal{B}}$ ist eine Topologie auf $X$ und wird als \textbf{die von} $\mathbf{\mathcal{B}}$ \textbf{erzeugte Topologie} \index{Topologie>von Basis erzeugte $\sim$}\textbf{auf }$\bm{X}$ bezeichnet.  \cite[\S 13]{munkres2000topology}

\end{defi}

\begin{defi}
Sei $(X,d)$ metrischer Raum, dann wird die von der Basis:
\begin{equation*}
\mathcal{B}=\left\{B_R(x): x\in X\text{ und }R>0 \right\}
\end{equation*}

erzeugte Topologie auf $X$ als die \textbf{von der Metrik} $\bm{d}$ \textbf{erzeugte Topologie}\index{Topologie>von Metrik erzeugte $\sim$} bezeichnet. Ist $X=\mathbb{R}^n$ und $d$ die euklidische Metrik, so hei\ss{}t diese Topologie \textbf{Standardtopologie auf }$\mathbb{R}^n$.\index{Topologie>Standard-$\sim$}
\end{defi}

\begin{defi}
Ein topologischer Raum $(X,\mathcal{T})$, dessen Topologie von einer Metrik $d$ induziert wird, hei\ss{}t \textbf{metrisierbar}.\index{topologischer Raum>metrisierbarer $\sim$}
\end{defi}

\begin{defi}
Sei $(X,\mathcal{T})$ ein topologischer Raum, und $S\subseteq X$ eine Teilmenge von X. Dann wird ein Mengensystem $\mathcal{T}_S\subseteq\mathcal{P}(X)$ wie folgt definiert:
\begin{align*}
U\in\mathcal{T}_S \Longleftrightarrow \exists V\in\mathcal{T} \text{ mit } U=V\cap S. 
\end{align*}

Das Mengensystem $\mathcal{T}_S$ ist eine Topologie auf $S$, und wird als die von $\bm{X}$ \textbf{auf} $\bm{S}$ \textbf{induzierte Unterraumtopologie}\index{Topologie>Unterraum-$\sim$} bezeichnet. Ein Element in $\mathcal{T}_S$ wird \textbf{in }$\bm{S}$\textbf{ relativ offene Menge} \index{Menge>relativ offene $\sim$} genannt. \cite[\S 16]{munkres2000topology}

\end{defi}

\begin{defi}
Sei $(X,\mathcal{T})$ ein topologischer Raum. Eine Teilmenge $S\subseteq X$ hei\ss{}t \textbf{relativ kompakt}, \index{Menge>relativ kompakte $\sim$} wenn $(\overline{S},\mathcal{T}_{\overline{S}})$ ein kompakter topologischer Raum ist.
\end{defi}

\begin{defi}
Ein topologischer Raum $(X,\mathcal{T})$ hei\ss{}t \textbf{Hausdorff-Raum}, \index{Hausdorff-Raum}falls f\"ur eine beliebige Auswahl von Punkten $p,q\in X$, mit $p\neq q$, eine offene Menge $U\in\mathcal{T}$ mit $p\in U$, und eine offene Menge $V\in\mathcal{T}$ existiert, so dass $U\cap V=\emptyset$ gilt.
\end{defi}

\begin{defi}
Ein topologischer Raum $(X,\mathcal{T})$ hei\ss{}t \textbf{separabel}, \index{topologischer Raum>separabler $\sim$} wenn es eine abz\"ahl"-bare Teilmenge $S\subseteq X$ gibt, so dass $\overline{S}=X$ gilt.
\end{defi}

Ist $U\in\mathcal{T}$, dann ist $U$, versehen mit der Unterraumtopologie, separabel, falls $(X,\mathcal{T})$ separabel ist. \cite[\S 30]{munkres2000topology}. Der Vektorraum $\mathbb{R}^n$, versehen mit der Standardtopologie, ist separabel, da die Menge $\mathbb{Q}^n$ abz\"ahlbar ist, und dicht in $\mathbb{R}^n$ liegt.

\begin{defi}
Ein topologischer Raum $(X,\mathcal{T})$ erf\"ullt das \textbf{zweite Abz\"ahlbarkeits"-axiom}, \index{zweites Abz\"ahlbarkeitsaxiom}falls es eine abz\"ahlbare Basis $\mathcal{B}$ gibt, welche die Topologie $\mathcal{T}$ erzeugt.
\end{defi}

Sowohl die Hausdorff-Eigenschaft, als auch das zweite Abz\"ahlbarkeitsaxiom \"ubertragen sich auf jeden beliebigen Unterraum von $(X,\mathcal{T})$ \cite[Theorem 17.11 und Theorem 30.2]{munkres2000topology}. Jeder metrisierbare topologische Raum erf\"ullt die Hausdorff-Eigenschaft. Ist der topologische Raum zus\"atzlich separabel, so erf\"ullt der Raum auch das zweite Abz\"ahl\-barkeitsaxiom. \cite[\S 21 und \S 30]{munkres2000topology} 

\begin{defi}
Eine Abbildung $F: X\longrightarrow Y$ zwischen zwei topologischen R\"aumen $(X,\mathcal{T}_X)$ und $(Y,\mathcal{T}_Y)$ hei\ss{}t \textbf{stetig}, \index{Stetigkeit} falls:
\begin{equation*}
U\in\mathcal{T}_Y \Longrightarrow F^{-1}(U)\in\mathcal{T}_X
\end{equation*}

gilt. Ist die Abbildung $F$ zus\"atzlich invertierbar, und ist die inverse Abbildung $F^{-1}:Y\longrightarrow X$ ebenfalls stetig, so hei\ss{}t diese Abbildung \textbf{Hom\"oomorphismus}. \index{Hom\"oomorphismus}Die topologischen R\"aume hei\ss{}en in diesem Fall \textbf{hom\"o{}omorph}. \index{hom\"o{}omorph|see{Hom\"oomorphismus}}
\end{defi}

\begin{defi}
Eine \textbf{topologische Mannigfaltigkeit} \index{Mannigfaltigkeit>topologische $\sim$}der Dimension $n\in\mathbb{N}$ ist ein topologischer Raum $(M,\mathcal{O})$, welcher die folgenden Eigenschaften erf\"ullt:
\begin{enumerate}[(i)]
\begin{item}
die Hausdorff-Eigenschaft
\end{item}
\begin{item}
das zweite Abz\"ahlbarkeitsaxiom
\end{item}
\end{enumerate}

Zus\"atzlich soll dieser topologische Raum \textbf{lokal euklidisch} \index{lokal euklidisch} sein. Das bedeutet: F\"ur jedes $p\in M$ existiert eine Umgebung $U\in\mathcal{O}$ von $p$, und eine in $\mathbb{R}^n$ offene Menge $V\subseteq\mathbb{R}^n$, sowie ein Hom\"oomorphismus $\varphi: U\longrightarrow V$. Hierbei werden die Unterraumtopologien verwendet.\\
Die Menge $U$ wird als \textbf{Koordinatenumgebung von} $\bm{p}$\index{Koordinaten>$\sim$-Umgebung} bezeichnet. Das Tupel $(U,\varphi)$ hei\ss{}t \textbf{Karte}. \index{Karte}Da $V=\varphi(U)$ gilt, wird $V$ in der Notation nicht direkt ber\"ucksichtigt. Falls $\varphi(p)=0$ gilt, dann hei\ss{}t $U$ \textbf{um} $\bm{p}$ \textbf{zentrierte Koordinatenumgebung}\index{Koordinaten>zentrierte $\sim$-Umgebung}, und falls $V=B_{R}(x)$, f\"ur ein $R>0$ und $x\in\mathbb{R}^n$, so hei\ss{}t U \textbf{Koordinatenball}. \index{Koordinaten> $\sim$-ball}
\end{defi}

Um die Dimension einer topologischen Mannigfaltigkeit $(M,\mathcal{O})$ auszudr\"ucken, wird auch kurz $\gls{dim(M)}=n$ geschrieben.

\begin{defi}
Sei $(M,\mathcal{O})$ ein topologische Mannigfaltigkeit, und seien $(U,\varphi), (V,\psi)$ Karten mit $U\cap V\neq\emptyset$. Dann hei\ss{}t die Abbildung:
\begin{equation*}
\psi\circ\varphi^{-1}: \varphi(U\cap V)\longrightarrow \psi(U\cap V)
\end{equation*}

\textbf{Koordinatentransformation}. \index{Koordinaten>$\sim$-Transformation}Zwei Karten $(U,\varphi)$, $(V,\psi)$ hei\ss{}en $\bm{C^{\infty}}$-\textbf{vertr\"aglich}, \index{$C^{\infty}$-vertr\"aglich} falls $U\cap V=\emptyset$ gilt, oder $U\cap V\neq\emptyset$ und die Koordinatentransformation $\psi\circ\varphi^{-1}$ ein $C^{\infty}$-\textbf{Diffeomorphismus}\index{$C^{\infty}$-Diffeomorphismus} ist, das bedeutet, dass $\psi\circ\varphi^{-1}$ invertierbar ist, und es ist sowohl $\psi\circ\varphi^{-1}\in C^{\infty}(\varphi(U\cap V),\psi(U\cap V))$, als auch $\varphi\circ\psi^{-1}\in C^{\infty}(\psi(U\cap V),\varphi(U\cap V))$ erf\"ullt.
\end{defi}

\begin{defi}
Sei $(M,\mathcal{O})$ eine topologische Mannigfaltigkeit. Eine Familie von Karten $\mathcal{A}=\{(U_i,\varphi_i)\}_{i\in I}$ hei\ss{}t  \textbf{Atlas}, \index{Atlas} falls:
\begin{equation*}
\bigcup_{i\in I}{U_i}=M
\end{equation*}

gilt. Sind zus\"atzlich, f\"ur jede Auswahl $i,j\in I$, die Karten $(U_i,\varphi_i)$ und $(U_j,\varphi_j)$ $C^{\infty}$-vertr\"aglich, so wird $\mathcal{A}$ als $\bm{C^{\infty}}$-\textbf{Atlas}\index{Atlas}\index{Atlas>$C^{\infty}$-$\sim$} auf $(M,\mathcal{O})$ bezeichnet.

\end{defi}

\begin{defi}
Sei $(M,\mathcal{O})$ eine topologische Mannigfaltigkeit. Ein $C^{\infty}$-Atlas $\mathcal{A}$ auf $(M,\mathcal{O})$ hei\ss{}t \textbf{maximaler} $\bm{C^{\infty}}$-\textbf{Atlas}\index{maximaler $C^{\infty}$-Atlas} auf $(M,\mathcal{O})$, wenn die folgende Aussage gilt:\\ 
Ist $(U,\varphi)$ eine Karte, die mit allen Karten $(V,\psi)\in\mathcal{A}$\ $C^{\infty}$-vertr\"aglich ist, dann gilt $(U,\varphi)\in\mathcal{A}$.\\ 
Ein solcher Atlas $\mathcal{A}$ hei\ss{}t \textbf{glatte Struktur}\index{Struktur>glatte $\sim$}, oder auch $\bm{C^{\infty}}$-\textbf{Struktur}\index{Struktur>$C^{\infty}$-$\sim$|see{glatte $\sim$}} auf $(M,\mathcal{O})$. 
\end{defi}

\end{defi}

In \cite[Lemma 1.10. (a)]{lee2003introduction} wird gezeigt, dass f\"ur jeden $C^{\infty}$-Atlas $\widehat{\mathcal{A}}$ genau ein maximaler $C^{\infty}$-Atlas $\overline{\mathcal{A}}$ existiert, der $\widehat{\mathcal{A}}$ enth\"alt.

\begin{defi}\label{defi:B6}
Eine \textbf{n}-\textbf{dimensionale glatte Mannigfaltigkeit}\index{Mannigfaltigkeit>glatte $\sim$} $(M,\mathcal{O},\mathcal{A})$, oder auch $\bm{C^{\infty}}$-\textbf{Mannigfaltigkeit}\index{Mannigfaltigkeit>$C^{\infty}$-$\sim$|see{glatte-$\sim$}}, ist eine topologische Mannigfaltigkeit $(M,\mathcal{O})$, zusammen mit einem maximalen $C^{\infty}$-Atlas $\mathcal{A}$.
\end{defi}

Ist aus dem Zusammenhang heraus erkennbar, welche Topologie $\mathcal{O}$, und welcher Atlas $\mathcal{A}$ auf $M$ definiert sind, oder sind Verwechselungen auszuschlie\ss{}en, so wird nur $\gls{M}$, oder um die Dimension hervorzuheben $M^n$, anstatt $(M,\mathcal{O},\mathcal{A})$ geschrieben. In diesem Zusammenhang wird mit \gls{Atlas} immer der Atlas auf $(M,\mathcal{O})$ bezeichnet.\\
Im Folgenden soll der Begriff der Produktmannigfaltigkeit definiert werden, dazu seien $\{(M_i,\mathcal{O}_i,\mathcal{A}_i)  \}_{i\in \{1,...,k\}}$ glatte Mannigfaltigkeiten, mit $dim(M_i)=n_i$ f\"ur alle $i\in\{1,...,k\}$. Es sei:
\begin{equation*}
M_{\times}:=\gls{prod_{i=1}^k {M_i}}
\end{equation*}

das kartesische Produkt der Mengen $M_1,...,M_k$, und es sei $\mathcal{O}_{\times}\subseteq \mathcal{P}\left(M_{\times}\right)$ die Produkttopologie auf der Menge $M_{\times}$ \cite[\S 19]{munkres2000topology}. Mit \cite[\S 21 und \S 30]{munkres2000topology} erf\"ullt der topologische Raum $\left(M_{\times},\mathcal{O}_{\times} \right)$ sowohl die Hausdorff-Eigenschaft, als auch das zweite Abz\"ahlbarkeitsaxiom. F\"ur ein $(p_1,...,p_k)\in M_{\times}$ k\"onnen, f\"ur $i\in\{1,...,k\}$, Karten\linebreak $(U_i,\varphi_i)\in\mathcal{A}_i$  gew\"ahlt werden, so dass $p_i\in U_i$ f\"ur alle $i\in \{1,...k\}$ gilt. Die Abbildung:
\begin{align}
\label{eq:B1}
\begin{split}
\varphi_{\times}: \prod_{i=1}^k {U_i}&\longrightarrow\mathbb{R}^{\sum_{i=1}^k{n_i}}\\
(x_1,...,x_k)&\mapsto\left(\varphi(x_1),..., \varphi(x_k) \right)
\end{split}
\end{align}

ist ein Hom\"oomorphismus auf das Bild:
\begin{equation*}
\varphi_{\times}\left(\prod_{i=1}^k {U_i} \right)=\prod_{i=1}^k {\varphi_i(U_i)}
\end{equation*}

Der topologische Raum $\left(M_{\times},\mathcal{O}_{\times} \right)$ ist also eine topologische Mannigfaltigkeit. Die Vereinigung aller Karten der Form \eqref{eq:B1} bildet einen $C^{\infty}-$Atlas \cite[Chapter 1]{lee2003introduction}. Mit $A_{\times}$ wird der maximale $C^{\infty}-$Atlas bezeichnet, der diesen Atlas enth\"alt. Es wird die folgende Bezeichnung festgelegt:

\begin{defi}\label{defi:B9}
Die glatte Mannigfaltigkeit $(M_{\times},\mathcal{O}_{\times},\mathcal{A}_{\times})$ wird als \textbf{glatte Produktmannigfaltigkeit}\index{Mannigfaltigkeit>glatte Produkt-$\sim$} bezeichnet.
\end{defi}

Der Konstruktion dieser Mannigfaltigkeit ist direkt zu entnehmen, dass $dim\ M_{\times}= \mathbb{R}^{\sum_{i=1}^k{n_i}}$ gilt.\\
F\"ur eine beliebige Menge $S$ wird die Abbildung $\gls{id_S}: S\longrightarrow S$ mit $id_S(x)=x$ f\"ur alle $x\in S$ definiert. Diese Abbildung wird auch als \textbf{Identit\"at auf} $\bm{S}$ bezeichnet.\index{Identit\"at auf S}

\begin{defi}
Sei $\Omega\subseteq\mathbb{R}^n$ offen und $\mathcal{A}$, der Atlas der die Abbildung $id_{\Omega}: \Omega\longrightarrow\Omega$ enth\"alt. Dann wird die glatte Struktur, die aus $\Omega$, versehen mit der Unterraumtopologie und dem Atlas $\mathcal{A}$ besteht, als \textbf{Standardstruktur}\index{Struktur>Standard-$\sim$} \textbf{auf} $\bm{\Omega}$ bezeichnet. 
\end{defi}

\begin{defi}\label{defB1}
Seien $(M,\mathcal{O}_M,\mathcal{A}_M)$ und $(N,\mathcal{O}_N,\mathcal{A}_N)$ glatte Mannigfaltigkeiten. Eine Abbildung $F: M\longrightarrow N$ hei\ss{}t \textbf{glatt}\index{Abbildung>glatte $\sim$}, falls es zu jedem $p\in M$ eine Karte $(U,\varphi)\in \mathcal{A}_M$ mit $p\in U$, und eine Karte $(V,\psi)\in\mathcal{A}_N$ mit $F(U)\subseteq V$ gibt, so dass die \textbf{Koordinatendarstellung}\index{Koordinaten>$\sim$-darstellung} von $F$ bez\"uglich $(U,\varphi)$ und $(V,\psi)$:
\begin{equation}
\label{eq:B2}
\gls{^{varphi}_{psi}F}:=\psi\circ F\circ \varphi^{-1}: \varphi(U)\longrightarrow \psi(V)
\end{equation}

eine glatte Abbildung im Sinne von Definition \ref{defi:A.1a} ist. In diesem Fall wird die Schreibweise $F\in \gls{C^{infty}(M,N)}$ verwendet. Im Spezialfall $N=\mathbb{R}$, wobei $\mathbb{R}$ die Standardstruktur tr\"agt, wird nur $F\in \gls{C^{infty}(M)}$ geschrieben. Hierf\"ur wird die Koordinatendarstellung in \eqref{eq:B2} mit \gls{^{varphi}F} notiert, falls $\psi=id_V$ ist.

\end{defi}

Die Menge $C^{\infty}(M,\mathbb{R}^q)$, mit der punktweisen Addition und der skalaren Multiplikation, ist ein $\mathbb{R}$-Vektorraum. Mit der punktweisen Multiplikation von Funktionen ist der Vektorraum $C^{\infty}(M)$ eine $\mathbb{R}$-Algebra. Verkn\"upfungen von glatten Abbildungen auf glatten Mannigfaltigkeiten sind glatt.\cite[Lemma 2.4]{lee2003introduction}

\begin{defi}

Sei $(M,\mathcal{O},\mathcal{A})$ eine glatte Mannigfaltigkeit. F\"ur ein $p\in M$ hei\ss{}t eine lineare Abbildung $X: C^{\infty}(M)\longrightarrow\mathbb{R}$ \textbf{Derivation}\index{Derivation} \textbf{in} $\bm{p}$, falls f\"ur alle $f,g\in C^{\infty}(M)$ die Produktregel:
\begin{equation*}
X(fg)=f(p)\cdot Xg+Xf\cdot g(p)
\end{equation*}

erf\"ullt ist. F\"ur zwei Derivationen $X,Y$ im Punkt p wird die Addition:
\begin{equation*}
(X+Y)f:= Xf+ Yf\hspace{0.5cm}\text{f\"ur }f\in C^{\infty}(M)
\end{equation*}

und die skalare Multiplikation:
\begin{equation*}
(\lambda X)f:= \lambda\cdot Xf\hspace{0.5cm}\text{f\"ur }\lambda\in\mathbb{R}\hspace{0.5cm}\text{und}\hspace{0.5cm}f\in C^{\infty}(M)
\end{equation*}

definiert. Der mit diesen Operationen definierte Vektorraum aller Derivationen im Punkt $p$ wird \textbf{Tangentialraum}\index{Tangentialraum} \textbf{von} $\bm{M}$ \textbf{im Punkt} $\bm{p}$ genannt, und mit $\gls{T_p M}$ bezeichnet. 

\end{defi}

\begin{defi}\label{defB2}

Seien $(M,\mathcal{O}_M,\mathcal{A}_M)$ und $(N,\mathcal{O}_N,\mathcal{A}_N)$ glatte Mannigfaltigkeiten sowie $F\in C^{\infty}(M,N)$. F\"ur ein $p\in M$ hei\ss{}t die Abbildung $F_{\ast}: T_p M\longrightarrow T_{F(p)}N$, die wie folgt definiert ist:
\begin{equation*}
(F_{\ast}X)f=X(f\circ F)\hspace{0.5cm}\text{f\"ur }f\in C^{\infty}(N)\text{ und }X\in T_p M
\end{equation*}

\textbf{Ableitung}\index{Ableitung} \textbf{von} $\bm{F}$ \textbf{im Punkt} $\bm{p}$.

\end{defi}

Die grundlegenden Eigenschaften der Ableitung, wie zum Beispiel die Linearit\"at, sind in \cite[Lemma 3.5]{lee2003introduction} aufgef\"uhrt.

\begin{defi}\label{defi:B10}

Sei $(M,\mathcal{O},\mathcal{A})$ eine glatte Mannigfaltigkeit, und sei $U\subseteq M$ eine offene Menge. Dann wird der topologische Raum $U$, zusammen mit der glatten Struktur:
\begin{equation*}
\left\{(V,\varphi)\in\mathcal{A}:  V\subseteq U\right\}
\end{equation*}
als \textbf{offene Untermannigfaltigkeit}\index{Mannigfaltigkeit>Unter-$\sim$>offene $\sim$} bezeichnet.

\end{defi}

F\"ur eine beliebige Menge $X$ und eine Teilmenge $S\subseteq X$ wird die Abbildung $i: S\longrightarrow X$, mit $i(x)=x$ f\"ur $x\in U$, als \textbf{Inklusion}\index{Inklusion} bezeichnet. In \cite[Lemma 3.7]{lee2003introduction} wird gezeigt, dass f\"ur eine offene Untermannigfaltigkeit $U\subseteq M$ und ein beliebiges $p\in U$ die Ableitung der Inklusion $i_{\ast}: T_p U\longrightarrow T_p M$ ein Isomorphismus ist. Diese Tatsache erlaubt die Identifikation des Tangentialraumes $T_p U$ mit dem Tangentialraum $T_p M$.

\begin{defi}
Sei $(M,\mathcal{O},\mathcal{A})$ eine glatte Mannigfaltigkeit der Dimension $n$, und sei $p\in M$, ferner sei $(U,\varphi)\in\mathcal{A}$ eine Karte mit $p\in U$. Dann wird f\"ur $i\in\{1,..,n\}$ der $\bm{i}$\textbf{-te} \textbf{Koordinatenvektor}\index{Koordinaten>Koordinaten.$\sim$} $\left. ^{\varphi}\frac{\partial}{\partial x^i}\right|_{p}\in T_p U$ folgenderma\ss{}en definiert:
\begin{equation*}
\gls{phi i Koord}(f):=\partial_i\, ^{\varphi}f(\varphi(p))\hspace{0.5cm}\text{f\"ur }f\in C^{\infty}(U)
\end{equation*}                                                                                     

In Anlehnung an die vorhergehende Bemerkung, wird der Tangentialvektor $i_{\ast}\left(\left. ^{\varphi}\frac{\partial}{\partial x^i}\right|_{p}\right)$ auch mit $\left. ^{\varphi}\frac{\partial}{\partial x^i}\right|_{p}$ bezeichnet.

\end{defi}

In \cite[Lemma 3.9]{lee2003introduction} wird gezeigt, dass die Menge der Koordinatenvektoren:
\begin{equation}\label{eq:B3}
\left\{\left. ^{\varphi}\frac{\partial}{\partial x^1}\right|_{p},...,\left. ^{\varphi}\frac{\partial}{\partial x^n}\right|_{p} \right\}\subseteq T_p U
\end{equation}

eine Basis des Tangentialraumes $T_p U$ bildet. Daraus folgt insbesondere, dass $T_p U$ beziehungsweise $T_p M$ ein Vektorraum der Dimension n ist. Der gem\"a\ss{} \cite[6.1]{fischer2008lineare} definierte Dualraum zu $T_p M$ wird mit \gls{Dualraum} bezeichnet. Die Elemente in $T_p^{\ast} M$ werden \textbf{Kovektoren}\index{Kovektor} genannt. Der Dualraum besitzt eine Basis:
\begin{equation}\label{eq:B4}
\left\{\left. ^{\varphi}dx^1\right|_{p},..., \left. ^{\varphi}dx^n\right|_{p}\right\}\subseteq T_p^{\ast} M
\end{equation}

welche \"uber die Bestimmungsgleichungen:
\begin{equation*}
\gls{i Dual}\left(\left. ^{\varphi}\frac{\partial}{\partial x^j}\right|_{p} \right)=\delta^i_j\hspace{1cm}\text{f\"ur alle }i,j\in \{1,...,n\}
\end{equation*}

definiert ist. Diese Basis wird als die zu \eqref{eq:B3}\index{Basis>duale $\sim$} \textbf{duale Basis} bezeichnet. In \cite[Chapter 3]{lee2003introduction} wird gezeigt, dass, in der Situation in Definition \ref{defB2}, die Ableitung der Abbildung $F$, in der Basis der Koordinatenvektoren zu $(U,\varphi)\in\mathcal{A}_M$, mit $p\in U$ und $(V,\psi)\in\mathcal{A}_N$ mit $F(U)\subseteq V$, von der Jacobimatrix der Abbildung $^{\varphi}_{\psi}F$ im Punkt $\varphi(p)$ repr\"asentiert wird. Konkret bedeutet das f\"ur $i\in\{1,...,dim(M)\}$:
\begin{equation}\label{eq:B6}
F_{\ast}\left. ^{\varphi}\frac{\partial}{\partial x^i}\right|_{p}=\sum_{j=1}^{dim(N)}{\partial_i\, ^{\varphi}_{\psi}F^j(\varphi(p))\left. ^{\psi}\frac{\partial}{\partial x^j}\right|_{F(p)}}
\end{equation}

\begin{defi}
Sei $F\in C^{\infty}(M,N)$ eine glatte Abbildung zwischen zwei glatten Mannigfaltigkeiten $(M,\mathcal{O}_M,\mathcal{A}_M)$ und $(N,\mathcal{O}_N,\mathcal{A}_N)$, dann wird der \textbf{Rang der Abbildung}\index{Rang einer Abbildung} $\bm{F}$ \textbf{im Punkt} $\bm{p}$ als der Rang der Ableitung $F_{\ast}: T_p M\longrightarrow T_{F(p)}N$ definiert.
\end{defi}

Hat die Abbildung $F$ in jedem Punkt $p\in M$ den Rang $dim(M)$, dann wird die Abbildung \textbf{Immersion}\index{Immersion} genannt. Falls $F$ eine Immersion ist und die Menge $M$ hom\"oomorph auf die Menge $F(M)$ abbildet, wobei $F(M)$ mit der Unterraumtopologie versehen wird, dann wird diese Abbildung \textbf{glatte Einbettung}\index{Einbettung>glatte $\sim$} oder $\bm{C^{\infty}}-$\textbf{Einbettung}\index{Einbettung>$C^{\infty}-$ $\sim$} genannt.

\begin{defi}
Sei $M$ ein topologischer Raum. Ein \textbf{Vektorb\"undel vom Rang k}\index{Vektorb\"undel} \textbf{\"uber/auf} $\bm{M}$  ist ein topologischer Raum $E$, zusammen mit einer surjektiven stetigen Abbildung $\pi: E\longrightarrow M$, wobei die folgenden Eigenschaften erf\"ullt sein sollen:

\begin{enumerate}[(i)]

\begin{item}
F\"ur jedes $p\in M$ hat die Menge $E_p:=\pi^{-1}(p)\subseteq E$, welche auch als \textbf{Faser}\index{Faser} \textbf{von} $\bm{E}$ \textbf{\"uber} $\bm{p}$ bezeichnet wird, die Gestalt eines k-dimensionalen $\mathbb{R}$-Vektorraumes.
\end{item}

\begin{item}
F\"ur jedes $p\in M$ existiert eine Umgebung $U$ von $p$ und ein Hom\"oomorphismus $\Phi: \pi^{-1}(U)\longrightarrow U\times \mathbb{R}^k$, so dass das folgende Diagramm kommutiert:

\begin{equation*}
\begin{xy}
\xymatrix{
         \pi^{-1}(U) \ar[rr]^{\Phi} \ar[rd]_{\pi}  &     &  U\times \mathbb{R}^k \ar[dl]^{\pi_1}  \\
                                &  U  &
     }
\end{xy}
\end{equation*}

\end{item}

Hierbei ist $\pi_1: U\times \mathbb{R}^k\longrightarrow U$ mit $\pi_1(x,y)=x$ f\"ur $(x,y)\in U\times \mathbb{R}^k$ definiert. Dar\"uber hinaus soll die Einschr\"ankung:
\begin{equation*}
\left.\Phi \right|_{E_p}: {E_p}\longrightarrow \{q\}\times \mathbb{R}^k \cong \mathbb{R}^k
\end{equation*}
ein linearer Isomorphismus sein.

\end{enumerate}

Die in (ii) definierte Abbildung $\Phi$ wird als \textbf{lokale Trivialisierung}\index{lokale Trivialisierung} \textbf{von} $\bm{E}$ \textbf{\"uber} $\bm{M}$ bezeichnet. Die Mannigfaltigkeit $M$ wird in diesem Zusammenhang auch als \textbf{Basis} \index{Basis>$\sim$ eines Vektorb\"undels}bezeichnet, w\"ahrend f\"ur $E$ die Bezeichnung \textbf{Totalraum}\index{Totalraum} gebr\"auchlich ist. Ein Vektorb\"undel wird kurz \gls{Vektorbundel} notiert. 
\end{defi}

Falls $(M,\mathcal{O}_M,\mathcal{A}_M)$ und $(E,\mathcal{O}_E,\mathcal{A}_E)$ glatte Mannigfaltigkeiten sind, die Abbildung $\pi: E\longrightarrow M$ eine glatte Abbildung ist, und die lokalen Trivialisierungen als Diffeomorphismen gew\"ahlt werden k\"onnen, so hei\ss{}t $E$ \textbf{glattes Vektorb\"undel}\index{Vektorb\"undel>glattes $\sim$} \"uber $M$.

\begin{defi}
Ein glattes Vektorb\"undel $(E,M,\pi)$ vom Rang $k$ \"uber einer glatten Mannigfaltigkeit $(M,\mathcal{O},\mathcal{A})$ hei\ss{}t \textbf{glatt trivialisierbar}\index{Vektorb\"undel>glatt trivialisierbares $\sim$}, falls eine glatte lokale Trivialisierung $\Phi: E\longrightarrow M\times \mathbb{R}^k$, also eine glatte globale Trivialisierung, existiert.
\end{defi}

\begin{defi}
Gegeben sei ein glattes Vektorb\"undel $(E,M,\pi)$ auf einer glatten Mannigfaltigkeit $(M,\mathcal{O},\mathcal{A})$. Eine stetige Abbildung $\sigma : M\longrightarrow E$ hei\ss{}t \textbf{Schnitt in}\index{Schnitt} $\bm{E}$, falls $\pi\circ\sigma=id_M $ gilt. Ist die Abbildung $\sigma$ glatt, so wird sie \textbf{glatter Schnitt in}\index{Schnitt>glatter $\sim$} $\bm{E}$ genannt.
\end{defi}

Im Folgenden soll der Begriff der eingebetteten Untermannigfaltigkeit definiert werden, hierzu sei $(M,\mathcal{O},\mathcal{A})$ eine n-dimensionale glatte Mannigfaltigkeit. Es sei eine Menge $S\subseteq M$ mit der folgenden Eigenschaft gegeben:\\
Es existiert ein $k\in \mathbb{N}$, f\"ur dass f\"ur alle $p\in M$ ein $(U_p,\varphi_p)\in\mathcal{A}$ mit $p\in U_p$ existiert, so dass:
\begin{equation*}
\varphi_p(S\cap U_p)=\left(\mathbb{R}^k\times\{0\}\right)\cap\varphi_p(U_p)
\end{equation*}

gilt. Dann ist f\"ur jedes $p\in M$ die Abbildung:
\begin{align}
\begin{split}\label{eq:B40}
\widehat{\varphi}_p: S\cap U_p&\longrightarrow V_p:=(pr_k\circ\varphi_p)(U_p)\subseteq\mathbb{R}^k\\
x&\mapsto (pr_k\circ\varphi_p)(x)
\end{split}
\end{align}

wobei:
\begin{align*}
pr_k:\mathbb{R}^n&\longrightarrow\mathbb{R}^k\\
(y^1,...,y^n)&\mapsto(y^1,...,y^k)
\end{align*}

ein Hom\"oomorphismus \cite[Theorem 8.2]{lee2003introduction}. Somit ist $(S,\mathcal{O}_S)$, mit der Unterraumtopologie $\mathcal{O}_S$, eine topologische Mannigfaltigkeit der Dimension $k$. Ferner ist der Atlas $\{(\widehat{\varphi}_p,S\cap U_p)\}_{p\in M}$, mit Karten der Form \eqref{eq:B40}, ein $C^{\infty}$-Atlas auf $(S,\mathcal{O}_S)$. Dies rechtfertigt die folgende Definition:

\begin{defi}\label{defi:B8}

Die glatte Mannigfaltigkeit $(S,\mathcal{O}_S,\mathcal{A}_S)$, wobei $\mathcal{A}_S$ der maximale Atlas ist, der die Familie von Karten $\{(\widehat{\varphi}_p,S \cap U_p)\}_{p\in M}$ enth"alt, wird als \textbf{eingebettete Untermannigfaltigkeit der Dimension}\index{Mannigfaltigkeit>Unter-$\sim$>eingebettete $\sim$} $\bm{k}$ bezeichnet.

\end{defi}

\begin{defi}

Gegeben sei eine n-dimensionale glatte Mannigfaltigkeit $(M,\mathcal{O},\mathcal{A})$ und ein $p\in M$. Ein \textbf{kovarianter 2-Tensor}\index{Tensor>kovarianter 2-$\sim$} \textbf{auf} $\bm{T_p M}$ ist eine bilineare Abbildung:
\begin{equation*}
T: T_p M \times T_p M\longrightarrow \mathbb{R}
\end{equation*}

Die Menge aller kovarianten $2-$Tensoren auf $T_p M$ wird mit \gls{Raum kov 2} bezeichnet. 

\end{defi}

Auf dieser Menge werden Operationen:
\begin{align*}
(a T)(X_1,X_2)&:= a\cdot T(X_1,X_2)& &X_1,X_2\in T_p M,\ a\in\mathbb{R}\\
(T+T')(X_1,X_2)&:= T(X_1,X_2)+T'(X_1,X_2)& & X_1,X_2\in T_p M
\end{align*}

definiert. Mit diesen Operationen wird die Menge $T^2(T_p M)$ zu einem $\mathbb{R}-$Vektorraum. Ein $T\in T^2(T_p M)$ hei\ss{}t \textbf{symmetrisch}\index{Tensor>kovarianter 2-$\sim$>symmetrischer $\sim$}, falls $T(X_1,X_2)=T(X_2,X_1)$ f\"ur  alle $X_1,X_2\in T_p M$ gilt. Die Menge aller symmetrischen kovarianten 2-Tensoren bildet einen Unterraum des Vektorraumes $T^2(T_p M)$, und wird mit \gls{Sigma^2(T_p M)} notiert. Ein $T\in T^2(T_p M)$ hei\ss{}t \textbf{positiv definit}, falls $T(X,X)> 0$ f\"ur alle $X\in T_p M\backslash\{ 0\}$\index{Tensor>kovarianter 2-$\sim$>positiv definiter $\sim$} gilt. F\"ur zwei Kovektoren $\omega,\eta \in T_p^{\ast} M$, wird das \textbf{Tensorprodukt}:\index{Tensor>$\sim$-Produkt}
\begin{align}\label{eq:B5}
\begin{split}
\gls{TenProd}: T_p M\times T_p M&\longrightarrow\mathbb{R}\\
(X,Y)&\mapsto\omega(X)\cdot \eta(Y)
\end{split}
\end{align}

definiert. Aus dieser Definition folgt direkt, dass $\omega\otimes\eta\in T^2(T_p M)$ gilt. In \cite[Proposition 11.2]{lee2003introduction} wird gezeigt, dass f\"ur ein $(U,\varphi)\in\mathcal{A}$ mit $p\in U$, die Menge:
\begin{equation*}
\left\{\left. ^{\varphi}dx^i\right|_{p}\otimes \left. ^{\varphi}dx^j\right|_{p}\right\}_{1\leq i, j\leq n}\subseteq T^2(T_p M)
\end{equation*}

eine Basis des Vektorraumes $T^2(T_p M)$ bildet. Daraus folgt insbesondere, dass die Dimension des Vektorraumes $T^2(T_p M)$ gleich $n^2$ ist. Neben dem in $\eqref{eq:B5}$ definierten Tensorprodukt, wird, f\"ur zwei Kovektoren $\omega,\eta\in T^{\ast}(T_p M)$, das \textbf{symmetrisierte Tensorprodukt}\index{Tensor>$\sim$-Produkt>symmetrisertes $\sim$}:
\begin{align*}
&\gls{symmTenProd}:  T_p M\times T_p M\longrightarrow\mathbb{R}\\
&\omega\eta:= \frac{1}{2}(\omega\otimes\eta+\eta\otimes\omega)
\end{align*}

verwendet. Die Menge:
\begin{equation*}
\left\{\left. ^{\varphi}dx^i\right|_{p} \left. ^{\varphi}dx^j\right|_{p}\right\}_{1\leq i\leq j\leq n}\subseteq T^2(T_p M)
\end{equation*}

bildet eine Basis des Vektorraumes $\Sigma^2(T_p M)$. Damit ist $\Sigma^2(T_p M)$ ein $\frac{n}{2}(n+1)$-dimen\-sionaler Unterraum von $T^2(T_p M)$.
Im Folgenden sei, f\"ur ein beliebiges Mengensystem $\left\{X_{\alpha} \right\}_{\alpha\in A}$, die \textbf{disjunkte Vereinigung}\index{disjunkte Vereinigung} \gls{disVerein} von $\left\{X_{\alpha} \right\}_{\alpha\in A}$ wie folgt definiert:
\begin{equation*}
\coprod_{\alpha\in A}{X_{\alpha}}:=\left\{(x,\alpha): \alpha\in A\text{ und } x\in X_{\alpha} \right\}
\end{equation*}

\begin{defi}
Sei $(M,\mathcal{A},\mathcal{O})$ eine glatte Mannigfaltigkeit, dann wird die Menge:
\begin{equation*}
\gls{BunKov}:=\coprod_{p\in M}{T^2(T_p M)}
\end{equation*}

als \textbf{B\"undel von kovarianten 2-Tensoren}\index{Tensor>B\"undel von kovarianten 2-$\sim$en} bezeichnet.

\end{defi}

Nun wird beschrieben, wie aus der Menge $T^2(M)$ ein glattes Vektorb\"undel konstruiert werden kann. Zun\"achst sei:
\begin{align*}
\pi: T^2(M)&\longrightarrow M\\
(u,p)&\mapsto p
\end{align*}

Es sei $(U,\varphi)\in\mathcal{A}$, dann wird eine Abbildung $\Phi_U$ wie folgt defniert:
\begin{align*}
\Phi_U: \pi^{-1}(U)&\longrightarrow \varphi(U)\times\mathbb{R}^{n^2}\\\
(T,p)&\mapsto\left(\varphi(p),(^{\varphi}T_{ij}(p))_{1\leq i,j\leq n}\right)
\end{align*}

wobei, unter Beachtung von \eqref{eq:B4}, die Koeffizientenmenge $(^{\varphi}T_{ij}(p))_{1\leq i,j\leq n}$ \"uber die Gleichung:
\begin{equation*}
T(p)= \ \gls{lok kov 2}(p) \left. ^{\varphi}dx^i\right|_{p}\otimes \left. ^{\varphi}dx^j\right|_{p}
\end{equation*}

eindeutig bestimmt ist. Mithilfe allgemeiner Resultate \"uber die Konstruktion von glatten Mannigfaltigkeiten \cite[Lemma 1.23]{lee2003introduction} und Verktorb\"undeln \cite[Lemma 5.5]{lee2003introduction}, wird aus diesen Karten eine Topologie und ein Atlas konstruiert, so dass $(T^2(M),M,\pi)$ ein glattes Vektorb\"undel vom Rang $n^2$ ist. Ein glatter Schnitt in diesem Vektorb\"undel wird \textbf{glattes kovariantes 2-Tensorfeld}\index{Tensor>glattes kovariantes 2-$\sim$feld} genannt. Die Menge aller glatten kovarianten 2-Tensorfelder wird mit \gls{mathcal{T}2M} bezeichnet. Mit punktweiser Addition und skalarer Multiplikation mit reellen Koeffizienten bildet diese Menge einen $\mathbb{R}$-Vektorraum.

\begin{defi}
Sei $M$ eine glatte Mannigfaltigkeit. Ein $T\in\mathcal{T}^2(M)$ hei\ss{}t \textbf{glattes symmetrisches kovariantes 2-Tensorfeld}\index{Tensor>glattes symmetrisches kovariantes 2-$\sim$-feld} auf $M$, falls f\"ur alle $p\in M$ der kovariante $2-$Tensor $T(p)$ symmetrisch ist.
\end{defi}

\begin{defi}\label{defi:B7}
Sei $(M,\mathcal{A},\mathcal{O})$ eine glatte Mannigfaltigkeit. Eine \textbf{Riemannsche Metrik}\index{Metrik>Riemannsche $\sim$} $\gls{g}\in \mathcal{T}^2(M)$ auf $M$ ist ein glattes symmetrisches kovariantes $2-$Tensorfeld auf $M$, mit der Eigenschaft, dass f\"ur alle $p\in M$ der kovariante $2-$Tensor $g(p)$ positiv definit ist.
\end{defi}

\begin{defi}
Eine \textbf{Riemannsche Mannigfaltigkeit}\index{Mannigfaltigkeit>Riemannsche $\sim$} \gls{(M,g)} ist eine glatte Mannigfaltigkeit $M$, zusammen mit einer Riemannschen Metrik $g\in \mathcal{T}^2(M)$.
\end{defi}

F\"ur ein $p\in M$ wird die Norm:
\begin{align*}
\gls{Riem Norm}: T_p M&\longrightarrow\mathbb{R}_{\geq 0}\\
X&\mapsto g(p)(X,X)
\end{align*}

als die von der \textbf{Riemannschen Metrik induzierte Norm auf}\index{Metrik>Riemannsche $\sim$>von $\sim$ induzierte Norm auf $T_p M$} $\bm{T_p M}$ bezeichnet.

Ein wichtiges Beispiel einer Riemannschen Mannigfaltigkeit ist der $\mathbb{R}^q$, versehen mit der Standardtopologie, der Standardstruktur und der \textbf{Standard-Metrik auf}\index{Metrik>Standard-$\sim$ auf $\mathbb{R}^q$} $\mathbb{R}^q$:
\begin{align*}
\gls{g^{can}}: M&\longrightarrow T^2(M)\\
p&\mapsto \sum_{i=1}^q{\left. ^{id}dx^i\right|_{p}^2}
\end{align*}

Hierbei ist $\left. ^{id}dx^i\right|_{p}^2:= \left. ^{id}dx^i\right|_{p} \left. ^{id}dx^i\right|_{p}$.

\begin{defi}
Gegeben seien glatte Mannigfaltigkeiten $M$ und $N$ sowie eine Abbildung $F\in C^{\infty}(M,N)$, dann wird die Abbildung:
\begin{align*}
\gls{Pullback}: T^2(T_{F(p)} N)&\longrightarrow T^2(T_p M)\\
(F^{\ast} S)(X_1,X_2)&:=S(F_{\ast} X_1, F_{\ast} X_2)\hspace{0.5cm}\text{f\"ur }S\in T^2(T_{F(p)} N)\text{ und }X_1,X_2\in T_p M
\end{align*}

\textbf{Pullback}\index{Pullback} genannt.

\end{defi}

Die wesentlichen Eigenschaften des Pullbacks\index{Pullback>$\sim$ eines glatten kovarianten 2-Tensorfeldes} sind in \cite[Proposition 11.8]{lee2003introduction} zusammengefasst. Der Begriff des Pullbacks wird auch f\"ur glatte Tensorfelder verwendet. Ein Pullback eines $\sigma\in\mathcal{T}^2(N)$ wird mit:
\begin{align*}
(F^{\ast}\sigma): M&\longrightarrow T^2(M)\\
p&\mapsto (F^{\ast}\sigma_{F(p)})
\end{align*}

definiert. Es gilt $F^{\ast}\sigma\in\mathcal{T}^2(M)$. F\"ur weitere Informationen diesbez\"uglich sei auf \cite[Proposition 11.9]{lee2003introduction} verwiesen.

\chapter{Hilfsresultate}
\thispagestyle{fancy}

In diesem Kapitel werden einige Lemmata bewiesen, die in den Kapiteln 4 bis 7 verwendet werden. Zuerst wird die Existenz, der in Lemma \ref{lem:4.2} verwendeten Immersion $F$, induktiv bewiesen.

\begin{lem}\label{lem:4.1}
F\"ur jedes $n\in\mathbb{N}\backslash\{0\}$ existiert eine Immersion $F^{(n)}\in C^{\infty}(\mathbb{R}^n,\mathbb{R}^{n+1})$, die in jedem Argument $2\pi$-periodisch ist, und f\"ur die gilt:
\begin{equation}\label{eq:4.16}
\partial_i F^{(n)}_l\equiv 0\hspace{0.5cm}\text{f\"ur alle }l\in\{1,...,n-1\}\text{ und }i\in\{l+1,...,n\}
\end{equation}
\end{lem}

\begin{bew}
F\"ur $n=1$ wird die Abbildung $F^{(1)}\in C^{\infty}(\mathbb{R},\mathbb{R}^{2})$ mit:
\begin{equation}\label{eq:4.83}
F^{(1)}(x):=
\begin{pmatrix}
\cos(x) \\
\sin(x)	
\end{pmatrix}
\end{equation}

definiert. Daraus wird die Abbildung $F^{(2)}\in C^{\infty}(\mathbb{R}^2,\mathbb{R}^{3})$ wie folgt konstruiert:
\begin{equation*}
F^{(2)}(x_1,x_2):=
\begin{pmatrix}
0 \\
F^{(1)}(x_2)
\end{pmatrix}
+\frac{1}{2}
\begin{pmatrix}
\cos(x_1) \\
\sin(x_1)\cdot F^{(1)}(x_2)
\end{pmatrix}
=
\begin{pmatrix}
\frac{1}{2}\cos(x_1) \\
\left( 1 +\frac{1}{2}\sin(x_1)\right)\cdot F^{(1)}(x_2)
\end{pmatrix}
\end{equation*}

Es gilt:
\begin{equation}\label{eq:4.84}
D F^{(2)}(x_1,x_2)=
\begin{pmatrix}
-\frac{1}{2}\sin(x_1) & 0 \\
\frac{1}{2}\cos(x_1) \cdot  F^{(1)}(x_2)  & \left( 1 +\frac{1}{2}\sin(x_1)\right)\cdot  \partial F^{(1)}(x_2) 
\end{pmatrix}
\end{equation}

Da $F^{(1)}(\mathbb{R})=\mathbb{S}^1$ gilt, erf\"ullt $F^{(2)}$ alle geforderten Eigenschaften. Insbesondere ist, wegen der Form der Jacobimatrix in \eqref{eq:4.84}, die Bedingung \eqref{eq:4.16} erf\"ullt. Nun sei $\epsilon\in\mathbb{R}_{>0}$, welches sp\"ater spezifiziert werden wird. Um $F^{(3)}$ aus $F^{(2)}$ zu konstruieren, wird vorher die Abbildung $\widehat{F}^{(2)}\in C^{\infty}(\mathbb{R}^2,\mathbb{R}^{3})$ mit:
\begin{equation*}
\widehat{F}^{(2)}(x_1,x_2):= \frac{\partial_1 F^{(2)}(x_1,x_2)\times \partial_2 F^{(2)}(x_1,x_2)}{\left|\partial_1 F^{(2)}(x_1,x_2)\times \partial_2 F^{(2)}(x_1,x_2)\right|_{\mathbb{R}^3}}\in \mathbb{S}^{2}\backslash lin (DF^{(2)}(x_1,x_2))
\end{equation*}

eingef\"uhrt. Aus \eqref{eq:4.84} folgt, unter Beachtung der Orthogonalit\"at der Spaltenvektoren:
\begin{align}\label{eq:4.alpha}
\begin{split}
\left|\partial_1 F^{(2)}(x_1,x_2)\times \partial_2 F^{(2)}(x_1,x_2)\right|_{\mathbb{R}^3}=&\left|\partial_1 F^{(2)}(x_1,x_2)\right|_{\mathbb{R}^3}\cdot \left|\partial_2 F^{(2)}(x_1,x_2)\right|_{\mathbb{R}^3}\\
=&\left(\frac{1}{2}+\frac{1}{4}\sin(x_1)  \right)=: f(x_1)
\end{split}
\end{align}

Nun wird die Abbildung $F^{(3)}\in C^{\infty}(\mathbb{R}^3,\mathbb{R}^{4})$ wie folgt definiert:
\begin{align*}
F^{(3)}(x_1,x_2,x_3):= &
\begin{pmatrix}
\epsilon \cos(x_1) \\
F^{(2)}(x_2,x_3) + \epsilon\sin(x_1)\cdot \widehat{F}^{(2)}(x_2,x_3)
\end{pmatrix}
\end{align*}

definiert. Es gilt:
\begin{align}\label{eq:4.85}
\begin{split}
&D F^{(3)}(x_1,x_2,x_3)=\\
&\begin{pmatrix}
-\epsilon\sin(x_1) & 0\\
 \epsilon\cos(x_1)\cdot \widehat{F}^{(2)}(x_2,x_3) & D F^{(2)}(x_2,x_3) +\epsilon \sin(x_1)\cdot D\widehat{F}^{(2)}(x_2,x_3)(x_2,x_3)
\end{pmatrix}
\end{split}
\end{align}

Ist $\epsilon\in\mathbb{R}_{>0}$ klein genug, so ist $F^{(3)}$ wegen \eqref{eq:4.85}, unter Verwendung der Periodizit\"at von $F^{(2)}$ und den Eigenschaften des Vektorproduktes, eine Immersion. Es wird noch gezeigt, dass die Bedingung $\partial_2 \widehat{F}^{(2)}_1\equiv 0$ erf\"ullt ist, woraus dann mit \eqref{eq:4.85} die Eigenschaft \eqref{eq:4.16} f\"ur $F^{(3)}$ folgt. Unter Anwendung der Produktregel gilt:
\begin{align}\label{eq:4.86}
\begin{split}
&\partial_2 \widehat{F}^{(2)}(x_1,x_2)\stackrel{\eqref{eq:4.alpha}}{=}f(x_1)\cdot\partial_2[\partial_1 F^{(2)}(x_1,x_2)\times \partial_2 F^{(2)}(x_1,x_2)]\\
=&f(x_1)\cdot \partial_2\partial_1 F^{(2)}(x_1,x_2)\times \partial_2 F^{(2)}(x_1,x_2)+f(x_1)\cdot\partial_1 F^{(2)}(x_1,x_2)\times \partial_2^2 F^{(2)}(x_1,x_2)
\end{split}
\end{align}

Mit \eqref{eq:4.84} gilt:
\begin{align*}
\partial_2 \partial_1 F^{(2)}(x_1,x_2)&=\frac{1}{2} 
\begin{pmatrix}
0 \\
\cos(x_1)\cdot\partial F^{(1)}(x_2)
\end{pmatrix} \\
\partial_2^2 F^{(2)}(x_1,x_2)&=
\begin{pmatrix}
0 \\
\left(1+ \frac{1}{2} \sin(x_1)\right)\cdot \partial^2 F^{(1)}(x_2)
\end{pmatrix}
\stackrel{\eqref{eq:4.83}}{=}-
\begin{pmatrix}
0 \\
\left(1+ \frac{1}{2} \sin(x_1)\right)\cdot  F^{(1)}(x_2)
\end{pmatrix}
\end{align*}

Dann ist mit \eqref{eq:4.86}:
\begin{align*}
\partial_2 \widehat{F}^{(2)}
=&\frac{f(x_1)}{2}
\begin{pmatrix}
\sin(x_1) \\
-\cos(x_1)\cdot  F^{(1)}(x_2)
\end{pmatrix}
\times 
\begin{pmatrix}
0 \\
\left(1+ \frac{1}{2} \sin(x_1)\right)\cdot  F^{(1)}(x_2)
\end{pmatrix}
\end{align*}

woraus $\partial_2 \widehat{F}^{(2)}_1\equiv 0$ folgt. Definiere nun f\"ur den Induktionschritt eine Abbildung $\widehat{F}^{(3)}\in C^{\infty}(\mathbb{R}^3,\mathbb{R}^4)$ mit:
\begin{equation*}
\widehat{F}^{(3)}(x_1,x_2,x_3):=
\begin{pmatrix}
\cos(x_1) \\
\sin (x_1)\cdot \widehat{F}^{(2)}(x_2,x_3)
\end{pmatrix}
\in \mathbb{S}^{3}
\end{equation*}

Falls $\epsilon\in\mathbb{R}_{>0}$ klein genug ist, dann gilt mit \eqref{eq:4.85}, wegen $|\widehat{F}^{(2)}|_{\mathbb{R}^3}\equiv 1$ und der Periodizit\"at von $F^{(2)}$ beziehungsweise $\widehat{F}^{(2)}$:
\begin{equation}\label{eq:4.87}
\widehat{F}^{(3)}(x_1,x_2,x_3)\notin lin (DF^{(3)}(x_1,x_2,x_3))
\end{equation}

f\"ur alle $(x_1,x_2,x_3)\in\mathbb{R}^3$. Ferner sei erw\"ahnt, dass $\widehat{F}^{(3)}$ die Eigenschaft \eqref{eq:4.16} erf\"ullt. Nun sei f\"ur $n\geq 3$ eine Abbildung $F^{(n)}\in C^{\infty}(\mathbb{R}^n,\mathbb{R}^{n+1})$ gegeben, die alle geforderten Eigenschaften erf\"ullt. Wegen \eqref{eq:4.87} kann  zun\"achst angenommen werden, dass  eine Abbildung $\widehat{F}^{(n)}\in C^{\infty}(\mathbb{R}^n,\mathbb{R}^{n+1})$ existiert, welche sowohl:
\begin{equation}\label{eq:Cbeta}
\widehat{F}^{(n)}(x)\in \mathbb{S}^n\backslash lin (DF^{(n)}(x))
\end{equation}

f\"ur alle $x\in\mathbb{R}^n$, als auch \eqref{eq:4.16} erf\"ullt. Dann wird f\"ur ein $\epsilon_n\in\mathbb{R}_{>0}$, welches wieder sp\"ater konkretisiert wird, die Abbildung $F^{(n+1)}\in C^{\infty}(\mathbb{R}^{n+1},\mathbb{R}^{n+2})$ folgenderma\ss{}en definiert:
\begin{equation*}
F^{(n+1)}(x_1,...,x_{n+1}):=
\begin{pmatrix}
\epsilon_n \cos(x_1) \\
F^{(n)}(x_2,...,x_{n+1})+\epsilon_n\sin(x_1)\cdot \widehat{F}^{(n)}(x_2,...,x_{n+1})
\end{pmatrix}
\end{equation*}

Dann ist:
\begin{align}\label{eq:4.88}
\begin{split}
&DF^{(n+1)}(x_1,...,x_{n+1})\\
=&
\begin{pmatrix}
-\epsilon_n \sin(x_1) & 0 \\
\epsilon_n\cos(x_1)\cdot \widehat{F}^{(n)}(x_2,...,x_{n+1}) & DF^{(n)}(x_2,...,x_{n+1}) +\epsilon_n\sin(x_1)\cdot D\widehat{F}^{(n)}(x_2,...,x_{n+1})
\end{pmatrix}
\end{split}
\end{align}

Durch eine hinreichend kleine Wahl von $\epsilon_n\in\mathbb{R}_{>0}$ wird mit \eqref{eq:Cbeta} gew\"ahrleistet, dass $F^{(n+1)}$ eine Immersion ist. F\"ur die Abbildung $\widehat{F}^{(n+1)}\in C^{\infty}(\mathbb{R}^{n+1},\mathbb{R}^{n+2})$ mit:
\begin{equation}\label{eq:4.89}
\widehat{F}^{(n+1)}(x_1,...,x_{n+1}):= 
\begin{pmatrix}
\cos(x_1) \\
\sin (x_1) \cdot \widehat{F}^{(n)}(x_2,...,x_{n+1})
\end{pmatrix}
\in \mathbb{S}^{n+1}
\end{equation}

kann dann angenommen werden, dass $\epsilon_n\in\mathbb{R}_{>0}$ klein genug ist, so dass:
\begin{equation}
\widehat{F}^{(n+1)}(x)\notin lin (DF^{(n+1)}(x))
\end{equation}

f\"ur alle $x\in \mathbb{R}^{n+1}$ gilt. Mit \eqref{eq:4.88} ist die Bedingung \eqref{eq:4.16} nach Induktionsvoraussetzung erf\"ullt.

\end{bew}

Die im folgenden Lemma notierten Abbildungen werden in Lemma \ref{lem:4.5}, speziell \eqref{eq:4.30}, daf\"ur verwendet, die Abbildung $v_1\in C^{\infty}(\mathbb{S}\times\overline{\mathbb{B}},\mathbb{R}^q)$ zu konstruieren.

\begin{lem}\label{lem:4.3}
Die Abbildungen $\alpha_1,\alpha_2\in C^{\infty}(\mathbb{S})$ mit:
\begin{align*}
\alpha_1(t)&:= \cos(t)-\frac{1}{\sqrt{3}}\cos(3t)\\
\alpha_2(t)&:= \sin(t)+\frac{1}{\sqrt{3}}\sin(3t)
\end{align*}

erf\"ullen die Eigenschaften:
\begin{equation}\label{eq:4.22}
\alpha_1'(t)\,\alpha_2''(t)-\alpha_2'(t)\,\alpha_1''(t)\neq 0 \hspace{0.5cm}\forall t\in\mathbb{S}
\end{equation}
\begin{equation}\label{eq:4.23}
\int_{-\pi}^{\pi}{\alpha_i'(t)\, \alpha_j(t) \ dt}=0 \hspace{0.5cm}\text{f\"ur }i,j\in\{1,2\}
\end{equation}
\begin{equation}\label{eq:4.24}
\int_{-\pi}^{\pi}{\alpha_i(t)\,\sqrt{\alpha_1'(t)^2+\alpha_2'(t)^2} \ dt}=0 \hspace{0.5cm}\text{f\"ur }i\in\{1,2\}
\end{equation}

Ferner ist die periodische Fortsetzung auf ganz $\mathbb{R}$ glatt.

\end{lem}

\begin{bew}
Es gilt f\"ur alle $t\in \mathbb{S}$:
\begin{align*}
&\alpha_1'(t)=-\sin(t)+\sqrt{3} \sin(3t) & & \alpha_2'(t)=\cos(t)+\sqrt{3} \cos(3t)\\
&\alpha_1''(t)=-\cos(t)+3\sqrt{3} \cos(3t) & & \alpha_2''(t)=-\sin(t)-3 \sqrt{3} \sin(3t)
\end{align*}

Dann ist:
\begin{align}\label{eq:4.21}
\notag&\alpha_1'(t)\,\alpha_2''(t)-\alpha_2'(t)\,\alpha_1''(t)\\
\notag=&(-\sin(t)+\sqrt{3} \sin(3t))(-\sin(t)-3 \sqrt{3} \sin(3t))\\
\notag&-(\cos(t)+\sqrt{3} \cos(3t))(-\cos(t)+3\sqrt{3} \cos(3t))\\
\notag=&\sin^2(t)+3\sqrt{3}\sin(t)\sin(3t)-\sqrt{3} \sin(3t)\sin(t)-9\sin^2(3t)\\
\notag&+\cos^2(t)-3\sqrt{3}\cos(t)\cos(3t)+\sqrt{3} \cos(3t)\cos(t)-9\cos^2(3t)\\
\notag=&2\sqrt{3}\sin(t)\sin(3t)-2\sqrt{3}\cos(t)\cos(3t)-8\\
=&2\sqrt{3}\,(\sin(t)\sin(3t)-\cos(t)\cos(3t))-8\\
\notag=&2\sqrt{3}\,\left[3\sin^2(t)-4\sin^4(t)-4\cos^4(t)+3\cos^2(t)\right]-8\\
\notag=&2\sqrt{3}\,\left[-4(\sin^4(t)+\cos^4(t))+3\right]-8\\
\notag=&-8\sqrt{3}\,(\sin^4(t)+\cos^4(t))+6\sqrt{3}-8
\end{align}

Die Funktion $t\mapsto \sin^4(t)+\cos^4(t)$ wird genauer untersucht. Es gilt:
\begin{align*}
&\frac{d}{dt} (\sin^4(t)+\cos^4(t))=4(\sin^3(t)\cos(t)-\cos^3(t)\sin(t))\\
=&4\sin(t)\cos(t)\cdot(\sin^2(t)-\cos^2(t))=4\sin(t)\cos(t)\cdot (2\sin^2(t)-1)\\
=&2\sin(2t)\cdot (2\sin^2(t)-1)=-2\sin(2t)\cos(2t)=-\sin(4t)
\end{align*}

Daraus ergibt sich:
\begin{equation*}
\min_{t\in\mathbb{S}}{(\sin^4(t)+\cos^4(t))}=\frac{1}{2}
\end{equation*}

Mit \eqref{eq:4.21} folgt dann:
\begin{equation*}
\alpha_1'(t)\,\alpha_2''(t)-\alpha_2'(t)\,\alpha_1''(t)\leq 2\sqrt{3}-8<0
\end{equation*}

womit \eqref{eq:4.22} gezeigt ist. Nun sei $i\in\{1,2\}$ dann ist:
\begin{equation*}
\int_{-\pi}^{\pi}{\alpha_i'(t)\,\alpha_i(t)\> dt}=\frac{1}{2}\int_{-\pi}^{\pi}{\frac{d}{dt}(\alpha_i^2)(t)\> dt}=\frac{\alpha_i^2(\pi)-\alpha_i^2(-\pi)}{2}=0
\end{equation*}

Weiterhin ist:
\begin{align*}
&\int_{-\pi}^{\pi}{\alpha_1'(t)\,\alpha_2(t)\> dt}=\int_{-\pi}^{\pi}{(-\sin(t)+\sqrt{3} \sin(3t))\cdot(\sin(t)+\frac{1}{\sqrt{3}}\sin(3t))  \> dt}\\
=&-\int_{-\pi}^{\pi}{\sin^2(t)  \> dt}+\left(\sqrt{3}-\frac{1}{\sqrt{3}}\right)\int_{-\pi}^{\pi}{\sin(t)\cdot \sin(3t) \> dt}+\int_{-\pi}^{\pi}{\sin^2(3t)  \> dt}\\
=&-\pi+0+\pi=0
\end{align*}

Daraus folgt mittels partieller Integration:
\begin{align*}
&\int_{-\pi}^{\pi}{\alpha_2'(t)\,\alpha_1(t)\> dt}=\alpha_2(\pi)\, \alpha_1(\pi)-\alpha_2(-\pi)\, \alpha_1(-\pi)-\int_{-\pi}^{\pi}{\alpha_2(t)\, \alpha_1'(t)\> dt}=0\\
\end{align*}

womit \eqref{eq:4.23} gezeigt ist. Es bleibt noch \eqref{eq:4.24} zu zeigen. Zun\"achst gilt:
\begin{align*}
&\alpha_1'(t)^2+\alpha_2'(t)^2=(-\sin(t)+\sqrt{3} \sin(3t))^2+(\cos(t)+\sqrt{3} \cos(3t))^2\\
=&\sin^2(t)-2\sqrt{3}\sin(t)\,\sin(3t)+3\sin^2(3t)+\cos^2(t)+2\sqrt{3}\cos(t)\,\cos(3t)+3\cos^2(3t)\\
=&-2\sqrt{3}\,(\sin(t)\,\sin(3t)-\cos(t)\,\cos(3t))+4\\
=&-2\sqrt{3}\,(3\sin^2(t)-4\sin^4(t)-4\cos^4(t)+3\cos^2(t))+4\\
=&-2\sqrt{3}\,(-4(\sin^4(t)+\cos^4(t))+3)+4\\
=&8\sqrt{3}\,(\sin^4(t)+\cos^4(t))-6\sqrt{3}+4
\end{align*}

Somit ist die Funktion $t\mapsto \sqrt{\alpha_1'(t)^2+\alpha_2'(t)^2}$ symmetrisch und es folgt, aus der Antisymmetrie der Abbildung $\alpha_2$, die Bedingung \eqref{eq:4.24} f\"ur den Fall $i=2$. Da $\alpha_1$ um $t=\frac{3}{2}\pi$ antisymmetrisch ist, folgt \eqref{eq:4.24} auch f\"ur den Fall $i=1$.

\end{bew}

Das folgende Lemma spielt bei der Betrachtung der Abbildung $F_{\epsilon,k}$ am Anfang von \autoref{sec:4.2} eine wichtige Rolle. Dort wird $\beta\in C^{\infty}(\mathbb{R})$ konkret definiert.

\begin{lem}\label{lem:C.3}
F\"ur jedes $s\in\mathbb{N}\backslash\{0\}$ gilt:
\begin{equation*}
\sup_{t\in\mathbb{R}}{|\partial^s \beta(t)|}<\infty
\end{equation*}

\end{lem}

\begin{bew}
Mit \eqref{eq:4.39} gilt:
\begin{equation*}
\sup_{t\in\mathbb{R}}{|\beta'(t)|}=\frac{1}{\min_{t\in\mathbb{S}}{|\varrho(t)|}}\stackrel{\eqref{eq:4.34}}{<}\infty
\end{equation*}

womit die Aussage f\"ur $s=1$ gezeigt ist. Nun sei die Aussage f\"ur $\{1,...,s\}\subseteq\mathbb{N}\backslash\{0\}$ richtig. Mit der Formel von Fa\`{a} di Bruno \cite{tortolini1855annali} gilt:
\begin{equation}\label{eq:4.42}
\partial^{s+1}(P\circ \beta)=\sum_{(k_1,...,k_{s+1})\in T_{s+1}}{\frac{n!}{k_1!\cdot...\cdot k_{s+1}!} \left(\partial^{k_1+...+k_{s+1}}P\circ \beta \right)\prod_{\substack{m=1\\k_m\geq 1}}^{s+1}\left(\frac{\partial^m \beta}{m!} \right)^{k_m}}
\end{equation}

mit $T_{s+1}:= \left\{(k_1,...,k_{s+1})\in\mathbb{N}^{s+1}:\sum_{j=1}^{s+1}{j k_j}=s+1\right\}$. Wegen $P\circ \beta=id_{\mathbb{R}}$ ist mit \eqref{eq:4.42}:
\begin{align*}
&0\equiv\partial^{s+1}(P\circ \beta)\\
=&\sum_{(k_1,...,k_{s+1})\in \widehat{T}_{s+1}}{\frac{n!}{k_1!\cdot...\cdot k_{s+1}!} \left(\partial^{k_1+...+k_{s+1}}P\circ \beta \right)\prod_{\substack{m=1\\k_m\geq 1}}^{s+1}\left(\frac{\partial^m \beta}{m!} \right)^{k_m}}+n!\cdot  P'\circ \beta\cdot \frac{\partial^{s+1}\beta}{(s+1)!}
\end{align*}

mit $\widehat{T}_{s+1}=T_{s+1}\backslash\{(0,...,0,1)\}$. Es folgt, wegen $P'(t)=\varrho(t)\neq 0$ f\"ur alle $t\in\mathbb{R}$:
\begin{align*}
\partial^{s+1}\beta=\frac{(s+1)!}{n!}\cdot \frac{-1}{ P'\circ \beta}\sum_{(k_1,...,k_{s+1})\in \widehat{T}_{s+1}}{\frac{n!}{k_1!\cdot...\cdot k_{s+1}!} \left(\partial^{k_1+...+k_{s+1}}P\circ \beta \right)\prod_{\substack{m=1\\k_m\geq 1}}^{s+1}\left(\frac{\partial^m \beta}{m!} \right)^{k_m}}
\end{align*}

Da die letzte Komponente eines Multiindex in $\widehat{T}_{s+1}$ stets $0$ ist, folgt die Behauptung aus der Induktionsvoraussetzung, unter Beachtung der Periodizit\"at von $\varrho$.

\end{bew}

Das folgende Lemma wird in \autoref{chap:Kap5} und \autoref{chap:Kap7} daf\"ur verwendet, um die Ableitungen der Approximationsfolgeglieder gleichm\"a\ss{}ig abzusch\"atzen.

\begin{lem}\label{lem:7.8}
Es sei $(a_k)_{k\in\mathbb{N}}\subseteq \mathbb{R}_{\geq 0}$ so dass f\"ur alle $k\in\mathbb{N}$ die Absch\"atzung:
\begin{equation*}
a_{k+1}\leq \frac{a_k+C}{2}
\end{equation*}

f\"ur eine feste, nicht von $k$ abh\"angige Konstante $C\in\mathbb{R}_{\geq 0}$ gilt, dann gilt f\"ur alle $k\in\mathbb{N}$ die Absch\"atzung:
\begin{equation*}
a_{k}\leq a_0+C
\end{equation*}

\end{lem}

\begin{bew}
F\"ur $k=0$ ist die Aussage richtig. Unter der Annahme, dass die Aussage f\"ur ein $k\in\mathbb{N}$ gilt, ist:
\begin{equation*}
a_{k+1}\leq \frac{a_0+2C}{2} \leq a_0+C
\end{equation*}

womit die Aussage bewiesen ist.

\end{bew}

\pagestyle{empty}

\setcounter{secnumdepth}{-1} 

\addcontentsline{toc}{chapter}{Notation}
\printglossaries

\bibliography{diplomarbeit}

\addcontentsline{toc}{chapter}{Index}
\printindex

\section*{Erklärung}		
\thispagestyle{empty}

Hiermit erkläre ich, dass ich die vorliegende Arbeit selbstständig und ohne Benutzung anderer als der angegebenen Quellen und Hilfsmittel angefertigt habe.
\newline\newline\newline
Magdeburg, den 19. November 2013 \hspace{6cm} Norman Zerg\"ange

\end{document}